\documentclass[10pt]{article}
\usepackage[left=0.75in,right=0.75in,top=0.75in,bottom=0.75in]{geometry} 
\usepackage{amsmath}        
\usepackage{graphicx}
\usepackage{amssymb}
\usepackage{amsthm}
\usepackage{epstopdf}
\usepackage{amsfonts}
\usepackage{hyperref}
\makeatletter
\def\amsbb{\use@mathgroup \M@U \symAMSb}
\makeatother
\usepackage{bbold}
\usepackage{xcolor}
\usepackage{booktabs}
\usepackage{bm}
\usepackage{placeins}
\usepackage[mathscr]{euscript}
\usepackage{subfig}
\usepackage{authblk}
\usepackage{thmtools}
\usepackage[]{algorithm2e}
\usepackage{epigraph}

\DeclareFontFamily{U}{mathx}{\hyphenchar\font45}
\DeclareFontShape{U}{mathx}{m}{n}{
      <5> <6> <7> <8> <9> <10>
      <10.95> <12> <14.4> <17.28> <20.74> <24.88>
      mathx10
      }{}
\DeclareSymbolFont{mathx}{U}{mathx}{m}{n}
\DeclareMathSymbol{\bigtimes}{1}{mathx}{"91}

\usepackage[sectionbib]{natbib}
\usepackage{chapterbib}

\usepackage{tikz}

\newtheorem{thm}{Theorem}
\newtheorem{cor}{Corollary}
\newtheorem{de}{Definition}
\newtheorem{re}{Remark}
\theoremstyle{lemma}

\newtheorem{prop}{Proposition}

\def\d{\delta}
\def\c{\circ}
\def\th{\theta}
\def\la{\lambda}
\def\E{\amsbb{E}}
\def\Var{\amsbb{V}ar}
\def\Cov{\amsbb{C}ov}
\def\C{\amsbb{C}ov}

\def\P{\amsbb{P}}
\def\R{\amsbb{R}}

\def\N{\amsbb{N}}
\def\Z{\amsbb{Z}}

\def\D{\mathrm{d}}

\newcommand\ind[1]{\amsbb{I}_{#1}}
\providecommand{\norm}[1]{\lVert#1\rVert}

\DeclareMathOperator*{\argmin}{arg\,min}

\DeclareMathOperator*{\Erf}{Erf}

\DeclareMathOperator*{\supp}{supp}

\DeclareMathOperator*{\for}{\quad\text{for}\quad}

\title{Uncertainty quantification by random measures and fields}

\author[1,2]{Caleb Deen Bastian\thanks{\href{mailto:cbastian@princeton.edu}{cbastian@princeton.edu}}}
\affil[1]{Program in Applied and Computational Mathematics, Princeton University, Princeton, NJ. 08544, USA}
\affil[2]{Massive Dynamics, Princeton, NJ. 08542, USA}
\author[1,3]{Herschel Rabitz}
\affil[3]{Department of Chemistry, Princeton University, Princeton, NJ. 08544, USA}
\date{\today}      

\begin{document}
%

\maketitle

\begin{abstract} We present a general framework for uncertainty quantification that is a mosaic of interconnected models. We define global first and second order structural and correlative sensitivity analyses for random counting measures acting on risk functionals of input-output maps. These are the ANOVA decomposition of the intensity measure and the decomposition of the random measure variance, each into subspaces. Orthogonal random measures furnish sensitivity distributions. We show that the random counting measure may be used to construct positive random fields, which admit decompositions of covariance and sensitivity indices and may be used to represent interacting particle systems. The first and second order global sensitivity analyses conveyed through random counting measures elucidate and integrate different notions of uncertainty quantification, and the global sensitivity analysis of random fields conveys the proportionate functional contributions to covariance. This framework complements others when used in conjunction with for instance algorithmic uncertainty and model selection uncertainty frameworks. 
\end{abstract}


\section{Introduction} Uncertainty quantification (UQ) is a fundamental area. There are multiple notions of UQ in various applications. Perhaps the most common notion is propagation of input uncertainty through a model in uncertainty propagation (UP), such as in global sensitivity analysis based on decomposition of variance. Another is uncertainty of estimators, a second-order analysis based on decomposition of variance. Another is focusing on entropy, instead of variance. Another is model uncertainty, such as that of Bayesian model averaging. 

We shall assume the reader has access to some independency of input-output random variables $(\mathbf{X},\mathbf{Y})=\{(X_i,Y_i)\}$ taking values in measurable space $(E\times F,\mathcal{E}\otimes\mathcal{F})$ with distribution $\nu\times Q$, where $\nu$ is a probability measure on $(E,\mathcal{E})$ and $Q$ is a transition probability kernel from $(E,\mathcal{E})$ into $(F,\mathcal{F})$. Oftentimes $Q$ is deterministic and specified through function $g:E\mapsto F$ as $Q(x,\cdot)=\delta_{g(x)}(\cdot)$. For example, $g$ could be a \emph{regressor} or a \emph{classifier} or some other function. Here we let $\mathcal{E}_+$ denote the collection of positive $\mathcal{E}$-measurable functions. Suppose we have some \emph{risk function} $f\in\mathcal{E}_+$ that is a function of $g$. The statistics of $f$ are readily computed as \begin{align*}\E f &=\nu f\\\Var f &= \nu f^2 - (\nu f)^2\end{align*} 

 Many times, one estimates these statistics using $n$ samples using the empirical distribution \[F_n(A) = \frac{1}{n}\sum_{i}^n\ind{A}(X_i)\for A\in\mathcal{E}\] in view of the fact that \[\lim_{n\rightarrow\infty}F_n f=\nu f\for f\in\mathcal{E}_+\] almost surely. Further for $f\in\mathcal{E}_+$ we have \begin{align*}\E F_n f &= \nu f\\\Var F_n f &= \frac{1}{n}\Var f\end{align*} and for $f,g\in\mathcal{E}_+$ we have \[\Cov(F_nf,F_ng) = \frac{1}{n}(\nu(fg)-\nu f\nu g)\] If $f$ and $g$ are disjoint, then the covariance simplifies to $\Cov(F_nf,F_ng)=-\frac{1}{n}\nu f\nu g$. The empirical distribution has negative covariance for disjoint functions and so its uncertainties are correlated.  The Monte Carlo estimator of the expected value with respect to the empirical distribution is called the bootstrap estimator \citep{bootstrap1}. As we shall see later in the article, the correlation destroys the probabilistic interpretation of the normalized variances.

Another area of focus is on $g$. In most settings, $E$ is finite with dimension $n$. For regression, a common tool is \emph{global sensitivity analysis} to understand the propagation of input uncertainty. This conducts a functional analysis of variance (ANOVA) analysis which conveys a decomposition of variance into subspaces. Functional HDMR has other remarkable properties, such as in many practical settings low-order expansions well approximate $g$ and the fact that the expansion is exact and finite. 

In this article we introduce a counting distribution $\kappa$, which in conjunction with $\nu\times Q$, forms a random counting measure. This framework integrates random measure and mean measure uncertainty quantification for risk functionals of the form $f=(g-\E g)^2\in\mathcal{E}_+$. We introduce a positive measurable mapping $k$ that, in conjunction with the random counting measure $(\kappa,\nu\times Q)$, forms a random field $(\kappa,\nu\times Q,k)$. Each aspect has a foundational theorem. As a `case-control' analysis, we compare the binomial process $nF_n$ (control) to the mixed binomial process (case) $N$ with $\kappa$ having positive variance throughout the article. Each analysis furnishes an ANOVA. The relations among these models are shown below in Figure~\ref{fig:relations}.  

Importantly, despite the generality of the RM-MM-RF-ANOVA UQ framework, the presentation here is not of a complete system---in this article we do not consider algorithmic uncertainty neither do we consider uncertainty of model selection, such as Bayesian model averaging; nor do we consider real-valued fields, as our random fields are necessarily positive being derived from (driven by) random counting measures. We believe the RM-MM-RF-ANOVA is coherent and complements other frameworks such as those mentioned but not contained.  

\begin{figure}[h]
\centering
\includegraphics[width=7in]{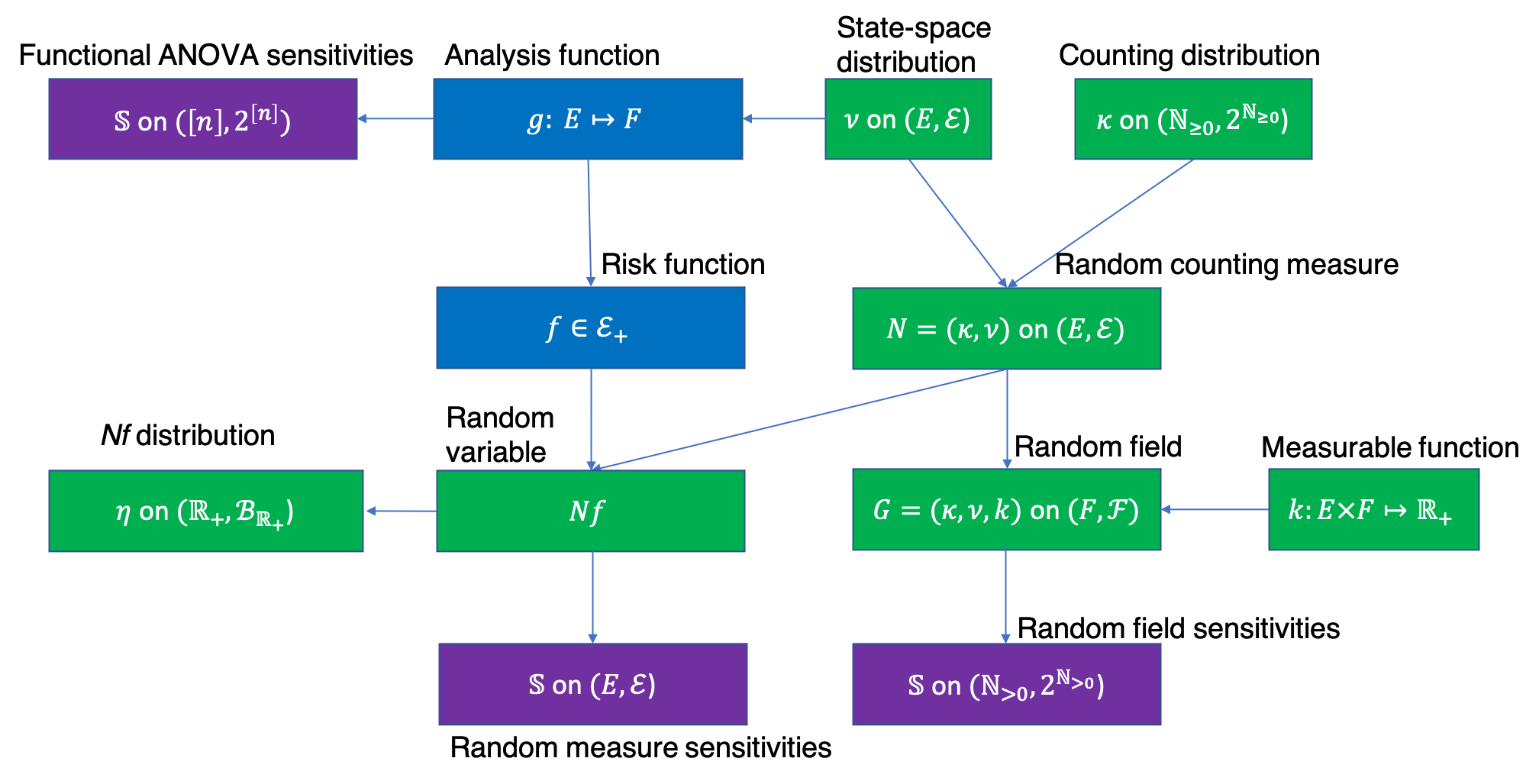}
\caption{Relations among mathematical models ({\color{green}\textbf{green}} = random variables, measures, and fields, distributions, and positive functions; {\color{blue}\textbf{blue}} = general and positive valued test functions; {\color{purple}\textbf{purple}} = sensitivities)}\label{fig:relations}
\end{figure}
 
 \FloatBarrier

In Section~\ref{sec:mixed} we give the formal backdrop in terms of the mixed binomial process, which generalizes the binomial process (bootstrap estimator). In Section~\ref{sec:hdmr} we give a brief description of functional ANOVA (high dimensional model representation) and global sensitivity analysis. In Section~\ref{sec:anova} we give ANOVA decompositions of the mean (intensity measure) and variance of the random counting measure into subspaces. In Section~\ref{sec:related} we discuss related work. In Section~\ref{sec:app} we apply this framework to risk functionals. In Section~\ref{sec:rf} we discuss random fields built from random measures and their ANOVA decompositions. In Section~\ref{sec:sim} we perform UQ comparing analysis and interpretation across different counting distributions and mean measures. In Section~\ref{sec:discuss} we end with discussion and conclusions. In the appendices, we provide additional examples. In Table~\ref{tab:0}, we show the examples, their locations in this article, and the principal ANOVA analyses conducted for the examples. 

\begin{table}[h!]
\begin{center}
\begin{tabular}{p{2in}cp{1.5in}lccc}
\toprule
Example Name & Section & Law $\nu$ & Dimension $n$ & MM & RM  & RF \\\midrule
Elementary symmetric polynomial &  \ref{sec:elem} & $\nu: \text{mean}\ne 0$, Bernoulli & $\N_{>0}$, 100, 1 & $\bullet$ &  $\bullet$ &\\
Ishigami function & \ref{sec:ishi} & $\text{Uniform}[-\pi,\pi]^3$ & 3 & $\bullet$ &  $\bullet$ &\\
Gaussian process regressor & \ref{sec:regress} & general, Ishigami & any, 3 & $\bullet$ &$\bullet$ &\\
$m$-class classifier & \ref{sec:class} & general & any &$\bullet$ & $\bullet$ &\\
Interacting particle systems &\ref{sec:sda}& general, Wiener  & $\infty$ &  & $\bullet$ & $\bullet$\\
Adaptive randomized controlled trials & \ref{ex:rct} & general & any &$\bullet$ &$\bullet$ &\\
Dynamic survival analysis of epidemics & \ref{ex:dsa} & SIR PT network & 2 & &$\bullet$ & $\bullet$\\\midrule
Symmetric polynomial with correlation & \ref{sec:poly} & Gaussian & 2 & $\bullet$  & $\bullet$ &\\
Graph property & \ref{sec:graph} & Erd\H{o}s-Reny\'{i}& $\N_{>0}$, 10 &$\bullet$  & $\bullet$ &\\
Ising model & \ref{sec:ising} &nearest-neighbor Ising & any & &  & $\bullet$ \\\bottomrule
\end{tabular}\caption{Article examples of RM-MM-RF ANOVAs (RM = random measure; MM = mean measure; RF = random field)}\label{tab:0}
\end{center}
\end{table}


\begin{table}[h!]
\begin{center}
\begin{tabular}{lp{4in}}
\toprule
Symbol &Definition \\\midrule
$\N_{\ge 0}$ & the set of the natural numbers $\{0,1,\dotsc\}$\\
$\R$ & the set of the real numbers $(-\infty,\infty)$\\
$\R_+$ & the set of the non-negative real numbers $[0,\infty)$\\
$(E,\mathcal{E})$ & measurable space\\
$\nu$ & probability measure on $(E,\mathcal{E})$\\
$\bm{X}$ & independency of random variables $\{X_i\}$ with law $\nu$ taking values in $(E,\mathcal{E})$\\
$(F,\mathcal{F})$ & measurable space\\
$(E\times F,\mathcal{E}\otimes\mathcal{F})$ & measurable product space\\
$Q$ & transitional probability kernel from $(E,\mathcal{E})$ into $(F,\mathcal{F})$\\
$\mu=\nu\times Q$ & probability measure $\mu(\D x, \D y)=\nu(\D x)Q(x,\D y)$ on $(E\times F,\mathcal{E}\otimes\mathcal{F})$\\
$\bm{Y}$ & independency of random variables $\{Y_i\}$ with law $Q(X_i,\cdot)$ taking values in $(F,\mathcal{F})$\\
$(\bm{X},\bm{Y})$ & independency of random variables $\{(X_i,Y_i)\}$ with law $\nu\times Q$ taking values in $(E\times F,\mathcal{E}\otimes\mathcal{F})$\\
$K$ & $\N_{\ge0}$-valued random variable\\
$\psi$ & probability generating function of $K$\\
$N=(\kappa,\nu)$ & random counting measure on $(E,\mathcal{E})$\\
$L$ & Laplace functional of $N$ \\
$\mathcal{E}_+$ & the collection of positive $\mathcal{E}$-measurable functions\\
$Nf$ & random variable in $(\R_+,\mathcal{B}_{\R_+})$ for function $f\in\mathcal{E}_+$\\
$F$ & Laplace transform of $Nf$\\
$\eta$ & distribution of $Nf$\\
$\eta_n$ & distribution of $Nf$ attained by maximum entropy with $n$ evaluations of $F$\\
$N=(N_t)_{t\in\R_+}$ & counting process of random measure $M=(\kappa,\nu)$ on $(\R_+,\mathcal{B}_{\R_+})$ \\
$\mathscr{F}=(\mathscr{F}_t)_{t\in\R_+}$ & filtration generated by the counting process \\
$PT(c,\d^2)$ & Poisson-type (PT) distribution with mean $c>0$ and variance $\d^2>0$ with family members binomial, Poisson, and negative binomial\\
$M=(K,\nu\times Q)$ & random counting measure on $(E\times F,\mathcal{E}\otimes\mathcal{F})$\\
$\delta_c$ & Dirac measure sitting at $c\in\N_{\ge0}$\\
$L^2(E,\mathcal{E},\nu)$ & space of square-integrable functions\\
$k$ & positive $\mathcal{E}\otimes\mathcal{F}$ measurable function\\
$G=(\kappa,\nu,k)$ & positive random field on $(F,\mathcal{F})$ formed by $N=(\kappa,\nu)$ and $k$\\
$U$ & mean of random field\\
$C$ & covariance of random field \\
\end{tabular}\caption{Symbols and definitions}\label{tab:00}
\end{center}
\end{table}

\FloatBarrier

\section{Background} In this section we give the mathematical backdrop of the article. In Section~\ref{sec:mixed} we discuss the mixed binomial process, whereas in Section~\ref{sec:hdmr} we describe global sensitivity analysis as conveyed through functional ANOVA or high dimensional model representation (HDMR).

\subsection{Mixed binomial process}\label{sec:mixed}
Let $(E,\mathcal{E})$ be a measure space and let $\nu$ be a probability measure on it. Let  $\mathbf{X}=\{X_i\}$ be an independency (collection) of (iid) $E$ valued random variables with law $\nu$. Let $K\sim\kappa$ be a $\N_{\ge 0}$-valued random variable independent of $\mathbf{X}$ with mean $c>0$ and variance $\d^2\ge0$. The \emph{mixed binomial process} is identified to the pair of deterministic probability measures $N=(\kappa,\nu)$ on $(E,\mathcal{E})$ through \emph{stone throwing construction}\citep{cinlar,bastian,rm} (STC) as \begin{equation}\label{eq:stc}N(A)=N\ind{A}=\int_{E}N(\D x)\ind{A}(x)\equiv\sum_{i}^{K}\ind{A}(X_i)\for A\in\mathcal{E}\end{equation} where $\ind{A}$ is a set function. We denote $\mathcal{E}_+$ the set of non-negative $\mathcal{E}$-measurable functions. $N$ is said to be \emph{orthogonal} if $c=\delta^2$: for disjoint $f,g\in\mathcal{E}_+$ we have $\Var(Nf+Ng)=\Var Nf+\Var Ng$.

We point out that the STC of the mixed binomial process is also known as a \emph{proper point process}. Similar ideas are of \emph{randomly-stopped sums}. Our interpretation and treatment here is from the random counting measure perspective. 

For $f\in\mathcal{E}_+$ we have \emph{mean} and \emph{variance} of the random variable $Nf$ \begin{align}\E Nf &= c\nu f\label{eq:mu0}\\\Var Nf &= c\nu f^2 + (\delta^2-c)(\nu f)^2\label{eq:var0}\end{align} The mean $c\nu$ is also known as the \emph{intensity measure} of $N$. For arbitrary $f,g\in\mathcal{E}_+$, we have \emph{covariance} \begin{equation}\label{eq:cov0}\Cov(Nf,Ng)=c\nu(fg)+(\delta^2-c)\nu f\nu g\end{equation} For disjoint partition $A,\dotsb,B$ of $E$, we have joint probability \begin{align*}
\P(N(A)=i,\ldots, N(B)=j)&=\P(N(A)=i,\ldots, N(B)=j|K=k)\P(K=k)\\
&=\frac{k!}{i!\cdots j!}\,\nu^i(A)\cdots \nu^j(B)\P(K=k).
\end{align*} 

Recall that the law of $N$ is uniquely determined by the \emph{Laplace functional} $L$ from $\mathcal{E}_+$ into $[0,1]$ \begin{equation}\label{eq:laplace} L(f)=\E e^{-Nf} = \E e^{-\int_EN(\D x)f(x)}= \psi(\nu e^{-f})\for f\in\mathcal{E}_+\end{equation} where $\psi$ is the \emph{probability generating function} (pgf) of $K$. The Laplace functional encodes all the information of $N$: its distribution, moments, etc. The \emph{moments} of $Nf$ (if they exist) can be attained from the Laplace functional \[ \E(Nf)^n = (-1)^n\lim_{q\downarrow0}\frac{\partial^n}{\partial q^n}L(qf)\quad\text{for}\quad n\in\N_{>0}\] The distribution of $Nf$, denoted by $\eta$, i.e. $\eta(\D x)=\P(Nf\in \D x)$, is encoded by the Laplace transform, which may be expressed in terms of the Laplace functional \begin{equation}F(\alpha) = \E e^{-\alpha Nf} =  \E e^{-N(\alpha f)}=L(\alpha f)\for \alpha\in\R_+\end{equation}

We can replace $\nu$ with $\nu\times Q$ in all expressions for the random measure $M=(\kappa,\nu\times Q)$ on $(E\times F,\mathcal{E}\otimes\mathcal{F})$. 

 Given the mixed binomial process $N=(\kappa,\nu)$ on $(E,\mathcal{E})$, a number of operations can be considered. One operation is \emph{restriction} (thinning) of $N$ to subspaces $A\subset E$ with mass $a=\nu(A)>0$, giving restricted random measure $N_A=(N\ind{A},\nu_A)$, where $\nu_A(\cdot)=\nu(A\cap\cdot)/\nu(A)$ is the restricted law. The law of $N_A$ is encoded in its Laplace functional \begin{equation}L_A(f) =\psi(a\nu_A e^{-f}+1-a)\for f\in\mathcal{E}_+\end{equation} Another is the \emph{image} (transformation) random measure $N\circ h^{-1}=(\kappa,\mu=\nu\circ h^{-1})$ on $(F,\mathcal{F})$ for some measurable transformation $h:E\mapsto F$ with Laplace functional defined through \[L(f)=\psi(\mu e^{-f})\for f\in\mathcal{F}_+\] In the article we use a random transformation (transition probability kernel) $Q$ where $\mu=\nu Q$.

For $(E,\mathcal{E})=(\R_+,\mathcal{B}_{\R_+})$, we can define a \emph{counting process} $N=(N_t)_{t\in\R_+}$ from the random measure $M=(\kappa,\nu)$ on $(E,\mathcal{E})$ as \[N_t = M([0,t])\for t\in\R_+\] i.e. a process with state-space $(\N_{\ge0},2^{\N_{\ge0}})$ whose every path $t\mapsto N_t$ starts with $N_0=0$ and is increasing and right-continuous with jumps of size one, with generated \emph{filtration} $\mathscr{F}=(\mathscr{F}_t)_{t\in\R_+}$ where \[\mathscr{F}_t = \sigma\{N_s: s\le t\} = \sigma\{M(A): A\in\mathcal{B}_{[0,t]}\}\] $N$ is automatically \emph{adapted} to $\mathscr{F}$, the filtration it generates. 


\subsection{High dimensional model representation, global sensitivity analysis}\label{sec:hdmr}
Suppose we have measurable space $(F,\mathcal{F})$ and some mapping $g:E\mapsto F$ that is measurable relative to $\mathcal{E}$ and $\mathcal{F}$. We assume that either $F=\R$ or $F$ is atomic and finite with dimension $m$. If the first case we assume $g\in L^2(E,\mathcal{E},\nu)$ and identify such $g$ as a \emph{regressor}. In the latter case, we identify $g$ to a \emph{classifier}. Hereafter we take $E$ as finite with dimension $n$. All such regressors may be decomposed into a hierarchy of \emph{high dimensional model representation} (HDMR) \citep{alis99} component functions \[g(x_1,\dotsb,x_n) = g_0 + \sum_ig_i(x_i) + \sum_{i<j}g_{ij}(x_i,x_j)+\dotsb+g_{1\dotsb n}(x_1,\dotsb,x_n)\] which convey a decomposition of variance \[\Var g = \sum_{u\subseteq\{1,\dotsb,n\}}\Var g_u + \sum_{u\ne v\subseteq\{1,\dotsb,n\}}\Cov(g_u,g_v)\] The component functions $\{g_u\}$ are hierarchically orthogonal under inner product $\langle\cdot,\cdot\rangle_{\nu}=\Cov(\cdot,\cdot)$ \citep{hooker}. This is also known as the \emph{functional ANOVA expansion} or the \emph{Sobol system} of $g$. When the inputs are independent $\nu=\prod_i\nu_i$, then the component functions are mutually orthogonal and may be recursively constructed as follows. Defining the operator \[\mathbf{M}^{i_1\dotsb i_l}g(x) =\int_{\{0,1\}^{\binom{n}{2}-l}}\prod_{j: j\notin\{i_1,\dotsb,i_l\}}\nu_j\{x_j\}g(x_1,\dotsb,x_{n}) \] the HDMR of $g$ for $\nu=\prod_i\nu_i$ is computed recursively, \begin{align*}g_0&= \mathbf{M}g(x)\\g_i(x_i) &=\mathbf{M}^ig(x) - g_0\for i\in\{1,\dotsb,n\}\\g_{ij}(x_i,x_j) &= \mathbf{M}^{ij}g(x) - g_i(x_i)-g_j(x_j)-g_0\for i<j\in\{1,\dotsb,n\}\\\vdots&\\g_{1\dotsb n}(x_1,\dotsb,x_n)&=g(x) - g_0-\sum_ig_i(x_i)-\sum_{i<j}g_{ij}(x_i,x_j)-\dotsb\end{align*} 

For the classifier $g$, we put $F=\{1,\dotsb,m\}$ and introduce additional mapping $h:F\mapsto\{0,1\}^m$ defined as \begin{equation}\label{eq:h}h\circ g(x) = (\ind{\{1\}},\dotsb,\ind{\{m\}})\circ g(x)\for x\in E\end{equation} Then the $m$ indicator dimensions of $h\circ g=(h_1\circ g,\dotsb,h_m\circ g)=(l_1,\dotsb,l_m)$ each belong to $L^2(E,\mathcal{E},\nu)$ and may be furnished with HDMRs.

HDMR analysis conveys a \emph{global sensitivity analysis} through the definition of structural and correlative sensitivity indices $\{(\amsbb{S}_u^a,\amsbb{S}_u^b): u\subseteq\{1,\dotsb,n\}\}$ \citep{li12,sobol93,sobol01}.

\begin{de}[Sensitivity indices of $g$] The structural sensitivity index of subspace $u\subseteq\{1,\dotsb,n\}$ is \begin{equation}\amsbb{S}_u^a = \frac{\Var g_u}{\Var g}\for u\subseteq\{1,\dotsb,n\}, |u|\ge1\end{equation} and the correlative sensitivity index is \begin{equation}\amsbb{S}_u^b = \sum_{v\subseteq\{1,\dotsb,n\}:v\ne u}\frac{\Cov(g_u,g_v)}{\Var g}\for u\subseteq\{1,\dotsb,n\}, |u|\ge1\end{equation} so that \begin{equation}1 = \sum_{u\subseteq\{1,\dotsb,n\}:|u|\ge1}(\amsbb{S}_u^a+\amsbb{S}_u^b) = \amsbb{S}^a+\amsbb{S}^b\end{equation}
\end{de}


\section{Select mixed binomial processes} 

To illustrate some of the expressive power of the mixed binomial process, in Section~\ref{sec:poisson} we give canonical examples of how the mixed binomial process may be used to construct fundamental random quantities. We give short and simple proofs. These are listed below in Table~\ref{tab:00}, grouped by the counting distribution.

\begin{table}[h!]
\begin{center}
\begin{tabular}{p{3in}ll}
\toprule
Object & $\kappa$ & Location \\\midrule
Binomial distribution & Dirac & Proposition~\ref{prop:bin}\\\midrule
Zero inflation & Binomial & Proposition~\ref{prop:inflate}\\\midrule
Watanabe's martingale theorem & Poisson & Proposition~\ref{prop:watanabe}\\
Additive (completely random) random measure & Poisson & Proposition~\ref{prop:additive}\\
Compound Poisson distribution & Poisson  & Remark~\ref{re:compound}\\
Negative binomial distribution & Poisson & Proposition~\ref{prop:negbin}\\
Gamma distribution & Poisson  & Proposition~\ref{prop:gamma} \\
Wiener process & Poisson & Proposition~\ref{prop:wiener}\\
Gaussian distribution & Poisson & Remark~\ref{re:gaussian}\\\midrule
Bayesian random measures & General & Remark~\ref{re:bayesian}\\
Mean-squared-error with Gaussian noise  & General & Proposition~\ref{prop:denoise}\\
Cluster process & General & Proposition~\ref{prop:cluster}\\\midrule
Poisson-type random measures & Poisson-type & Proposition~\ref{thm:ptrms}\\\bottomrule
\end{tabular}\caption{Mixed binomial process examples}\label{tab:00}
\end{center}
\end{table}

\subsection{Select examples}\label{sec:poisson}  

We start with the simplest mixed binomial process, the binomial process, and use this to construct the binomial distribution as a sum of Bernoulli random variables. 

\begin{prop}[Binomial]\label{prop:bin}Let $N=(\kappa,\nu)$ be a Dirac random measure on $(E,\mathcal{E})=(\{0,1\},2^{\{0,1\}})$ where $\kappa=\text{Dirac}(n)$ for $n\in\N_{\ge0}$, put $\nu=\text{Bernoulli}(p)$ for $p\in[0,1]$ and take $f(x)=x$. Then \[Nf\stackrel{D}{=}\text{Binomial}(n,p)\]  
\end{prop}
\begin{proof} The Laplace functional of $N$ is\[\E e^{-Nf}= (\nu e^{-f})^n\for f\in\mathcal{E}_+\] The result follows from \[\nu e^{-\alpha f} = 1-p+p\beta\for\alpha\in\R_+, \quad\beta=e^{-\alpha}\in[0,1]\] so that \[F(\alpha) = L(\alpha f) = (1-p+p\beta)^n\for\alpha\in\R_+,\quad \beta=e^{-\alpha}\in[0,1]\] which is the pgf of the binomial distribution.
\end{proof}
%

The following proposition provides an interpretation of the binomial random measure as a sum of zero-inflated non-negative random variables. 

\begin{prop}[Zero inflation]\label{prop:inflate} Let $N=(\kappa,\nu)$ be a binomial random measure on $(E,\mathcal{E})$ formed by $\mathbf{X}$, where $\kappa=\text{Binomial}(n,p)$. Then for $f\in\mathcal{E}_+$ the random variable $Nf$ is the sum of $n$ iid zero-inflated random variables $\mathbf{Z}=\{Z_i\}$, each with $\P(Z_i=0)=1-p$ and $\P(Z_i=f(X_i))=p$,  that is, \[Nf = \sum_{i}^K f(X_i) = \sum_{i}^n Z_i\] 
  \end{prop}
\begin{proof}
$N$ has Laplace functional \[L(f)=\E e^{-Nf} = (1-p+p\nu e^{-f})^n\for f\in\mathcal{E}_+\] and for $f\in\mathcal{E}_+$ Laplace transform \[F(\alpha) = L(\alpha f) = (1-p + p \nu e^{-\alpha f})^n\for \alpha\in\R_+\] The inner quantity can be written \[1-p + p \nu e^{-\alpha f} = (1-p)\nu e^{\alpha 0} + p \nu e^{-\alpha f}=(1-p)F_0(\alpha)+p F_f(\alpha)\] which is the mixture of a random variable taking value zero with probability $1-p$ with Laplace transform $F_0$ and a non-negative random variable with Laplace transform $F_f$ with probability $p$, so $Nf$ is the sum of $n$ iid zero-inflated random variables as claimed. 
\end{proof} 

We give a connection to the theory of martingales with Watanabe's theorem. 

\begin{prop}[Watanabe's Martingale Theorem]\label{prop:watanabe} Let $M=(\kappa,\nu)$ be a random measure on $(E,\mathcal{E})=(\R_+,\mathcal{B}_{\R_+})$ and $\kappa$ with mean $c\in(0,\infty)$. Let $N=(N_t)_{t\in\R_+}$ be the counting process where $N_t = M([0,t])$ with generated filtration $\mathscr{F}=(\mathscr{F}_t)_{t\in\R_+}$, where $\mathscr{F}_t = \sigma\{M(A): A\in\mathcal{B}_{[0,t]}\}$. Then $N$ is a Poisson process with rate function $c_\bullet=c\nu([0,\bullet])$ with respect (adapted) to $\mathscr{F}$ iff \[N_t-c_t\for t\in\R_+\] is an $\mathscr{F}$-martingale.
\end{prop}

We construct additive (completely random) measures.

\begin{samepage}
 \begin{prop}[Additive random measures]\label{prop:additive} Let $N=(\kappa,\nu)$ be a Poisson random measure on $(E\times\R_+,\mathcal{E}\otimes\mathcal{B}_{\R_+})$. Then \[L(A) = \int_{A\times\R_+}N(\D x,\D z)z\for A\in\mathcal{E}\] defines an additive random measure $L$ on $(E,\mathcal{E})$, where $L(A),\dotsb,L(B)$ are independent for all choices of finitely many disjoint $A,\dotsb,B$ in $\mathcal{E}$. The law of $L(A)$ is defined through the Laplace transform for $A$ in $\mathcal{E}$ as \[\E e^{-\alpha L(A)} = \exp_-{c\int_{A\times\R_+}\nu(\D x, \D z)(1-e^{-\alpha z})}\for \alpha\in\R_+\]
 \end{prop}
 \end{samepage}
 \begin{proof} By Fubini's theorem, $L$ is a random measure. $L(A)$ is determined from the trace of $N$ on $A\times\R_+$. The independence follows from the splitting property of Poisson, so $L$ is additive. The Laplace transform follows from the Laplace functional, noting that $\alpha L(A) = Nf$ for $f(x,z)=\alpha z\ind{A}(x)$. Then we have \[\nu e^{-f} = \int_{E\times\R_+}\nu(\D x,\D z)e^{-\alpha z\ind{A}(x)} = \int_{A\times\R_+}\nu(\D x,\D z)e^{-\alpha z}+\int_{A^c\times\R_+}\nu(\D x,\D z) \] so that \[\nu e^{-f}-1 = \nu(e^{-f}-1) =\int_{A\times\R_+}\nu(\D x,\D z)(e^{-\alpha z}-1)\] \end{proof}
  
   \begin{re}[Integer-valued] Suppose Poisson $N=(\kappa,\nu)$ is on $(E\times\N_{\ge0},\mathcal{E}\otimes 2^{\N_{\ge0}})$. We define $\beta L(A) = Nf$ for $f(x,z)=\beta z\ind{A}(x)$ where $\beta=-\log(\alpha), \alpha\in[0,1]$. Then \[\E e^{-\beta L(A)} = \exp_-{c\int_{A\times\N_{\ge0}}\nu(\D x, \D z)(1-\alpha^z)}\for \alpha\in[0,1], \quad \beta=-\log(\alpha)\in\R_+\]
 \end{re}
 
 These results allow us to define the compound Poisson random measure.  
 \begin{samepage}
\begin{re}[Compound Poisson]\label{re:compound} Let $N$ be Poisson and consider $\nu(\D x,\D z)=\mu(\D x)\eta(\D z)$. Then $F(\alpha)=\int_{\R_+}\eta(\D z)e^{-\alpha z}$ is the Laplace transform of $\eta$ and \[\int_{A\times\R_+}\nu(\D x,\D z)(1-e^{-\alpha z})=\mu(A)\left(1-F(\alpha)\right)\] so the Laplace transform of $L(A)$ is given by \[\E e^{-\alpha L(A)} = \exp_-c\mu(A)(1-F(\alpha))\for \alpha\in\R_+\] which is the compound Poisson distribution with mean $c\mu(A)$ and independent marking distribution $\eta$ with Laplace transform $F$. $L$ is additive. Similarly, for $\eta$ on $(\N_{\ge0},2^{\N_{\ge0}})$, we have \[\E e^{-\beta L(A)} = \exp_-c\mu(A)(1-\phi(\alpha))\for \alpha\in[0,1], \quad\beta=-\log(\alpha)\]  where $\phi(\alpha)=\sum_{z\in\N_{\ge0}}\eta\{z\}\alpha^z$ is the probability generating function of $\eta$.
\end{re}
\end{samepage}

\begin{prop}[Negative binomial]\label{prop:negbin} Consider the set-up of Remark~\ref{re:compound}. Let $\kappa=\text{Poisson}(r\log(\frac{1}{1-p}))$ for $r\in(0,\infty)$ and $p\in(0,1)$, take $\eta=\text{Logarithmic}(p)$ on $(\N_{>0},2^{\N_{>0}})$ with pgf \[\phi(t) = \frac{\log(1-pt)}{\log(1-p)},\] and put $f(x,y)=\ind{A}(x)y$ for $A\in\mathcal{E}$. Then \[L(A)=Nf\stackrel{D}{=}Z\sim\text{Negative-binomial}(r\mu(A),p)\] \end{prop}
\begin{proof}We have \[\E e^{-\beta L(A)} =  \exp_-c\mu(A)(1-\phi(\alpha)) = \exp_-r\log(\frac{1}{1-p})\mu(A)(1-\frac{\log(1-p\alpha)}{\log(1-p)}) = \left(\frac{1-p}{1-p\alpha}\right)^{r\mu(A)}\]\end{proof} 
 
 Next we construct the gamma distribution. 
  \begin{samepage}
 \begin{prop}[Gamma]\label{prop:gamma} Let $N=(\kappa,\nu)$ be a Poisson random measure on $([0,t]\times(0,1],\mathcal{B}_{[0,t]}\otimes\mathcal{B}_{(0,1]})$ formed by independency $(\mathbf{X},\mathbf{Y})$. Let $\nu=\text{Uniform}[0,t]\times\text{Uniform}(0,1]$ and put $f(x,y)=\frac{1}{\lambda}e^{-x}\log 1/y$ for $\lambda\in(0,\infty)$ fixed. Then the random variable $Nf$ formed as \[Nf = \sum_i^K \frac{1}{\lambda}e^{-X_i}\log 1/Y_i\] has Laplace transform \[F(\alpha) = L(\alpha f) = \E e^{-N(\alpha f)}= \left(\frac{\lambda + e^{-t}\alpha}{\lambda+\alpha}\right)^{\frac{c}{t}}\xrightarrow[c/t\rightarrow d\in(0,\infty)]{c,t\rightarrow\infty}\left(\frac{\lambda}{\lambda+\alpha}\right)^d\for \alpha\in\R_+\] where the limit is the Laplace transform of a gamma random variable with shape parameter $d$ and rate parameter $\lambda$, that is \[\lim_{\substack{c,t\rightarrow\infty\\c/t\rightarrow d}}Nf\stackrel{D}{=}Z\sim\text{Gamma}(d,\lambda)\] 
\end{prop}\end{samepage}
\begin{proof} The result follows from noting that \[\nu e^{-\alpha f} = \frac{1}{t}\log \left(\frac{\lambda+e^{-t}\alpha}{\lambda+\alpha}\right)+1\] and that the uniform random variable $U\sim\text{Uniform}(0,1]$ random variable forms an exponential random variable $W\sim\text{Exponential}(\lambda)$ as $W=\frac{1}{\lambda}\log 1/U$.
\end{proof}
This shows that the gamma distribution can be interpreted as the limit of a Poisson random measure involving uniform random variables. This idea also shows up in our next example, where we construct the Wiener process.

\begin{samepage}
\begin{prop}[Wiener process]\label{prop:wiener} Let $N=(\kappa,\nu)$ be a Poisson random measure on $([0,T],\mathcal{B}_{[0,T]})$ with $\nu=\text{Uniform}[0,T]$. Define \[N_t = \frac{N([0,t]) - tc/T}{\sqrt{c/T}}\for t\in[0,T]\] Then \begin{enumerate}\item[(i)] $N_t$ is a L\'{e}vy process \item[(ii)] $N_t$ has limit \[\lim_{c\rightarrow\infty}N_t\stackrel{D}{=}W_t\for t\in[0,T]\] where $W_t$ is the Wiener process \item[(iii)] $N_t$ has mean, variance, and covariance \begin{align*}\E N_t &=0\\\Var N_t &=t\\\Cov(N_s,N_t) &= s\wedge t\end{align*} \end{enumerate} 
\end{prop}
\end{samepage}
\begin{proof} (1) (a) Note that $N_t$ is right-continuous left-limited, inherited from $N([0,t])$, starting with $N_0=N([0,0])=0$. (b) Next note that \[N_{t+u}-N_t = \frac{(N([0,t+u])-N([0,t]) + c((t+u) - t)/T}{\sqrt{c/T}}=\frac{N((t,t+u]) - cu/T}{\sqrt{c/T}}\]  By the splitting property of Poisson, $N((t,t+u])$ and is independent of $N((0,t])$ and hence $N_{t+u}-N_t$ is independent of $N_t$. (c) Moreover, $N((t,t+u])\sim\text{Poisson}(cu/T)$ has the same distribution as $N((0,u])$ for all $t\in\R_+$, so $N_{t+u}-N_t$ has the same distribution as $N_u$. By (a), (b), (c), $N_t$ is a L\'{e}vy process. 

(2)  Because $N_t$ is a L\'{e}vy process, so is $W_t$. The law of $W_t$ is therefore specified through the characteristic function. Hence it suffices to show the claimed limit through characteristic functions. The characteristic function of $aN([0,t])=\frac{1}{\sqrt{c/T}}N([0,t])$ is \[\varphi_{aN([0,t])}(r)=\sum_{k\ge0}e^{irka}(ct/T)^ke^{-ct/T}/k! = \exp{\frac{ct}{T}(e^{ira}-1)}=\varphi_{N([0,t])}(ar)\for r\in\R\] The characteristic function of $N_t$ is given by \[\varphi_{N_t}(r)=\varphi_{aN([0,t])}(r)\varphi_{\delta_{-t\sqrt{c/T}}}(r) = \exp\left({\frac{ct}{T}(e^{ir/\sqrt{c/T}}-1)-irt/\sqrt{c/T}}\right)\for r\in\R\] Then we have that \[\lim_{c\rightarrow\infty}\varphi_{N_t}(r) = e^{-r^2t/2} = \varphi_{W_t}(r)\for r\in\R\] which is the characteristic function of a L\'{e}vy process where $W_t$ is a Gaussian random variable with mean zero and variance $t$, i.e. the Wiener process.

(3) The mean, variance, and covariance follow from rescaling  of the Poisson random measure by the factor $1/\sqrt{c/T}$, where $\E N([0,t])= ct/T$ and $\Cov(N_s,N_t)=\frac{T}{c}c\nu([0,s]\cap[0,t])=s\wedge t$.

\end{proof}
Proposition~\ref{prop:wiener} allows us to construct generic one-dimensional Gaussian random variables.
  
\begin{re}[Gaussian random variables]\label{re:gaussian} By Proposition~\ref{prop:wiener} the Poisson distribution may be used to construct $\text{Gaussian}(\mu,\sigma^2)$ random variable $Z$: For $c\in(0,\infty)$, take $K\sim\text{Poisson}(c\sigma^2)$ and put $N = \mu + (K-c\sigma^2)/\sqrt{c}$. Then $\lim_{c\rightarrow\infty}N\stackrel{D}{=}Z\sim\text{Gaussian}(\mu,\sigma^2)$. 
\end{re}

%
We give a brief remark on Bayesian random measures. 
\begin{re}[Bayesian random measures]\label{re:bayesian} Let $N=(\kappa,\nu)$ be a random measure on parameter space $(E,\mathcal{E})$ where $\nu$ is a prior distribution, i.e. $\nu(x)=\P(x)$. We call $N$ a \emph{prior} random measure. Let $f\in\mathcal{E}_+$ be a measurable function defined as \[f(x) = L(x|\mathfrak{X})\] where $L$ is a \emph{likelihood} function and $\mathfrak{X}$ is some data.  We have \emph{evidence} \[\P(\mathfrak{X}) = \nu f = \int_E\nu(\D x)L(x|\mathfrak{X})=\frac{1}{c}\E Nf\] The image random measure $M=N\circ f^{-1}=(\kappa,\eta=\nu\circ f^{-1})$ is the \emph{posterior} random measure with density \[\eta(\D x) = (\nu\circ f^{-1})(\D x)= \P(x\in \D x|\mathfrak{X}) =\frac{L(x|\mathfrak{X})\nu(\D x)}{\P(\mathfrak{X})}\] 
\end{re}

Another idea is experimental uncertainty. We give a proposition that shows how to attain the denoised mean-squared-error when the contaminating noise is independent additive Gaussian.  

\begin{samepage}
\begin{prop}[Denoising]\label{prop:denoise}Let $(\mathbf{X},\mathbf{Y})$ an independency distributed according to $\nu\times Q$ on $(E\times\R,\mathcal{E}\otimes\mathcal{B}_\R)$, and let $\mathbf{Z}$ be independent $\text{Gaussian}(0,\sigma^2)$ noise on $(\R,\mathcal{B}_{\R})$, where $\sigma^2\in\R_+$. The triple $(\mathbf{X},\mathbf{Y},\mathbf{Z})$ forms a random measure $N=(\kappa,\mu=\nu\times Q\times\xi)$ on $(E\times\R\times\R,\mathcal{E}\otimes\mathcal{B}_\R\otimes\mathcal{B}_\R)$. Consider $f(x,y,z) = (g(x) +z - y)^2$ where $g\in L^2(E,\mathcal{E},\nu)$ and define $h(x,y)=(g(x)-y)^2$. Then the random variable $Nf$ formed as \[Nf = \sum_i^K(g(X_i)+Z_i-Y_i)^2\] has mean and variance \begin{align*}\E N f &= c(\sigma^2+(\nu\times Q)h)\\\Var Nf &=  c\left(3\sigma^4+(\nu\times Q)\left(6\sigma^2h+h^2\right)\right) + (\delta^2-c)(\sigma^2+(\nu\times Q)h)^2\end{align*} and Laplace transform \[F(\alpha) =  \psi(\mu e^{-\alpha f}) = \psi(F_f(\alpha))\for\alpha\in\R_+\] where $F_f$ is the Laplace transform of $f$ given by \[F_f(\alpha)=\frac{1}{\sqrt{1+2\sigma^2\alpha}}(\nu\times Q)e^{-\alpha h/(1+2\sigma^2\alpha)}=\frac{1}{\sqrt{1+2\sigma^2\alpha}}F_h(\frac{\alpha}{1+2\sigma^2\alpha})\for\alpha\in\R_+\] and $F_h$ is the Laplace transform of $h$ given by \[F_h(\alpha)=(\nu\times Q)e^{-\alpha h}\for \alpha\in\R_+\]
\end{prop}
\end{samepage}
\begin{proof} The result follows from \begin{align*}\mu f &= \sigma^2+(\nu\times Q)h\\\mu f^2 &=3\sigma^4+(\nu\times Q)\left(6\sigma^2h+h^2\right)\end{align*} and  \[\xi e^{-\alpha f} = \frac{1}{\sqrt{1+2\sigma^2\alpha}}e^{-\alpha h/(1+2\sigma^2\alpha)}\] so that \[\mu e^{-\alpha f} = (\nu\times Q)(\xi e^{-\alpha f}) = \frac{1}{\sqrt{1+2\sigma^2\alpha}}(\nu\times Q)e^{-\widetilde{\alpha} h} = \frac{1}{\sqrt{1+2\sigma^2\alpha}} F_h(\widetilde{\alpha})\] where $\widetilde{\alpha} = \alpha/(1+2\sigma^2\alpha)\in[0,\frac{1}{2\sigma^2})$ and $F_h(\alpha)=(\nu\times Q)e^{-\alpha h}$ is the Laplace transform of $h$.
\end{proof}

\begin{re}[Chi-square] Consider the set-up of Proposition~\ref{prop:denoise}. Take $\kappa=\text{Dirac}(c)$ for $c\in\N_{>0}$  assume that $\sigma^2=1$ and $h=0$. Then \[F(\alpha) = \frac{1}{(1+2\alpha)^{c/2}}\] is the Laplace transform of a Chi-square distribution with $c$ degrees of freedom.
\end{re}

Proposition~\ref{prop:denoise} shows that mean-squared-error may by de-noised for independent additive Gaussian noise. The interpretation of $f$ may be either the model containing error, i.e. $g(x)+z$, or the truth containing noise $z+y$. The result shows that the Laplace transforms of $f$ and $h$ are related for known $\sigma^2$: knowledge of $F_f$ provides knowledge of $F_h$ and vice versa. For Gaussian noise, this conveys a deconvolution of the density of $Nf$ (observed) to attain the density of $Mh$ (unobserved), where $M=(\kappa,\nu\times Q)$ is the random measure on $(E\times\R,\mathcal{E}\otimes\mathcal{B}_{\R})$. 

\begin{samepage}
\begin{re}[Example] Consider the set-up of Proposition~\ref{prop:denoise}. Consider $g(x)=x$ on $E=[0,1]$ with $\nu=\text{Uniform}[0,1]$. Let $Q(x,\cdot)=\delta_{1/2}(\cdot)$ and take $\sigma^2=1$. Put $f(x,y,z)=(x+z-y)^2$ and $h(x,y)=(x-y)^2$. Then \[F_f(\alpha) = (\nu\times Q\times \xi) e^{-\alpha f} = \sqrt{\frac{\pi}{\alpha}}\Erf\left(\frac{1}{2}\sqrt{\frac{\alpha}{1+2\alpha}}\right)\] and \[ F_h(\frac{\alpha}{1+2\alpha})=\sqrt{1+2\alpha}F_f(\alpha) =\sqrt{\frac{\pi}{\bullet}}\Erf\left(\frac{\sqrt{\bullet}}{2}\right)\circ\frac{\alpha}{1+2\alpha} \] The laws of $f$ and $h$ can be determined by the inverse Laplace transforms, which give \[(\mu\circ f^{-1})(\D x) = \D x\frac{1}{2 \sqrt{x}}\left(\Erf\left(\frac{1+2 \sqrt{x}}{2 \sqrt{2}}\right)+\Erf\left(\frac{1-2 \sqrt{x}}{2 \sqrt{2}}\right)\right)\for x\in\R_+\] and  \[(\nu\times Q)\circ h^{-1}(\D x) = \D x\frac{1}{\sqrt{x}}\for x\in(0,1/4)\] 
\end{re}
\end{samepage}

Another common usage of the mixed binomial process is construction of cluster processes. This is a nested mixed binomial process, where the mark-space forms a mixed binomial process. We describe the general case below. 

\begin{samepage}
\begin{prop}[Cluster process]\label{prop:cluster} Let $N=(\kappa,\nu)$ be a random measure on $(E,\mathcal{E})$ where $\kappa$ has mean $c$ and variance $\delta^2$. Consider transition probability kernel $Q$ from $(E,\mathcal{E})$ into the space of random counting measures $(F,\mathcal{F})$, defined as \[Q(x,\D y) =\P(M_x\in \D y)\] where $M_x=(\xi_x,\zeta_x)$ is a random counting measure on $(E,\mathcal{E})$, $\xi_x$ is a $\N_{\ge0}$-valued distribution conditioned on $x$ with mean $c_x$ and variance $\delta_x^2$, and $\zeta_x$ is a conditional distribution on $(E,\mathcal{E})$. The random measure $M=(\kappa,\nu\times Q)$ on $(E\times F,\mathcal{E}\otimes\mathcal{F})$ is called a cluster process and is formed as \[Mf = \sum_i^K f\circ (X_i,M_{X_i})\for f\in(\mathcal{E}\otimes\mathcal{F})_+\] Consider $f(x,y)=y(E)$. Then $Mf$ is the total number of points of the cluster process with mean and variance \begin{align*}\E Mf &=c\nu c_\bullet\\\Var Mf &=c\nu(c_\bullet^2+\delta^2_\bullet)+(\delta^2-c)(\nu c_\bullet)^2\end{align*} and Laplace transform \[F(\alpha) = \psi(\nu\phi_\bullet(\beta))\for \alpha\in\R_+,\, \beta = e^{-\alpha}\in[0,1]\] where $\psi$ is the pgf of $\kappa$ and $\phi_x$ is the pgf of $\xi_x$. Now consider $f(x,y)=y(A)$ for $A\subset E$, where $Mf$ is the number of points of the cluster process located in $A$. Then $Mf$ has mean and variance \begin{align*}\E Mf &=  c\nu(c_\bullet a_\bullet)\\\Var Mf &= c\nu((c_\bullet^2+\delta_\bullet^2)a_\bullet)+(\delta^2-c)\nu^2(c_\bullet a_\bullet)\end{align*} where $a_x=\zeta_x(A)$ and Laplace transform \[F(\alpha) =  \psi(\nu\phi_\bullet^A(\beta))\for \alpha\in\R_+,\, \beta = e^{-\alpha}\in[0,1]\] where $\phi_x^A(t)=\phi_x(a_xt+1-a_x)$ is the restricted pgf of $\xi_x$ to $A\subset E$.
\end{prop}
\end{samepage}

\begin{proof} For $f(x,y)=y(E)$, we have \[\E Mf = c(\nu\times Q)f = c\int_E\nu(\D x)\E M_x(E) = c\nu c_\bullet\] and \begin{align*}\Var Mf &= c(\nu\times Q)f^2 + (\delta^2-c)((\nu\times Q)f)^2 \\&=c\int_E\nu(\D x)\E M_x^2(E) +  (\delta^2-c)(\nu c_\bullet)^2\\&=c\nu(\E^2 M_x(E) + \Var M_x(E)) + (\delta^2-c)(\nu c_\bullet)^2 \\ &= c\nu(c^2_\bullet+\delta_\bullet^2) + (\delta^2-c)(\nu c_\bullet)^2\end{align*} The Laplace transform follows from \[F(\alpha)=\psi((\nu\times Q) e^{-\alpha f}) = \psi(\int_E\nu(\D x)\E e^{-\alpha M_x(E)}) = \psi(\nu\phi_\bullet(\beta))\] where $\beta=e^{-\alpha}$.

For $f(x,y)=y(A)$, we have \[\E Mf = c(\nu\times Q)f = c\int_E\nu(\D x)\E M_x(A) = c\nu(c_\bullet a_\bullet)\] and \begin{align*}\Var Mf &= c(\nu\times Q)f^2 + (\delta^2-c)((\nu\times Q)f)^2\\&=c\nu(\Var M_x(A)+\E^2 M_x(A))+ (\delta^2-c)\nu^2(c_\bullet a_\bullet)\\&=c\nu(c_\bullet a_\bullet+(\delta^2_\bullet-c_\bullet)a_\bullet + c^2_\bullet a_\bullet) + (\delta^2-c)\nu^2(c_\bullet a_\bullet)\\&=c\nu((c_\bullet^2+\delta_\bullet^2)a_\bullet)+(\delta^2-c)\nu^2(c_\bullet a_\bullet)\end{align*} The Laplace transform follows from  \[F(\alpha) = \psi((\nu\times Q)e^{-\alpha f}) =\psi(\nu \E e^{-\alpha M_\bullet(A)}) =\psi(\nu\phi_\bullet^A(\beta))\for \alpha\in\R_+\]

\end{proof}

In Table~\ref{tab:cluster} we list some common cluster processes for choices of $\kappa$ and $\xi_\bullet$. 

\begin{table}[h!]
\begin{center}
\begin{tabular}{lll}
\toprule
Cluster process & $\kappa$ & $\xi_\bullet$ \\\midrule
Independent cluster process & general & general \\
Cox cluster process & general & Poisson \\
Poisson cluster process & Poisson & general \\
Neyman-Scott process & Poisson & Poisson \\\bottomrule
\end{tabular}\caption{Cluster processes}\label{tab:cluster}
\end{center}
\end{table}

\begin{samepage}
\begin{re}[Example] Consider the set-up of Proposition~\ref{prop:cluster}. Let $(E,\mathcal{E})=(\R^2,\mathcal{B}_{\R^2})$.  Let $\kappa$ be Poisson or Dirac. Let $\xi$ be independent of $\nu$ and Poisson. Let $\zeta_\bullet$ be Gaussian with mean $\bullet$ and variance $\sigma^2\in\{0.01^2,0.1^2\}$. Consider $\nu=\text{Uniform}[0,1]^2$. Let $\E\kappa=20$ and $\E\xi=10$. Consider $A=(-\infty,\frac{1}{2}]\times(-\infty,\frac{1}{2}]$ and note that $\nu a_\bullet = 1/4$. We sample $Mf$ for the two $f$'s for Poisson and Dirac. These correspond to the Neyman-Scott and Cox cluster processes respectively. In Table~\ref{tab:cluster} we show the means and variances for the analytic and estimated from $10^5$ random realizations. The estimated values are very close to the true values. In Figure~\ref{fig:cluster} we show the distribution of $Mf$ for the two $f$'s as estimated for $\sigma^2\in\{0.01^2,0.1^2\}$. Dirac possesses `spikes' on $A$ for $\sigma^2=0.01^2$, corresponding to the low variance of the clusters (the spikes occurring at integers at roughly multiples of $\E\xi$), whereas Poisson is invariant to the choice of $\sigma^2$, performing smoothening. 

\begin{table}[h!]
\begin{center}
\begin{tabular}{c|cccc|cccc}
\cline{2-9}
&\multicolumn{4}{c}{$f(x,y)=y(E)$}&\multicolumn{4}{c}{$f(x,y)=y(A)$}\\\cline{2-9}
&$\E Mf $ & $\widetilde{\E}Mf$ & $\Var Mf $ & $\widetilde{\Var}Mf$ & $\E Mf $ & $\widetilde{\E}Mf$ & $\Var Mf $ & $\widetilde{\Var}Mf$ \\\midrule
Poisson &200&200.1&2\,200&2\,202.3&50&50.0&550&535.1\\
Dirac &200&200.3&2\,000&2\,007.4&50&50.2&500&495.6\\
\bottomrule
\end{tabular}\caption{Cluster random measure simulation results from $10^5$ samples for $\sigma^2=0.01^2$}\label{tab:cluster}
\end{center}
\end{table}

\begin{figure}[h!]
\centering
\begingroup
\captionsetup[subfigure]{width=3in,font=normalsize}
\subfloat[$\sigma^2=0.01^2$, $f(x,y)=y(A)$\label{fig:clusterA}]{\includegraphics[width=3.5in]{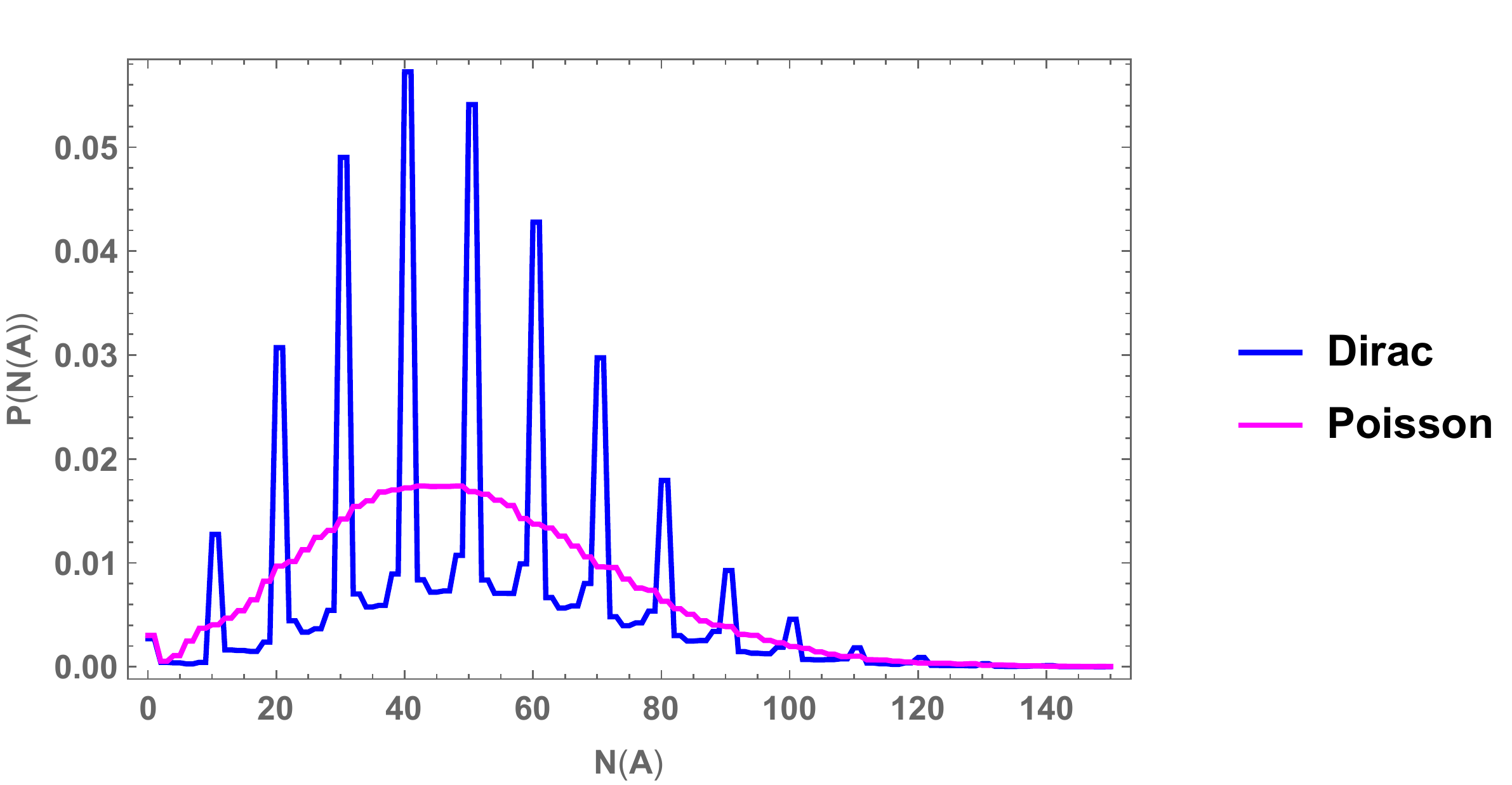}}
\subfloat[$\sigma^2=0.01^2$, $f(x,y)=y(E)$\label{fig:clusterB}]{\includegraphics[width=3.5in]{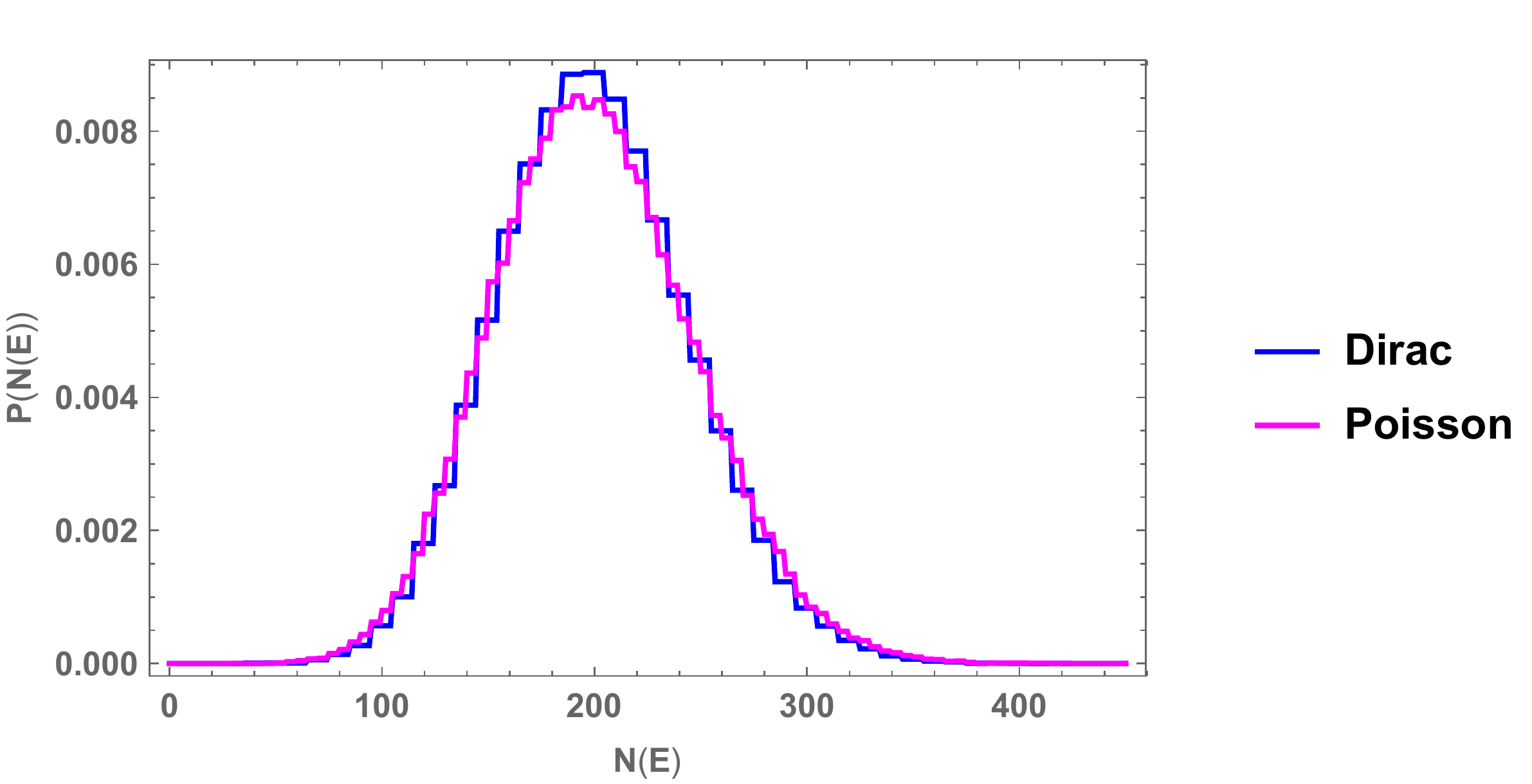}}\\
\subfloat[$\sigma^2=0.1^2$, $f(x,y)=y(A)$\label{fig:clusterA1}]{\includegraphics[width=3.5in]{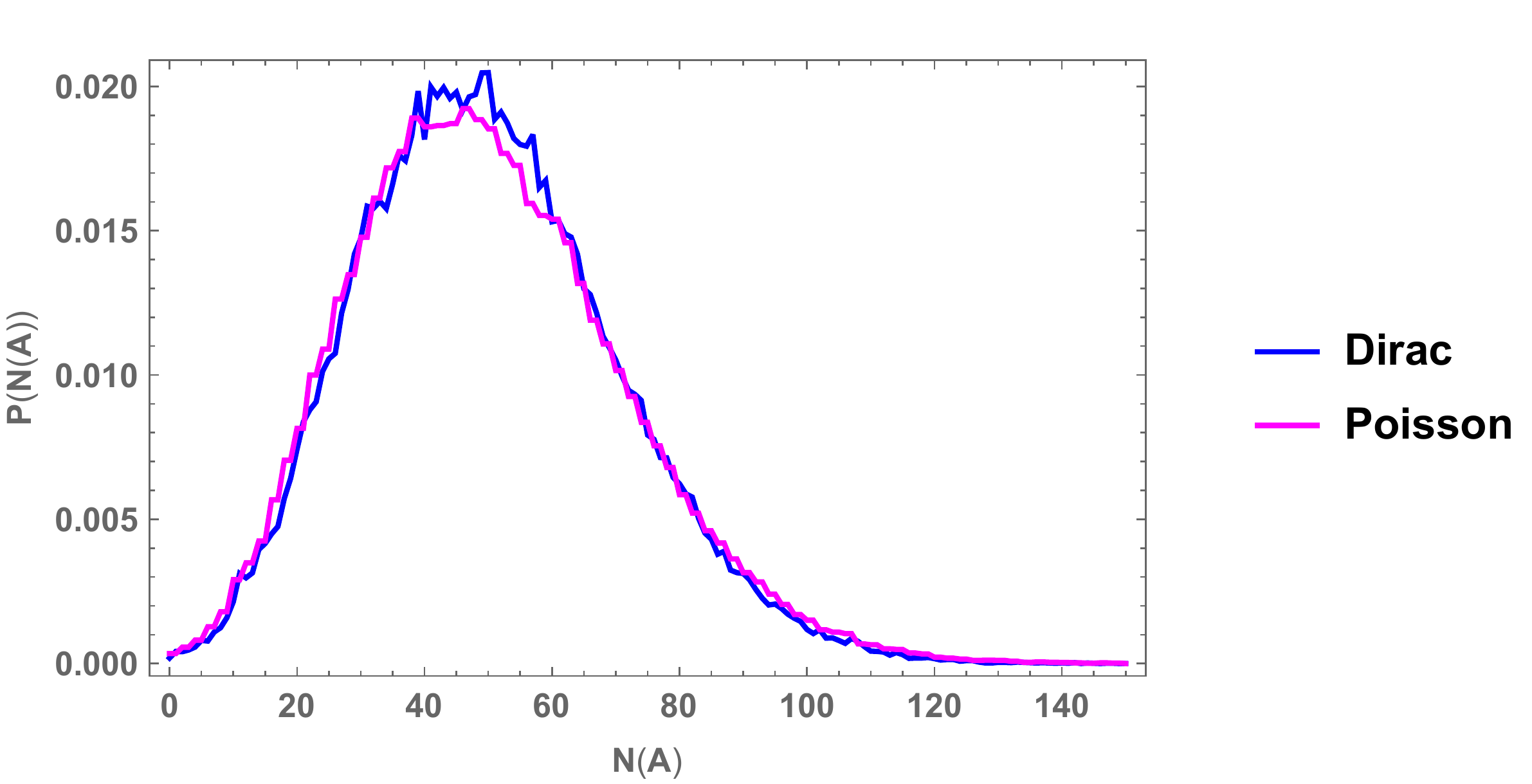}}
\subfloat[$\sigma^2=0.1^2$, $f(x,y)=y(E)$\label{fig:clusterB1}]{\includegraphics[width=3.5in]{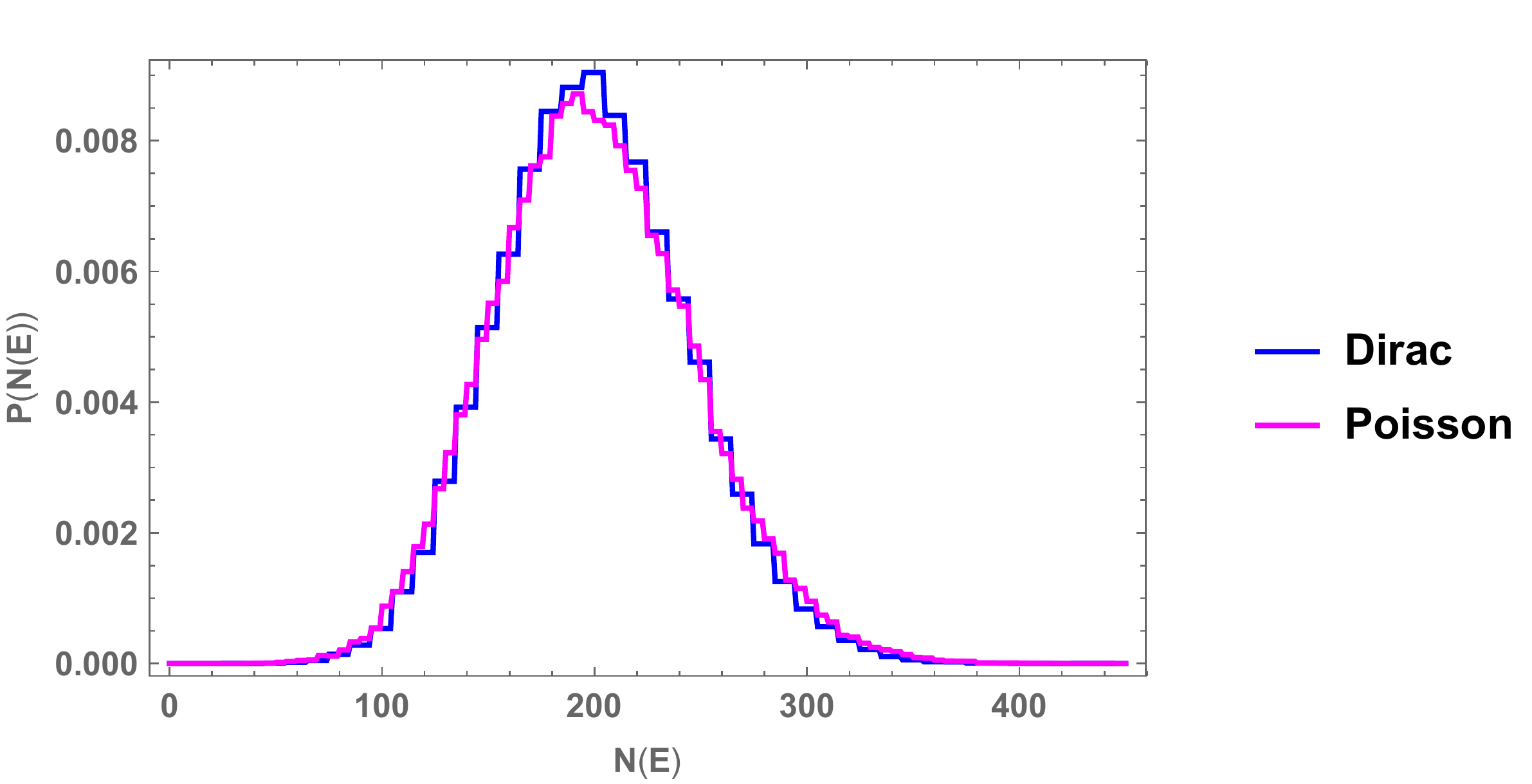}}\\
\endgroup
\caption{Simulated distribution of cluster process $Mf$}\label{fig:cluster}
\end{figure}

\FloatBarrier
 
\end{re}
\end{samepage}

\subsection{Poisson-type random measures} 

The class of mixed binomial processes where $\kappa$ is Poisson, negative binomial or binomial are called \emph{Poisson-type} (PT) random measures, with corresponding PT members $\kappa$. The Poisson-type (PT) family distributions are \[K\sim\kappa=\begin{cases}\text{Binomial}(n,p)\in\{0,\dotsc,n\}&\text{for}\quad n\in\N_{>0}, p\in(0,1)\\\text{Poisson}(c)\in\N_{\ge0}&\text{for}\quad c\in(0,\infty)\\\text{Negative-binomial}(r,p)\in\N_{\ge0}&\text{for}\quad r\in\N_{>0}, p\in(0,1)\end{cases}\] The degenerate member of the PT family is the Dirac $\delta_c=\lim_{\substack{p\rightarrow1\\n\rightarrow c}}\text{Binomial}(n,p)$. The following result establishes the existence and uniqueness of Poisson-type (PT) random measures as random counting measures closed under restriction to (thinning by) subspaces.

\begin{prop}[Existence and uniqueness of PT random measures]\label{thm:ptrms} Assume that $K\sim \kappa_\theta$ where pgf $\psi_\theta$ belongs to the canonical non-negative power series family of distributions with parameter $\theta\in\Theta$ and $\{0,1\}\subset \supp(K)$. Consider the random measure $N=(\kappa_\theta,\nu)$ on the  space $(E, \cal E)$ and assume that $\nu$ is diffuse. Then for any $A\subset E$ with $\nu(A)=a>0$ there exists a mapping $h_a:\Theta\rightarrow \Theta$ such that the restricted random measure is $N_A=(\kappa_{h_a(\theta)}, \nu_A)$, that is, \begin{equation}\label{eq:4a}\E e^{-N_Af}=\psi_{h_a(\theta)}(\nu_A e^{-f})\quad \text{for}\quad f\in {\cal E}_+\end{equation}  iff $K$ is Poisson, negative binomial or binomial. 
\end{prop}
\begin{proof}Theorem 3 \citep{bastian}\end{proof}

In Table~\ref{tab:can} we give Poisson-type distributions and their canonical parameters, pgfs, and ``bone'' mappings.

\begin{table}
\begin{center}
\begin{tabular}{l c c c }
\toprule
Name &Parameter $\th$ & $\psi_\th(t)$& $h_a(\theta)$ \\\midrule
Poisson & $\la$ & $\exp[{\th(t-1)}]$ &$a\th$ \\
Bernoulli & $p/(1-p)$ & $(1+\th t)/(1+\th)$& $a\th/(1+(1-a)\th)$\\
Geometric & $p$ & $(1-\th)/(1-t\th)$& $a\th/(1-(1-a)\th)$\\\bottomrule
\end{tabular}\caption{Poisson-type distributions with corresponding pgfs and mappings of their canonical parameters}\label{tab:can}
\end{center}
\end{table}

\FloatBarrier

\begin{re}[Closure for atomic] The sufficiency part of Proposition~\ref{thm:ptrms} holds for atomic $\nu$, whereas the necessity part (uniqueness) fails. This means that the PT random measures are closed under restriction for all measures.
\end{re}

\begin{re}[Bone mapping] The mapping $h_a$ is a called a ``bone'' mapping, for which the pgf satisfies \[\psi_\theta(at+1-a) = \psi_{h_a(\theta)}(t)\] 
\end{re}

\begin{re}[Consistency at limit] For the PT family, \[\lim_{c\rightarrow\infty}\P(K>k)=1\quad\text{for all}\quad k\] (sending $r,n\rightarrow\infty$ for negative binomial and binomial respectively). Hence $K=+\infty$ almost surely for $c=+\infty$ for the PT family. 
\end{re}

The rescalings for $A\subseteq E$ with $a=\nu(A)\in(0,1]$ are given by \[N(A)=N\ind{A}=K_A\sim\begin{cases}\text{Binomial}(n,ap) &\text{if }K\sim\text{Binomial}(n,p)\\\text{Poisson}(ac)&\text{if }K\sim\text{Poisson}(c)\\\text{NegativeBinomial}(r,ap/(1-(1-a)p))&\text{if }K\sim\text{Negative-binomial}(r,p)\end{cases}\] The mean of $Nf$ for $f\in\mathcal{E}_+$ is \begin{align*}a&=\nu f\\\E Nf&=\begin{cases}anp &\text{if }K\sim\text{Binomial}(n,p)\\ ac&\text{if }K\sim\text{Poisson}(c)\\ arp/(1-p)&\text{if }K\sim\text{Negative-binomial}(r,p)\end{cases}\end{align*} and the variance is \begin{align*}a&=\nu f\\b&=\nu f^2\\\Var Nf&=\begin{cases}np(b-pa^2) &\text{if }K\sim\text{Binomial}(n,p)\\ cb&\text{if }K\sim\text{Poisson}(c)\\ (\frac{rp}{1-p})(b+(\frac{p}{1-p})a^2) &\text{if }K\sim\text{Negative-binomial}(r,p)\end{cases}\end{align*}

For arbitrary $f,g\in\mathcal{E}_+$, we have covariance \begin{align*}a&=\nu f\\b&=\nu g\\d&=\nu(fg)\\\C(Nf,Ng) &= \begin{cases}np(d-pab)&\text{if }K\sim\text{Binomial}(n,p)\\ cd &\text{if }K\sim\text{Poisson}(c)\\(\frac{rp}{1-p})(d+(\frac{p}{1-p})ab)&\text{if }K\sim\text{Negative-binomial}(r,p)\end{cases}\in\R\end{align*} 

\subsection{Orthogonal die random measures} 



The set of permissible supports for the discrete uniform distribution is defined as \[A\equiv\{(m,n): \text{integers }0\le m\le n\text{ except }m=n=0\}\] Consider the discrete uniform family of distributions $\mathcal{K}$ \begin{equation}\label{eq:U}\mathcal{K}=\{\kappa_{mn}=\text{Uniform}\{m,m+1\dotsc,n-1,n\}: (m,n)\in A\}\end{equation} $\kappa_{mn}$ has mean $c=(m+n)/2>0$ and variance $\delta^2=((n-m+1)^2-1)/12\ge0$. Random realizations of $\kappa_{mn}$ may be thought of in terms of rolling fair dice, where the number of sides of the dice is equal to $n-m+1$. 

\begin{de}[Orthogonal die] We call $\kappa_{mn}$ an \emph{orthogonal die} if $c=\delta^2$. In turn we call $N=(\kappa_{mn},\nu)$ an orthogonal die random measure.\end{de}


We give the existence and uniqueness of orthogonal dice. There are an infinite number of orthogonal dice, whose numbers of sides are coprime to two and three.

\begin{samepage}
\begin{prop}[Existence and uniqueness of orthogonal rectangular dice]\label{thm:1} Orthogonal dice $\kappa_{mn}$ with support $\{m,m+1\dotsc,n-1,n\}$ for integers $0\le m<n$ are completely enumerated by the collection \begin{equation}\label{eq:K}\mathcal{S}=\{\kappa_{mn}: (m,n)\in A\}\subset\mathcal{K}\end{equation} where \[A=\{(m,n): k=1,2,4,5,7,8,\dotsb, m=(k^2-1)/3, n=2k+m+2\}\] with $|\mathcal{S}|=\infty$ and $m/n\rightarrow1$ as $m,n\rightarrow\infty$. Moreover, for each $\kappa_{mn}\in\mathcal{S}$, the integer $n-m+1$ is a product of one or more primes each having value equal to or greater than five. \end{prop} 
\end{samepage}\begin{proof} Theorem 3 \citep{bastian2020orthogonal}\end{proof}

The first 15 orthogonal dice are shown below in Table~\ref{tab:dice}. We have the following corollary that shows that the orthogonal dice span the naturals.

\begin{cor}[Naturals]\[\bigcup_{\kappa_{mn}\in\mathcal{S}}\{m,\dotsb,n\} = \N_{\ge0}\]\end{cor}\begin{proof} The result follows from putting $m(k) = (k^2-1)/3$ and $n(k)=2k+m(k)+2$ and noting that $m(k+1)<m(k+2)<n(k)$ for all $k\ge1$, so that starting with $m(1)=0$ and $n(1)=4$, all the integers are enumerated at least once for the given sequence of $k$.
\end{proof} 

The orthogonal dice are special because, in addition to being orthogonal, they possess a Poisson-limit theorem (PLT): thinning and ranging over the support recapitulates the Poisson random measure. 

\FloatBarrier

\begin{prop}[Orthogonal die PLT]\label{thm:main2}The sequence of thinned discrete uniform random measures $(N_{m,n}^a)$ converges in distribution (converges weakly) to the Poisson random measure $N$, that is, \[\lim_{\substack{m,n\rightarrow\infty\\m/n\rightarrow1\\a\rightarrow0\\na\rightarrow c}}\E e^{-N_{m,n}^af}=\E e^{-Nf}\quad\text{for}\quad f\in\mathcal{E}_+\]
\end{prop}\begin{proof}Theorem 2 \citep{bastian2020orthogonal} \end{proof}

\begin{table}[h!]
\begin{center}
\begin{tabular}{c|cc|cc}
\toprule
$k$ &$m$ & $n$ &$c$ & $n-m+1$ \\\midrule
 1& 0 &4 &2 &   5\\
2& 1&7 & 4&   7\\
4 & 5&15 &10& 11\\ 
5& 8&20 &14& 13\\
7&16&32 &24& 17\\
8&21&39 &30& 19\\
10&33&55 &44& 23\\
11&40&64 &52& 25\\
13&56&84 &70& 29\\
14&65&95 &80& 31\\
16&85&119 &102& 35\\
17&96&132 &114 & 37\\
19&120&160 &140& 41\\
20&133&175 &154& 43\\
22&161&207 &184& 47\\
\bottomrule
\end{tabular}\caption{First 15 orthogonal dice}\label{tab:dice}
\end{center}
\end{table} 

\FloatBarrier

\section{Random measure ANOVA of risk functionals}\label{sec:anova}

In Section~\ref{sec:var} we decompose the variance \eqref{eq:var0}. In Section~\ref{sec:mean} we decompose the mean \eqref{eq:mu0}. 

\subsection{Variance decomposition (RM-ANOVA)}\label{sec:var}
  
In the following result, we have that the variance of the full model is decomposed into the variance-covariance structure of the partition, obtained through the variance-covariance formulas. Based on the class of mixed binomial processes we refer to the decomposition of random measure variance in terms of partitions the \emph{random measure analysis of variance} (RM-ANOVA).

\begin{thm}[RM-ANOVA] Let $N=(\kappa,\nu)$ be a random measure on $(E,\mathcal{E})$. Let $f\in\mathcal{E}_+$ and consider disjoint partition $\{A,\dotsb,B\}$ of $E$. Define $f_a = f\ind{A}, \dotsb, f_b = f\ind{B}$. Then
 \begin{align}\Var Nf &= \sum_{D\in\{A,\dotsb,B\}} \Var Nf_d + \sum_{D_i\ne D_j\in\{A,\dotsb,B\}} \Cov(Nf_{d_i},Nf_{d_j})\\&=\sum_{D\in\{A,\dotsb,B\}}(c\nu f_d^2+(\delta^2-c)(\nu f_d)^2) + \sum_{D_i\ne D_j\in\{A,\dotsb,B\}} (\delta^2-c)\nu f_{d_i}\nu f_{d_j}\nonumber\end{align} 
\end{thm}
\begin{proof} Noting that $f=f\ind{A}+\dotsb+f\ind{B}$ and that disjointness implies that $\Cov(Nf_a,Nf_b)=(\delta^2-c)\nu f_a\nu f_b$, we have the decomposition using the variance and covariance formulas of the mixed binomial process.
\end{proof}

If we normalize by the overall variance, then we retrieve a kind of sensitivity analysis.

\begin{de}[Sensitivity indices of $Nf$] For the random measure $N$ and relative to the disjoint partition $\{A,\dotsb,B\}$ of $E$, the structural sensitivity index of $Nf$ is defined as \begin{equation}\amsbb{S}_d^a\equiv\frac{\Var Nf_d}{\Var Nf}\for D\in\{A,\dotsb,B\}\end{equation} and correlative sensitivity index is defined as \begin{equation}\amsbb{S}^b_d\equiv\sum_{D_i\in\{A,\dotsb,B\}: D_i\ne D}\frac{\Cov(Nf_d,Nf_{d_i})}{\Var Nf}\for D\in\{A,\dotsb,B\}\end{equation} where \begin{equation}1 = \sum_{D\in\{A,\dotsb,B\}}(\amsbb{S}^a_d+\amsbb{S}^b_d)= \amsbb{S}^a+\amsbb{S}^b\end{equation}  This gives the set of sensitivity indices $\{(\amsbb{S}^a_d,\amsbb{S}^b_d): D\in\{A,\dotsb,B\}\}$. 
\end{de}

These sensitivity indices indicate the contributions to variance of the random measure in the partitions. If $c=\delta^2$ for $\kappa$, such as with Poisson \citep{Bastian:2020aa} or an orthogonal die \citep{bastian2020orthogonal}, then $\amsbb{S}^a=1$ and $\amsbb{S}^b=0$, so $(\amsbb{S}^a_d)$ is a probability vector, conveying a distribution of uncertainty on the partition. Otherwise, the vector is positively or negatively defective and hence looses a probabilistic interpretation. 


A common setting is $\kappa=\delta_n$ for some $n\in\N_{>0}$ in the binomial process $N$. For the binomial process, \begin{equation}\label{eq:binsa}\amsbb{S}^a = \sum_{D\in\{A,\dotsb,B\}}\frac{\Var Nf_d}{\Var Nf}=\sum_{D\in\{A,\dotsb,B\}}\frac{\Var f_d}{\Var f}>1\end{equation} and  \begin{equation}\label{eq:binsb}\amsbb{S}^b = \sum_{D\in\{A,\dotsb,B\}}\sum_{D_i\in\{A,\dotsb,B\}: D_i\ne D}\frac{-\nu f_d \nu f_{d_i}}{\Var f}<1\end{equation} 
 
\begin{cor}[Sensitivity probability measure]For orthogonal $N$ and $f\in\mathcal{E}_+$, the sensitivity probability measure $\amsbb{S}$ on $(E,\mathcal{E})$ is given by \begin{equation}\label{eq:sdx}\amsbb{S}(\D x) = \frac{\nu(\D x)f^2(x)}{\nu f^2}\end{equation} so that \[\amsbb{S}_d\equiv\int_D\amsbb{S}(\D x)=\amsbb{S}(D)\for D\in\mathcal{E}\] For variable subset $u\subseteq\{1,\dotsb\}$, putting $E_{-u}\equiv\bigtimes_{i\in\{1,\dotsb\}: i\notin u}E_i$, we have \begin{equation}\label{eq:sdxi}\amsbb{S}_u(\D x_u) \equiv \int_{E_{-u}}\amsbb{S}(\D x)\end{equation} 
\end{cor}

\subsection{Mean decomposition (MM-ANOVA)}\label{sec:mean} In the following result we have a decomposition of the mean measure.

\begin{thm}[MM-ANOVA]Let $N=(\kappa,\nu)$ be a random measure on $(E,\mathcal{E})$. Let $g\in L^2(E,\mathcal{E},\nu)$ and define $f=(g-\nu g)^2\in\mathcal{E}_+$. Suppose $E$ has dimension $n$. Then \begin{equation}\E Nf =c\sum_{u,v\subseteq\{1,\dotsb,n\}}\Cov(g_u,g_v)\end{equation} where the $\{g_u\}$ are the HDMR component functions of $g$. 
\end{thm}\begin{proof} Note that $\nu f=\nu(g-\E g)^2 =\Var g$. The HDMR decomposition of $g
\in L^2(E,\mathcal{E},\nu)$ into component functions $\{g_u\}$ and the associated variance decomposition finishes the proof.
\end{proof}

If $\nu$ has a product form $\nu=\prod_i\nu_i$, then the component functions are mutually orthogonal and \begin{equation}\Var g= \sum_i\Var g_i+\sum_{i<j}\Var g_{ij} + \dotsb + \Var g_{1\dotsb n}\end{equation} We refer to the decomposition of the mean measure using functional ANOVA (HDMR) as \emph{mean measure ANOVA} (MM-ANOVA). This gives us the set of sensitivity indices $\{(\amsbb{S}^a_u,\amsbb{S}^b_u): u\subseteq\{1,\dotsb,n\}\}$. When $\nu=\prod_i\nu_i$, then similar to RM-ANOVA, we can define entropy of $f$ on $E$ through the structural sensitivity indices \begin{equation}H(f) = \sum_{u\subseteq\{1,\dotsb,n\}}-\amsbb{S}^a_u\log\amsbb{S}^a_u\end{equation}

\subsection{Combined (RM-MM-ANOVA)} We combine the ANOVAs of the previous sections. Let $\{A,\dotsb,B\}$ be a disjoint partition of $E$. Let $g\in L^2(E,\mathcal{E},\nu)$. Take $f=(g-\E g)^2\in\mathcal{E}_+$. Let $f_a=f\ind{A}, \dotsb, f_b=f\ind{B}$. Let $N=(\kappa,\nu)$ be the random measure on $(E,\mathcal{E})$ formed by STC from independency $\mathbf{X}=\{X_i\}$ with mean and variance \begin{align}\E Nf &= c\nu f = c\sum_{u,v\subseteq\{1,\dotsb,n\}}\Cov(g_u,g_v)\\\Var Nf &= c\nu f^2 + (\delta^2-c)(\nu f)^2=\sum_{D\in\{A,\dotsb,B\}} \Var Nf_d + \sum_{D_i\ne D_j\in\{A,\dotsb,B\}} \Cov(Nf_{d_i},Nf_{d_j})\end{align} The variance is finite if $g\in L^4(E,\mathcal{E},\nu)$. We have respective sensitivity indices $\{(\amsbb{S}^a_u,\amsbb{S}^b_u): u\subseteq\{1,\dotsb,n\}\}$ and $\{(\amsbb{S}^a_d,\amsbb{S}^b_d)\}$. For orthogonal $N$, we have sensitivity distribution $\amsbb{S}$ on $(E,\mathcal{E})$.

\subsection{Counting distributions} The key quantity of $\kappa$ is $\delta^2-c$. This determines the correlative structure of $N$. When $\delta^2-c=0$, $N$ is orthogonal. In Table~\ref{tab:1} we give some counting distributions and indicate properties. We give pairs of counting distributions across negative, zero, and positive correlation. We indicate whether or not the random measures are closed under restriction to subspaces, as well as some limiting relations. As the orthogonal $N$ are special and we have two choices, we suggest that the orthogonal die be used whenever there is a theoretical reason for bounded support, such as a finite bound on the number of points, otherwise Poisson should be used. The orthogonal die random measure has a Poisson limit, so in many problems the distinction is meaningless.

\begin{table}[h!]
\begin{center}
\begin{tabular}{lclcl}
\toprule
Name & Support & $\delta^2-c$ & Closure & Limit(s)\\\midrule
Dirac$(c)$& $\{c\}$ & $-c$ & No &  \\
Binomial$(n,p)$ & $\{0,\dotsb,n\}$ & $-np^2$ & Yes & Dirac, Poisson \\
Poisson$(c)$ & $\N_{\ge0}$ & 0 & Yes &\\
Orthogonal-die$(m,n)$ & $\{m,\dotsb,n\}$ & 0 & No & Poisson \\
Negative-binomial$(r,p)$ & $\N_{\ge0}$ & $+r(\frac{p}{1-p})^2$ & Yes & Poisson\\
Zeta$(s)$ & $\N_{>0}$ & $+\frac{\zeta (s-1) \zeta (s+1)-\zeta (s) (\zeta (s)+\zeta (s+1))}{\zeta (s+1)^2}$ for $s>2$ & No &\\
\bottomrule
\end{tabular}\caption{Counting measures}\label{tab:1}
\end{center}
\end{table}

\FloatBarrier

\subsection{Related work}\label{sec:related}In this section we describe related work: other work that solves the same problems with different methods, work that uses the same methods to solve different problems, work that is similar to ours that solves similar problems, and a discussion of related problem domains. 

\paragraph{Connection to functional HDMR} An alternative approach to UQ of risk functionals is representation using functional HDMR \citep{alis99}. Taking $(E,\mathcal{E})$ as a function space, functional ANOVA of \emph{functionals} has been considered in functional HDMR and functional-cut-HDMR, which project the risk functional $f\in\mathcal{E}_+$ into functional subspaces and decomposes $\Var f$ . For Dirac $N$, we have $\Var Nf = c\Var f$, and hence, assuming $\nu=\prod_i\nu_i$, functional HDMR may be used to decompose the variance into functional subspaces \[\Var f = \sum_i\Var f_i + \sum_{i<j}\Var f_{ij} + \dotsb + \Var f_{1\dotsb n}\] This is distinct from RM-ANOVA, which decomposes $\Var Nf$ on a partition of disjoint subspaces. 

\paragraph{Other UQ methods for $g$}We describe some other UQ methods for studying the structure of $g$ in terms of the input coordinates. Besides functional ANOVA, other approaches towards gaining functional UQ of $g$ are variable importance and dependence measures defined by partial dependence \citep{friedman2001}, derivative-based global sensitivity indices \citep{dgsm1}, entropy-based methods \citep{mi,ent}, functional principal component analysis \citep{fpca}, and polynomial chaos expansions \citep{pc}. For example, for derivative-based GSA, we have functionals \[\eta_i = \int_E\left(\frac{\partial g(x_1,\dotsc,x_n)}{\partial x_i}\right)^2\nu(\D x)\for i\in\{1,\dotsb,n\}\] with the relation that \[\sum_{u\subseteq i}\amsbb{S}_u^a\le \frac{\eta_i}{\pi^2\Var g}\] 





\paragraph{Other applications of mixed binomial processes} 
 Mixed binomial processes are very general and include many well-known processes such as the binomial and Poisson processes. For example, the Poisson random measure is related to the structure of L\'{e}vy processes, Markov jump processes, and the excursions of Brownian motion, and the Poisson random measure is prototypical, enabling construction of a variety of more evolved processes.
 
 \paragraph{Other random measures} The random counting measures considered here are $\N_{\ge0}$ valued. For real-valued random measures $N$, the law of $N$ is encoded in the characteristic functional $\E e^{iNf}$ for $f\in\{\mathcal{E}-\text{measurable functions}\}$. A prototypical real-valued additive Gaussian random measure is Wiener.
 
 \begin{samepage}
 \begin{prop}[Wiener measure]\label{prop:wienerm} Consider the $\R$-valued Wiener measure $W_t$ on $(\R_+,\mathcal{B}_{\R_+})$ defined as \[W_t f = \int_{[0,t]}f(s)\D W_s\for f\in L^2([0,t],\mathcal{B}_{[0,t]},\text{Leb})\] with characteristic functional \[\E e^{i W_t f} = \exp_-\frac{1}{2}\int_{[0,t]}f^2(s)\D s\for f\in L^2([0,t],\mathcal{B}_{[0,t]},\text{Leb})\] Then $W_t$ is additive and for $f\in L^2([0,t],\mathcal{B}_{[0,t]},\text{Leb})$ the random variable $W_tf$ is Gaussian with characteristic function \[\E e^{irW_tf}= \exp_-\frac{1}{2} r^2 \int_{[0,t]}f^2(s)\D s\for r\in\R\] mean and variance \begin{align*}\E W_t f &= 0\\\Var W_t f &= \int_{[0,t]}f^2(s)\D s\end{align*} and for $f,g\in L^2([0,t],\mathcal{B}_{[0,t]},\text{Leb})$ covariance \[\Cov(W_sf,W_tg) = \int_{[0,s\wedge t]}f(u)g(u)\D u\] 
 \end{prop}\end{samepage}
 
 \begin{proof}For all choices of disjoint $f,\dotsb, g$ in $L^2([0,t],\mathcal{B}_{[0,t]},\text{Leb})$, e.g. $fg=0$, additivity follows from \[\E e^{iW_t(f+\dotsb+g)} = \exp_-\frac{1}{2}\int_{[0,t]}(f(s)+\dotsb+g(s))^2\D s = \exp_-\frac{1}{2}\int_{[0,t]}(f^2(s)+\dotsb+g^2(s))\D s = \E e^{iW_t f}\dotsb\E e^{i W_t g}\] so $W_tf,\dotsb,W_tg$ are independent random variables. Similarly independence holds for orthogonal $f,\dotsb,g$ in $L^2([0,t],\mathcal{B}_{[0,t]},\text{Leb})$. The characteristic function follows from the characteristic functional as $\E e^{irW_tf}=\E e^{iW_t(rf)}$, which shows $W_tf\sim\text{Gaussian}(0,\int_{[0,t]}f^2(s)\D s)$. The covariance follows from It\^{o} isometry \[\Cov(W_sf,W_tg)=\E W_sfW_tg=\E\int_{[0,s]}f(u)\D W_u\int_{[0,t]}g(u)\D W_u=\int_{[0,s\wedge t]}f(u)g(u)\D u\] 
 \end{proof}

\paragraph{Other applications of HDMR} Consider $n$ variate function $g(x_1,\dotsc,x_n)$. Defining the $T$-order HDMR of $g$ as  \[g^T(x_1,\dotsc,x_n)=\sum_{u\subseteq\{1,\dotsb,n\}: |u|\le T}g_u(x_u)\] in many practical problems for $T\ll n$ we have $g^T\simeq g$, so $g^T$ is a reduced-order representation of $g$. Oftentimes evaluating $g$ may be expensive, whereas evaluation of the reduced-order model $g^T$ is typically fast and efficient. 


\paragraph{Connection to bootstrap} 
Random measure uncertainty quantification is commonly practiced through use of the binomial process $N=nF_n$ (Dirac $\kappa=\delta_n$) and its empirical distribution $F_n=\frac{1}{n}N$ in the bootstrap estimator \citep{bootstrap1}. The binomial process is a degenerate mixed binomial process, so bootstrap analyses are contained within the framework of mixed binomial processes. 



\section{Risk functionals of input-output models}\label{sec:app}

A common scenario in data science is an independency of input-output data $(\mathbf{X},\mathbf{Y})=\{(X_i,Y_i)\}$ taking values in $(E\times F,\mathcal{E}\otimes\mathcal{F})$ with distribution $\mu=\nu\times Q$ where $Q$ is a transition probability kernel from $(E,\mathcal{E})$ into $(F,\mathcal{F})$, i.e. $\mu(\D x, \D y)=\nu(\D x)Q(x,\D y)$. Let $M=(\kappa,\nu\times Q)$ be a random counting measure on $(E\times F,\mathcal{E}\otimes\mathcal{F})$.

\subsection{Regression}\label{sec:regress}
For regression, consider $Mf_\theta$ for $f_\theta(x,y)=(y-g_\theta(x))^2$, where $g_\theta$ is some regressor with parameters $\theta\in\Theta$ and $F=\R$ \[Mf_\theta = \sum_{i}^K(Y_i-g_\theta(X_i))^2\] with mean and variance \begin{align*}\E Mf_\theta &= c\mu f_\theta\\\Var Mf_\theta &= c\mu f_\theta^2 + (\delta^2-c)(\mu f_\theta)^2\end{align*} Note that \[\mu f_\theta = \int_{E\times F}\mu(\D x,\D y)f_\theta(x,y)=\int_E\nu(\D x)\int_FQ(x,\D y)f_\theta(x,y)\] We have risk $R(\theta)=\mu f_\theta =\frac{1}{c}\E Mf_\theta$ so the risk estimate is attained as \[\argmin_{\theta\in\Theta}R(\theta)\] Now consider a finite disjoint partition $\{A,\dotsb,B\}$ of $E\times F$ and define $f_\theta^a(x,y) = f_\theta(x,y)\ind{A}((x,y))$ and $f_\theta^b(x,y) = f_\theta(x,y)\ind{B}((x,y))$. This gives \[\Cov(Mf_\theta^a,Mf_\theta^b)=(\delta^2-c)\mu f_\theta^a\mu f_\theta^b\] Therefore we have the decomposition \begin{align}\Var Mf_\theta &= \sum_{D\in\{A,\dotsb,B\}} \Var Mf_\theta^d + \sum_{D_i\ne D_j\in\{A,\dotsb,B\}} \Cov(Mf_\theta^{d_i},Mf_\theta^{d_j})\end{align}This gives us the set of sensitivity indices $\{(\amsbb{S}^a_d,\amsbb{S}^b_d): D\in\{A,\dotsb,B\}\}$. For orthogonal $N$ the sensitivity density is defined as \[\amsbb{S}^a_u(\D x_u) = \frac{1}{(\nu\times Q)f_\theta^2}\int_{E_{-u}}\nu(\D x)\int_F Q(x,\D y)f_\theta^2(x,y)\for x_u\in E_u\] We can also define sensitivity densities on the output and input-output spaces as  \[\amsbb{S}^a(\D y) = \frac{1}{(\nu\times Q)f_\theta^2}\int_{E}\nu(\D x)Q(x,\D y)f_\theta^2(x,y)\for y\in F\] and  \[\amsbb{S}^a_u(\D x_u,\D y) = \frac{1}{(\nu\times Q)f_\theta^2}\int_{E_{-u}}\nu(\D x)Q(x,\D y)f_\theta^2(x,y)\for (x_u,y)\in E_u\times F \] 

For the decomposition of risk $R(\theta)$, we take $E$ with $n$ dimensions, set $c_0=\nu g_\theta\in\R$ and put $F=\{c_0\}$ with $Q(x,\cdot)=\ind{F}(x)$. Therefore \begin{equation}R(\theta)=\mu f_\theta = \mu(g_\theta-c_0)^2 =\Var g_\theta = \sum_{u,v\subseteq\{1,\dotsb,n\}}\Cov(g_u,g_v)\end{equation} This gives us the set of sensitivity indices $\{(\amsbb{S}^a_u,\amsbb{S}^b_u): u\subseteq\{1,\dotsb,n\}\}$. 


\subsection{Classification}\label{sec:class}

For classification, consider $Mf_\theta$ for $f_\theta(x,y)=\ind{}(y\ne g_\theta(x))$ where $g_\theta$ is some classifier and $F$ is countable. We have that \begin{equation}\mu f_\theta = \int_{E\times F}\mu(\D x,\D y)f_\theta(x,y)=\int_E\nu(\D x)\sum_{y\in F}Q(x,\{y\})f_\theta(x,y)=\P(y\ne g_\theta(x))\end{equation} Risk is similarly defined. Consider the partition of $\{A,\dotsb,B\}$ of $E\times F$ and define $f_\theta^a(x,y) = f_\theta(x,y)\ind{A}((x,y))$ so that we have the vector $(Mf_\theta^a,\dotsb,Mf_\theta^b)$. The $f_\theta^a,\dotsb,f_\theta^b$ are disjoint. We have $\Cov(Mf_\theta^a,Mf_\theta^b)=(\delta^2-c)\mu f_\theta^a\mu f_\theta^b$. In a similar manner to regression and noting that $f_\theta^2=f_\theta$, we can attain a partition of the variance $\Var Mf_\theta$ with sensitivity indices $\{(\amsbb{S}^a_d,\amsbb{S}^b_d):D\in\{A,\dotsb,B\}\}$. For orthogonal $N$, we have the sensitivity density \[\amsbb{S}^a_u(\D x_u) = \frac{1}{(\nu\times Q)f_\theta}\int_{E_{-u}}\nu(\D x)\sum_{y\in F}Q(x,\{y\})f_\theta(x,y)\for x_u\in E_u\] We also can define a sensitivity density on $F$ as \[\amsbb{S}^a\{y\} = \frac{1}{(\nu\times Q)f_\theta}\int_{E}\nu(\D x)Q(x,\{y\})f_\theta(x,y)\for y\in F\] and density on $E_u\times F$ as \[\amsbb{S}^a(\D x_u,y) = \frac{1}{(\nu\times Q)f_\theta}\int_{E_{-u}}\nu(\D x)Q(x,\{y\})f_\theta(x,y)\for (x_u,y)\in E_u\times F\]


For Dirac $\kappa$ we have \[\Var Mf_\theta = c\mu f_\theta(1-\mu f_\theta)\] and sensitivity indices \[\amsbb{S}_d^a = \frac{\Var Mf_\theta^d}{\Var Mf_\theta} = \frac{\mu f_\theta^d(1-\mu f_\theta^d)}{\mu f_\theta(1-\mu f_\theta)}\for D\in\{A,\dotsb,B\}\] and \[\amsbb{S}_d^b = \sum_{D_i\in\{A,\dotsb,B\}: D_i\ne D}\frac{\Cov(Mf_\theta^d,Mf_\theta^{d_i})}{\Var Mf_\theta} = \sum_{D_i\in\{A,\dotsb,B\}: D_i\ne D}-\frac{\mu f_\theta^d\mu f_\theta^{d_i}}{\mu f_\theta(1-\mu f_\theta)}\for D\in\{A,\dotsb,B\}\]

For Poisson or orthogonal die $\kappa$, we have $\delta^2=c$ and \begin{equation}\Var Mf_\theta = c\nu f_\theta^d\end{equation} and \begin{equation}\amsbb{S}_d^a = \frac{\mu f_\theta^d}{\mu f_\theta}\for D\in\{A,\dotsb,B\}\end{equation} The sensitivity distribution is given by \[\amsbb{S}(\D x,\D y) = \mu(\D x,\D y) f_\theta(x,y) / \mu f_\theta\]

\section{Random field ANOVA}\label{sec:rf} We discuss a connection between STC of mixed binomial processes and positive random fields, touching fPCA and random field theories. 

\begin{samepage}
\begin{prop}[Positive random field]\label{prop:1} Let $N=(\kappa,\nu)$ be a random measure on $(E,\mathcal{E})$ formed by independency $\mathbf{X}=\{X_i\}$. Let $(F,\mathcal{F})$ be a measurable space and let $k:E\times F\mapsto\R_+$ be $\mathcal{E}\otimes\mathcal{F}$-measurable. Then \[G(y) = \int_E N(\D x)k(x,y)=\sum_i^Kk(X_i,y)\for y\in F\] defines a \emph{positive random field} $G=\{G(y): y\in F\}$ on $(F,\mathcal{F})$ with law specified by the Laplace transform \[\E e^{-\alpha G} = \psi(\nu e^{-\int_F\alpha(\D y)k(\cdot,y)})\quad\text{for every finite measure }\alpha\text{ on }(F,\mathcal{F})\] and, putting $f_y(\cdot)=k(\cdot,y)\in\mathcal{E}_+$ so that $G(y)=Nf_y$ for $y\in F$, with mean and covariance \begin{align*}U(y)&=\E G(y)=c\nu f_y\for y\in F\\C(y,z)&=\Cov(G(y),G(z))=c\nu(f_yf_z)+(\delta^2-c)\nu f_y \nu f_z\for y,z\in F\end{align*} 
\end{prop}
\end{samepage}
\begin{proof} The law of $G$ is specified by the finite-dimensional distributions of $(G(y_1),\dotsb,G(y_n))$ for $n\ge1$ and $y_1,\dotsb,y_n\in F$ or equivalently by the Laplace transform $\E e^{-\alpha G} = \E e^{-\int_F\alpha(\D y)G(y)}$ for all finite measures $\alpha$ on $(F,\mathcal{F})$. Given $N$ and $k$, the Laplace transform of $G$ follows from the Laplace functional of $N$ for test function $f_y(\cdot)=k(\cdot,y)\in\mathcal{E}_+$. The mean and covariance follow from $N$.\end{proof}


\begin{re}[Densities]The density of $Nf_y=G(y)$, $\eta_y$, defines a transition probability kernel $\eta_y(\cdot)=Q(y,\cdot)$ from $(F,\mathcal{F})$ into $(E,\mathcal{E})$. Let $\alpha$ be a probability measure on $(F,\mathcal{F})$. Then  $\alpha\times Q$ is the joint distribution on $(F\times E,\mathcal{F}\otimes\mathcal{E})$  \[(\alpha\times Q)f = \int_F\alpha(\D y)\int_EQ(y,\D x)f(x,y)\for f\in(\mathcal{F}\otimes\mathcal{E})_+\] and $\alpha Q$ is the marginal distribution of $\alpha\times Q$ on $(E,\mathcal{E})$ \[(\alpha Q)f = \int_F\alpha(\D y)\int_EQ(y,\D x)f(x)\for f\in\mathcal{E}_+\]
\end{re}

\begin{figure}[h!]
\centering
\begingroup
\captionsetup[subfigure]{width=3in,font=normalsize}
\subfloat[Orthogonal $N$\label{fig:kernel1}]{\includegraphics[width=3.5in]{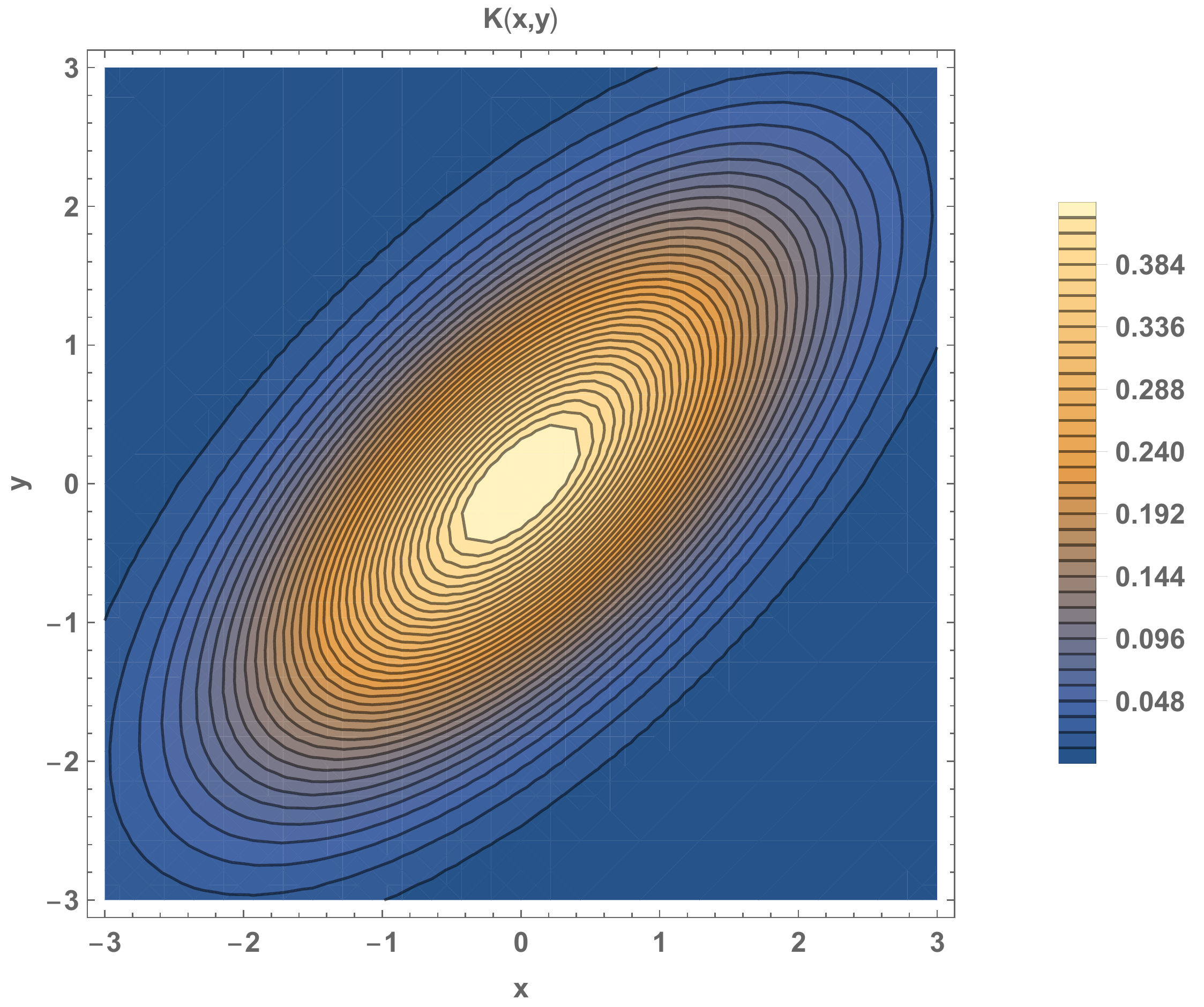}}
\subfloat[Dirac $N$\label{fig:kernel3}]{\includegraphics[width=3.5in]{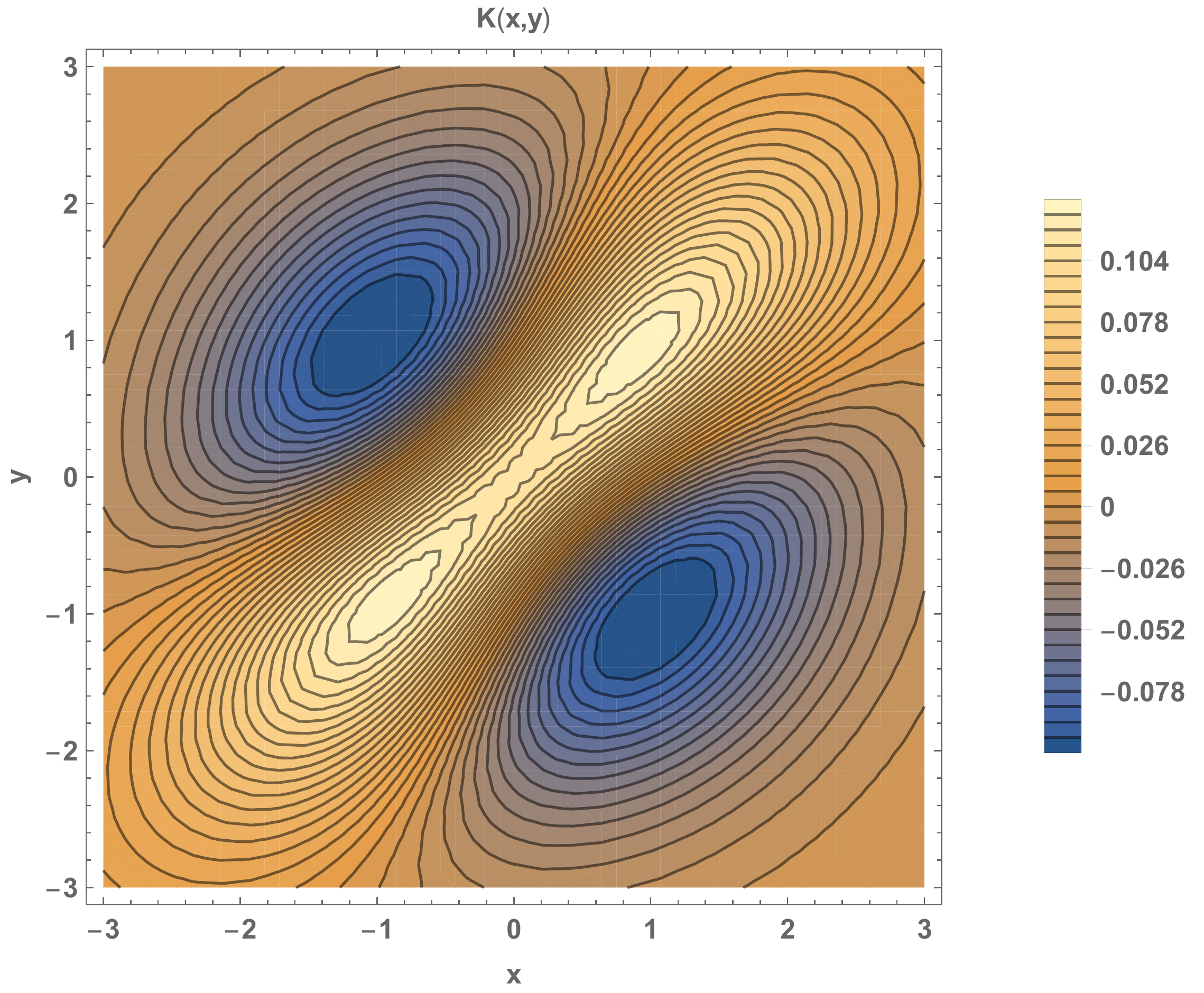}}\\
\endgroup
\caption{Random field kernel $C$ using radial basis kernel $k$ with $c=\gamma=1$ for orthogonal and Dirac $N$}\label{fig:kernel}
\end{figure}

We give an example for Proposition~\ref{prop:1}: consider $(F,\mathcal{F})=(E,\mathcal{E})$ where $E$ has dimension one, $\nu$ as standard Gaussian, and the radial basis kernel $k(x,y)=e^{-\gamma(x-y)^2}$ so that \[\nu(f_yf_z) = \frac{e^{-\frac{\gamma}{1+4\gamma} \left(y^2+2 \gamma  (y-z)^2+z^2\right)}}{\sqrt{1+4 \gamma}}\] and \[\nu f_y = \frac{1}{\sqrt{1+2\gamma}}e^{-\frac{\gamma}{1+2\gamma}y^2}\] In Figure~\ref{fig:kernel} we show $C$ on $[-3,3]\times[-3,3]$ for the radial basis kernel with $c=\gamma=1$ for orthogonal and Dirac $N$. The random field covariance kernels are very different: in Figure~\ref{fig:kernel1} orthogonal has positive correlation and a single mode located at $y=z=0$, whereas in Figure~\ref{fig:kernel3} Dirac has a four correlative modes, two positive and two negative, around $y,z\approx \pm 1$. 

\FloatBarrier

\begin{re}[Random fields on graphs]\label{re:graph} Let $(F,\mathcal{F})$ be indexed by the vertices of some graph space $(V,E)$ with distribution $\alpha$. Examples of $\alpha$ include the Markov random field and conditional random field. 
\end{re}

Next we use the Karhunen-Lo\`{e}ve theorem to get the eigensystem expansion of $C$. We note that KL theorem is a type of polynomial chaos expansion, where the basis is optimal in $L^2$ with respect to $\mu$. The normalized eigenvalues reveal the contributions to the variance. 

\begin{samepage}
\begin{prop}[fPCA expansion]\label{prop:fpca}
Assume $C$ is continuous. Then $C$ can be decomposed \[C(y,z)=\sum_{k\ge1}\lambda_k\varphi_k(y)\varphi_k(z)\] where $\{\lambda_k\}$ and $\{\varphi_k\}$ are the eigenvalues and orthonormal eigenfunctions of the operator of the form \begin{align*}f&\mapsto\int_F\mu(\D y)K(y,z)f(y)\\ L^2(F,\mathcal{F},\mu)&\mapsto L^2(F,\mathcal{F},\mu)\end{align*} and the centered process $G-U$ can be expressed as \[G(y)-U(y) = \sum_{k\ge1}\xi_k\varphi_k(y)\] where \[\xi_k = \int_F\mu(\D y)(G(y)-U(y))\varphi_k(y)\] with $\E\xi_k=0$, $\Var\xi_k=\lambda_k$, and $\Cov(\xi_k,\xi_l)=0$ for $k\ne l$. 
\end{prop}
\end{samepage}

\begin{de}[RF-ANOVA] Consider the random field $G$ on $(F,\mathcal{F})$ formed from random measure $N=(\kappa,\nu)$ on $(E,\mathcal{E})$ and function $k:E\times F\mapsto\R_+$. Putting $\lambda = \sum_i\lambda_i$, the sensitivity indices of $G$ are the normalized variances \[\amsbb{S}_i = \frac{\lambda_i}{\lambda}\for i\in\N_{>0}\]
\end{de}


\begin{re}[Discrete approximation]The continuous eigenfunctions may be difficult to compute. In practice a finite number of eigenvalues and eigenfunctions can be attained to arbitrary precision using the eigensystem decomposition of the matrix $(C(y,z): y,z\in\{a,\dotsb,b\})$. \end{re}

\subsection{Related work} Random fields are ubiquitous and are tantamount to stochastic processes. There are many random field models. This section shows construction of general positive random fields $G$ with law determined by the Laplace transform, derived from the Laplace functional of $N$. Poisson random fields are canonical, although the terminology in the literature is conflicting, as `Poisson random field' can sometimes mean either Poisson random measure or some other Poisson-derived quantity. 

The law of real-valued random field $G$ on $(F,\mathcal{F})$ is determined by the characteristic function $\E e^{irG}$ for functions $r$ such that $\int_F r(y)G(y)\D y<\infty$. A common choice is Gaussian. 


\section{Density identification through inversion of Laplace transforms} Let $N=(\kappa,\nu)$ be a random measure on $(E,\mathcal{E})$. For $f\in\mathcal{E}_+$, the random variable $Nf$ has Laplace transform $F(\alpha)$ and distribution $\eta$. In principle $\eta$ is attained from the inverse Laplace transform of $F(\alpha)$ but in practice this can be difficult. 

\begin{re}[Indicators] For indicator functions $f=\ind{A}$, $A\in\mathcal{E}$, the distribution of $N(A)$, denoted $\kappa_A$, is encoded by the pgf $\psi_A(t)=\psi(at+1-a)$ \[\kappa_A\{k\} = \P(N(A)=k)  = \psi_A^{(k)}(0) / k!\for k\in\N_{\ge0}\]
\end{re}

We consider $f\in\mathcal{E}_+$ such that $Nf\in\R_+$, whose law $\eta$ is encoded in the Laplace transform. We can attain $\eta$ through the method of maximum entropy using knowledge of some number of generalized moments of $Nf$, described in for instance \cite{maxent}, which are evaluations of the Laplace transform. We describe the basic set-up and application to random measures.

The basic problem of interest is the truncated moment problem \[\text{identify }\eta\text{ such that } F(\alpha_i)= f_i\in[0,1]\for i=1,\dotsb,n\] We introduce the change of variables $Nf\mapsto y=e^{-Nf}\in[0,1]$. Then \[F(\alpha) = \E y^{\alpha}=\psi(\nu e^{-\alpha f})\] is the generalized moment of $y$, where $y$ has distribution $\mu$.  It turns that out that the collection \[\{f_i=\E y^{\alpha_i}: i\in\N_{>0}, \alpha_1>\alpha_2>\dotsb, \sum_{i}\alpha_i=\infty, \lim_{i\rightarrow\infty}\alpha_i=0\}\] uniquely determines $\mu$. Given $\alpha_1>\dotsb>\alpha_n$, i.e. the first $n$ terms, the density $\mu(\D x)=\P(y\in\D x)$ on $[0,1]$ can be attained using the method of maximum entropy as $\mu_n$, that is, we can use maximum entropy to invert the Laplace transform evaluated on a set of points. The maximum entropy distribution of $\mu$ given $n$ generalized moment constraints is \begin{equation}\mu_n(\D y) = \D y\frac{1}{Z(\lambda_1,\dotsb,\lambda_n)}\exp_-(\lambda_1 y^{\alpha_1}+\dotsb+\lambda_n y^{\alpha_n})\for y\in[0,1]\end{equation} where \[Z(\lambda_1,\dotsb,\lambda_n) = \int_{[0,1]}\D y\exp_-(\lambda_1 y^{\alpha_1}+\dotsb+\lambda_n y^{\alpha_n})\] and $\lambda_1,\dotsb,\lambda_n$ in $\R$ are attained as \begin{equation}(\lambda_1,\dotsb,\lambda_n)=\argmin_{(\lambda_1^*,\dotsb,\lambda_n^*)\in\R^n}\log Z(\lambda_1^*,\dotsb,\lambda_n^*) + \lambda_1^*f_1+\dotsb+\lambda_n^*f_n\end{equation} where the function being minimized is strictly convex. 


Then we can attain the maximum-entropy-derived density of $Nf$ as \begin{equation}\eta_n(\D x)=\P(Nf\in\D x)=e^{-x}(\mu_n\circ e^{-x})(\D x)\end{equation} on $(\R_+,\mathcal{B}_{\R_+})$. 

If $f$ takes large values and the integral must be computed numerically, then the calculation of $F(\alpha)$ can introduce numerical problems. To circumvent this, we can introduce the normalized function $f^*=f/C$, where $C>0$, such as $C=\nu f$. Then the density of $y=e^{-Nf^*}$ is attained as $\mu_n$, followed by that of $Nf$, given by \[\eta_n(\D x)=\P(Nf\in\D x)=\frac{1}{C}e^{-x/C}(\mu_n\circ e^{-x/C})(\D x)\]


\paragraph{Specific calculation} Consider the gamma distribution $\text{Gamma}(d,\lambda)$, where $d$ is the shape parameter and $\lambda$ is the rate parameter. We have \[F(\alpha) = \left(\frac{\lambda}{\lambda+\alpha}\right)^d\] We draw $\alpha_i$ for $i=1,\dotsb,10$ according to $\alpha_i\sim\text{Exponential}(1)$. Then we compute $\{F(\alpha_i)\}$ for $d=2$ and $\lambda=1$. Finally we estimate $\lambda_1,\dotsb,\lambda_{10}$. We show the maximum entropy density $\eta_{10}(\D x)$ below in Figure~\ref{fig:laplace} in comparison to the underlying density of $\text{Gamma}(2,1)$. The reconstruction is very good for a modest number ($n=10$) of samples.   

\begin{figure}[h]
\centering
\includegraphics[width=5in]{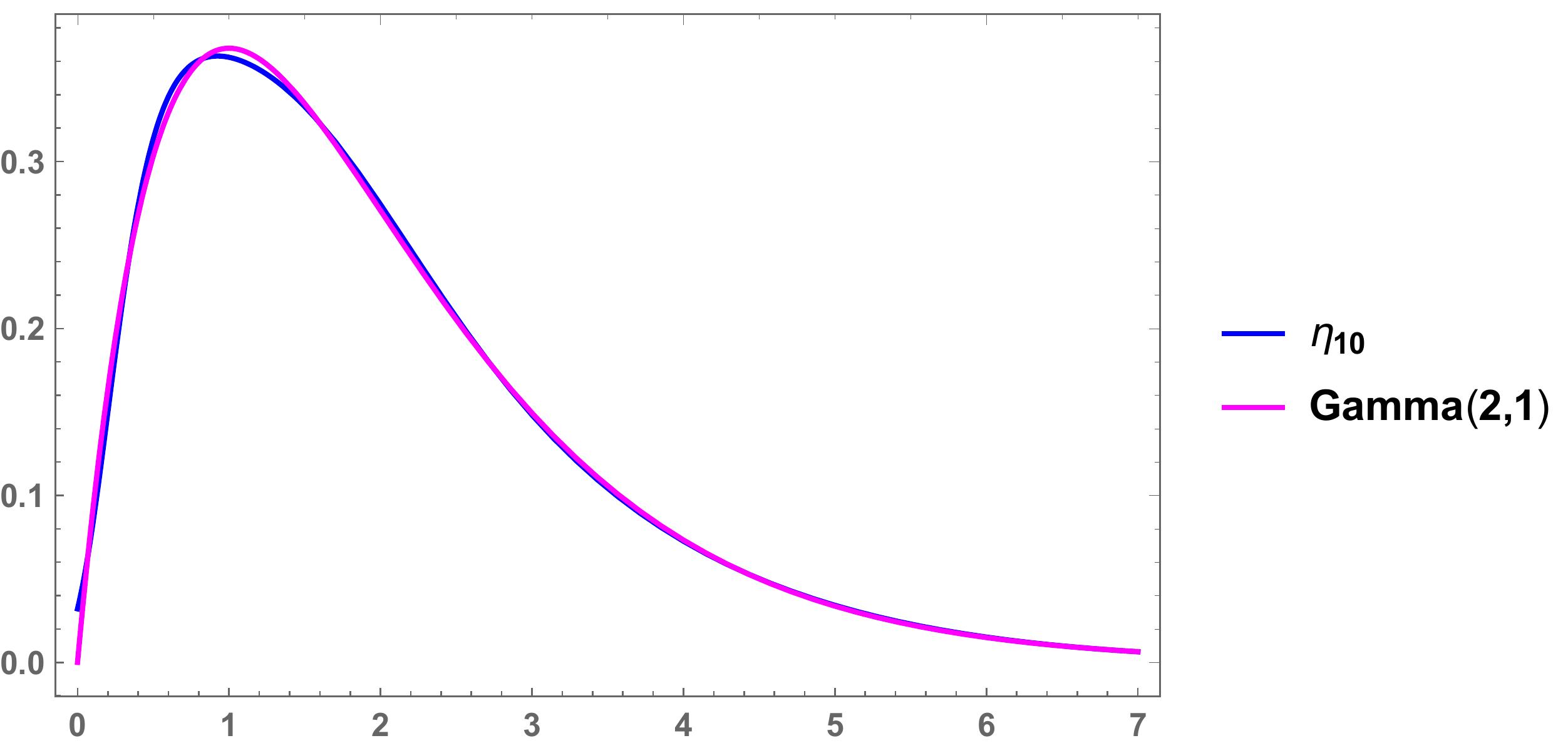}
\caption{Density of $Nf$, $\eta_n$ for $n=10$, based on maximum entropy for gamma random variable $\text{Gamma}(2,1)$}\label{fig:laplace}
\end{figure}

\FloatBarrier

\section{Examples}\label{sec:sim} Here we analyze in detail two canonical functions using RM-MM-ANOVA, an elementary symmetric polynomial of $n$ variables in Section~\ref{sec:elem} and a non-linear function commonly used in sensitivity analysis called the Ishigami function of three variables in Section~\ref{sec:ishi}. We compute the RM-MM-ANOVA decompositions for these problems. We describe application of the RM-MM-ANOVA framework to analyzing the risk functionals of regressors in Section~\ref{sec:regress} and classifiers in Section~\ref{sec:class} and to adaptive randomized controlled trials in Section~\ref{ex:rct} and to dynamic survival analysis of epidemics in Section~\ref{ex:dsa}. We provide two additional analyses of functions in the appendices: a symmetric bivariate polynomial with correlated inputs in Section~\ref{sec:poly} and a graph property in Section~\ref{sec:graph}. In Section~\ref{sec:sda} we discuss a random field model for interacting particle systems, and in Section~\ref{sec:ising} in appendices we have a short application on the Ising model. 

\subsection{Elementary symmetric polynomial}\label{sec:elem} As referenced in Table~\ref{tab:0}, we consider a real-valued test function $g$ and its HDMR and global sensitivity analysis. Consider the elementary symmetric polynomial \[g(x_1,\dotsc,x_n) = \prod_{i=1}^n x_i,\;\;\;\text{iid}\,x,\;\;\;\rho\equiv\sigma/\mu\ne 0\] with \begin{align*}\E g &=\mu^n\\\Var g &=\mu^{2n}((1+\rho^2)^n-1)\end{align*} where $\nu g^2 =\mu^{2n}(1+\rho^2)^n$.

For example, for $x_i\sim\text{Uniform}(0,1)$, we have \[\P(y=g(x_1,\dotsb,x_n)\in\D y) = \frac{(-1)^{n-1}}{(n-1)!}\log^{n-1}(y)\D y\for y\in(0,1)\]


The component functions of $g$ are linear combinations of elementary symmetric polynomials
\begin{align*}
g_0 &=  \mu^n\\
g_i(x_i) &= \mu^{n-1}x_i - g_0\\
g_{ij}(x_i,x_j)&= \mu^{n-2}x_ix_j - g_i(x_i) - g_j(x_j) - g_0\\
\vdots&.
\end{align*}
The partition of variance is given by \[\Var g =  \sum_i \Var g_i + \sum_{i<j}\Var g_{ij} + \dotsc + \Var g_{1\dotsb n},\]
where
\begin{align*}
\Var g_i &=\mu^{2n}\rho^2\\
\Var g_{ij} &=\mu^{2n}\rho^4\\
\vdots&\\
\Var g_{1\dotsb n}&=\mu^{2n}\rho^{2k}.
\end{align*} The sensitivity indices satisfy \[\sum_{\substack{k\\i_1<\dotsb<i_k}}\amsbb{S}_{i_1\dotsb i_k} = 1\] and at each order follow
\begin{align*}
\amsbb{S}^a_{u:|u|=k}\equiv\sum_{i_1<\dotsb<i_k}\amsbb{S}_{i_1\dotsb i_k} &= \frac{ \binom{n}{k}\rho^{2 k}}{\left(\rho^2+1\right)^n-1}
\end{align*}
where $\sum_{k=1}^n \amsbb{S}^a_{u:|u|=k} = 1$. Therefore $P=(\amsbb{S}^a_{u:|u|=k})$ forms a distribution on $\{1,\dotsc,n\}$. We have \begin{align*}\E P &= \frac{n \rho^2 \left(\rho ^2+1\right)^{n-1}}{\left(\rho ^2+1\right)^n-1}\\\Var P &= \frac{n \rho ^2 \left(\rho ^2+1\right)^{n-2} \left(\left(\rho ^2+1\right)^n-n \rho ^2-1\right)}{\left(\left(\rho ^2+1\right)^n-1\right)^2}\\H&=\sum_{k}^n-\amsbb{S}^a_{u:|u|=k}\log \amsbb{S}^a_{u:|u|=k}\end{align*} 

Below in Figure~\ref{fig:poly} we plot the sensitivity indices by subspace $(\amsbb{S}_{u:|u|=k}: k=1,\dotsc,n)$ for varying $\rho$ for $n=100$ and the entropy as a function of $\rho$.

\begin{figure}[h!]
\centering
\begingroup
\captionsetup[subfigure]{width=3in,font=normalsize}
\subfloat[Sensitivity distribution]{\includegraphics[width=4in]{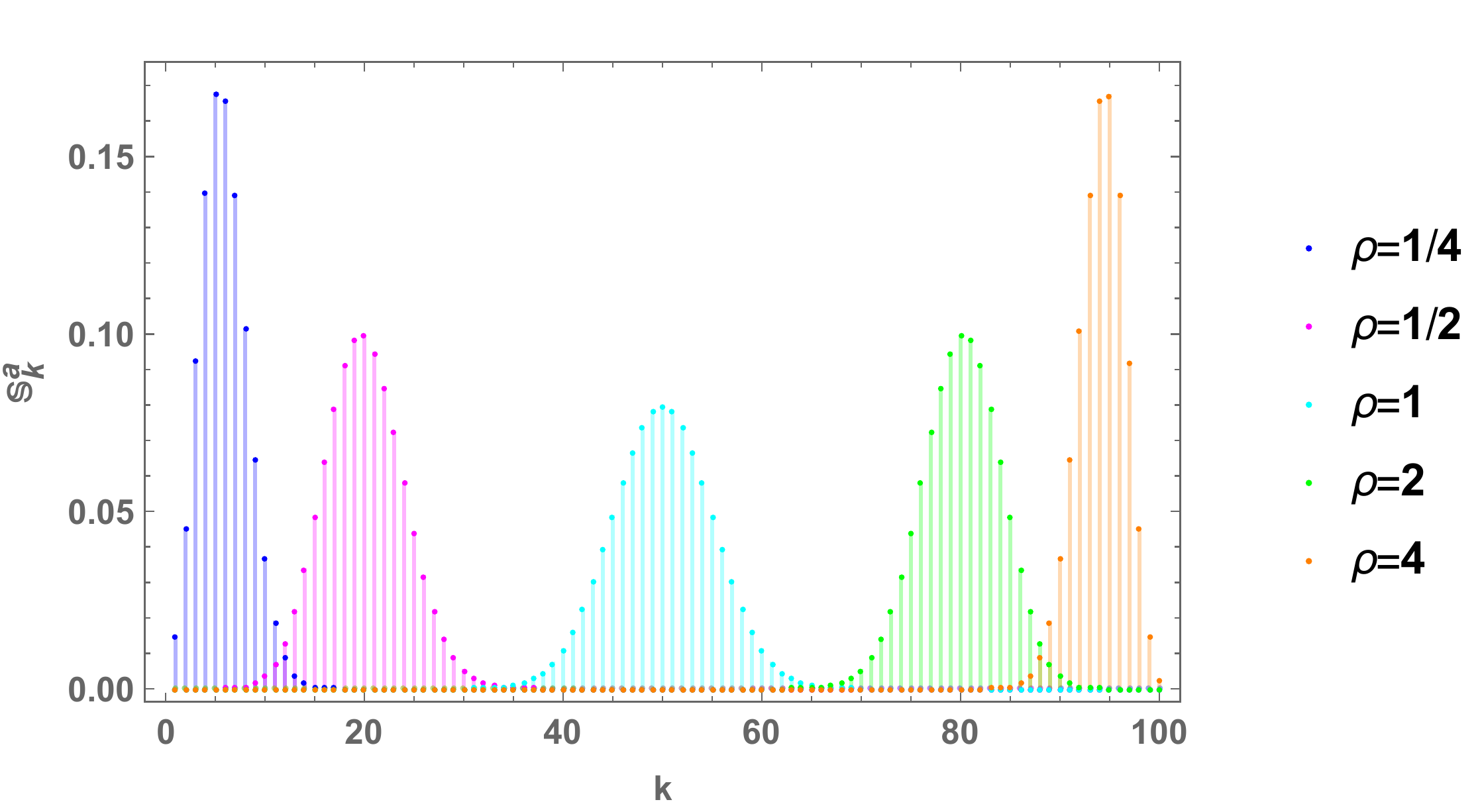}}
\subfloat[Entropy]{\includegraphics[width=3in]{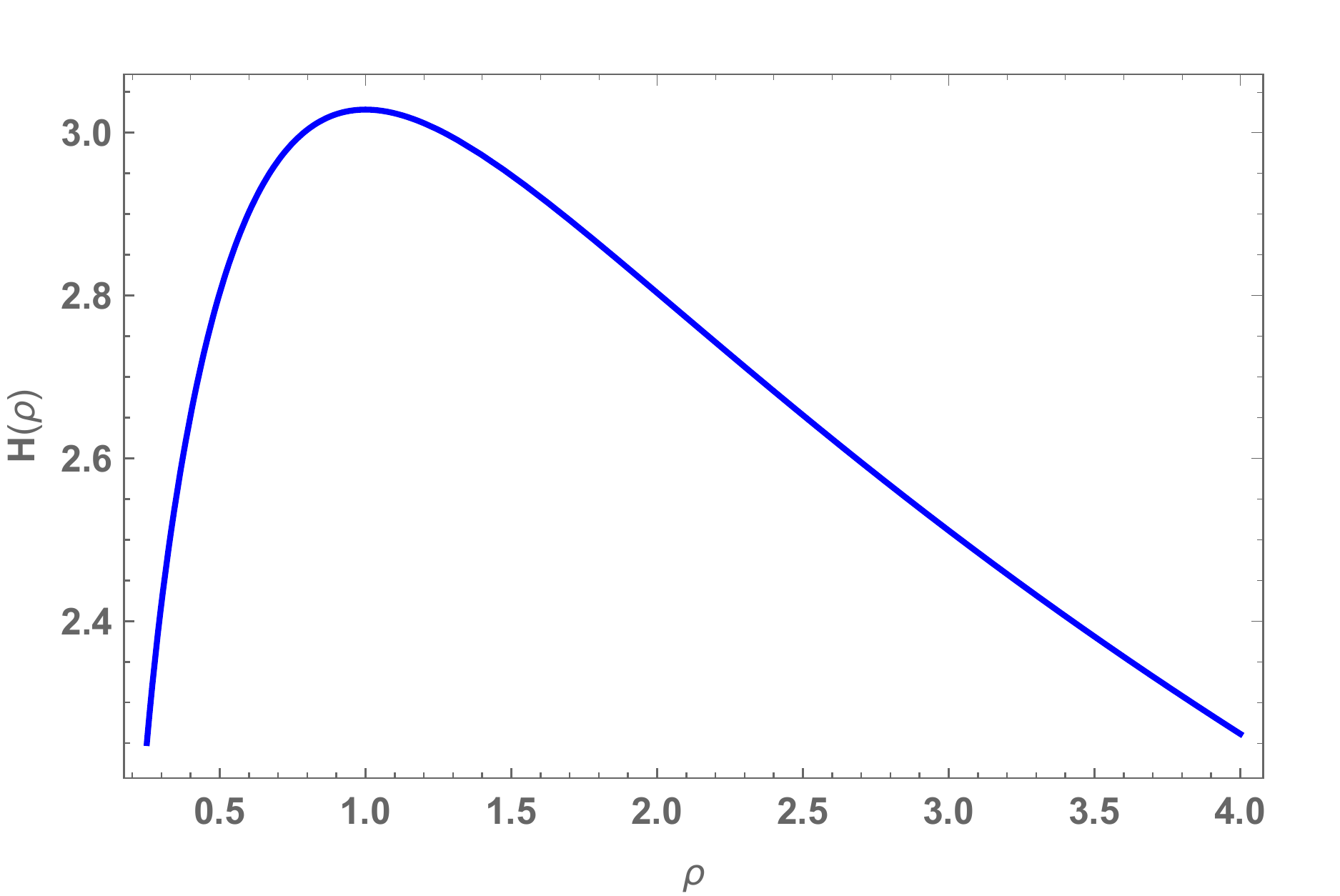}}\\
\endgroup
\caption{HDMR sensitivity indices of elementary symmetric polynomial of $n=100$ variables in subspace dimension and entropy for varying coefficient of variation $\rho$}\label{fig:poly}
\end{figure}

\FloatBarrier

\emph{These findings show that the effective dimension of the fractional contribution of the $k$-dimensional subspaces to the variance is regulated by the coefficient of variation $\rho$ and that the distribution on subspaces is not uniform.} If $\rho=1$, then \[\E P = \frac{n}{2-2^{-n+1}}\simeq\frac{n}{2}\] if $\rho<1$, then $\E P<n/2$; if $\rho>1$, then $\E P>n/2$. This is an elegant result that the mean measure of squared error risk functionals of monomials admits such a decomposition.  



Now consider \[f=(g-\E g)^2 =\left(\prod_{i=1}^n x_i - \mu^n\right)^2\in\mathcal{E}_+ \] so that we have the decomposition into various subspaces \begin{equation}\nu f = \Var g = \sum_{u\subseteq\{1,\dotsb,n\}}\Var g_u = \sum_{u\subseteq\{1,\dotsb,n\}:|u|\ge1} \mu^{2n}\rho^{2|u|}=\sum_{k=1}^n\binom{n}{k}\mu^{2n}\rho^{2k}\end{equation}

Consider random measure $N=(\kappa,\nu)$ on $(E,\mathcal{E})$. To define the Laplace functional $L$, note that \[\nu e^{-f} = \int_E\nu(\D x)e^{-(\prod_ix_i-\mu^n)^2}\] This may be calculated exactly in some cases, i.e. for $n=1$ and $E=[0,1]$ with uniform $\nu$ we have $\nu e^{-f} = \sqrt{\pi}\Erf(\frac{1}{2})$. Similarly for the Laplace transform we have $\nu e^{-\alpha f} = \sqrt{\frac{\pi}{\alpha}}\Erf\left(\frac{\sqrt{\alpha }}{2}\right)$.

 For $f\in\mathcal{E}_+$ the random variable $Nf$ has mean and variance \begin{align}\E Nf &= c\nu f =c\sum_{u\subseteq\{1,\dotsb,n\}} \mu^{2n}\rho^{2|u|}=c\sum_{k=1}^n\binom{n}{k}\mu^{2n}\rho^{2k}\\\Var Nf &=c\nu f^2 +(\delta^2-c)(\nu f)^2\end{align}  For the variance, we attain $\nu f^2$, which is the fourth central moment of $g$, \begin{align*}\nu f^2 &= \nu(g-\E g)^4 = \nu g^4 -4\E g\nu g^3+6(\E g)^2\nu g^2 -3(\E g)^4\end{align*} We have $\nu(\D x) g^j(x) = \prod_i \nu_i(\D x_i)x_i^j$ for $j=1,2,3,4$. For disjoint partition $\{A,\dotsb,B\}$ of $E$, we put $f_a=f\ind{A}, \dotsb, f_b=f\ind{B}$. We have \[\Cov(Nf_a,Nf_b)=(\delta^2-c)\nu f_a\nu f_b\] and \begin{align*}\Var Nf &= \sum_{D\in\{A,\dotsb,B\}} \Var Nf_d + \sum_{D_i\ne D_j\in\{A,\dotsb,B\}} \Cov(Nf_{d_i},Nf_{d_j})\\&= \sum_{D\in\{A,\dotsb,B\}} (c\nu f_d^2 +(\delta^2-c)(\nu f_d)^2) + \sum_{D_i\ne D_j\in\{A,\dotsb,B\}}(\delta^2-c)\nu f_{d_i}\nu f_{d_j} \end{align*} For orthogonal $N$ the sensitivity distribution in each coordinate is defined as \[\amsbb{S}^a_i(\D x_i) = \int_{E_1\times\dotsb\times E_{i-1}\times E_{i+1}\times\dotsb\times E_n}\frac{\nu(\D x)f^2(x)}{\nu f^2}\]

\begin{figure}[h!]
\centering
\begingroup
\captionsetup[subfigure]{width=3in,font=normalsize}
\subfloat[Orthogonal $N$: $\nu f_a^2$, $\nu f_b^2$, and $\nu f^2$ as a function of $p$\label{fig:1a}]{\includegraphics[width=3.5in]{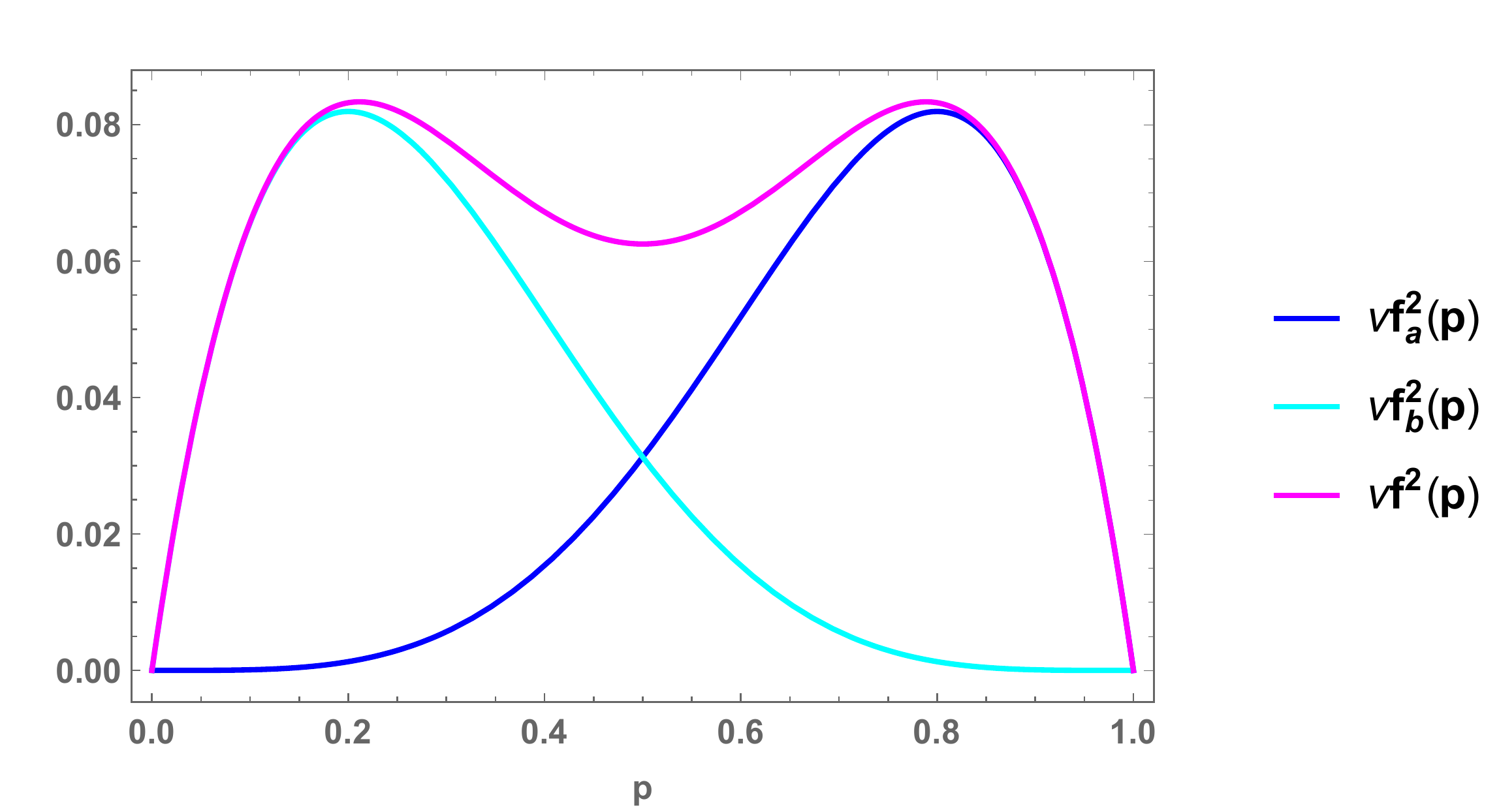}}
\subfloat[Dirac $N$: $\Var f_a$, $\Var f_b$, and $\Var f$ as a function of $p$\label{fig:1b}]{\includegraphics[width=3.5in]{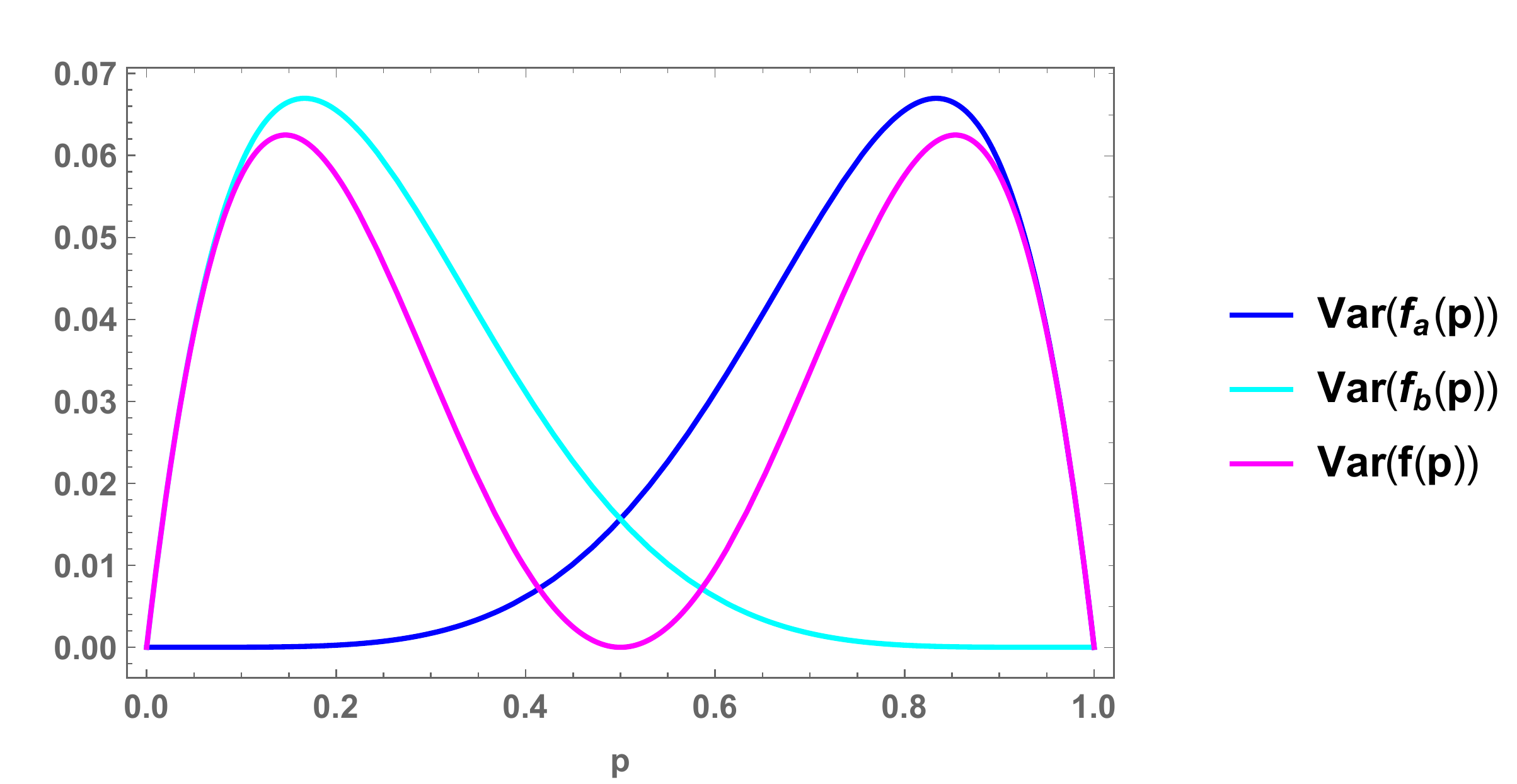}}\\
\subfloat[Orthogonal $N$: $\amsbb{S}_a^a$ and $\amsbb{S}_b^a$ as a function of $p$ \label{fig:1c}]{\includegraphics[width=3.5in]{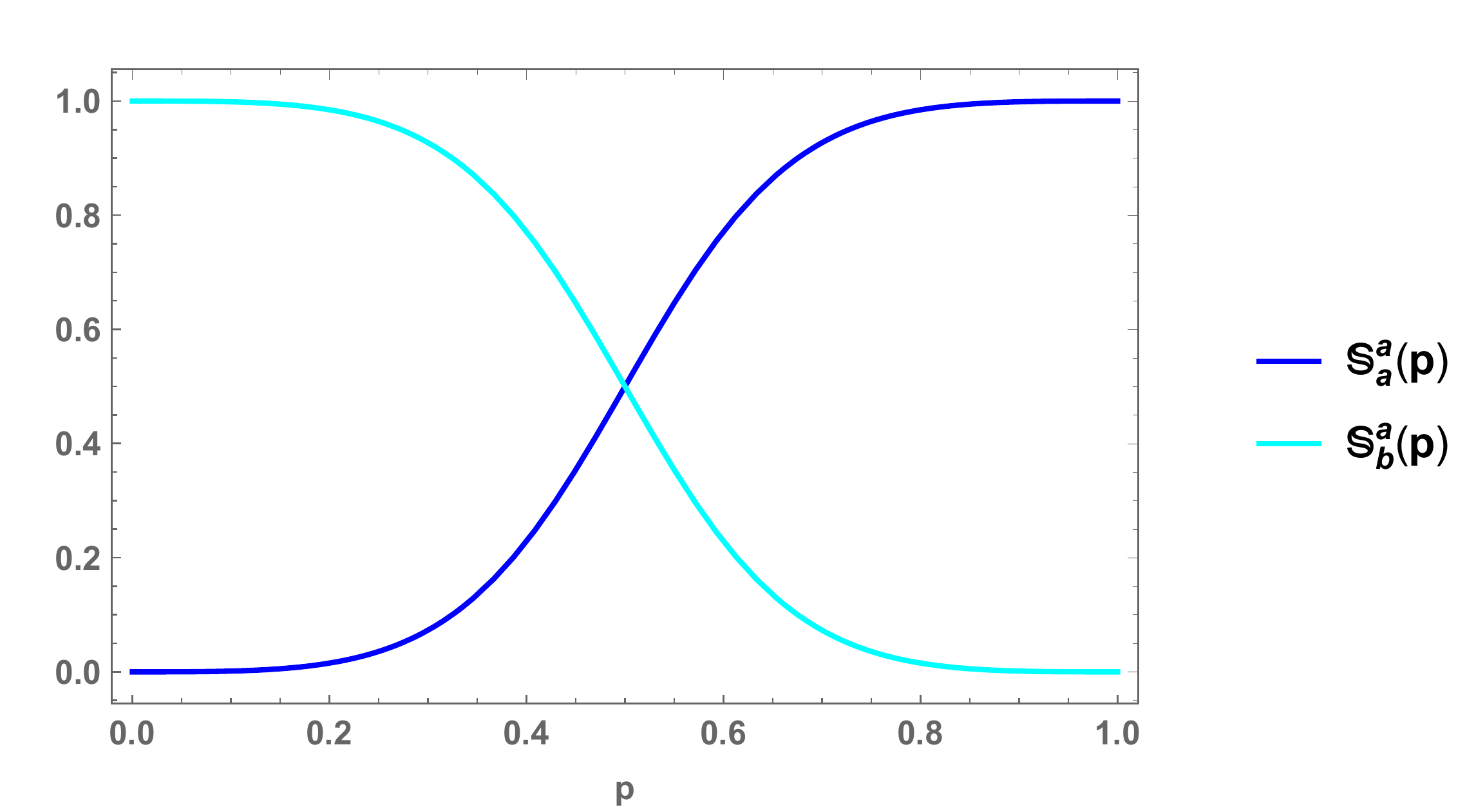}}
\subfloat[Dirac $N$: $\amsbb{S}_a^a$ and $\amsbb{S}_b^a$ as a function of $p$ \label{fig:1d}]{\includegraphics[width=3.5in]{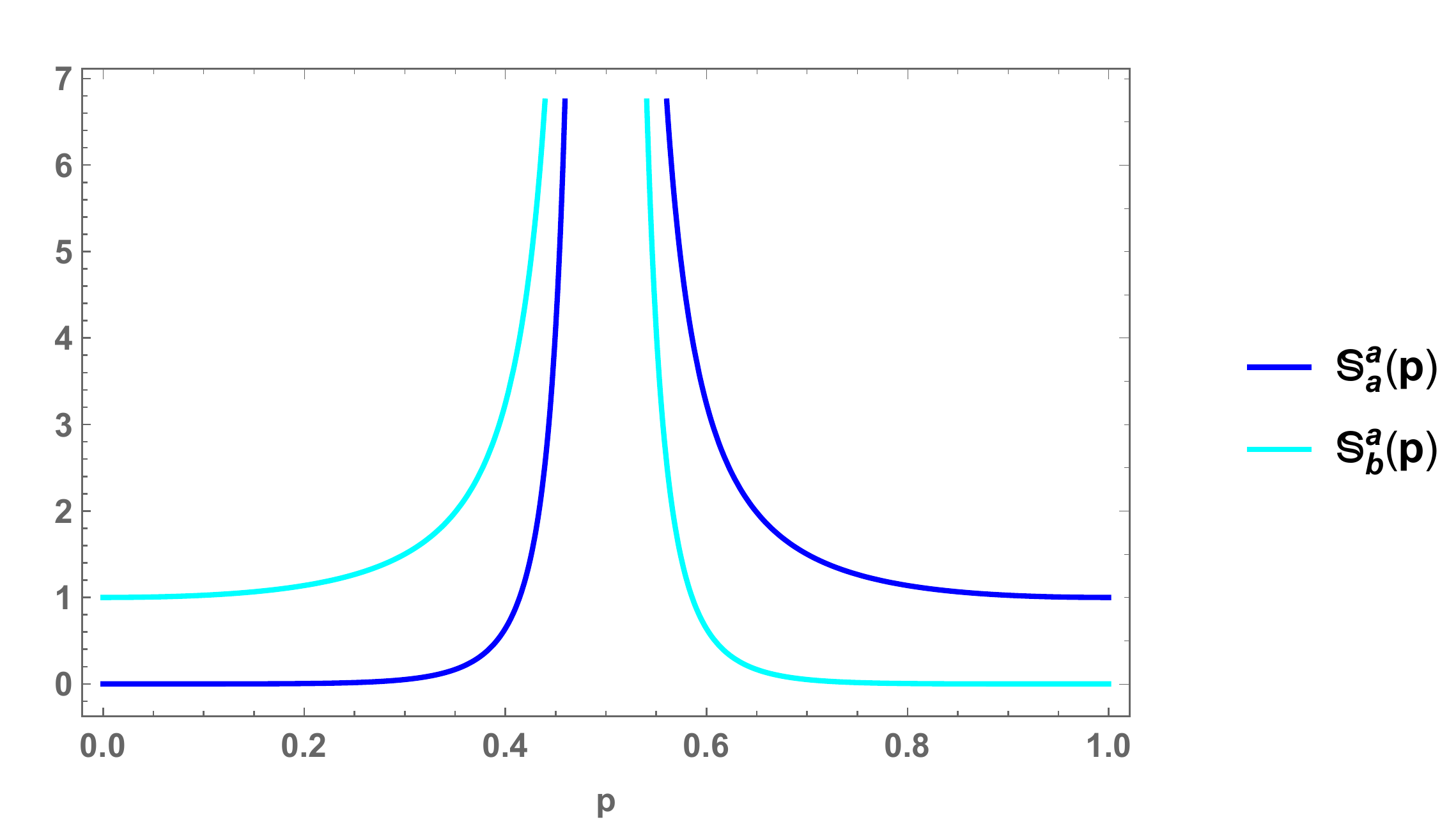}}\\
\subfloat[$H(\amsbb{S}_i^a)$ as a function of $p$\label{fig:1e}]{\includegraphics[width=3.5in]{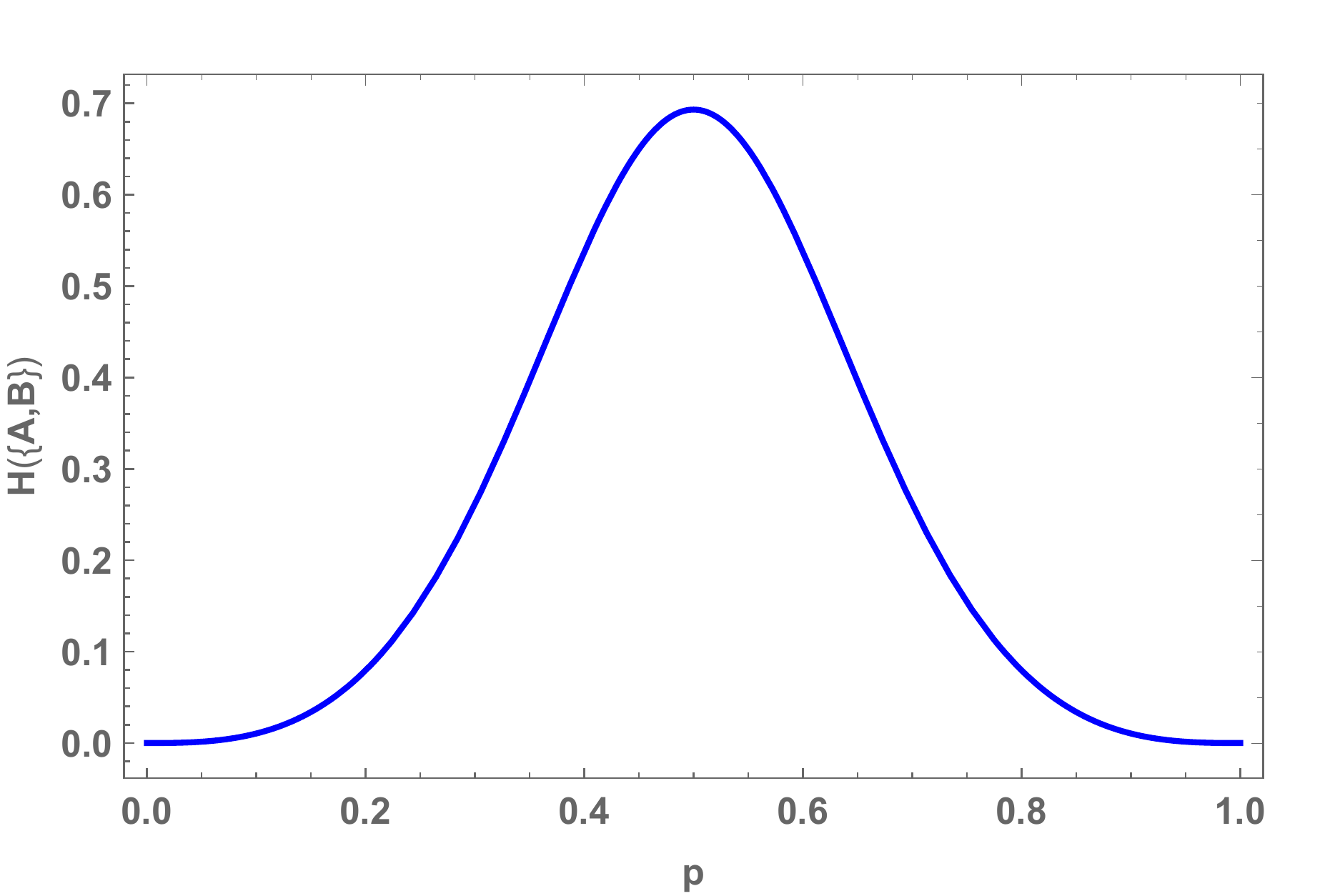}}\\
\endgroup
\caption{Random measure variance restrictions, sensitivity indices, and entropy for orthogonal and Dirac $N$, Bernoulli$(p)$, and univariate polynomial as a function of $p$}\label{fig:1}
\end{figure}

\paragraph{Specific calculations} For an orthogonal die or Poisson $\kappa$, we have $c=\delta^2$, so structural sensitivity indices are given by \begin{equation}\amsbb{S}^a_d = \frac{\nu f_d^2}{\nu f^2}\for D\in\{A,\dotsb,B\}\end{equation} where \begin{align*}\nu f_d^2 &= \nu f^2\ind{D}=\int_D\nu(\D x)(g(x)-\E g)^4\end{align*} Suppose $n=1$ so that $g(x)= x$ and let $\nu$ on $E=\{0,1\}$ be Bernoulli$(p)$ with mean $p\in[0,1]$. To define the Laplace functional, we have \[\nu e^{-f} = e^{-p^2}(1-p) + e^{-(1-p)^2}p\]Consider partition into the atoms $A=\{0\}$ and $B=\{1\}$. Then we have that \begin{align*}\nu f_a^2 &=p^4(1-p) \\\nu f_b^2&=  (1-p)^4p\\\nu f^2 &=p^4(1-p) + (1-p)^4p\\\amsbb{S}_a^a &= \frac{p^3}{3 p^2-3 p+1}\\\amsbb{S}_b^a &= \frac{(1-p)^3}{3 p^2-3 p+1}\\H(\amsbb{S}^a)&=-\amsbb{S}_a^a\log\amsbb{S}_a^a - \amsbb{S}^a_b\log\amsbb{S}^a_b\end{align*} These quantities are visualized in Figure~\ref{fig:1}. \emph{Figure~\ref{fig:1a} shows the random measure uncertainty is high for $x=0$ when $p$ is near one and similarly for $x=1$ when $p$ is near zero, whereas the second moment $\nu f^2$ shows interesting behavior of an inverted double-well (bistable) potential, with local maxima at $p=\frac{1}{6}(3\pm\sqrt{3})$ and local minima at $p=\frac{1}{2}$.} The potential is the uncertainty (second moment) of a point as a function of $p$ in the parameter space (unit interval). Moreover, the entropy of $N$ is shown in Figure~\ref{fig:1e}, which is symmetric and unimodal with maximum at $p=1/2$. Being the only non-trivial partition, $\{A,B\}$ is the maximum entropy partition of $E$. 

For Dirac we have as mentioned in \eqref{eq:binsa} and \eqref{eq:binsb} the sensitivity indices and variance terms \begin{align*}\nu f_a &=p^2(1-p)\\\nu f_b &=(1-p)^2p\\\Var f_a &= p^5(1-p)\\\Var f_b &= (1-p)^5p\\\Var f &= p(1-p)(1-2p)^2\\\amsbb{S}^a_a &=\frac{\Var f_a}{\Var f}=\frac{p^4}{(1-2p)^2}\\\amsbb{S}^a_b &=\frac{\Var f_b}{\Var f}=\frac{(1-p)^4}{(1-2p)^2}\\\amsbb{S}^b_a &= \amsbb{S}^b_b=-\frac{\nu f_a \nu f_{b}}{\Var f}=-\frac{p^2(1-p)^2}{(1-2p)^2}\\\end{align*} \emph{We see that the sensitivity indices for Dirac all have singularities at $p=1/2$, shown in Figure~\ref{fig:1d}.} The correlative terms are negative, so the structural terms do not form a probability vector but rather have positive defective mass. In Figure~\ref{fig:1b} we plot $\Var f_a$, $\Var f_b$, and $\Var f$ as a function of $p$. These are similar to $\nu f_a^2$, $\nu f_b^2$, and $\nu f^2$ although steeper and with the overall variance less than the sum of the variances of the two restrictions. In fact the overall variance is zero at $p=1/2$.

\FloatBarrier

\clearpage

\subsection{Ishigami function}\label{sec:ishi}

As referenced in Table~\ref{tab:0}, we consider the Ishigami function, a non-linear continuous function \[ g(x_1,x_2,x_3) = \sin x_1 + a \sin^2 x_2 + b x_3^4\sin x_1\for(x_1,x_2,x_3)\in[-\pi,\pi]^3\] with independent uniformly distributed coordinates and parameters $a,b\in\R_+$. The mean and variance are \begin{align*}\E g &= \frac{a}{2}\\\Var g &= \frac{1}{8}a^2+\frac{\pi ^8}{18}b^2+\frac{\pi ^4}{5}b+\frac{1}{2}\end{align*} We organize the HDMR results into Table~\ref{tab:2}.

\begin{table}[h!]
\begin{center}
\begin{tabular}{cccc}
\toprule
Subspace $u$ & $g_u$ & $\Var g_u$ & $\amsbb{S}_u^a$\\\midrule
$\{1\}$ & $(1+\frac{b\pi^4}{5})\sin x_1$ & $\frac{1}{50} \left(\pi ^4 b+5\right)^2$ & $\frac{36 \left(\pi ^4 b+5\right)^2}{5 \left(45 a^2+20 \pi ^8 b^2+72 \pi ^4 b+180\right)}$\\
$\{2\}$ & $-\frac{a}{2}\cos 2x_2$ & $\frac{a^2}{8}$ & $\frac{45 a^2}{45 a^2+20 \pi ^8 b^2+72 \pi ^4 b+180}$\\
$\{1,3\}$ & $b(x_3^4-\frac{\pi^4}{5})\sin x_1$ & $\frac{8 \pi ^8 b^2}{225}$ & $\frac{64 \pi ^8 b^2}{5 \left(45 a^2+20 \pi ^8 b^2+72 \pi ^4 b+180\right)}$\\\bottomrule
\end{tabular}\caption{HDMR of Ishigami function}\label{tab:2}
\end{center}
\end{table}


%
%

\FloatBarrier 

Let $f=(g-\E g)^2$. We have \begin{equation}\nu f =\Var g = \Var g_1+ \Var g_2 + \Var g_{13}\end{equation} For subspace $A\subset E$, we put $f_a=f\ind{A}$ and have \begin{equation}\nu f_a = \int_A\nu(\D x)f(x)\end{equation} Next we calculate $\nu f^2$, which is the fourth central moment with value \begin{equation}\nu f^2 = \frac{1}{24} \pi ^8 \left(a^2+6\right) b^2+\frac{3}{20} \pi ^4 \left(a^2+2\right) b+\frac{3}{128} \left(a^4+16 a^2+16\right)+\frac{3 \pi ^{16} b^4}{136}+\frac{3 \pi ^{12} b^3}{26}\end{equation} The restriction to subspace $A\subset E$ is similarly given as \begin{equation}\nu f_a^2 = \int_A\nu(\D x)f^2(x)\end{equation}

Consider the random measure $N=(\kappa,\nu)$ on $(E,\mathcal{E})$. To define the Laplace functional, we put \[\nu e^{-f} = \frac{1}{(2\pi)^3}\int_{[-\pi,\pi]^3}e^{-(\sin x_1 + a \sin^2 x_2 + b x_3^4\sin x_1-\frac{a}{2})^2}\D x_1\D x_2 \D x_3\] The Laplace transform is $F(\alpha)=\psi(\nu e^{-\alpha f})$. For $f\in\mathcal{E}_+$ the random variable $Nf$ has mean and variance \begin{align}\E Nf &= c\nu f =c(\Var g_1+\Var g_2+\Var g_{13})\\\Var Nf &=c\nu f^2 + (\delta^2-c)(\nu f)^2\end{align} 

The density of $Nf$, $\eta$, can be attained from evaluations of the Laplace transform $F$ using maximum entropy.  We take $\kappa=\text{Poisson}(100)$ and $\alpha_i\sim\text{Exponential}(1)$ for $i=1,\dotsb,n$ where $n=10$. We have $a=7$ and $b=0.1$. For the integral $\nu e^{-\alpha f}$, the variance of $g$ is large, which introduces numerical difficulties in computing the integral, so we define $f^*=f/\nu f$, compute $F^*(\alpha_i)$ for $i=1,\dotsb,n$, attain the maximum entropy distribution of $Nf^*\mapsto e^{-Nf^*}$ as $\mu_n$, and finally attain the maximum entropy distribution of $Nf$, $\eta_{n}$, as \[\eta_n(\D x) = \frac{1}{\nu f} e^{-x/\nu f}(\mu_n\circ e^{-x/\nu f})(\D x)\] In Figure~\ref{fig:laplacew} we show the density. In Table~\ref{tab:ishi}, we show the statistics compared to the exact values. The agreement is very close. 


\begin{table}[h!]
\begin{center}
\begin{tabular}{ccccc}
\toprule
$\E Nf$  &$\E\eta_{10}$ & $\Var Nf$ & $\Var\eta_{10}$\\\midrule
1\,384.5 & 1\,387.0 & 67\,223.4 & 66\,419.3\\\bottomrule
\end{tabular}\caption{Mean and variance of $Nf$ and $\eta_n$ where $\kappa=\text{Poisson}(100)$}\label{tab:ishi}
\end{center}
\end{table} 

\begin{figure}[h]
\centering
\includegraphics[width=4in]{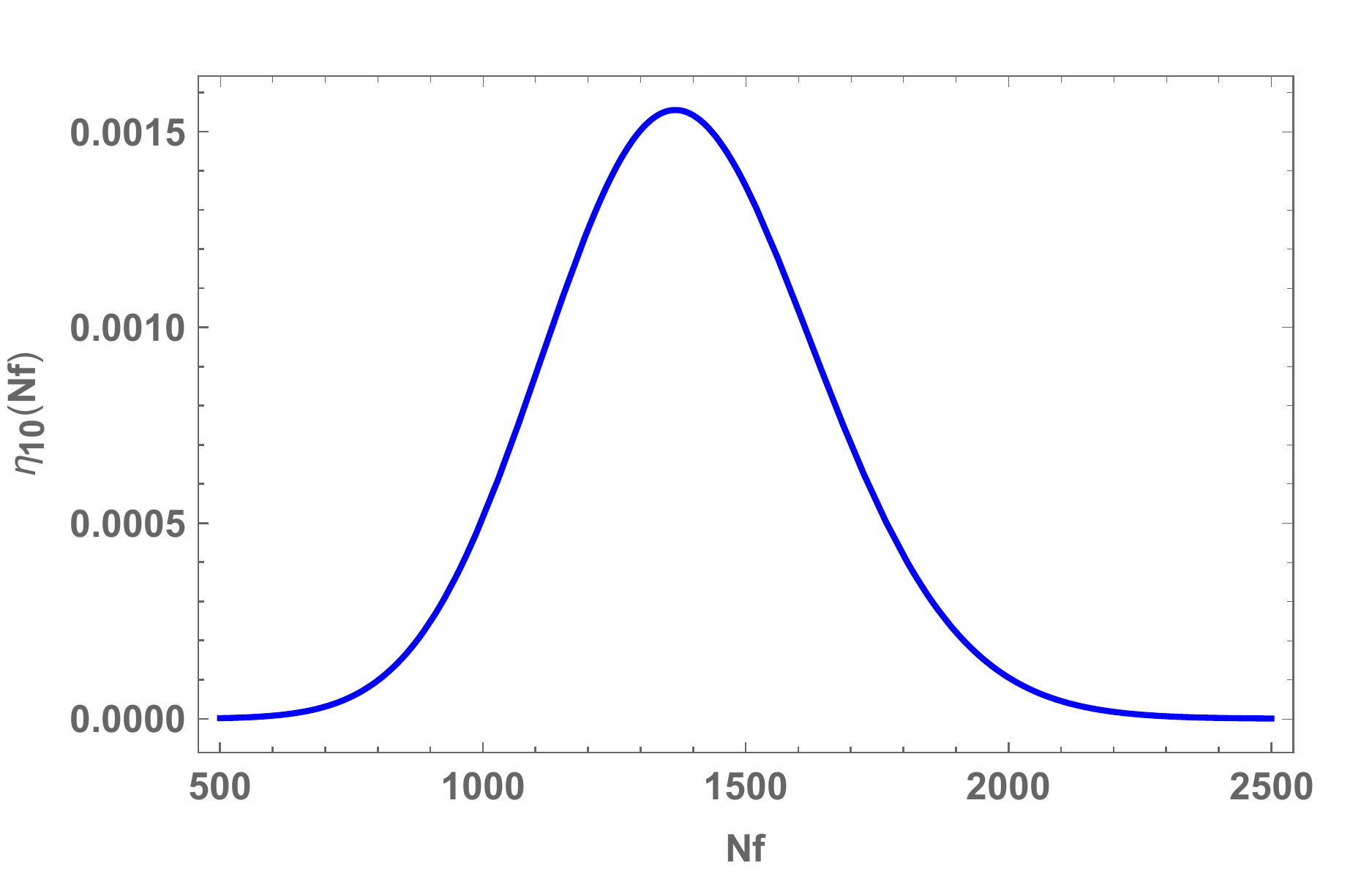}
\caption{Density of $Nf$, $\eta_n$ for $n=10$, for $\kappa=\text{Poisson}(100)$ }\label{fig:laplaceishi}
\end{figure}

\FloatBarrier

 For disjoint partition $\{A,\dotsb,B\}$ of $E$, we put $f_a=f\ind{A}, \dotsb, f_b=f\ind{B}$. Then we have \begin{align*}\Var Nf &= \sum_{D\in\{A,\dotsb,B\}} \Var Nf_d + \sum_{D_i\ne D_j\in\{A,\dotsb,B\}} \Cov(Nf_{d_i},Nf_{d_j})\\&= \sum_{D\in\{A,\dotsb,B\}}(c\nu f_d^2 + (\delta^2-c)(\nu f_d)^2) + \sum_{D_i\ne D_j\in\{A,\dotsb,B\}}(\delta^2-c)\nu f_{d_i}\nu f_{d_j} \end{align*} 

\begin{figure}[h!]
\centering
\begingroup
\captionsetup[subfigure]{width=3in,font=normalsize}
\subfloat[$(\nu f_d^2)$ for coordinate 1\label{fig:3a}]{\includegraphics[width=3.5in]{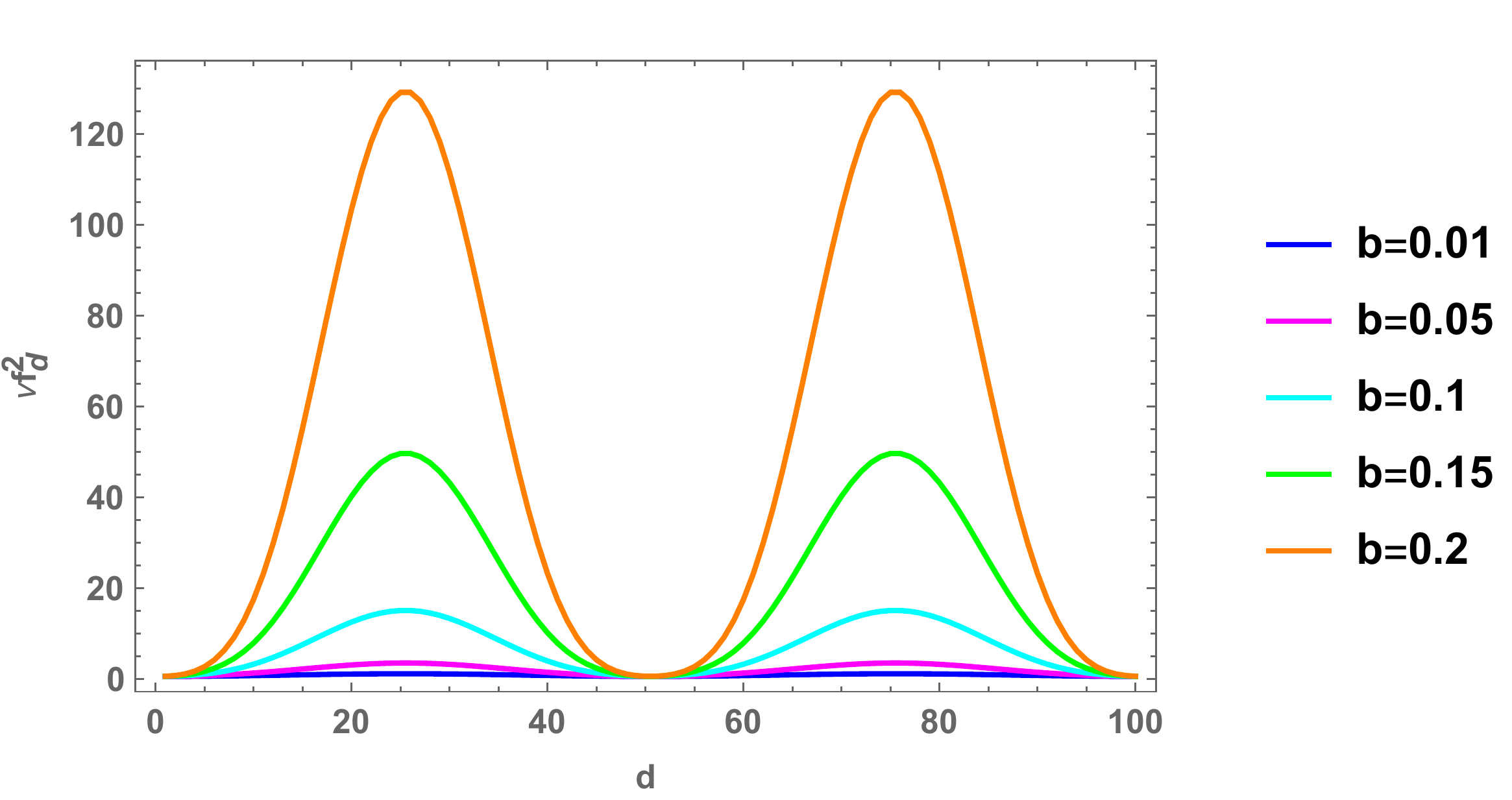}}
\subfloat[$(\amsbb{S}_d^a)$ for coordinate 1\label{fig:3b}]{\includegraphics[width=3.5in]{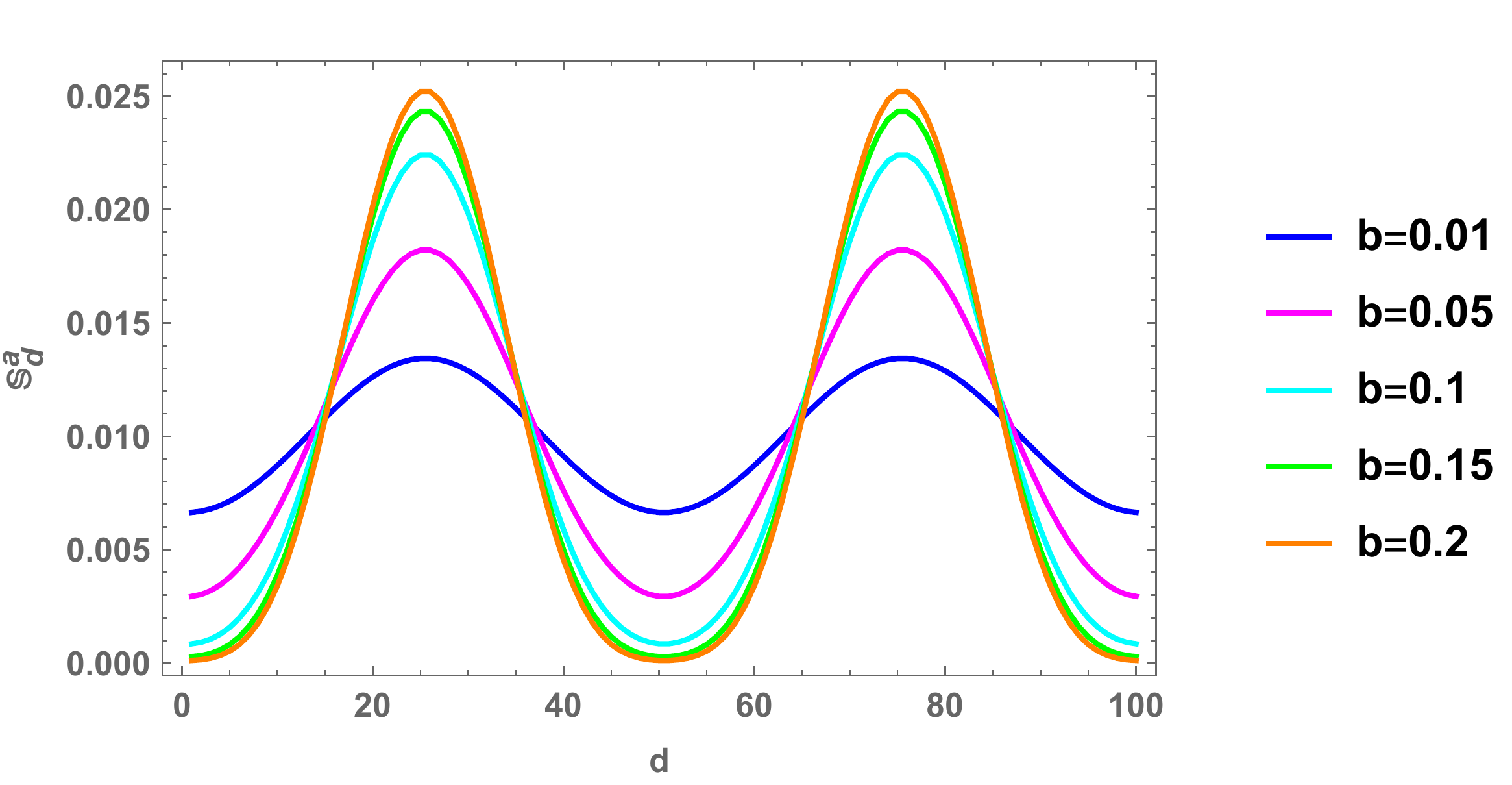}}\\
\subfloat[$(\nu f_d^2)$ for coordinate 2\label{fig:3c}]{\includegraphics[width=3.5in]{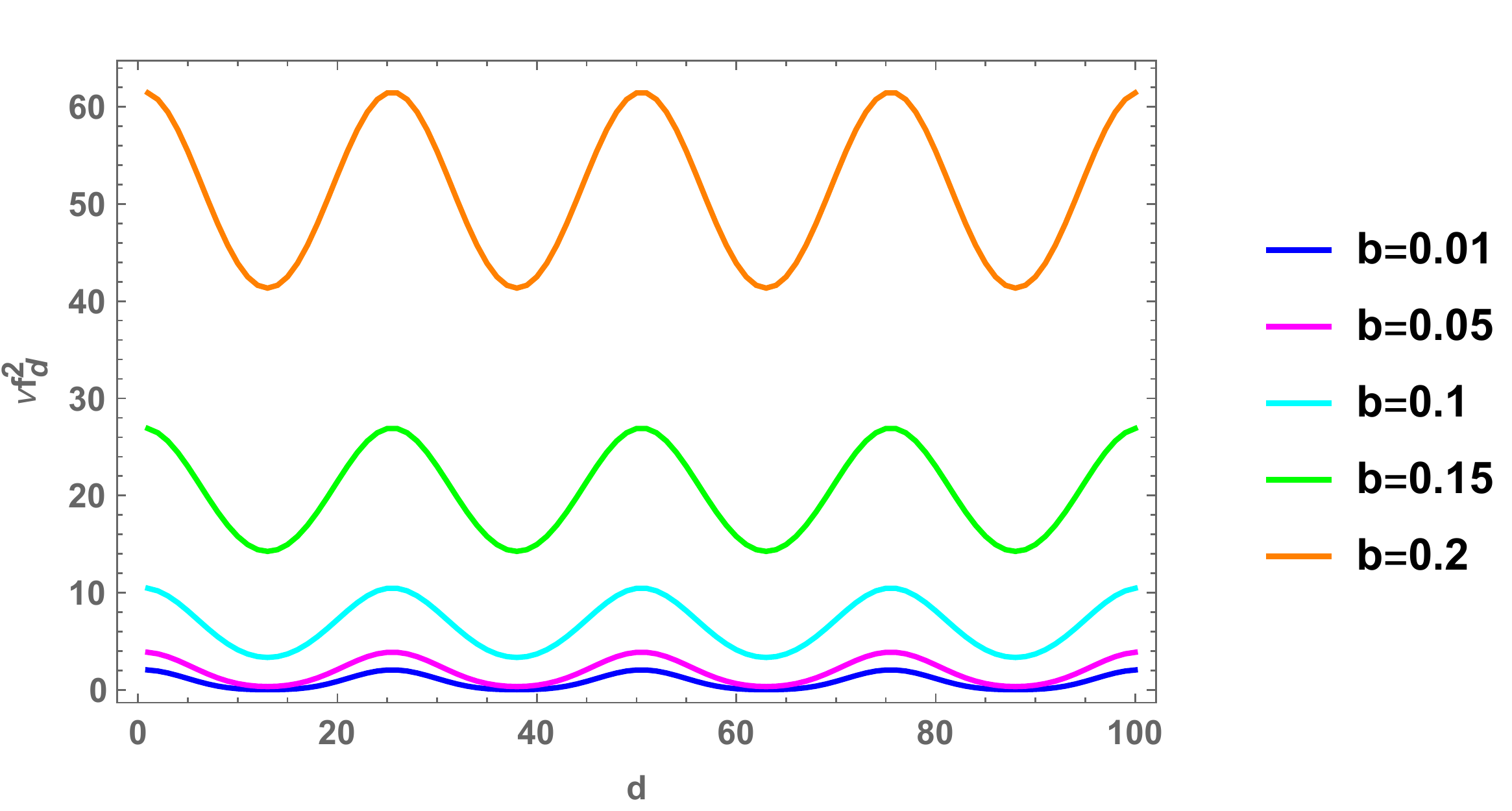}}
\subfloat[$(\amsbb{S}_d^a)$ for coordinate 2\label{fig:3d}]{\includegraphics[width=3.5in]{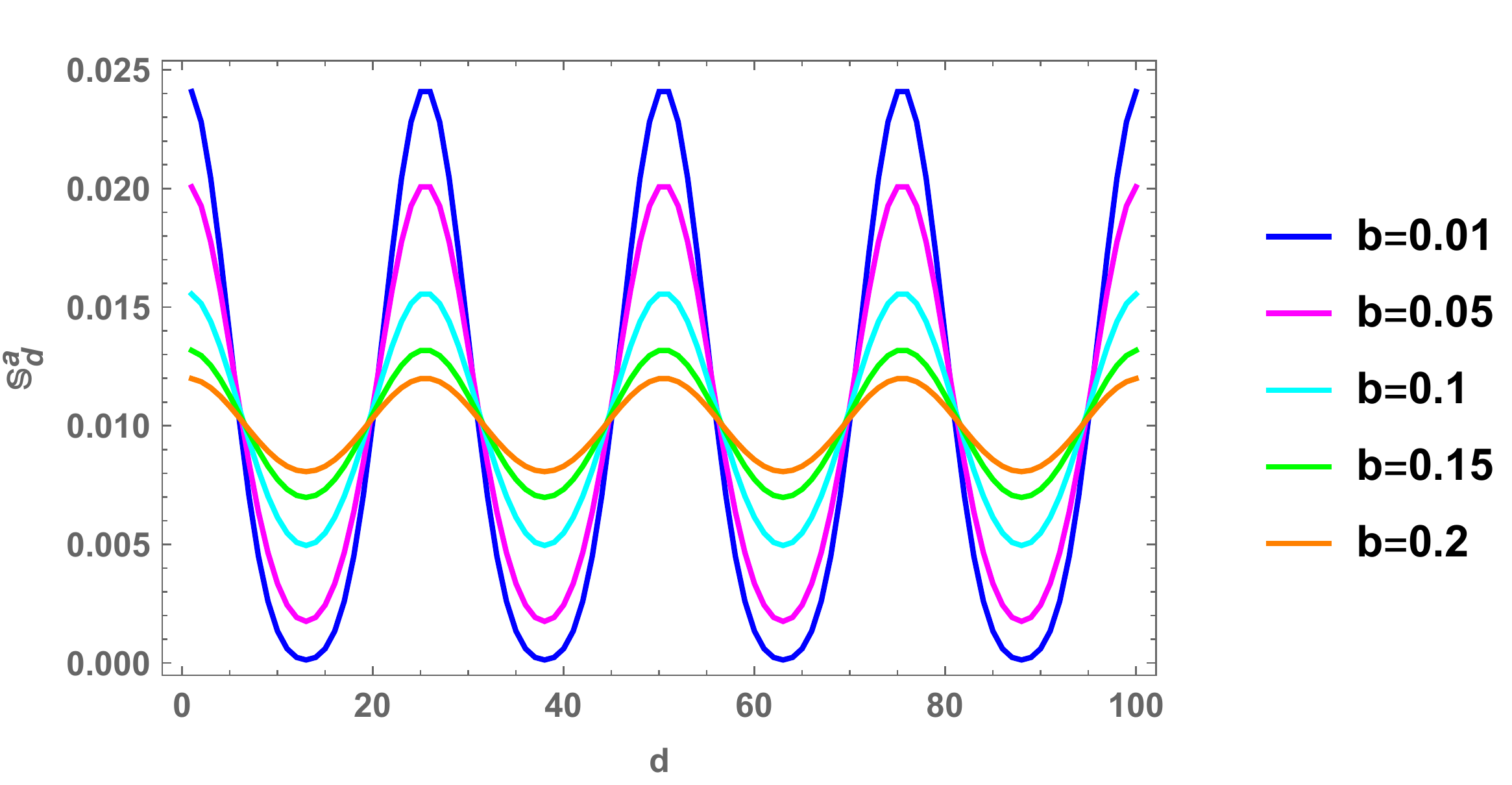}}\\
\subfloat[$(\nu f_d^2)$ for coordinate 3\label{fig:3e}]{\includegraphics[width=3.5in]{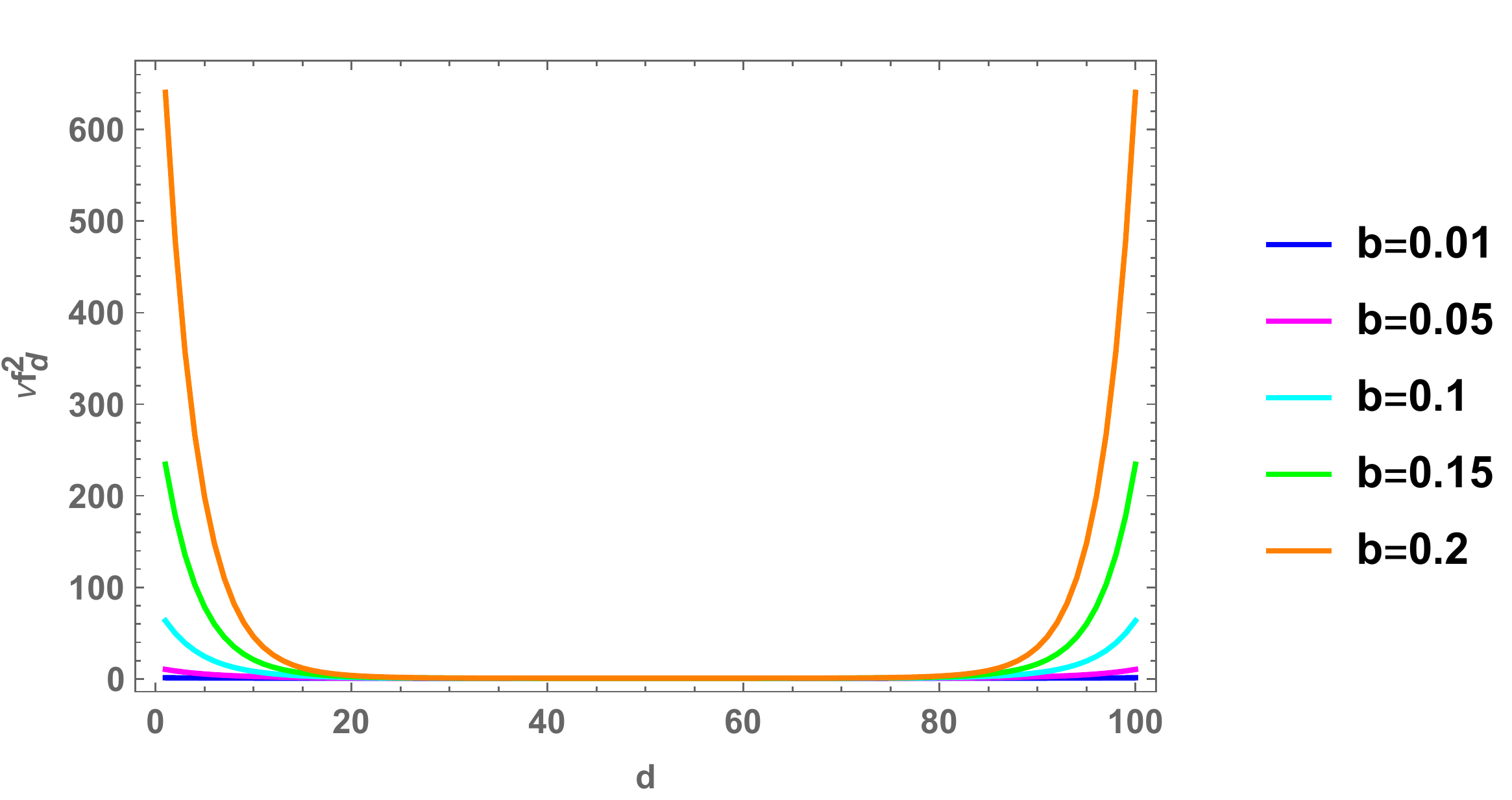}}
\subfloat[$(\amsbb{S}_d^a)$ for coordinate 3\label{fig:3f}]{\includegraphics[width=3.5in]{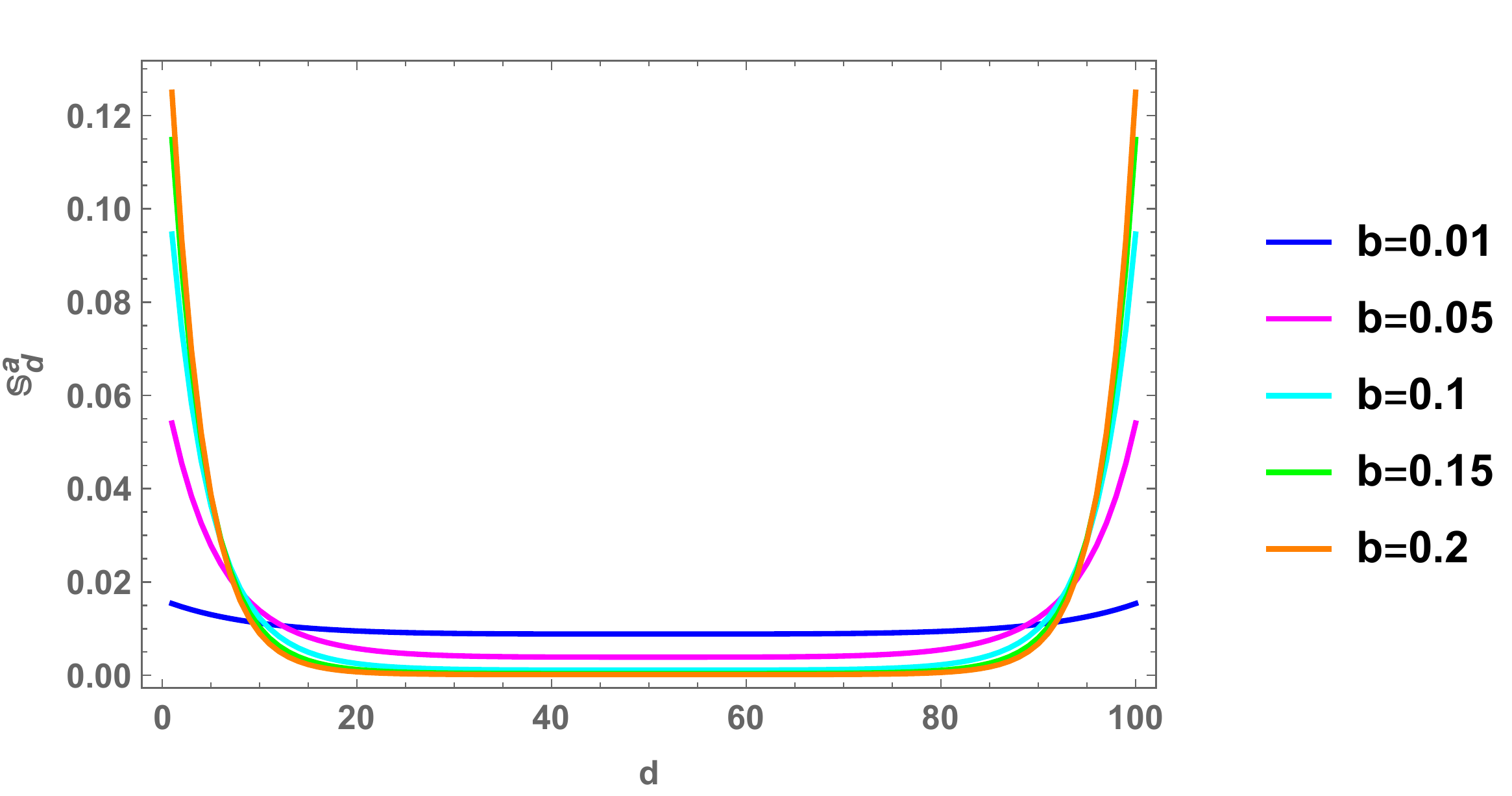}}\\
\endgroup
\caption{Ishigami function for orthogonal $N$: second moments $(\nu f_d^2)$ and sensitivity indices $(\amsbb{S}_d^a)$ for partitions by coordinates into $100$ intervals for orthogonal $N$ for $a=7$ and $b\in\{0.01,0.05,0.1,0.15,0.2\}$}\label{fig:3}
\end{figure}

For an orthogonal die or Poisson $\kappa$, we have $c=\delta^2$, so the structural sensitivity indices are given by \begin{equation}\amsbb{S}^a_d = \frac{\nu f_d^2}{\nu f^2}\for D\in\{A,\dotsb,B\}\end{equation} where \[\nu f_d^2 = \nu f^2\ind{D}=\int_D\nu(\D x)(g(x)-\E g)^4\] and its correlative are zero $\amsbb{S}^b_d=0$. As such, $(\amsbb{S}_d^a)$ forms the distribution $P$ on $\{A,\dotsc,B\}$. The sensitivity index density of \eqref{eq:sdxi} is given by \begin{align*}\xi&=2\pi\nu f^2=2 \pi  \left(\frac{3 a^4}{128}+\frac{1}{120} a^2 \left(5 \pi ^8 b^2+18 \pi ^4 b+45\right)+\frac{3 \pi ^{16} b^4}{136}+\frac{3 \pi ^{12} b^3}{26}+\frac{\pi ^8 b^2}{4}+\frac{3 \pi ^4 b}{10}+\frac{3}{8}\right)\\\amsbb{S}^a(\D x_1) &= \frac{1}{\xi}\left(\alpha + \beta \sin^2(x_1) + \gamma\sin^4(x_1)\right)\D x_1\\\alpha&=\frac{3 a^4}{128}\\\beta&=\frac{1}{60} a^2 \left(5 \pi ^8 b^2+18 \pi ^4 b+45\right)\\\gamma&=\frac{\pi ^{16} b^4}{17}+\frac{4 \pi ^{12} b^3}{13}+\frac{2 \pi ^8 b^2}{3}+\frac{4 \pi ^4 b}{5}+1\\\amsbb{S}^a(\D x_2)&=\frac{1}{\xi}\left(\alpha + \beta\cos^2(2x_2)+\gamma\cos^4(2x_2)\right)\D x_2\\\alpha&=\frac{3 \pi ^{16} b^4}{136}+\frac{3 \pi ^{12} b^3}{26}+\frac{\pi ^8 b^2}{4}+\frac{3 \pi ^4 b}{10}+\frac{3}{8}\\\beta&=\frac{1}{60} a^2 \left(5 \pi ^8 b^2+18 \pi ^4 b+45\right)\\\gamma&=\frac{a^4}{16}\\\amsbb{S}^a(\D x_3)&=\frac{1}{\xi}\left(\frac{3}{128} \left(a^4+16 \left(a b x_3^4+a\right)^2+16 \left(b x_3^4+1\right)^4\right)\right)\D x_3\end{align*} We set $a=7$ and vary $b\in\{0.01,0.05,0.1,0.15,0.2\}$. We form the following partitions by subdividing each coordinate interval into $m=100$ disjoint equisized intervals. This gives three partitions of size 100. For the first partition $\{A_i\}$, we take $A_i=[-\pi+2\pi i/m,-\pi+2\pi (i+1)/m]\times[-\pi,\pi]\times[-\pi,\pi]$ for $i=0,\dotsb,m-1$. In Figure~\ref{fig:3}, for each partition by a variable, we show the second moments and sensitivity indices across the intervals to reveal the uncertainties of the random measure estimator. The uncertainties are markedly different across the three coordinates: the uncertainties of the first two coordinates are periodic, with the second having higher frequency than the first, whereas the third coordinate has a `bathtub' appearance. In Figure~\ref{fig:4}, we show the first-order component functions and the HDMR-derived entropy as a function of $b$. We can see that the entropy has a maximum around $b\approx 0.13$.  We plot the entropy of the distribution defined by the random measure structural sensitivity indices as a function of $b$ for fixed $a=7$, corresponding to the ranges of $b$ used in Figure~\ref{fig:3}. All three exhibit monotonic behavior over the range of $b$, $[0.01,0.2]$: Coordinate 1 decreases, coordinate 2 increases, and coordinate 3 decreases.


\FloatBarrier

For Dirac $\kappa$, we have $\delta^2=0$ so that the structural index is \begin{equation}\amsbb{S}_d^a = \frac{\Var f_d}{\Var f}\for D\in\{A,\dotsb,B\}\end{equation} and the correlative index is \begin{equation}\amsbb{S}_d^b = \sum_{D_i\in\{A,\dotsb,B\}:D_i\ne D}-\frac{\nu f_{d}\nu f_{d_i}}{\Var f}\for D\in\{A,\dotsb,B\}\end{equation} We note that \begin{align*}\amsbb{S}_a&=\sum_{D\in\{A,\dotsb,B\}}\amsbb{S}_d^a > 1\\\amsbb{S}^b&=\sum_{D\in\{A,\dotsb,B\}}\amsbb{S}_d^b < 0\\\amsbb{S}^a+\amsbb{S}^b&=1\end{align*} 

For Dirac we see that the structural terms are defined in terms of the variances of the functions, whereas for the orthogonal random measures Poisson and orthogonal die the structural terms are defined in terms of the second moments. 

In Figure~\ref{fig:4a} we show the decomposition of the variances and structural sensitivity indices $(\amsbb{S}_d^a)$ across the partition for the coordinates for Dirac. The sensitivity indices exhibit distinct behavior in comparison to the sensitivity indices for orthogonal $N$ in Figure~\ref{fig:3}: the sensitivity indices of orthogonal $N$ are more stable, whereas the sensitivity indices of Dirac have more dramatic expression. The second moments and variances are very similar for both orthogonal and Dirac $N$, owing to the small contribution of the first moment of $f$ squared.

\subsection{Regressor}\label{sec:regress}As referenced in Table~\ref{tab:0}, here we have some dataset $\mathfrak{D}=(\mathfrak{X},\mathfrak{Y})=\{(x_i,y_i): i=1,\dotsc,n\}$ containing $n$ data points in $E\times F$ with induced distribution $\mu=\nu\times Q$ \[\mu f = \frac{1}{n}\sum_{i=1}^n f\circ(x_i,y_i)\for f\in(\mathcal{E}\otimes\mathcal{F})_+\] Let $(\mathbf{X},\mathbf{Y})=\{(X_i,Y_i)\}$ be the independency of values distributed according to $\mu=\nu\times Q$ which forms the random measure $M=(\kappa,\mu=\nu\times Q)$ on $(E\times F,\mathcal{E}\otimes\mathcal{F})$. We consider $g$ as in the data science model as a regressor and put $f(x,y)=(g(x)-y)^2\in(\mathcal{E}\otimes\mathcal{F})_+$. The Laplace functional of $M$ for $f$ is defined through \[\mu e^{-f} = \int_{E\times F}\mu(\D x,\D y) e^{-(g(x)-y)^2}\]

The HDMR of $g$ may be approximated with respect to the empirical random measure using tree functions, such as with Random Forest or the gradient boosting regressor. This gives an approximate MM-ANOVA decomposition of the intensity measure into subspaces. The ``true'' distribution is indicated by $\widetilde{\mu}=\widetilde{\nu}\times \widetilde{Q}$. In general $\widetilde{\nu}\ne\prod\widetilde{\nu}_i$, so the HDMR approximation must enforce hierarchical orthogonality with respect to the empirical measure $\nu$, such as through the QR decomposition. The sensitivity indices are $\{(\amsbb{S}_u^a,\amsbb{S}_u^b): u\subseteq\{1,\dotsb,n\}\}$. Examples of such constructions for the Ishigami function using tree models are found in Section 5 and the appendices of \cite{bastian2018}.

For the RM-ANOVA decomposition of the random measure variance into disjoint subspaces $\{A,\dotsb,B\}$, we identify the Voronoi cell partition through K-means clustering with some number of clusters and an appropriate metric, although any unsupervised clustering algorithm will do. The idea of K-means is to identify the partition that minimizes intracluster variances. The sensitivity indices are attained as $\{(\amsbb{S}_d^a,\amsbb{S}_d^b): D\in\{A,\dotsb,B\}\}$. Thus the sensitivity indices reveal the random measure uncertainty of the risk function in the various disjoint volume elements. 

Let $g_\theta$ be a Gaussian process with mean $U_\theta$ and covariance $K_\theta$, where $E=\R^n$ and $\mathcal{E}=\mathcal{B}_{\R^n}$, that is, $(g_\theta(x_1),\dotsb,g_\theta(x_m))$ is Gaussian for every $x_1,\dotsb,x_m\in E$ and $m\ge1$. The Gaussian process is trained under an assumption of independent mean-zero additive Gaussian noise with variance $\sigma^2$. For the point $x\in E$, we have mean \[\E g_\theta(x) =U_\theta(x)  = K_\theta(x,\mathfrak{X})(K_\theta(\mathfrak{X},\mathfrak{X})+\sigma^2 I)^{-1}\mathfrak{Y}\for x\in E\] and variance \[\Var g_\theta(x) = K_\theta(x,x)- K_\theta(x,\mathfrak{X})(K_\theta(\mathfrak{X},\mathfrak{X})+\sigma^2 I)^{-1}K_\theta(\mathfrak{X},x)\for x\in E\] Putting $\alpha = (K_\theta(\mathfrak{X},\mathfrak{X})+\sigma^2 I)^{-1}\mathfrak{Y}$, we have the representer theorem \[U_\theta(x) = \sum_{i}^n \alpha_i K_\theta(x,X_i)\] Now define $f_\theta\in(\mathcal{E}\otimes\mathcal{F})_+$ as $f_\theta(x,y)=(U_\theta(x)-y)^2$. The Laplace transform is \[F(\alpha)=\psi(\mu e^{-\alpha f_\theta})\]

 The random measure has mean and variance \begin{align*}\E Mf_\theta &= c\widetilde{\mu} f_\theta\\\Var Mf_\theta&=c\widetilde{\mu} f_\theta^2 + (\delta^2-c)(\widetilde{\mu} f_\theta)^2\end{align*}

\begin{figure}[h!]
\centering
\begingroup
\captionsetup[subfigure]{width=3in,font=normalsize}
\subfloat[HDMR Component function $g_1$\label{fig:3ab4}]{\includegraphics[width=3.5in]{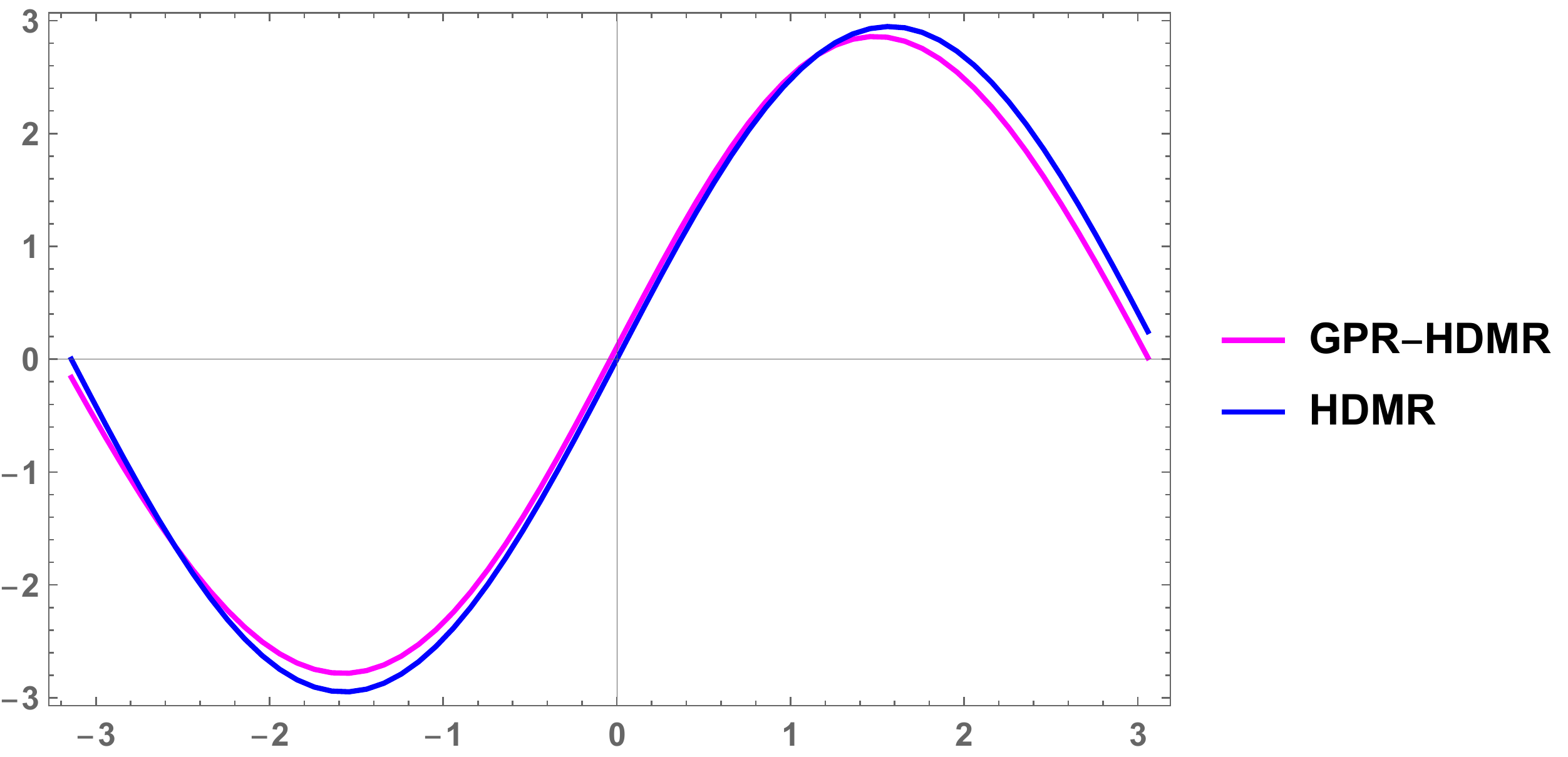}}
\subfloat[HDMR Component function $g_2$\label{fig:3ac2}]{\includegraphics[width=3.5in]{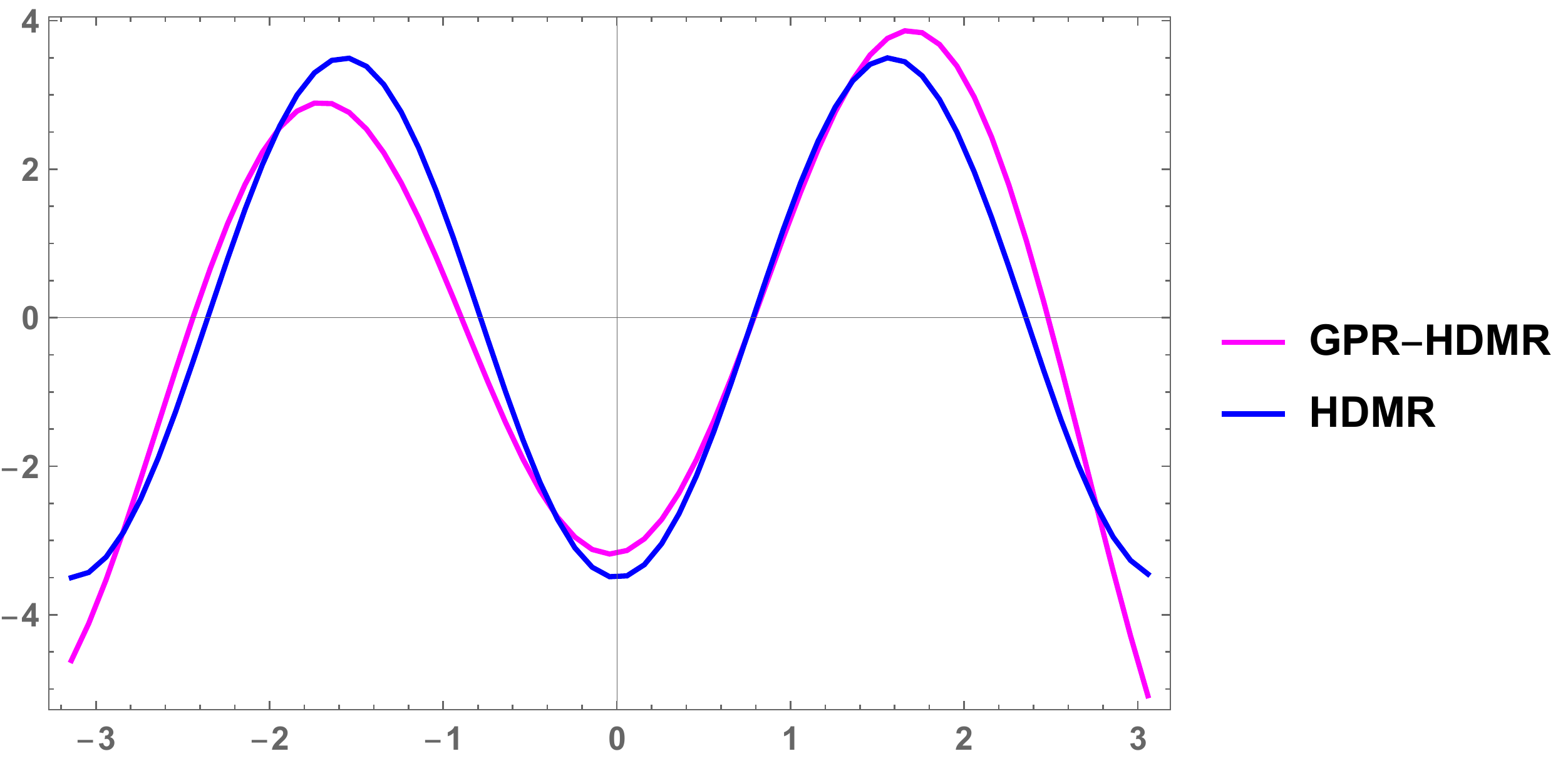}}\\
\subfloat[Orthogonal random measure sensitivity densities of $f=(g_\theta-\E g_\theta)^2$ for coordinate 1 \label{fig:3ad}]{\includegraphics[width=3.5in]{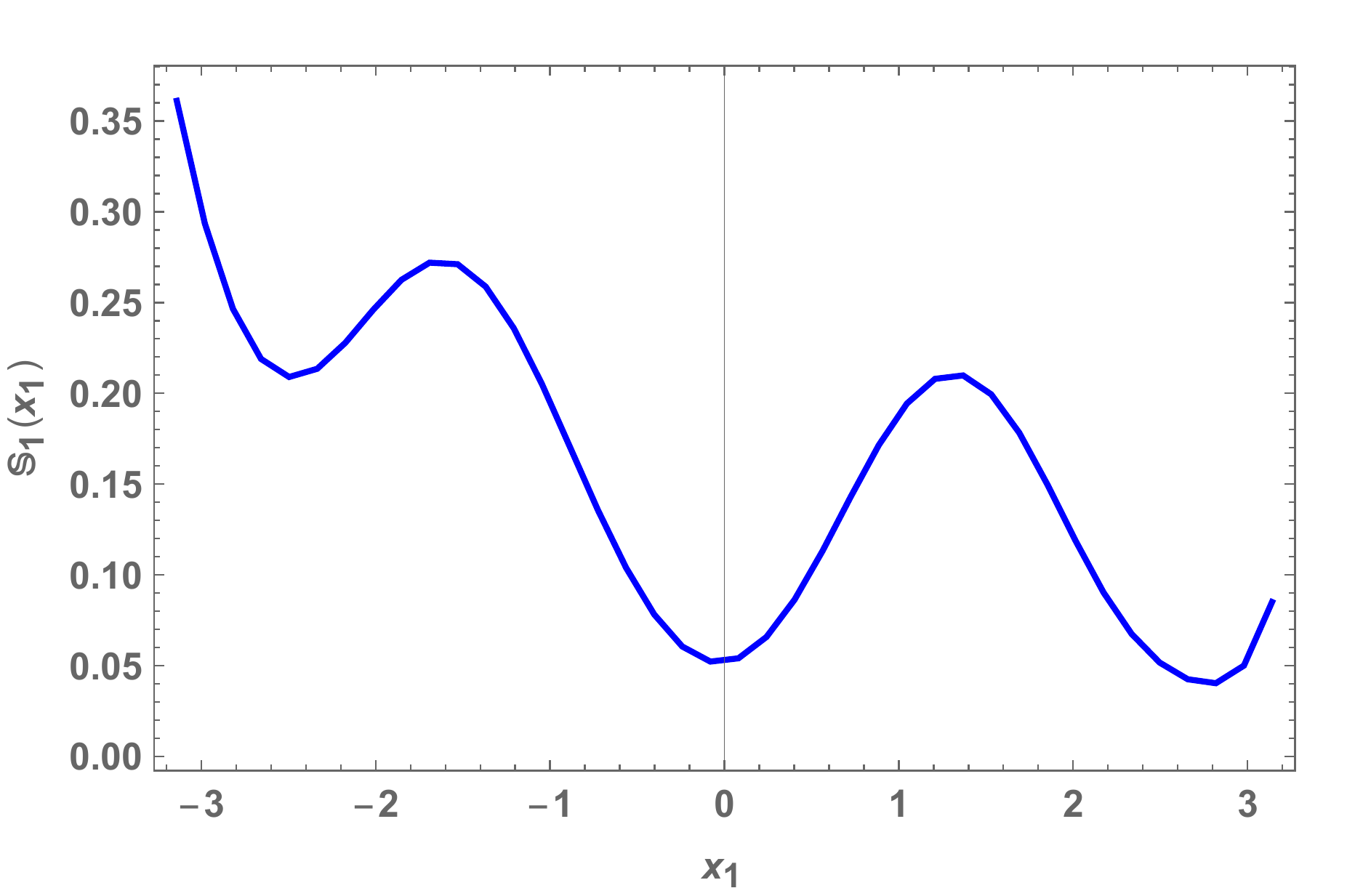}}
\subfloat[Orthogonal random measure sensitivity densities of $f=(g_\theta-\E g_\theta)^2$ for coordinate 2 \label{fig:3ae}]{\includegraphics[width=3.5in]{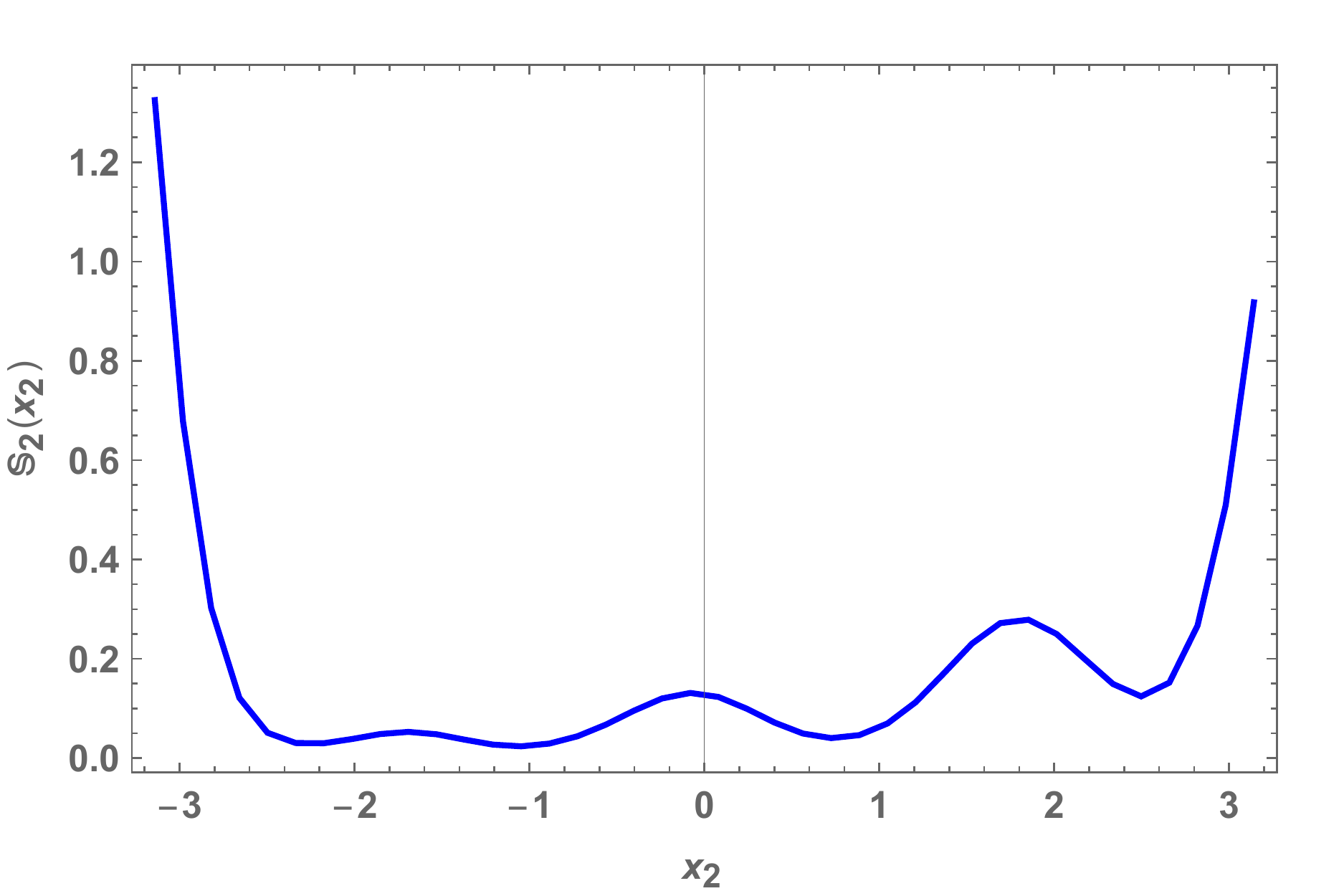}}\\
\subfloat[Orthogonal random measure sensitivity densities of $f=(g_\theta-\E g_\theta)^2$ for coordinate 3 \label{fig:3af}]{\includegraphics[width=3.5in]{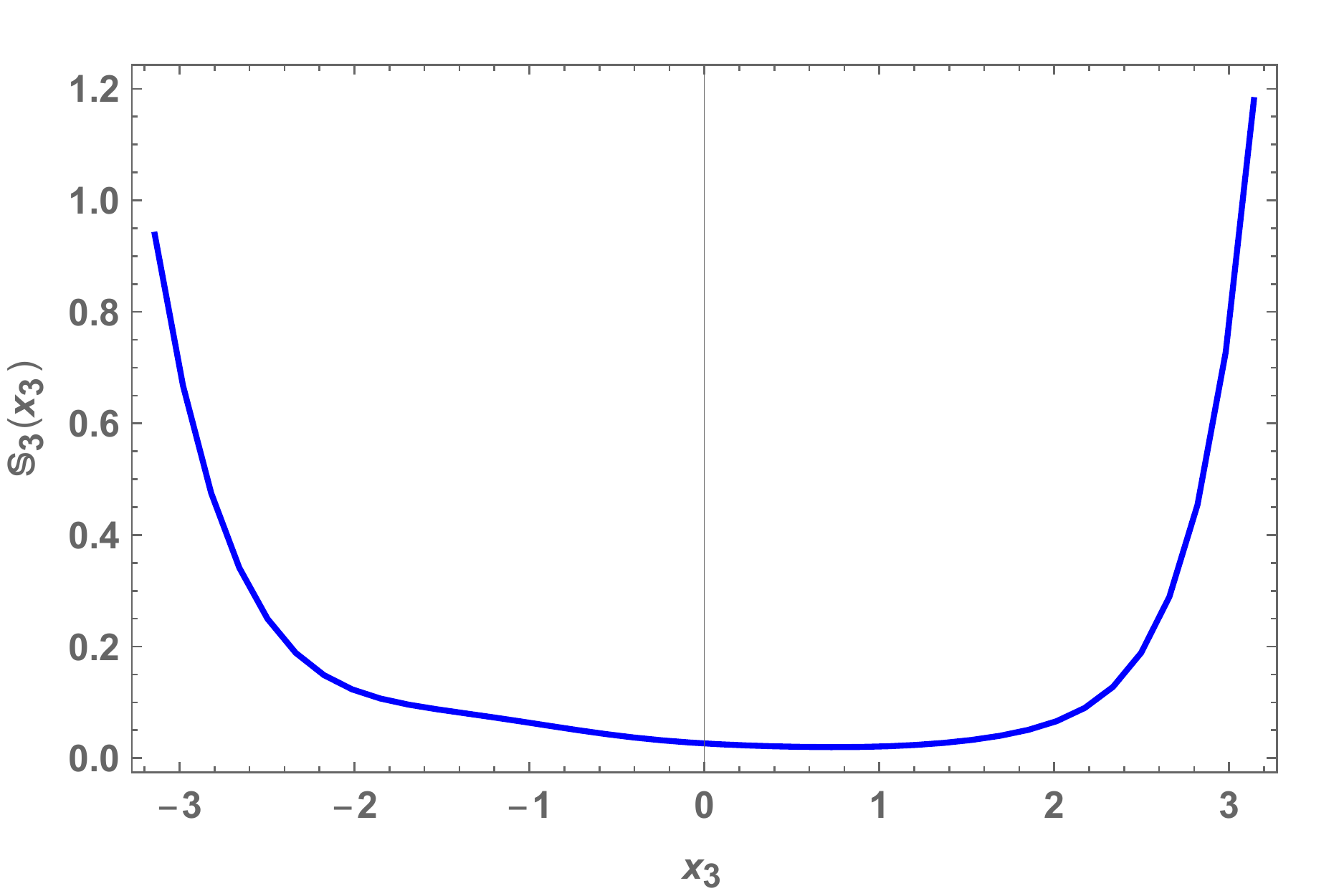}}
\subfloat[Orthogonal random measure sensitivity densities of MSE $f=(g_\theta-g)^2$ \label{fig:3aa5}]{\includegraphics[width=3.5in]{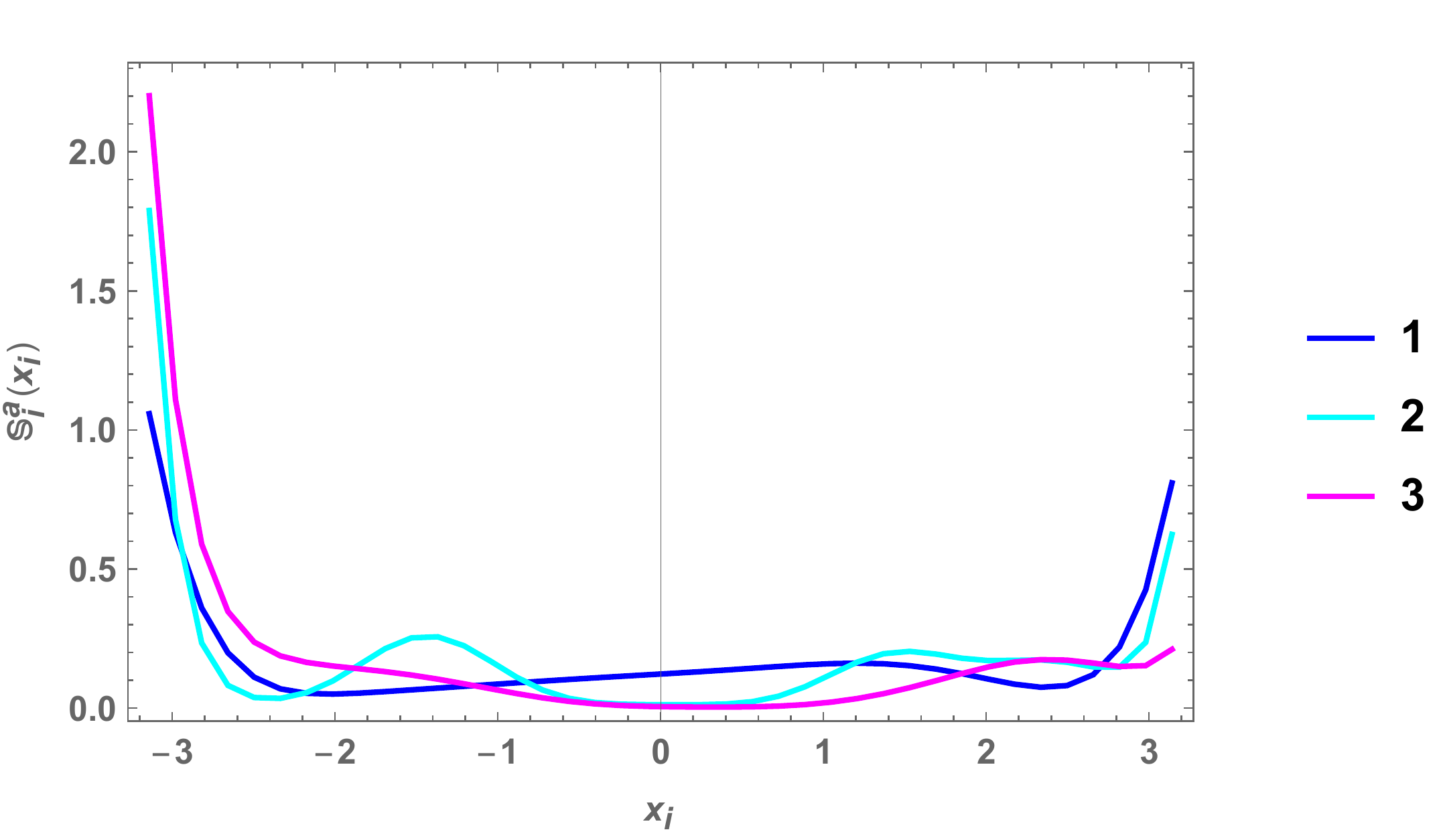}}\\
\endgroup
\caption{Gaussian process regressor of the Ishigami for $a=7$ and $b=0.1$ using radial basis kernel with $\gamma=0.1$ and $n=100$ samples}\label{fig:gp}
\end{figure}

For a specific calculation, consider the Ishigami function $g$ on $E=[-\pi,\pi]^3$. We construct a Gaussian process regressor with $n=100$ random input-output samples and noise $\sigma^2=0$ and the radial basis function with $\gamma=0.1$. In Table~\ref{tab:2} we show basic diagnostics of the model, including its mean-squared-error (MSE) and the coefficient of determination (COD) (each integrated with respect to the true density to eliminate estimation error). The regressor achieves approximately 72\% reconstruction fidelity. We set $Q(x,\cdot)=\delta_{g(x)}(\cdot)$. For orthogonal $M=(\kappa,\widetilde{\mu})$, we compute the first-order HDMR component functions, the sensitivity distributions $\amsbb{S}_1$, $\amsbb{S}_2$, and $\amsbb{S}_3$ for MSE and for variance and show these in Figure~\ref{fig:gp}. The sensitivity distributions are random, based on the random sample of size 100, and so for independent datasets $\mathfrak{D}$ the distributions are different. In Figure~\ref{fig:3aa5} the sensitivities of MSE of the random measure in the coordinates can be seen to coincide with the errors in the respective component function estimates in Figures~\ref{fig:3ab4} and \ref{fig:3ac2}. In Figures~\ref{fig:3ad}, \ref{fig:3ae}, and \ref{fig:3af} the random measure sensitivity indices of $f=(g_\theta-\E g_\theta)^2$ are shown in the coordinates. These are noisy versions of the true sensitivity curves seen in the earlier section on the Ishigami function.

\begin{table}[h!]
\begin{center}
\begin{tabular}{cccccc}
\toprule
$\E g$ & $\E g_\theta$ & $\Var g$ & $\Var g_\theta$ & $\widetilde{\nu}(g_\theta-g)^2$ (MSE) & $1-\frac{\widetilde{\nu}(g_\theta-g)^2}{\Var g}$ (COD)\\\midrule
3.5 & 3.8617 & 13.8446  &15.6533 & 3.8581 & 0.7213\\\bottomrule
\end{tabular}\caption{Gaussian process regressor performance measures for Ishigami function with $a=7$, $b=0.1$, $n=100$ samples, and zero additive error $\sigma^2=0$}\label{tab:2}
\end{center}
\end{table} 

\FloatBarrier

\subsection{Classifier}\label{sec:class} As referenced in Table~\ref{tab:0}, we suppose the same setup as the Regressor example in terms of a dataset, which forms the random measure $M=(\kappa,\mu=\nu\times Q)$ on $(E\times F,\mathcal{E}\otimes\mathcal{F})$, for a classifier. We assume $M$ is orthogonal. 

Define risk function $f(x,y)=\ind{}(y\ne g(x))$ for $(x,y)\in E\times F$. The Laplace functional of $f(x,y)=\ind{}(y\ne g(x))$ is defined through \[\mu e^{-f} = \int_{E\times F}\mu(\D x,\D y) e^{-\ind{}(y\ne g(x))}\] We are interested in $\mu f=\P(y\ne g(x))$ and $\mu f_d=\P(y\ne g(x)|(x,y)\in D)$ so that we may attain the structural sensitivity indices $(\amsbb{S}_d^a)$ for some partition $\{A,\dotsb,B\}$ of $E\times F$. We have $\mu f_a = \P(y\ne g(x)|(x,y)\in A), \dotsb, \mu f_b = \P(y\ne g(x)|(x,y)\in B)$.

Suppose we have a binary classifier $g$ where $F=\{0,1\}$. For partitions, consider a partition on the output space, $A=E\times\{0\}$ and $B=E\times\{1\}$ so \begin{align*}\mu f_a &= \P(y\ne g(x)|y=0)=\int_E\nu(\D x)Q(x,\{0\})f(x,0)\\\mu f_b &= \P(y\ne g(x)|y=1)=\int_E\nu(\D x)Q(x,\{1\})f(x,1)\\\mu f &= \P(y\ne g(x))=\mu f_a+\mu f_b\end{align*} Interpreting the class 1 as the positive, then the sensitivity indices are \begin{align*}\amsbb{S}^a_a &= \frac{\P(y\ne g(x)|y=0)}{\P(y\ne g(x))}=\frac{\E FP}{\E FP+\E FN}\\\amsbb{S}^a_b &= \frac{\P(y\ne g(x)|y=1)}{\P(y\ne g(x))}=\frac{\E FN}{\E FP+\E FN}\end{align*} where $\E FP$ and $\E FN$ denote the expected number of false positives and negatives. The orthogonal RM-ANOVA necessarily strictly measures the false positives and false negatives in assessing estimator uncertainty for the risk function $f(x,y)=\ind{}(y\ne g(x))$. Maximum entropy has $\E FP=\E FN$, i.e. balanced error across classes. 


To get the usual set of HDMR sensitivity indices, define $l=h\circ g$ per \eqref{eq:h}, giving for each $x\in E$ a binary vector $l(x)=(l_1(x),l_2(x))$. Let \[\nu l = (\P(g(x)=0),\P(g(x)=1))\] Then the probability for each $y\in F$ may be decomposed in terms of HDMR, giving the sensitivity indices $\{(\amsbb{S}_u^a(y),\amsbb{S}_u^b(y)): u\subseteq\{1,\dotsb,n\}\}$. This gives a decomposition of the image intensity measure in terms of a partition on the elements of $F$.

These expressions readily generalize to the multiclass setting $F=\{y_1,\dotsb,y_m\}$. Class $y_i\in F$ has mean \[\mu f_i = \P(y\ne g(x)|y=y_i)=\int_E\nu(\D x)Q(x,\{y_i\})f(x,y_i)\] and sensitivity index \[\amsbb{S}_i^a = \frac{\P(y\ne g(x)|y=y_i)}{\P(y\ne g(x))} = \frac{\E FP_i+\E FN_i}{\sum_{y_j\in F}(\E FP_j + \E FN_j)}\] where $\E FP_i$ and $\E FN_i$ are the expected number of false positives and negatives. Maximum entropy has $\E FP_i + \E FN_i=\E FP_j + \E FN_j$ for $y_i,y_j\in F$, i.e. balanced errors across the classes. Putting $\E F_1 = \E FP_1+\E FN_1$, the maximum entropy manifold of $(\E FP_i, \E FN_i)$ of any class $i$ is $\{(x,y)\in\R_+^2: x+y=\E F_1\}$.

\subsection{Space domain awareness of interacting particle systems}\label{sec:sda}\epigraph{\emph{Astronomy compels the soul to look upwards and leads us from this world to another.}}{Plato}

As referenced in Table~\ref{tab:0}, here we describe a mathematical model for space domain awareness based on random counting measures and random fields, where the points of the system are space objects, such as spacecraft, asteroids, and so on. The random field model is built using the random measure. We describe a random measure as a superposition of three (spatio-temporal) counting processes: Poisson, binomial, and negative binomial. Together, these processes can exhibit any covariance structure. We show that this modeling framework may be used to represent, for instance, spatiotemporal interaction probabilities through the random field model, provisioning an interacting particle system model.

\subsubsection{Superposition random measures} Consider random measure $N=(\kappa,\nu)$ on $(E,\mathcal{E})$. We think of the state-space as a function space $(E,\mathcal{E})=(F,\mathcal{F})^{\R_+}$. Here we take element $X_i\sim\nu$ as a space object, such as a spacecraft or asteroid, and $F=\R^3$. For the counting distribution we consider the convolution $\kappa=\kappa_-*\kappa_0*\kappa_+$, where $\kappa_-$ is binomial, $\kappa_0$ is Poisson, and $\kappa_+$ is negative binomial, i.e. $K=K_-+K_0+K_+$, with mean $c=c_-+c_0+c_+$ and variance $\delta^2 = \delta_-^2 + \delta_0^2 + \delta_+^2$. The pgf of $K$ is a product of the pgfs $\psi_\theta = \psi_{\theta_-}\psi_{\theta_0}\psi_{\theta_+}$ where we put $\theta=(\theta_-,\theta_0,\theta_+)$. The random measure $N$ is formed through STC as \[N(A) = \sum_i^{K_-+K_0+K_+}\ind{A}(X_i) \for A\in\mathcal{E}\] The Laplace functional of $N$ is given by \[L(f) = \psi_\theta(\nu e^{-f})\for f\in\mathcal{E}_+\] The random measure $N$ is orthogonal if $c=\delta^2$. This is equivalent to requiring $(\delta^2_--c_-)+(\delta^2_+-c_+)=0$ or $\delta^2_-+\delta^2_+=c_-+c_+$. The sensitivity density of $Nf$ for $f\in\mathcal{E}_+$ and orthogonal $N$ is \[\amsbb{S}(\D x) = \frac{\nu(\D x)f^2(x)}{\nu f^2}\for x\in E\]

\subsubsection{Image random measures} Consider the mapping $h:E\mapsto F$ as $h(w)=w(t)$ for $w\in E$ and $t\in\R_+$ fixed. We consider the image measure $M_t=N\circ h^{-1}=(\kappa,\mu_t=\nu\circ h^{-1})$ on $(F,\mathcal{F})$, where \[\mu_t f = (\nu\circ h^{-1})f = \nu(f\circ h)\for f\in\mathcal{F}_+\] with Laplace functional \[L(f) = \psi_\theta(\mu_t e^{-f})\for f\in\mathcal{F}_+\] formed by STC as \[M_tf = \sum_i^Kf\circ X_i(t)\for f\in\mathcal{F}_+\] The sensitivity density of $M_tf$ for $f\in\mathcal{F}_+$ and orthogonal $M_t$ is \[\amsbb{S}(\D x) = \frac{\mu_t(\D x)f^2(x)}{\mu_t f^2}\for x\in F\] 


\subsubsection{Restricted random measures} Now consider the restriction to subspace $A\subset F$, such as some domain of interest. For the restricted random measure $M_A=(M\ind{A},\mu_A)$ on $(A\cap F, \mathcal{F}_A)$, where $\mu_A(\cdot)=\mu(A\cap\cdot)/\mu(A)$ and $\mathcal{F}_A = \{A\cap B: B\in\mathcal{F}\}$, the pgf of the counting distribution is given by \begin{align*}\psi_A(t)&=\psi_\theta(at+1-a)\\&=\psi_{\theta_-}(at+1-a)\psi_{\theta_0}(at+1-a)\psi_{\theta_+}(at+1-a)\\&=\psi_{h_a(\theta_-)}(t)\psi_{h_a(\theta_0)}(t)\psi_{h_a(\theta_+)}(t)\\&=\psi_{h_a(\theta)}(t)\end{align*} where $h_a(\theta)$ is the vector of bone mappings $(h_a(\theta_-),h_a(\theta_0),h_a(\theta_+))$. The Laplace functional is \[L_A(f) = \psi_{h_a(\theta)}(\mu_A e^{-f})\for f\in\mathcal{F}_+\] The sensitivity density of $M_Af$ for $f\in\mathcal{F}_+$ and orthogonal $M_A$ is \[\amsbb{S}_a(\D x) = \frac{\mu_A(\D x)f^2(x)}{\mu_A f^2}\for x\in A\]

\subsubsection{Interacting random fields} Consider the image random measure $M_t=(\kappa,\mu_t)$ on $(F,\mathcal{F})$. Now let $k:F\times F\mapsto\R_+$ be $\mathcal{F}\otimes\mathcal{F}$ measurable and define the random field \[G_t(y) = \int_FM_t(\D x)k(x,y)\for y\in F\]  with Laplace transform \[\E e^{-\alpha G_t} = \psi_\theta(\mu_t e^{-\int_F\alpha(\D y)k(\cdot,y)})\quad\text{for every finite measure }\alpha\text{ on }(F,\mathcal{F})\] Put $f_y(\cdot)=k(\cdot,y)\in\mathcal{F}_+$, so that we have the mean and covariance \begin{align*}U_t(y)&=\E G_t(y)=c\mu_t f_y\for y\in F\\C_t(y,z)&=\Cov(G_t(y),G_t(z))=c\mu_t(f_yf_z)+(\delta^2-c)\mu_t f_y \mu_t f_z\for y,z\in F\end{align*} For example let $k$ be the radial basis function $k(y,z)=e^{-\gamma\norm{y-z}^2}$, which we interpret as an interaction probability as a function of distance. Then $G_t(y)$ has the interpretation of the number of interactions of space objects at some point $y\in F$ at time $t\in\R_+$ with mean number of interactions $U_t(y)$ and across locations $z\in F$ interaction number covariance $C_t(y,z)$. 



\subsubsection{Specific calculation} Let $\nu$ be Wiener so that the image measure is given by \[\mu_t(A) = \int_\R\D x\frac{1}{\sqrt{2\pi t}}e^{-\frac{x^2}{2t}}\ind{A}(x)\for A\in\mathcal{B}_{\R}\] This is a prototypical disequilibrium Gaussian process. 

Suppose each particle pays ``rent'' based on its location according to $g=\cos^2$. The random variable $M_tg$ is total amount of rent paid by all the particles and has Laplace transform \[F(\alpha) = \psi_\theta(\mu_t e^{-\alpha g})\for \alpha\in\R_+\] For simplicity, we assume for the moment that all the particles belong to the Poisson random measure. 

Then putting $f=(g-\mu_t g)^2\in\mathcal{E}_+$ and noting that $g^{-1}(y) = \{\pm\arccos(\pm\sqrt{y})+2\pi k: k\in\Z\}$ and $\frac{\D}{\D y}g^{-1}(y) = \frac{1}{2\sqrt{y(1-y)}}$, we have \begin{align*}\mu_t\circ g^{-1}(\D y) &= \sum_{k=-\infty}^\infty \D y\frac{1}{2\sqrt{2\pi t y(1-y)}}\sum_{\pm}e^{-\frac{(\pm\arccos(\pm\sqrt{y})+2\pi k)^2}{2t}}\\\mu_t g &= e^{-t} \cosh (t)\xrightarrow{t\rightarrow\infty}\frac{1}{2}\\\mu_t f&= \frac{1}{8} e^{-8 t} \left(e^{4 t}-1\right)^2\xrightarrow{t\rightarrow\infty}\frac{1}{8}\\\mu_t f^2 &= \frac{1}{4} e^{-16 t} \sinh ^4(2 t) (4 \sinh (4 t)+\sinh (8 t)+8 \cosh (4 t)+2 \cosh (8 t)+5)\xrightarrow{t\rightarrow\infty}\frac{3}{128}\\\amsbb{S}(\D x)&=\frac{\mu_t(\D x)f^2(x)}{\mu_t f^2}\end{align*} Hence the statistics all admit stationary solutions. 

\begin{figure}[h]
\centering
\includegraphics[width=4in]{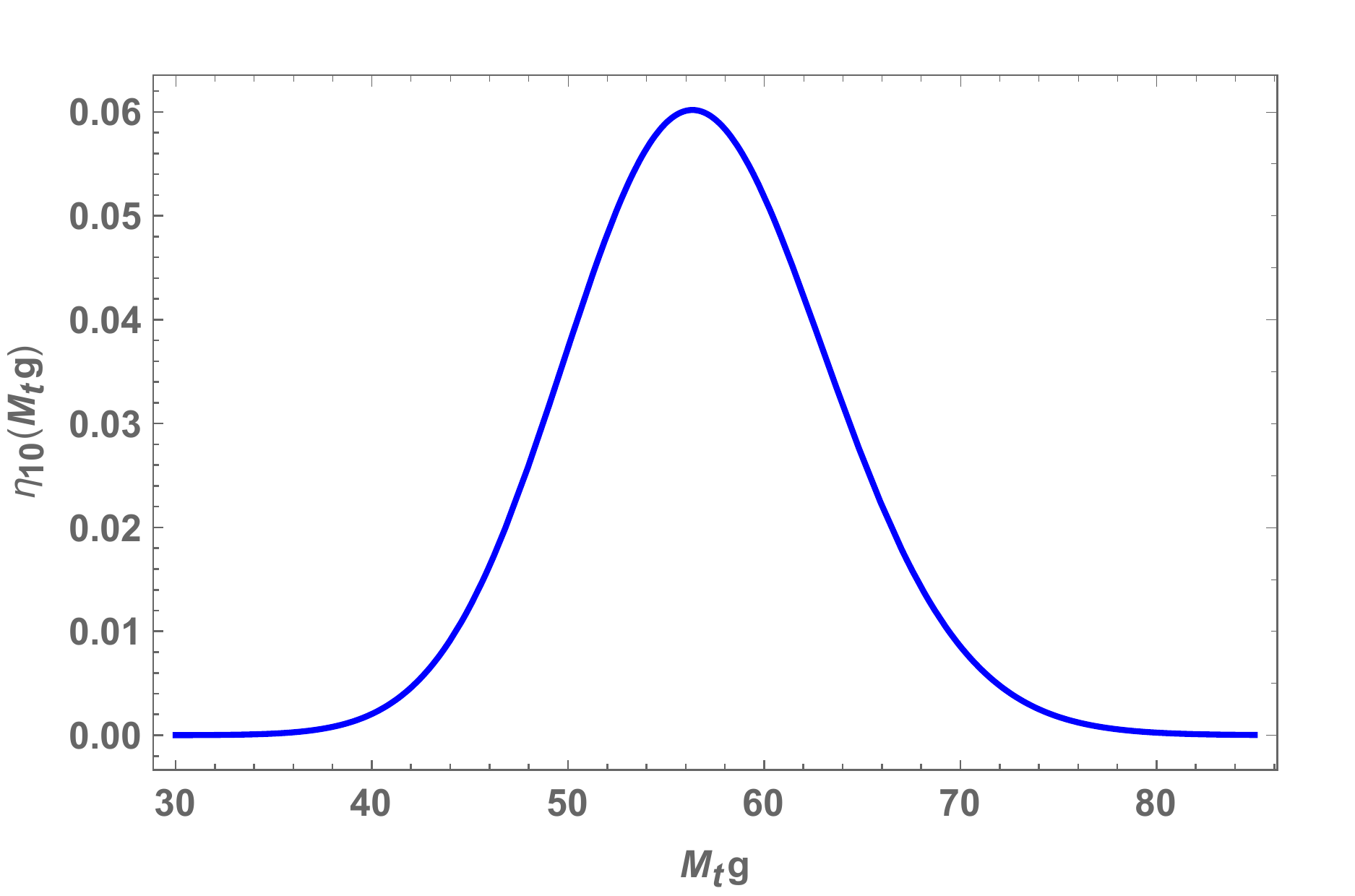}
\caption{Maximum entropy distribution $\eta_n$ for $n=10$ for $M_tg$}\label{fig:laplacew}
\end{figure}

We identify the density of $M_tg$, denoted by $\eta_{n}$, where $\kappa=\text{Poisson}(100)$, by evaluating $F(\alpha_i)$ for $i=1,\dotsb,n$ with $\alpha_i\sim\text{Exponential}(1)$ and attaining the density by maximum entropy. For $n=10$ we show the density in Figure~\ref{fig:laplacew}. We show the statistics compared to the exact values. The agreement is very close. 

\begin{table}[h!]
\begin{center}
\begin{tabular}{ccccc}
\toprule
$\E M_tg$  &$\E\eta_{10}$ & $\Var M_tg$ & $\Var\eta_{10}$\\\midrule
56.7668 & 56.7667 & 44.2710 & 44.2624\\\bottomrule
\end{tabular}\caption{Mean and variance of $M_tg$ and $\eta_n$ where $\kappa=\text{Poisson}(100)$}\label{tab:3}
\end{center}
\end{table} 

\FloatBarrier

We take $k:E\times E\mapsto\R_+$ as the radial basis kernel for interaction probability. Putting $f_y(\cdot)=k(\cdot,y)\in\mathcal{E}_+$, this gives expected interaction probability \[U_t(y)=\mu_t f_y = \frac{e^{-\frac{\gamma  y^2}{2 \gamma  t+1}}}{\sqrt{2 \gamma t +1}}\xrightarrow{t\rightarrow\infty}0\] and covariance term \[\mu_t(f_yf_z) = \frac{e^{-\frac{\gamma  \left(2 \gamma  t (y-z)^2+y^2+z^2\right)}{4 \gamma  t+1}}}{\sqrt{4 \gamma t +1}}\] We have that $U_t$ is maximized at time $\hat{t}$ with value $U_{\hat{t}}$ \begin{align*}\hat{t}&=\max\{y^2-\frac{1}{2\gamma},0\}\\U_{\hat{t}}&=\begin{cases}\frac{1}{|y|\sqrt{2e\gamma}}& \hat{t}>0\\e^{-\gamma y^2} & \hat{t}=0 \end{cases} \end{align*} These results show that the time of maximum interaction probability increases quadratically in interaction distance $\mathcal{O}(y^2)$, whereas the maximum interaction probability decreases with interaction distance $\mathcal{O}(1/|y|)$.  

We have covariance \[C_t(y,z) = \frac{c}{\sqrt{4 \gamma t +1}}e^{-\frac{\gamma  \left(2 \gamma  t (y-z)^2+y^2+z^2\right)}{4 \gamma  t+1}} + \frac{\delta^2-c}{2 \gamma t +1}e^{-\frac{\gamma (y^2+z^2)}{2 \gamma  t+1}}\] 

A higher-order interaction random field can be defined through defining $k: E\times E^n\mapsto\R_+$. For example, suppose $k$ is defined through a correlated bivariate radial basis function \[k(x,(y,z)) = \exp_- \frac{1}{2(1-\rho ^2)}\left((\frac{x-y}{\sigma_y})^2-2 \rho (\frac{x-y}{\sigma_y}) (\frac{x-z}{\sigma_z})+ (\frac{x-z}{\sigma_z})^2\right)\] Putting $f_{yz}(\cdot)=k(\cdot,(y,z))\in\mathcal{E}_+$, we have mean \[U_t(y,z)=\mu_tf_{yz} = \frac{\exp_- \left(\frac{t (y-z)^2+\sigma_y^2 z^2+\sigma_z^2 y^2-2 \rho  \sigma_y \sigma_z y z}{2 \left(1-\rho ^2\right) \sigma_y^2 \sigma_z^2+2 t \left(\sigma_y^2+\sigma_z^2-2 \rho  \sigma_y \sigma_z\right)}\right)}{\sqrt{1+t\frac{\left(\sigma_y^2+\sigma_z^2-2 \rho  \sigma_y \sigma_z\right)}{\left(1-\rho ^2\right) \sigma_y^2 \sigma_z^2}}}\xrightarrow{t\rightarrow\infty}0\] This reveals the expected probability of interaction of a particle with two correlated points.

\begin{figure}[h!]
\centering
\begingroup
\captionsetup[subfigure]{width=3in,font=normalsize}
\subfloat[Random measure sensitivity density $\amsbb{S}(\D x)$ for $t=1$\label{fig:wienera}]{\includegraphics[width=3.5in]{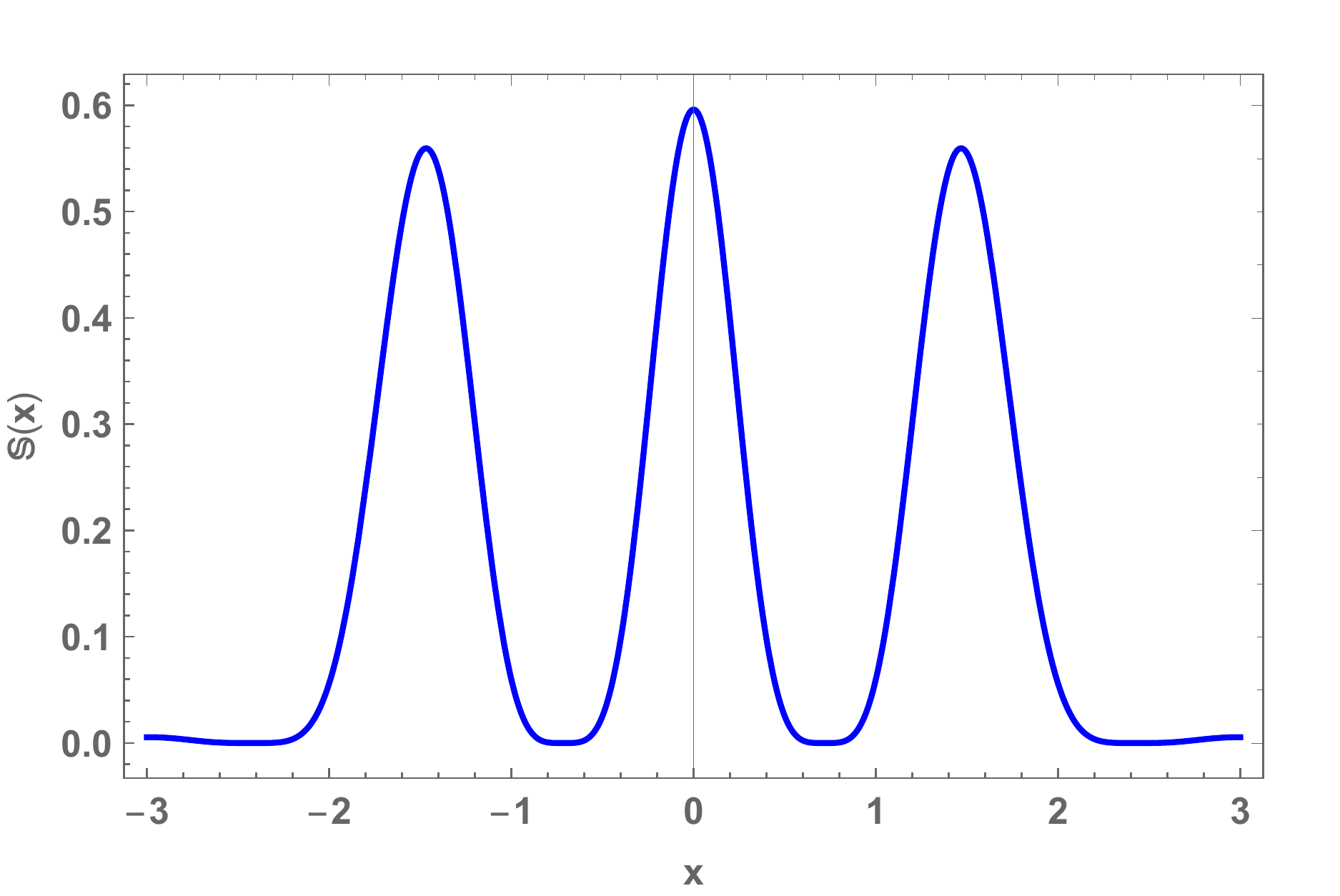}}
\subfloat[Random measure sensitivity density $\amsbb{S}(\D x)$ for $t=5$\label{fig:wienerb}]{\includegraphics[width=3.5in]{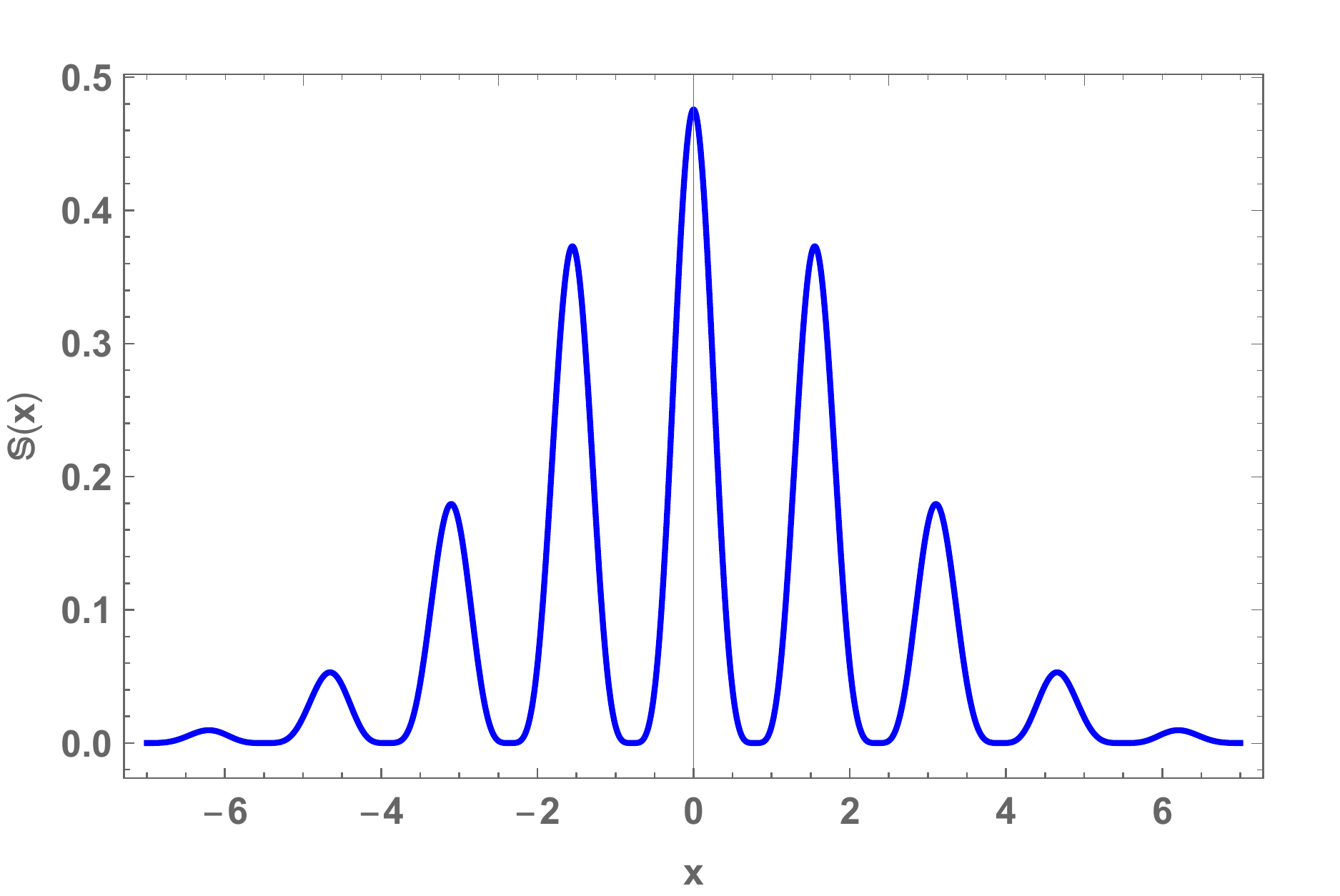}}\\
\subfloat[Orthogonal $M_t$: Random field sensitivity indices $\{\amsbb{S}_i\}$\label{fig:3aa}]{\includegraphics[width=3.5in]{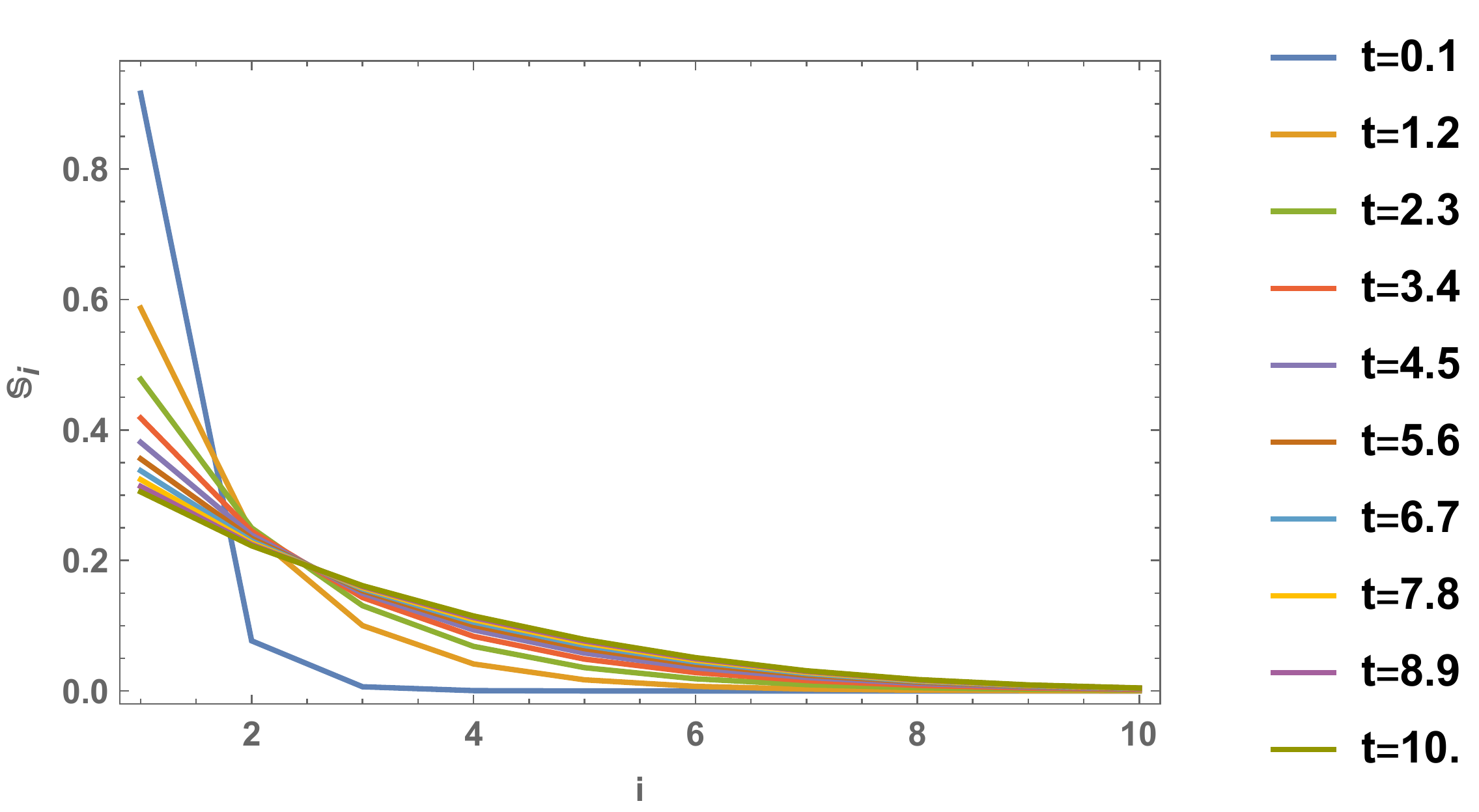}}
\subfloat[Orthogonal $M_t$: effective dimension of random field as a function of time $t$\label{fig:3ab}]{\includegraphics[width=3.5in]{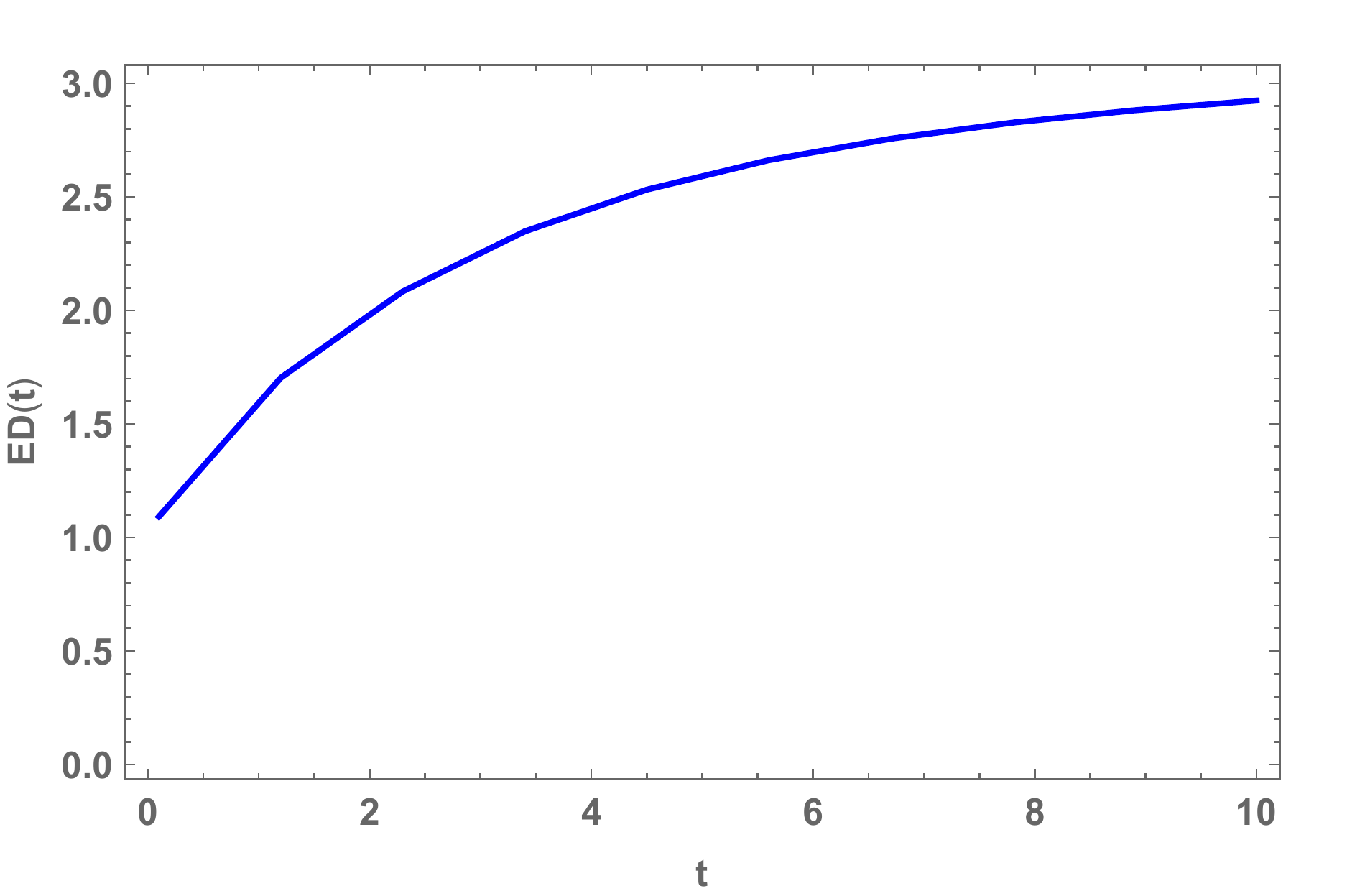}}\\
\subfloat[Orthogonal $M_t$: Random field covariance reconstruction $\lambda_1\varphi_1^\intercal\varphi_1$ at $t=1$\label{fig:3ab2}]{\includegraphics[width=3.5in]{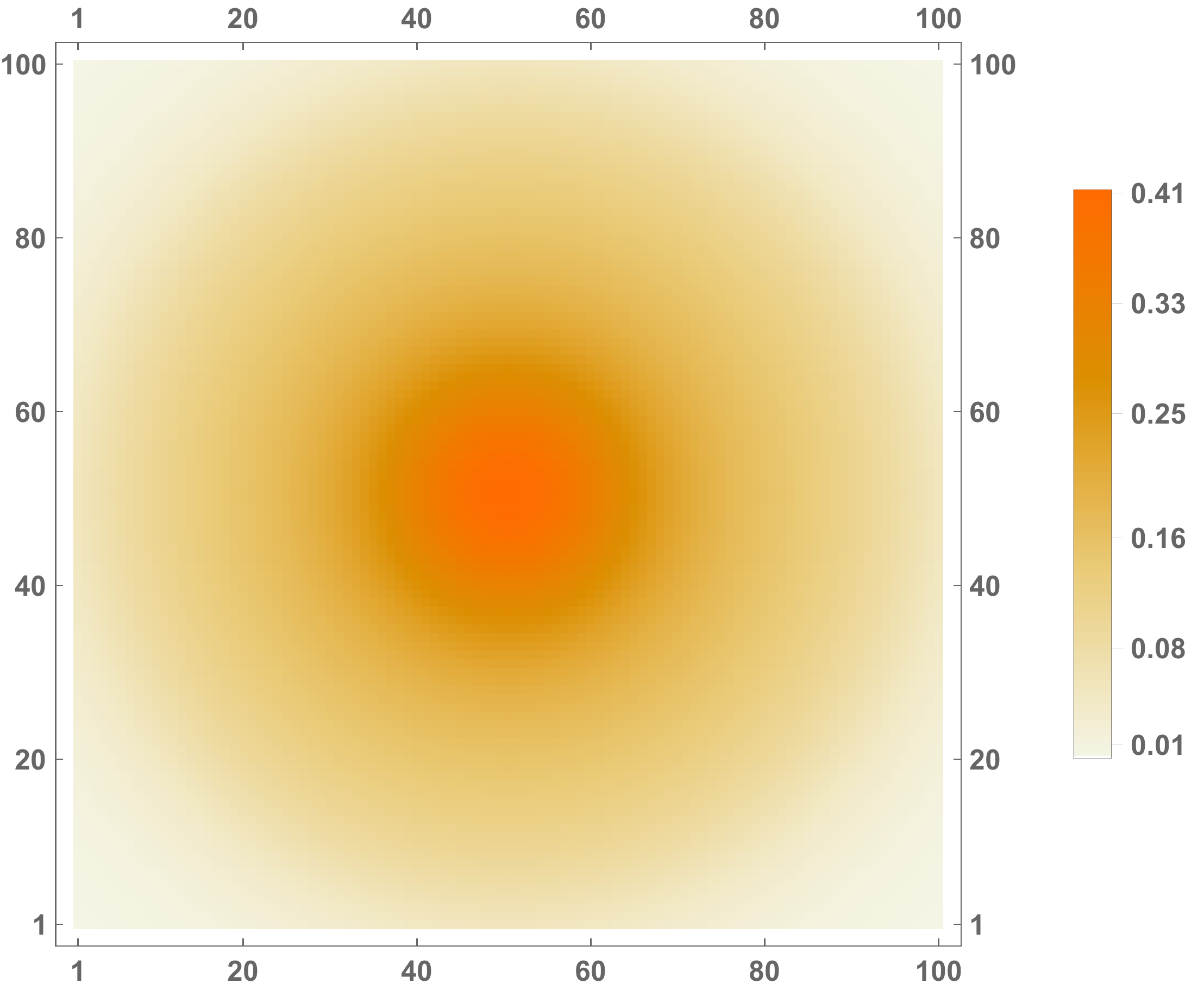}}
\subfloat[Dirac $M_t$: Random field covariance reconstruction $\lambda_1\varphi_1^\intercal\varphi_1$ at $t=1$\label{fig:3ab3}]{\includegraphics[width=3.5in]{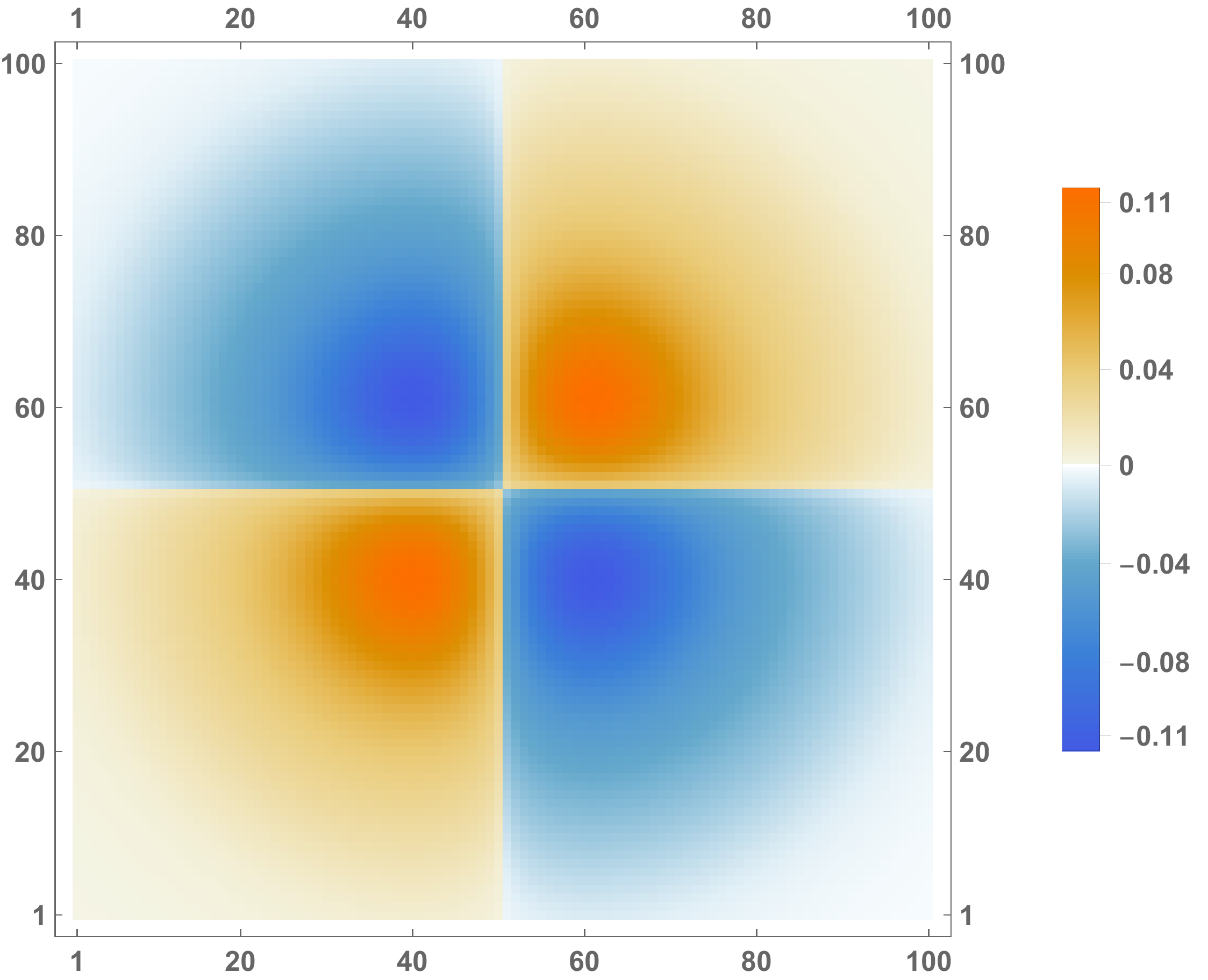}}\\
\endgroup
\caption{Wiener interaction probability analysis}\label{fig:wiener}
\end{figure}

In Figures~\ref{fig:wienera} and \ref{fig:wienerb}, we show the sensitivity indices for the random measure for $M_t f$ for $t\in\{1,5\}$. They are markedly different for the times and are each multimodal. This follows from the exponential dampening of periodic $g$. In Figure~\ref{fig:3aa}, we show the random field sensitivity indices, estimated from the eigenvalues of the covariance matrix on a grid on $[-5,5]\times[-5,5]$, each coordinate containing 100 equispaced values, and we take $c=\delta^2=\gamma=1$. As shown in Figure~\ref{fig:3ab}, as time increases, the effective dimension increases. Figures~\ref{fig:3ab2} and \ref{fig:3ab3} show the covariance reconstruction using the principal eigenvector and eigenvalue $\lambda_1\varphi_1^\intercal\varphi$ at $t=1$ for orthogonal and Dirac $M_t$. These are first-order approximations to the kernels of Figures~\ref{fig:kernel1} and \ref{fig:kernel3}. 

\paragraph{Specific calculation: drift and initial condition} Suppose we have $N=(\kappa,\nu)$ on $(E,\mathcal{E})$ where $\nu$ is Wiener. Now introduce an initial condition $X_0\in\R$ and drift constant $D\in\R$ with joint distribution $\eta$. The triple $(\mathbf{X},\mathbf{X}_0,\mathbf{D})$ forms the random measure $M=(\kappa,\mu=\nu\times\eta)$. Now consider the function $h(X,X_0,D) = X(t) + X_0 + Dt$ for $t\in\R_+$ fixed. Then the image measure $\mu_t=\mu\circ h^{-1}$ is given by \[\mu_t(A) = \int_{\R\times\R}\eta(\D x_0, \D D)\int_\R\D x\frac{1}{\sqrt{2\pi t}}e^{-\frac{(x-x_0-Dt)^2}{2t}}\ind{A}(x)\for A\in\mathcal{B}_{\R}\] If $\eta$ is Gaussian with hyperparameters for $X_0$ as $\mu_0$ and $\sigma_0^2$ and for $D$ as $\mu_D$ and $\sigma_D^2$ and correlation $\rho\in[-1,1]$, then \[\mu_t(A) = \int_\R\D x \frac{1}{\sqrt{2\pi(\sigma_0^2+t+\sigma_D^2 t^2+2 \sigma_0\sigma_D t \rho)}}\exp_-\frac{(x-\mu_0-\mu_D t)^2}{2(\sigma_0^2+t+\sigma_D^2 t^2+2 \sigma_0\sigma_D t \rho)}\ind{A}(x)\] so $\mu_t=\text{Gaussian}(\mu_0+\mu_D t,\sigma_0^2+t+\sigma_D^2 t^2+2 \sigma_0\sigma_D t \rho)$. The random field mean and covariance are similar to the case without drift, with additional parameters. 

For interactions such as collisions which generate debris fields, the interaction can define a transition into a cluster process, of the non-vaporized components and their drifts, with initial locations at the location of the object undergoing spallation. 



\FloatBarrier


\subsubsection{Identification} We recount the quantities that comprise the random measure and field models. Both are described through the tuple $\Xi=(\kappa,\mu_t,k)$. For each such tuple, there exists an indexed collection of restrictions $\Upsilon(\Xi)=\{(\kappa_A,\mu_{tA},k): A\subseteq F\}$. 

\subsubsection{Generalizations} All manner of generalizations can be had. We can assume the spacecraft are marked with properties, such as flag (country or company), and we can introduce additional complexity into the function $k$ by allowing interaction probability to depend on additional dimensions and to have different interaction types. Interactions could be collisions (the collision \emph{probability} resulting from the uncertainty of the particle's position at some time), the experience of space weather events (having spatiotemporal uncertainty), confrontation / dispute (where the particles are combatants), and so on. We can also consider interacting random fields for various restrictions of the system into subspaces. For example we could form restrictions through some clustering / partition of the domain and of space object marks / properties. With humanity's increasing expression of space-faring activities, e.g. \cite{spacex}, the models here can serve as canonical representations of the information it generates, i.e. the filtration $\mathscr{F}=(\mathscr{F}_t)_{t\in\R_+}$ generated by $M$, \[\mathscr{F}_t = \sigma\{M_s: s\le t\}.\]


\subsection{Adaptive randomized controlled trials}\label{ex:rct}

As referenced in Table~\ref{tab:0}, we discuss a framework for adaptive randomized controlled trials using random measures that mitigates the phenomenon of under or over-powered trials at final analysis time that is based on superpositions of orthogonal die random measures. The injection of counting noise by the orthogonal dice into trial design does not negatively impact type I or II error probabilities and eliminates effect correlation across groups caused by a fixed trial size. The normalized uncertainties (variances) may be interpreted as probabilities, which provides insight into the distribution of uncertainty of the random measure across groups.

\subsubsection{Randomized controlled trials}
A randomized controlled trial is a random realization of a random counting measure. Let $N=(\kappa,\nu)$ be a random counting measure on a clinical trial design space $(E,\mathcal{E})$, where $\kappa$ is a counting distribution with mean $c$ and variance $\delta^2$ and $\nu$ is a probability measure on $(E,\mathcal{E})$. In the simplest case, $E=\{C,T\}$ contains control and treatment groups with equal probability, i.e. $\nu\{C\}=\nu\{T\}=1/2$. Enrollees are materialized $\mathbf{X}=\{X_i\}$ and measurements are attained from the enrollees $\mathbf{Y}=\{Y_i\}$ in measurement space $(F,\mathcal{F})$ according to transition probability kernel $Q$ from $(E,\mathcal{E})$ into $(F,\mathcal{F})$, i.e. $Y_i\sim Q(X_i,\cdot)$. For simplicity, we assume $(F,\mathcal{F})=(\R_+,\mathcal{B}_{\R_+})$. The random measure $M=(\kappa,\nu\times Q)$ on $(E\times F,\mathcal{E}\otimes\mathcal{F})$ is formed from $(\mathbf{X},\mathbf{Y})$ by stone throwing construction as \[Mf = \sum_i^K f\circ (X_i,Y_i)\for f\in(\mathcal{E}\otimes\mathcal{F})_+\] 

Put $f_T(x,y)=\ind{\{T\}}(x)y$ and $f_C(x,y)=\ind{\{C\}}(x)y$. These are disjoint functions. In general, the random variables $Mf_T$ and $Mf_C$ are correlated, where \[\Cov(Mf_T,Mf_C) = (\delta^2-c)(\nu\times Q)f_T(\nu\times Q)f_C\] with mean \[\E Mf_T = c(\nu\times Q)f_T\] and variance \[\Var Mf_T = c(\nu\times Q)f_T^2 + (\delta^2-c)((\nu\times Q)f_T)^2\] The variance of the random measure describes its uncertainty. We have \[(\nu\times Q)f_T = \int_E\nu(\D x)\ind{\{T\}}(x)\int_FQ(x,\D y) y=\frac{1}{2}\int_F Q(\{T\},\D y)y=\frac{1}{2}c_T\] and \[(\nu\times Q)f_T^2 = \int_E\nu(\D x)\ind{\{T\}}(x)\int_FQ(x,\D y) y^2=\frac{1}{2}\int_F Q(\{T\},\D y)y^2=\frac{1}{2}(c_T^2+\delta_T^2)\]

The null hypothesis is that $(\nu\times Q)f_T=(\nu\times Q)f_C$ and depends on the mean measure. This is equivalent to $\E Mf_T = \E Mf_C$. The typical setting is $\kappa=\text{Dirac}(c)$ for $c\in\N_{>0}$. This random measure has minimum variance with $\delta^2=0$. The covariance is negative due to the fixed sample size, \[\Cov(Mf_T,Mf_C) = -\frac{c}{4}c_Tc_C\] and \[\Cov(\frac{1}{c}Mf_T,\frac{1}{c}Mf_C) = -\frac{1}{4c}c_Tc_C\] which decays to zero as $c\rightarrow\infty$. 

\begin{figure}[h]
\centering
\includegraphics[width=4in]{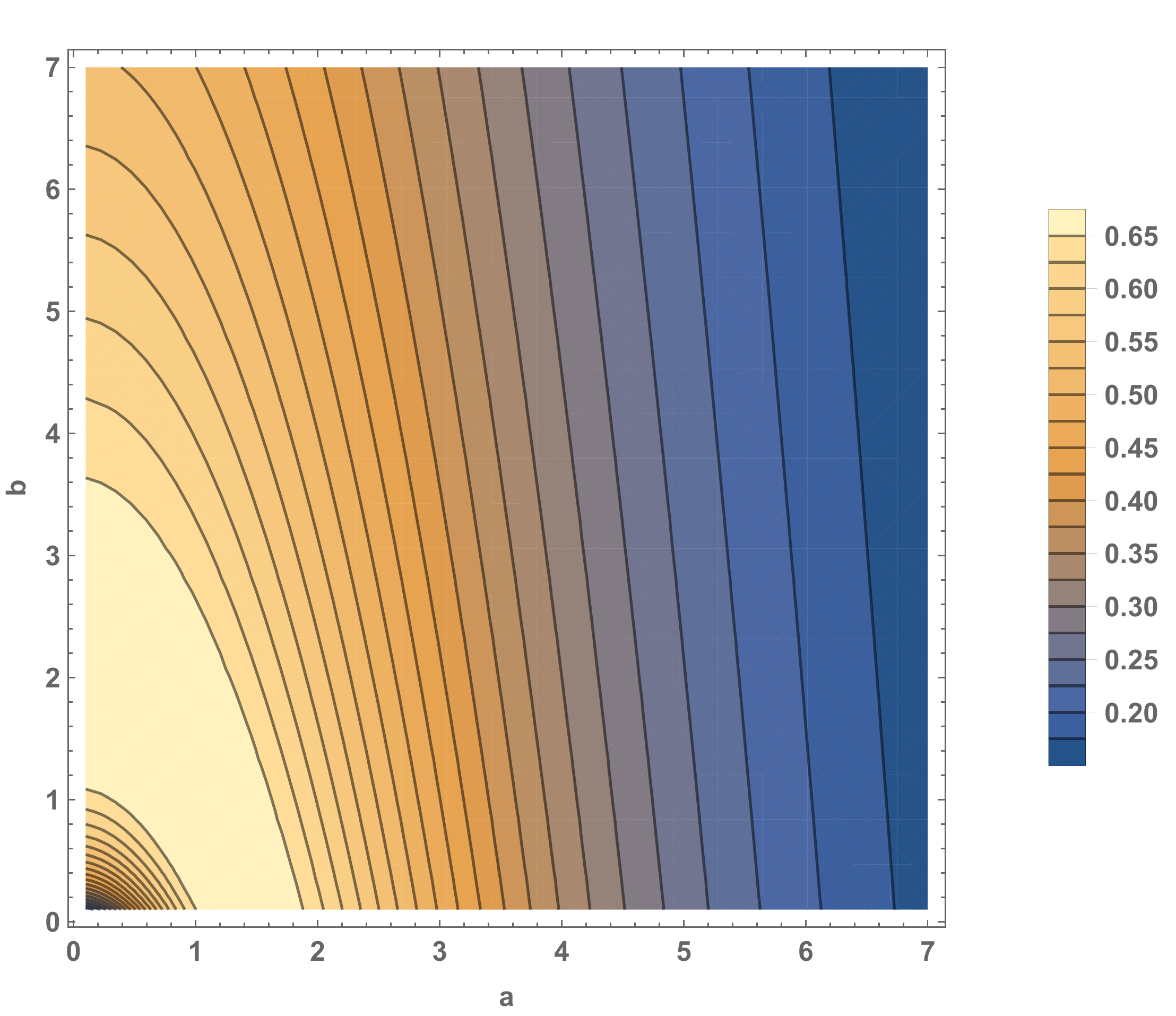}
\caption{Entropy of $\amsbb{S}$ for a randomized controlled trial with non-negative measurements using an orthogonal random measure}\label{fig:entcont}
\end{figure}

When $\kappa$ is orthogonal, the covariance is zero for any mean sample size $c\in(0,\infty)$, the variance is given by \[\Var Mf_T = c(\nu\times Q)f_T^2\] and the normalized variances form a sensitivity distribution \[\amsbb{S}_T = \frac{(\nu\times Q)f_T^2}{(\nu\times Q)f^2} = \frac{c_T^2+\delta_T^2}{c_T^2+\delta^2_T + c_C^2 + \delta_C^2}\] where $\amsbb{S}_T+\amsbb{S}_C=1$. We take $c_T = a c_C$ and $\delta_T^2 = b\delta_C^2$ for $a,b\in(0,\infty)$ and plot $H(\amsbb{S})$ as a function of $a$ and $b$ in Figure~\ref{fig:entcont} where $c_C=\delta_C^2=1$. 

\FloatBarrier

\begin{figure}[h]
\centering
\includegraphics[width=4in]{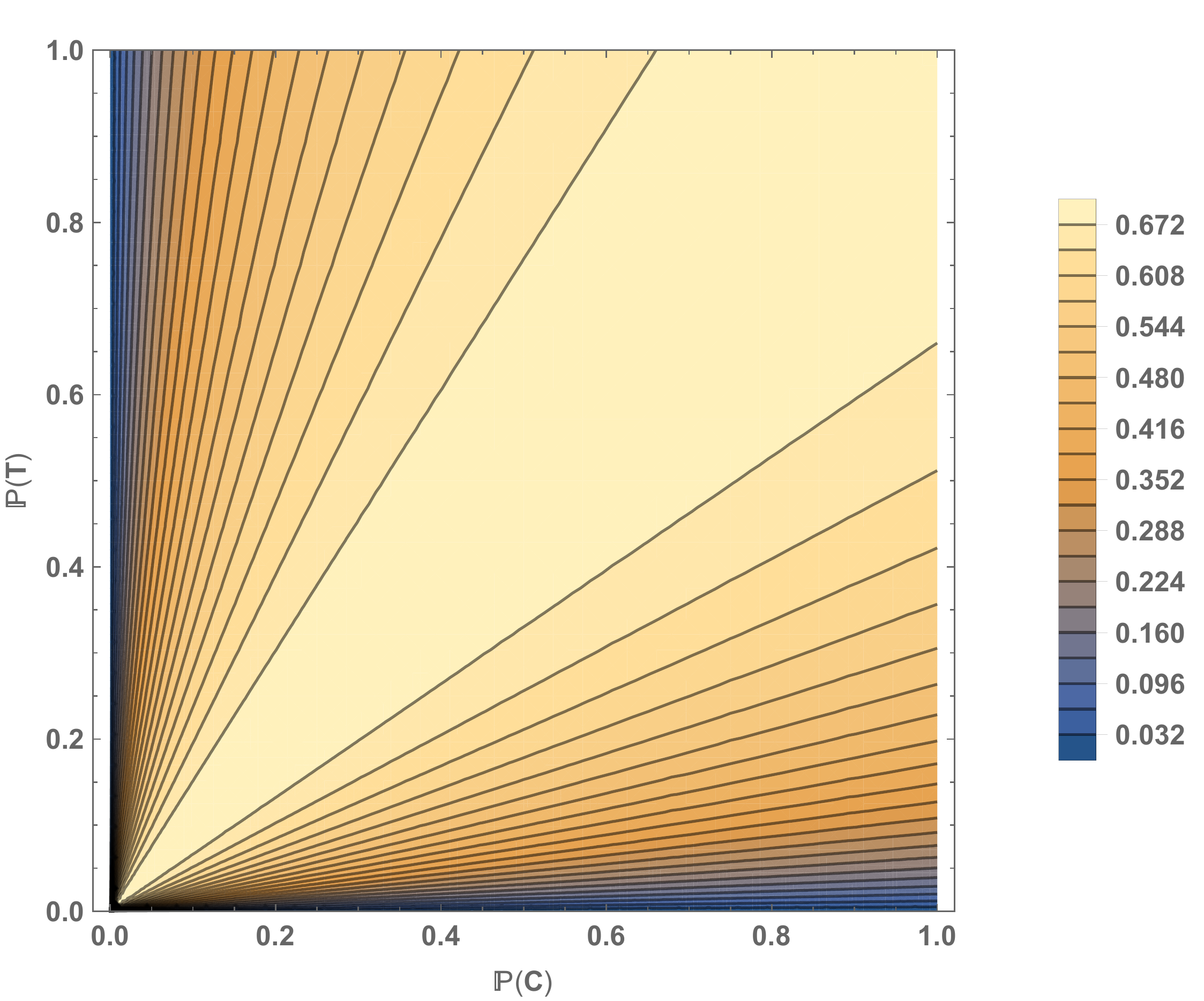}
\caption{Entropy of $\amsbb{S}$ for a randomized controlled trial testing vaccine efficacy using an orthogonal random measure}\label{fig:enttrial}
\end{figure}

For example, consider a randomized controlled trial testing vaccine efficacy. The measurement is an indicator of infection  of the underlying disease being vaccinated against over the period of the trial. Then $c_T = \P(T)$ and $c_C=\P(C)$ are the probabilities of infection based on group assignment and the variance is \[\Var Mf_T = c\P(T) + (\delta^2-c)\P^2(T)\] and the covariance is \[\Cov(Mf_T,Mf_C) = (\delta^2-c)\P(T)\P(C)\] For Dirac, this is \[\Var Mf_T = c\P(T)(1-\P(T))\] and \[\Cov(Mf_T,Mf_C) = -c\P(T)\P(C)\] For orthogonal $\kappa$ the variance is \[\Var Mf_T = c\P(T)\] and the sensitivity index in $T$ is \[\amsbb{S}_T = \frac{\P(T)}{\P(T)+\P(C)}\] When $\P(T)\ll\P(C)$, then the uncertainty of the effects is dominated by the infections of the control group, $\amsbb{S}_T\ll\amsbb{S}_C$. In Figure~\ref{fig:enttrial} we show the entropy of $\amsbb{S}$ in $\P(C)$ and $\P(T)$, i.e. \[H(\amsbb{S}) = -\amsbb{S}_C\log(\amsbb{S}_C)-\amsbb{S}_T\log(\amsbb{S}_T)\] Clinical trials of efficacious vaccines against infectious diseases have low entropy of uncertainty. Consider the recent Moderna and Pfizer vaccines for COVID-19 \citep{moderna,pfizer}. For Moderna, there are approximately $n=30\,400$ enrollees, where infection was recorded in 5 cases in the vaccine group and 90 cases in the control group, giving entropy $H(\amsbb{S})\approx0.206$. For Pfizer, there are approximately $n=44\,000$ enrollees, with 8 cases in the vaccine group and 162 cases in the control group, giving entropy $H(\amsbb{S})\approx0.190$. The theoretical largest entropy is $\log(2)\approx0.693$, whereas the minimum entropy is zero.  

\FloatBarrier




\subsubsection{Adaptation using superpositions of orthogonal dice}
We take $M$ to be a superposition of independent random measures that are constructed sequentially. Here we consider for simplicity two random measures with superposition counting law \[\kappa = \kappa_1*\kappa_2\] Because the random measures are independent, the superposition law of the random measure is described through the product of the Laplace functionals of the respective random measures. If these random measures are orthogonal, then the superposition random measure is orthogonal. 

Let $c_1$ be the initial total size informed by power analysis: that is, the minimum number of enrollees at some effect size (with some assumed type I and II error probabilities). We fractionate the design by some factor $a$, such as $a=1/2$. We take $K_1\sim\kappa_1$ as the first orthogonal die with support $\supp(K_1)=\{m,\dotsb,n\}$ with $m\ge ac_1$. $K_1$ is the number of enrollees in the first random measure. Then, based on the collection of the data and observation of effects in this first batch, power analysis with the empirical effect gives the imputed total number of new enrollees for the desired significance level, which we denote as $c_2$. If $c_2>K_1$, then additional data is generated, through new independent enrollees, to prevent an underpowered design. In this case, to construct $\kappa_2$, we choose the first orthogonal die with lower support $m\ge c_2-K_1$. Then $K_2\sim\kappa_2$ is sampled and the enrollees materialized. Then we have that $K_1+K_2\ge c_2$, so we have achieved or exceeded our type I and II error probabilities. Because $\kappa_1$ and $\kappa_2$ are independent and orthogonal, then $\kappa$ is orthogonal. In analysis of effects across treatment and control groups, the group effects are approximately compound Poisson random variables due to weak convergence of orthogonal die random measure to the Poisson random measure.  

The first 15 orthogonal dice are shown in Table~\ref{tab:dice}. There are an infinite number of orthogonal dice, and their support size $n-m+1$ is coprime to 2 and 3, i.e. their sizes enumerate all primes greater than or equal to 5 and their products. Please see \cite{bastian2020orthogonal} for more information about the orthogonal dice and weak convergence of the orthogonal die random measure to the Poisson random measure.

\FloatBarrier

\subsubsection{Functional information} Here we expand the randomized controlled trial space to incorporate a mark-space. We suppose each enrollee has some independent mark-space $(H,\mathcal{H})$, containing properties such as demographics, genome, contact network properties, etc., with distribution $\eta$. Now let $Q$ be a transition kernel from $(E\times H,\mathcal{E}\otimes\mathcal{H})$ into $(F,\mathcal{F})=(\R_+,\mathcal{B}_{\R_+})$. The random measure $M=(\kappa,\nu\times\eta\times Q)$ on $(E\times H\times F, \mathcal{E}\otimes\mathcal{H}\otimes\mathcal{F})$ is formed by $(\mathbf{X},\mathbf{Y},\mathbf{Z})$ through STC.

Consider again the set-up of a randomized controlled trial testing vaccine efficacy. We retrieve the existing analysis with $f_T(x,y,z) = \ind{\{T\}}(x)z$ and $f_C(x,y,z)=\ind{\{C\}}(x)z$, e.g., the mean is  \begin{align*}\E Mf_T &= c(\nu\times\eta\times Q)f_T \\&= c\int_E\nu(\D x)\ind{\{T\}}(x)\int_H\eta(\D y)\int_F Q((x,y),\D z) z \\&=\frac{c}{2}\int_H\eta(\D y)\int_F Q((\{T\},y),\D z) z \\&= \frac{c}{2}\int_H\eta(\D y)\P(T,y)\\&=\frac{c}{2}\P(T) \end{align*} Now suppose we have a partition $A,\dotsb,B$ of $H$. We define $f_T^A(x,y,z) = \ind{\{T\}}(x)\ind{A}(y)z$ and $f_C^A(x,y,z) = \ind{\{C\}}(x)\ind{A}(y)z$. The mean is \[\E Mf_T^A = \frac{c}{2}\P_A(T)\] For orthogonal $N$, we have sensitivity indices for restrictions to $A\subset H$ as \begin{align*}\amsbb{S}^A_T &= \frac{\P_A(T)}{\P_A(T)+\P_A(C)}\\\amsbb{S}^A_C&= \frac{\P_A(C)}{\P_A(T)+\P_A(C)}\end{align*} Restrictions could be various subpopulations of interest, such as indications of pregnancy or race.

Suppose we have some classifier $g:E\times H\mapsto\{0,1\}$ that predicts infection, trained on input-output samples $\{((X_i,Y_i),Z_i)\}$. Let $f(x,y,z) = (g(x,y)-z)^2$. Let $Q((x,y),\cdot)=\delta_{\E g}(\cdot)$ so that $f(x,y,z)=(g(x,y)-\E g)^2$. Then \[\E Mf = c(\nu\times\eta\times Q)f = c\Var g \] The HDMR of $g$ attains the hierarchy of component functions in $E\times H$. The classifier $g$ has real-world application to pharmaceutical development: the United States Food and Drug Administration (FDA) has an  Accelerated Approval program for ``serious conditions that fill an unmet medical need'' \citep{fda}. If (a) the classifier is trained on independent input-output samples, such as from some other trial or design, and if (b) the relation between the marks and the measurement has a strong clinical foundation, and if (c) the classifier is highly predictive, then in some cases the predicted outcomes may be used as surrogates for the actual outcomes to attain provisional acceptance until the actual outcomes may be measured. Similar hold for regressors and continuous measurements. 

\subsection{Dynamic survival analysis of epidemics}\label{ex:dsa} We develop a model generalizing SIR epidemic models to Poisson-type network models \citep{jacobsen2018large} with Poisson-type random measures \citep{bastian} called \emph{dynamic survival analysis} in order to capture the network structure of contacts and quarantine and to provide uncertainty quantification of infected and recovered individuals. This framework may be interpreted as an uncertainty quantification of epidemics based on compartmental models and may be instantiated for various particle systems of persons. 

\subsubsection{Binomial random measure}
Assume the independency of individuals $\mathbf{U}$ is surveyed for symptoms of infectious disease and forms a binomial random measure $N=(K,\nu\times Q )$ on the space $(E,\cal {E})$, where $E=\{(x,y): 0<x<y\}$, $\kappa=\text{Binomial}(n,p)$, and $p>0$ (specified below). Each individual $U_i=(X_i,Y_i)$ is described by a pair of infection and recovery times, distributed according to $\nu\times Q$. 
 
 \subsubsection{Mean measure through dynamic survival analysis}
To describe the relevant mean law $\nu\times Q$, consider a dynamic survival analysis (DSA) model describing the evolution of proportions of susceptible ($S$) infectious ($I$) and removed ($R$). Let $\mathcal{G}(\psi,n)$ be a CM random graph with Poisson-type (PT) degree distribution probability generating function (pgf) $\psi$ and population size $n$ where half-edges (stubs) are paired uniformly at random. The PT pgfs satisfy \[\partial\psi(z) = \mu\,\psi^\kappa(z).\] The nodes of the graph are individuals. Each node is labeled with an infectious status in $\{S,I,R\}$. Initially, $S$-type nodes have all unpaired half-edges, and infection occurs through random pairing with an infected half-edge. $I$-type nodes may infect ($I$), recover ($R$), or drop ($D$). Let $L_t(i)\in\{S,I,R\}$ be the evolution of the $i$-th node label. An ODE system may be derived for the time-evolution of $(S,I,R)$ using the law of large numbers, where $\theta=(\beta,\gamma,\rho,\kappa,\mu,\delta)$ and initial condition $(1,\rho,\mu\rho)$ 
\begin{align}
\dot{S}_t &= -\beta D_t S_t\label{eq:S2}\\
\dot{I}_t &= \beta D_t S_t - \gamma I_t\nonumber\\
\dot{D}_t &= \beta(1-\kappa)D_t^2+(\kappa\mu\beta S_t^{2\kappa-1}-\beta-\gamma-\delta)D_t\label{eq:D2}
\end{align} 
$\kappa$ is the average network density defined as \[\kappa=\frac{\partial^2\psi(1)}{(\partial\psi(1))^2}=\frac{\mu^{ex}}{\mu}\] where $\mu^{ex}$ is the average excess degree and $\mu$ is the average degree. $D_t$ describes the number of SI edges divided by the number of susceptible nodes and is called the average density of infection. We have basic reproduction number for the DSA model on $\mathcal{G}(\psi,n)$ \begin{equation}\mathcal{R}_0 = \frac{\kappa\beta\mu}{\beta+\gamma+\delta}.\end{equation} Let $\tau_\infty$ be the solution to \begin{align}\tau_\infty &= 1 - e^{-\mathcal{R}_0(\tau_\infty+\rho)}\quad\text{if }\kappa=1\label{eq:tau}\\\frac{\kappa}{\kappa-1}((1-\tau_\infty)^{1-\kappa}-1)&=\mathcal{R}_0((1-\tau_\infty)^\kappa-(1+\rho))\quad\text{if }\kappa\ne 1\nonumber\end{align} and note that $S_\infty=1-\tau_\infty$. To complete the definition of the binomial random measure, we set $p=\tau_\infty$. We think of $S_t$ as an improper survival function with improper density $-\dot{S}_t$. It is improper because \[\int_0^\infty -\dot{S}_t=1-S_\infty = \tau_\infty < 1\] Hence, the proper density and mean measure of the binomial random measure is given by \begin{equation}\nu(x) = -\dot{S}_x/\tau_\infty\end{equation} In view of Proposition~\ref{prop:inflate}, the binomial random measure may be interpreted as a collection of $n$ zero-inflated infection and recovery times. The DSA model parameters are shown below in Table~\ref{tab:dsa}. 

\begin{table}[h!]
\caption{DSA model parameters}
\label{tab:dsa}
\begin{center}
\begin{tabular}{cp{3in}c}
\toprule
Parameter & Description & Domain \\\midrule
$\beta$ & infection rate & $(0,\infty)$\\
$\gamma$ & recovery rate & $(0,\infty)$\\
$\rho$ & initial fraction of infected & $(0,1)$\\
$1+\rho$ & scaled initial population & $(1,2)$\\
$1/\rho$ & initial susceptible population size & $(1,\infty)$\\
$\tau_T$ & probability of a randomly selected susceptible individual being infected by time $T$; proportion of infected by time $T$ in the epidemic in infinite population; relative epidemic size at time $T$ & $(0,1)$\\
$\mu^{ex}$ & mean excess degree of contact network & $(0,\infty)$\\
$\mu$ & mean degree of contact network& $(0,\infty)$\\
$\kappa$ & average network density; $\kappa=\mu^{ex}/\mu$ & $(0,\infty)$\\
$\delta$ &quarantine and SD rate & $\R_+$\\
$\mathcal{R}_0$ & basic reproduction number; $\kappa\beta\mu/(\beta+\gamma+\delta)$ & $(1,\infty)$\\
 \bottomrule
\end{tabular}
\end{center}
\end{table}
\FloatBarrier

The restriction of $\nu$ to $A=(0,T]$ (or, thinning of $\nu$ by $A$) is given by $\nu_A$ as \begin{equation}\label{eq:nuA}\nu_A(x) = \nu(x)\ind{\{x<T\}}(x)/\nu(A)=-\dot{S}_x\ind{\{x<T\}}(x)/\tau_T\end{equation} The restriction of $\nu$ to $B=(T,U]$ is given by $\nu_B$ as \begin{equation}\label{eq:nuA}\nu_B(x) = \nu(x)\ind{\{T<x<U\}}(x)/\nu(B)=-\dot{S}_x\ind{\{T<x<U\}}(x)/(\tau_{U}-\tau_{T})\end{equation}

Consider the transition probability kernel from $(F,\mathcal{F})$ into $(F,\mathcal{F})$ defined by $Q$ as the shifted (Lebesgue) exponential distribution with parameter $\gamma$. \begin{equation}\label{eq:Q}Q(x,y)\sim\text{Exponential}(\gamma)\ind{x<y}(y).\end{equation} Another variant is to translate the density with parameter $\varepsilon\in\R$ where $\varepsilon+y>0$, with kernel \begin{equation}\label{eq:Qe}Q_\varepsilon(x,y)\sim\text{Exponential}(\gamma)\ind{x<y+\varepsilon}(y+\varepsilon).\end{equation} 

The law $\nu\times Q$ is defined on $(E,\mathcal{E})$ where \begin{equation}\label{eq:E}E=\{(x,y): 0<x<y<\infty\}\end{equation} and $\mathcal{E}=\mathcal{B}_{E}$. The marginal law $\nu Q$ is on $(F,\mathcal{F})$. The restricted law $\nu_A\times Q$ is on $(A\times F,\mathcal{B}_A\otimes\mathcal{F})$, and the restricted marginal law \[g(B) = (\nu_A Q)(B)\quad\text{for}\quad B\in\mathcal{F}\] is on $(F,\mathcal{F})$. If the recovery times are also restricted by $T$, then we have conditional law \[g_A(B) = \frac{g(B)}{g(A)}\quad\text{for}\quad B\in\mathcal{B}_A\]

Let $T\in\R_+$ be the final time of the epidemic with $\tau_T$ and number of infections $k_T$. The population size is $n=k_T/\tau_T$. The number of susceptibles at time $T$ is $s_T = n-k_T$. The total number of infected at the end of the epidemic is $k_\infty = \tau_\infty n$, and the total number of susceptibles at the end of the epidemic is $s_\infty=(1-\tau_\infty)n$. These quantities are summarized in Table~\ref{tab:dsaq}.

\begin{table}[h!]
\caption{DSA model quantities}
\label{tab:dsaq}
\begin{center}
\begin{tabular}{cp{3in}cc}
\toprule
Parameter & Description & Observed & Value \\\midrule
$k_T$ & count of infected by time $T$ & Yes & $k_T$\\
$n$ & population size at time $T$ (number of susceptibles and infected) & No & $k_T/\tau_T$\\
$s_T$ & count of susceptibles by time $T$ & No & $n-k_T$\\
$k_\infty$ & count of infected by end of epidemic & No & $\tau_\infty n$\\
$s_\infty$ & count of susceptibles by end of epidemic & No & $(1-\tau_\infty)n$\\
 \bottomrule
\end{tabular}
\end{center}
\end{table}
\FloatBarrier

\subsubsection{Reduced system}\label{sec:reduce}
We divide the third equation \eqref{eq:D2} by the first \eqref{eq:S2}, solve for $D_t$ as a function of $S_t$, then plug the result into \eqref{eq:S2}. This gives a reduced system with only one equation describing the decay of susceptibles \begin{equation}\label{eq:S3}-\dot{S}_t = \frac{a}{1-\kappa} S_t(1-S_t^{\kappa-1})+b(1-S_t^\kappa)S_t^\kappa+cS_t^\kappa\end{equation} with $S_0=1$ where \begin{align*}a&=\beta+\gamma+\delta\\b&=\beta\mu\\c&=\beta\mu\rho\end{align*} and \begin{align*}\mathcal{R}_0&=\kappa b/a\\\rho&=c/b.\end{align*} 


For $\kappa=1$ (Poisson $\psi$) the ODE system becomes 
\begin{align}
\dot{S}_t &= -\beta D_t S_t\label{eq:S}\\
\dot{I}_t &= \beta D_t S_t - \gamma I_t\nonumber\\
\dot{D}_t &= (\mu\beta S_t-\beta-\gamma-\delta)D_t\label{eq:D}
\end{align} 
Letting $a=\beta+\gamma+\delta$, we divide the third equation \eqref{eq:D} by the first \eqref{eq:S}, giving
\begin{align*}
\frac{\D D_t}{\D S_t} &= -\mu+\frac{a}{\beta S_t}\\
D_0(1) &= \mu\rho
\end{align*} 
which may be solved as \begin{equation}\label{eq:DS}D_t(S_t) = \frac{\tilde{\gamma}}{\beta}\log (S_t)+  \mu(1+  \rho  - S_t).\end{equation} Then we plug \eqref{eq:DS} into the first equation \eqref{eq:S} to give \begin{equation}\label{eq:SS}-\dot{S}_t = (\beta+\gamma+\delta)S_t\log (S_t)+  \beta\mu(S_t-S_t^2)+ \beta\mu\rho S_t.\end{equation} Note that equation \eqref{eq:SS} may be obtained from \eqref{eq:S3} by taking the limit $\kappa\rightarrow1$. Consistent with Section~\ref{sec:reduce}, \eqref{eq:SS} suggests the parameters \begin{align*}a&=\beta+\gamma+\delta\\b&=\beta\mu\\c&=\beta\mu\rho\end{align*} with \begin{align*}\mathcal{R}_0&=b/a\\\rho&=c/b.\end{align*} Note that \[S_\infty=-\frac{a}{b} W\left(-\frac{b}{a} e^{-\frac{b}{a}-\frac{c}{a}}\right)=-\frac{1}{\mathcal{R}_0} W(-\mathcal{R}_0 e^{-\mathcal{R}_0(1+\rho)})=1-\tau_\infty\] where $W$ is the product-log function. 

%

\subsubsection{Time-series analysis} Let $f\in\mathcal{E}_+$ be a non-negative time-series, i.e. $f(x)\in\R_+$, such as a stock market index or some other population measure, e.g. happiness \citep{happiness}. The random variable $Nf$ is the integral of $f$ with respect to the epidemic, i.e. \[Nf = \int_{\R_+}N(\D x)f(x)\]  The mean and variance are \begin{align*}\E Nf &= n\tau_\infty\nu f\\\Var Nf &= n\tau_\infty(\nu f^2 - \tau_\infty(\nu f)^2)\end{align*} Consider the restriction to $A=(T,U]$ with $\nu(A)=\tau_U-\tau_T$ and put $g(x)=f(x)\ind{A}(x)$. Then the mean and variance are \begin{align*}\E Ng &=\E N_A f = n\tau_\infty(\tau_U-\tau_T)\nu_A f\\\Var Ng &=\Var N_Af = n\tau_\infty(\tau_U-\tau_T)(\nu_A f^2-(\tau_U-\tau_T)(\nu_A f)^2)\end{align*}

\subsubsection{RM-ANOVA}\label{sec:dsa-rm} The mean and variance of the random measure $N$ are \begin{align*}\E N(E) &=  n\tau_\infty\\\Var N(E) &= n\tau_\infty(1-\tau_\infty)\end{align*} For the sets $A=(0,T]$ and $B=(T,\infty]$, we have \begin{align*}\E N(A) &= n\tau_T\\\E N(B) &= n(\tau_\infty-\tau_T)\\\Var N(A) &= n\tau_T(1 - \tau_T)\\\Var N(B) &= n(\tau_\infty-\tau_T)(1+\tau_T-\tau_\infty)\\\Cov(N(A),N(B)) &= -n\tau_T(\tau_\infty-\tau_T)\end{align*} The structural and correlative sensitivity indices are \begin{align*}\amsbb{S}_a^a &= \frac{\tau_T(1-\tau_T)}{\tau_\infty(1-\tau_\infty)}\\\amsbb{S}_b^a &= \frac{(\tau_\infty-\tau_T)(1+\tau_T-\tau_\infty)}{\tau_\infty(1-\tau_\infty)}\\\amsbb{S}_a^b &=\amsbb{S}^b_b = -\frac{\tau_T(\tau_\infty-\tau_T)}{\tau_\infty(1-\tau_\infty)}\end{align*} and are free of $n$.

\subsubsection{RF-ANOVA} Let the random measure $N=(\kappa,\nu)$ as in Section~\ref{sec:dsa-rm} and consider restricted random measure $N_A=(\kappa_A,\nu_A)$ to timeset $A=(0,T]$ with $\nu(A)=a>0$ with counting variable $K_A\sim\kappa_A=\text{Binomial}(n,\tau_T=a\tau_\infty)$. Let $(F,\mathcal{F})$ be a measurable space, such as a stock market indices or sentiment measures over time, i.e. $F=\R_+^{\R_+}$. Let $k:E\times F\mapsto\R_+$ be a $\mathcal{E}\otimes\mathcal{F}$-measurable function. Putting $f_y(\cdot)=k(\cdot,y)\in\mathcal{E}_+$ for $y\in F$, the random field $G$ is formed as \[G(y)=N_Af_y = \sum_i^{K_A}k(X_i,y)\for y\in F\] with  expected value \[U(y)=n\tau_T\nu_A(f_y)\for y\in F\] and covariance \[C(y,y') =n\tau_T(\nu_A(f_yf_{y'}) -\tau_T\nu_A(f_y) \nu_A(f_{y'}))\for y,y'\in F\]  For example we can put \[k(x,y) = \int_{\R_+}\nu_{(0,x]}(\D u)y(u)\] so that \[\nu_A(f_y) = \int_{\R_+}\nu_{A}(\D x)\int_{\R_+}\nu_{(0,x]}(\D u)y(u)\] and \begin{align*}\nu_A(f_yf_{y'}) &= \int_{\R_+}\nu_A(\D x)(\int_{\R_+}\nu_{(0,x]}(\D u)y(u))(\int_{\R_+}\nu_{(0,x]}(\D u)y'(u))\\&=\int_{\R_+}\nu_A(\D x)\int_{\R_+}\nu_{(0,x]}(\D u)y(u)y'(u)\end{align*} Another choice could be as simple as \[k(x,y) = y(x)\] so that  \[\nu_A(f_y) = \int_{\R_+}\nu_{A}(\D x)y(x)\] and \begin{align*}\nu_A(f_yf_{y'}) &= \int_{\R_+}\nu_A(\D x)y(x)y'(x)\end{align*}


\section{Discussion and conclusions}\label{sec:discuss} We describe a general integrated framework RM-MM-RF-ANOVA for uncertainty quantification, which furnishes first and second-order statistics and enables construction of general positive random fields. The first-order MM-ANOVA analysis is a functional ANOVA decomposition of the intensity measure of the squared loss functional into subspaces, whereas the second-order RM-ANOVA analysis is decomposition of random measure variance into subspaces. For orthogonal random measures, RM-ANOVA furnishes a sensitivity distribution which may be attained in all marginals. Positive random fields may be constructed from the random measures, and RF-ANOVA decomposes the field. 


This UQ framework is based on the variance-covariance structure of the mixed binomial process. As such it generalizes the binomial process, also known as the empirical random measure or the Dirac random counting measure, and its bootstrap estimator. The generalization is through a degree-of-freedom gain with the instantiation of a counting distribution $\kappa$. We show that for the squared loss function, the mean (intensity) measure may be decomposed using HDMR. Other loss functions may be used, such as absolute loss, thereby forsaking functional ANOVA for representation of the intensity measure. The random field model may be used to construct interacting particle systems. The method of maximum entropy allows the density of the random variable $Nf$ (formed by the random measure $N$ and test function $f\in\mathcal{E}_+$) to be attained from a collection of evaluations of the Laplace transform. This gives a general, scalable method to attain the density of $Nf$. 

The examples illuminate some interesting findings and are suggestive. We find that the effective dimension of the HDMR of the symmetric polynomial is regulated by the coefficient of variation $\rho$ of the input variables, with maximum entropy at $\rho=1$. Its HDMR is dense in subspaces. For a univariate polynomial with Bernoulli input, the orthogonal random measure variance exhibits strong dependence on $p$, with maximum entropy at $p=1/2$. Moreover we find the second moment, corresponding to the variance of the orthogonal random measure, to have the structure of an inverted double-well potential. We also see that the structural sensitivity indices for the binomial process (Dirac random counting measure) each possess a singularity. In the second example, we find that the Ishigami function has a non-monotone entropy profile in the second parameter for $a=7$ and that the orthogonal random measure uncertainty on partitions by each of the coordinates are markedly different across the coordinates: two are periodic while the other has a ``bath-tub'' appearance. Their entropies are monotone. In the third and fourth examples, we describe applications to regression and classification. In the regression example, we derive the uncertainty measures for a Gaussian process regressor trained on a small number of samples ($n=100$) from the Ishigami function. Despite the small sample size, the GPR is able to accurately reconstruct the HDMR component functions, and the sensitivity density is noisy.  For classification, we find that the orthogonal random measure sensitivity indices depend on the number of incorrect classifications. Maximum entropy of the orthogonal random measure is achieved with equal errors across classes. 

The example of interacting particle systems illustrates how spatiotemporal interaction probabilities may be attained from traffic flow random measure models. For the specific example, the sensitivity distribution exhibits multi-modal behavior, and the random field enables construction of the covariance and its reconstruction through the truncated eigensystem. The example of adaptive randomized controlled trials shows that the injection of counting noise into trial design does not negatively impact type I or II error probabilities and confers an orthogonal design, which enables the sensitivity distribution to be attained in the effects. Adaptation through superpositions of orthogonal dice mitigates the phenomenon of under- or over-powered trials at final analysis time. Dynamic survival analysis of epidemics has the structure of a binomial random measure, whose mean measure is encoded by a differential equation of a compartmental model on a Poisson-type network. Random field models are developed for stock market indicators, enabling computation of the effect of epidemic dynamics on market values. 

Across the examples, we calculate virtually all the quantities analytically and suggest these may be useful and canonical for UQ by random measures. These examples also underscore the significance of the orthogonal mixed binomial processes, that is, those random counting measures whose counting distributions have mean equals variance. For orthogonal $N$, the structural sensitivity indices define a probability measure on partitions, conferring a probabilistic interpretation of uncertainty and enabling calculation of entropy. Therefore we suggest UQ practitioners to employ orthogonal random measures, such as Poisson or the orthogonal dice, as a substitute for empirical random measures, in order to gain probabilistic insight into the UQ exercise. The examples indeed show that, at the cost of increasing the uncertainty of the random measure through orthogonal $N$ ($\kappa$ having positive variance equal to mean), one gets a return of an interpretation of the normalized random measure variances as probabilities.



  \bibliographystyle{apalike}

\newpage
\appendix

\section{Supplementary figures} 

\subsection{Ishigami function}
\begin{figure}[h!]
\centering
\begingroup
\captionsetup[subfigure]{width=3in,font=normalsize}
\subfloat[Component function $g_1$ for $a=7$, $b=0.1$]{\includegraphics[width=3.5in]{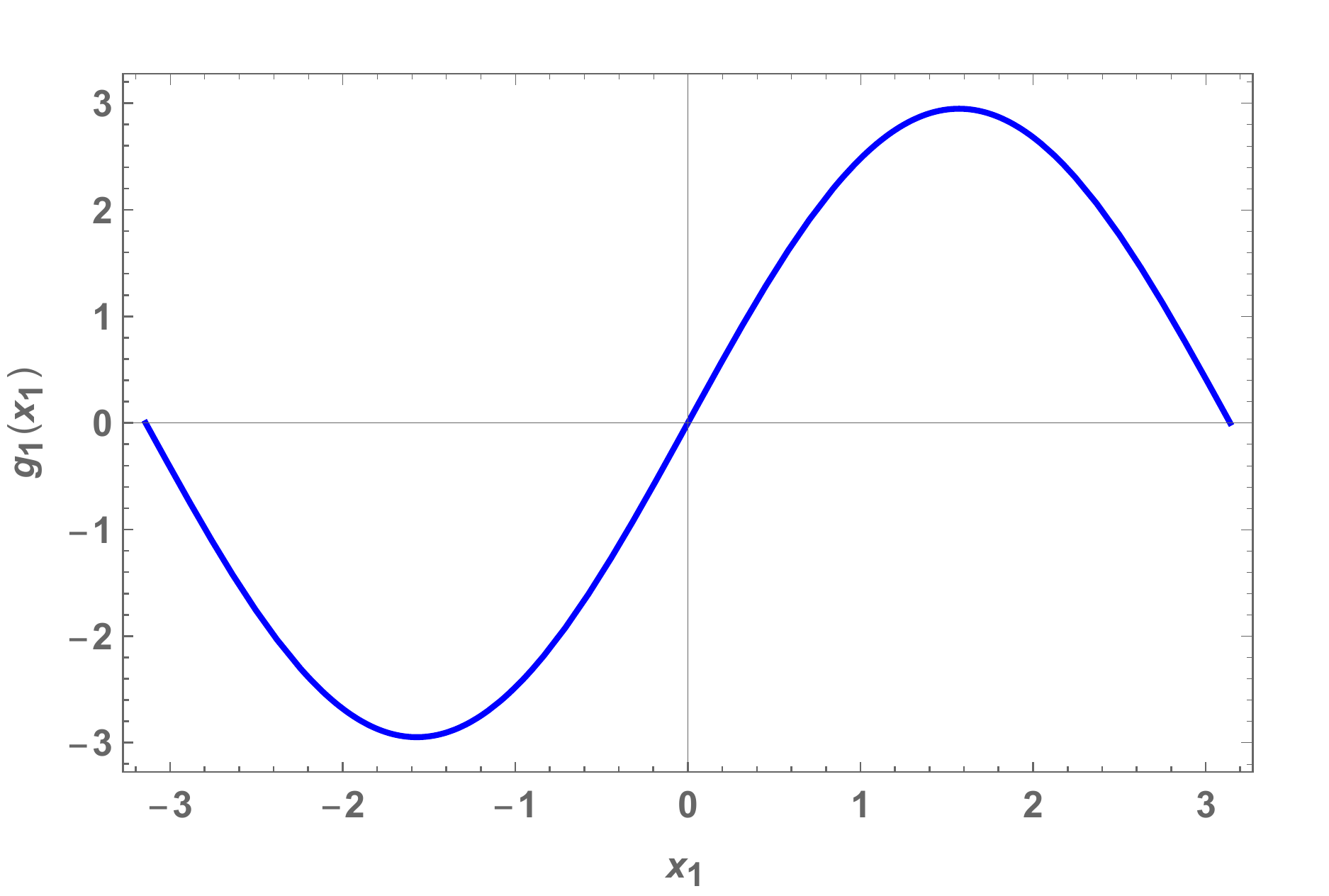}}
\subfloat[Random measure entropy coordinate 1]{\includegraphics[width=3.5in]{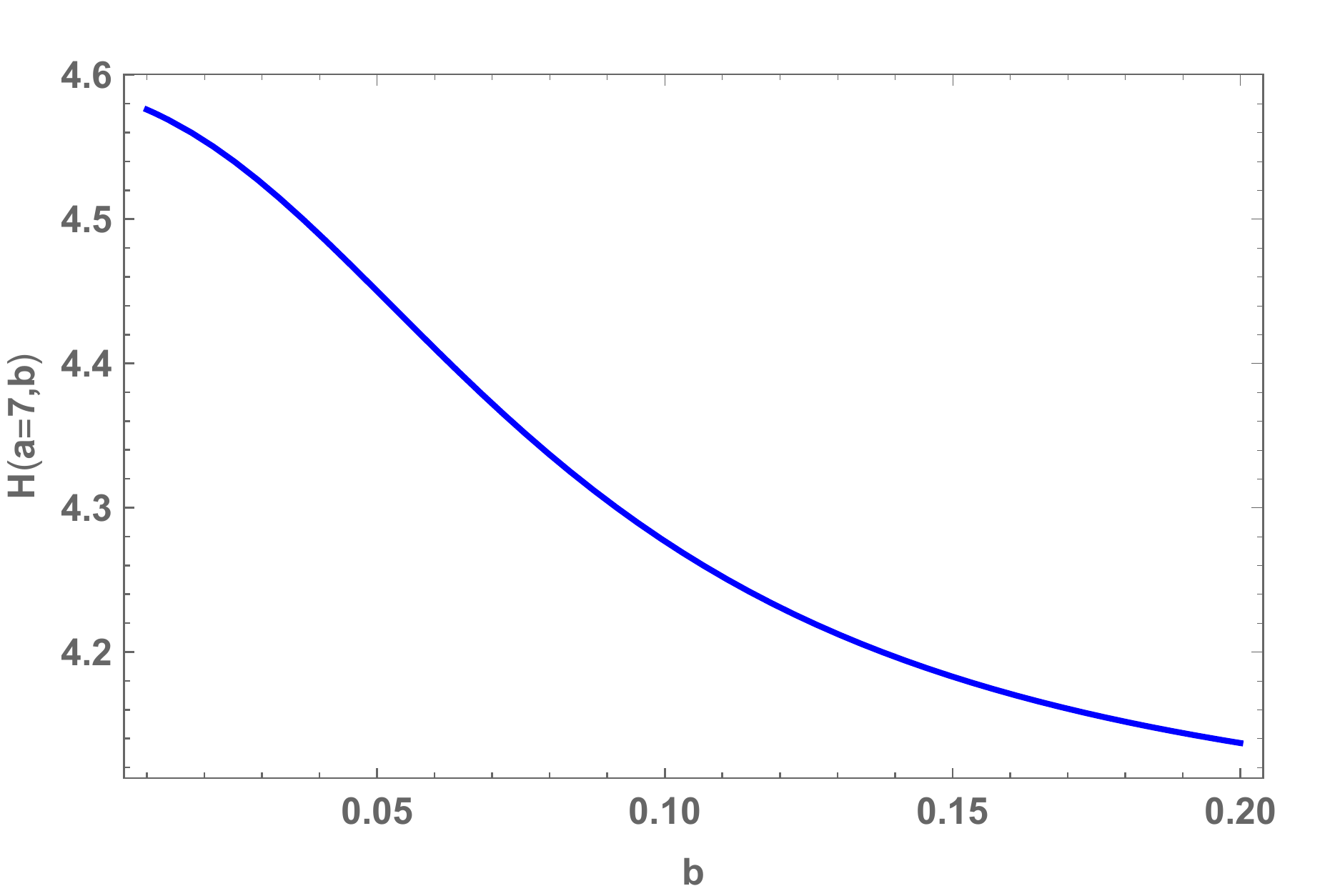}}\\
\subfloat[Component function $g_2$ for $a=7$, $b=0.1$]{\includegraphics[width=3.5in]{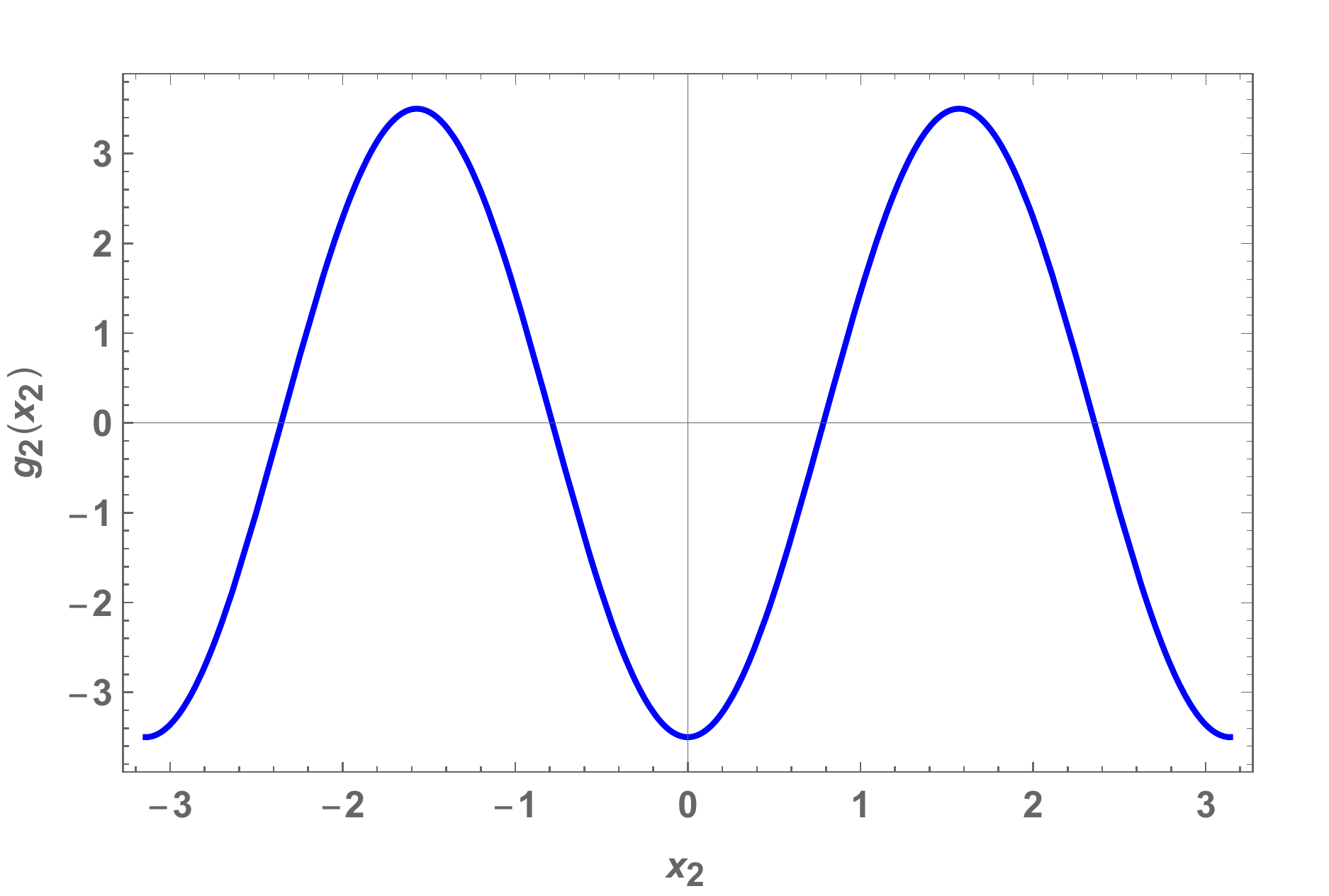}}
\subfloat[Random measure entropy coordinate 2]{\includegraphics[width=3.5in]{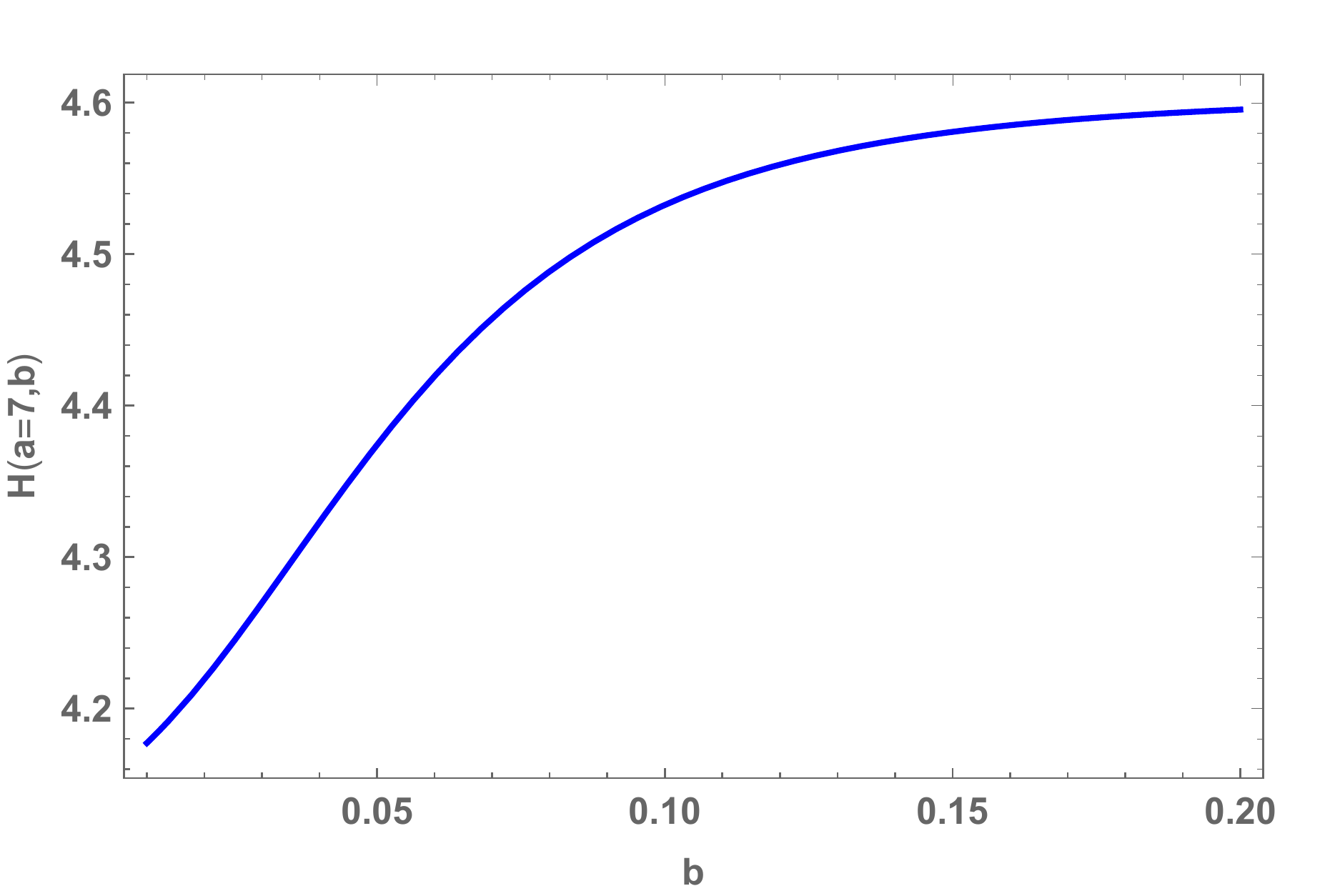}}\\
\subfloat[HDMR entropy $H(a=7,b)$ as function of $b$]{\includegraphics[width=3.5in]{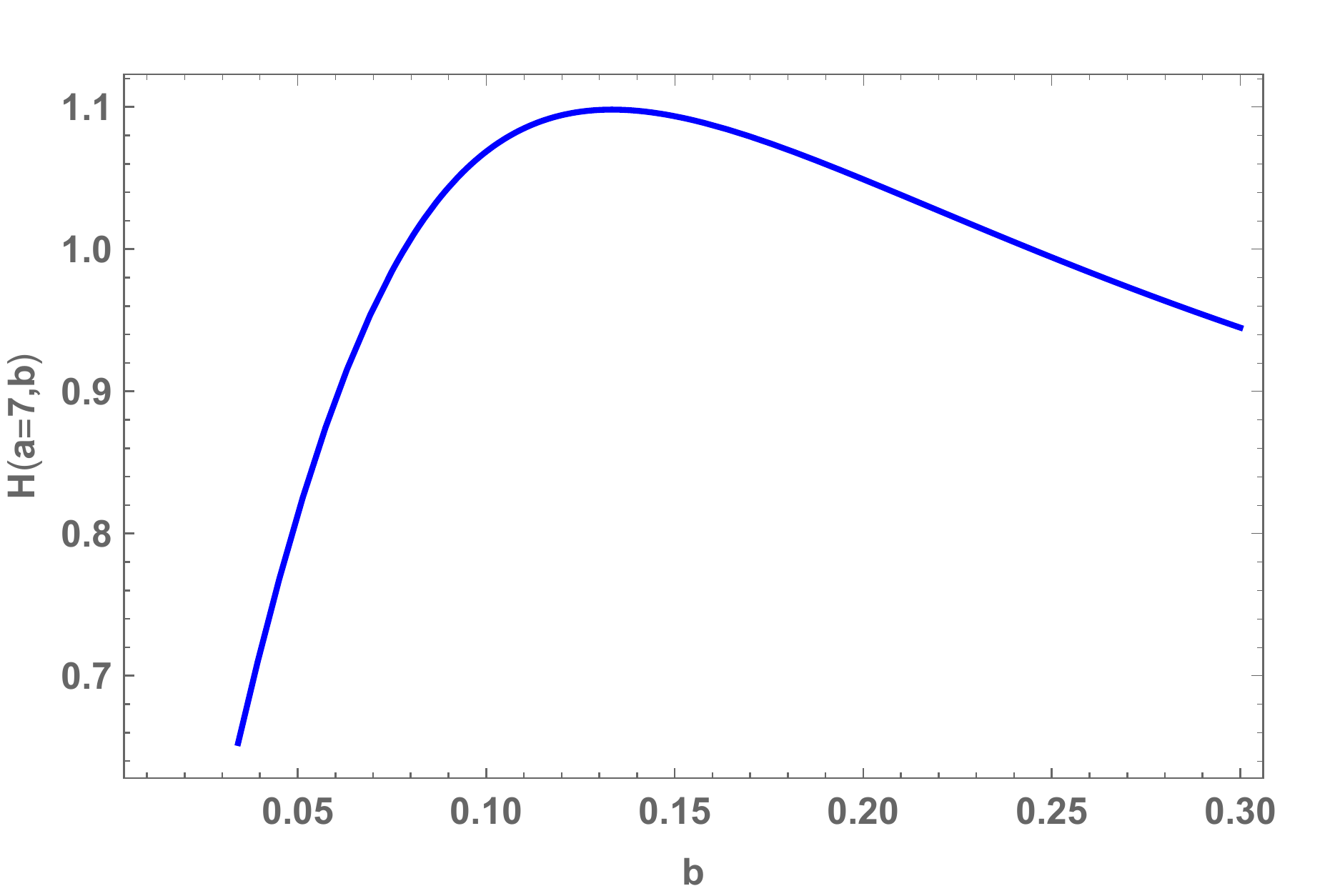}}
\subfloat[Random measure entropy coordinate 3]{\includegraphics[width=3.5in]{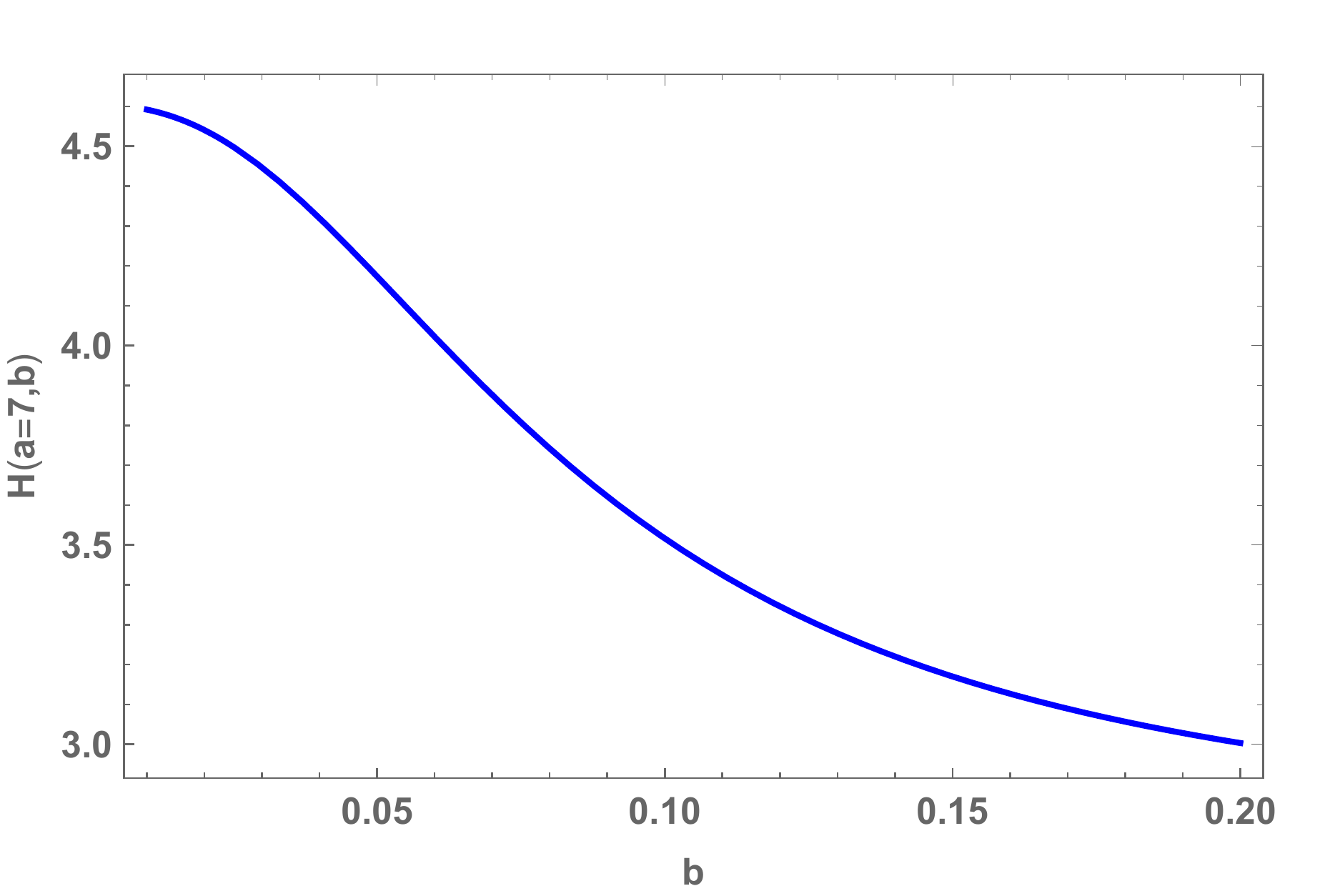}}
\endgroup
\caption{Ishigami first-order component functions $g_1$ and $g_2$ for $a=7$ and $b=0.1$; Ishigami function orthogonal random measure entropy for partitions by coordinates as a function of $b$; HDMR entropy as a function of $b$ for $a=7$}\label{fig:4}
\end{figure}

\begin{figure}[h!]
\centering
\begingroup
\captionsetup[subfigure]{width=3in,font=normalsize}
\subfloat[$(\Var f_d)$ for coordinate 1]{\includegraphics[width=3.5in]{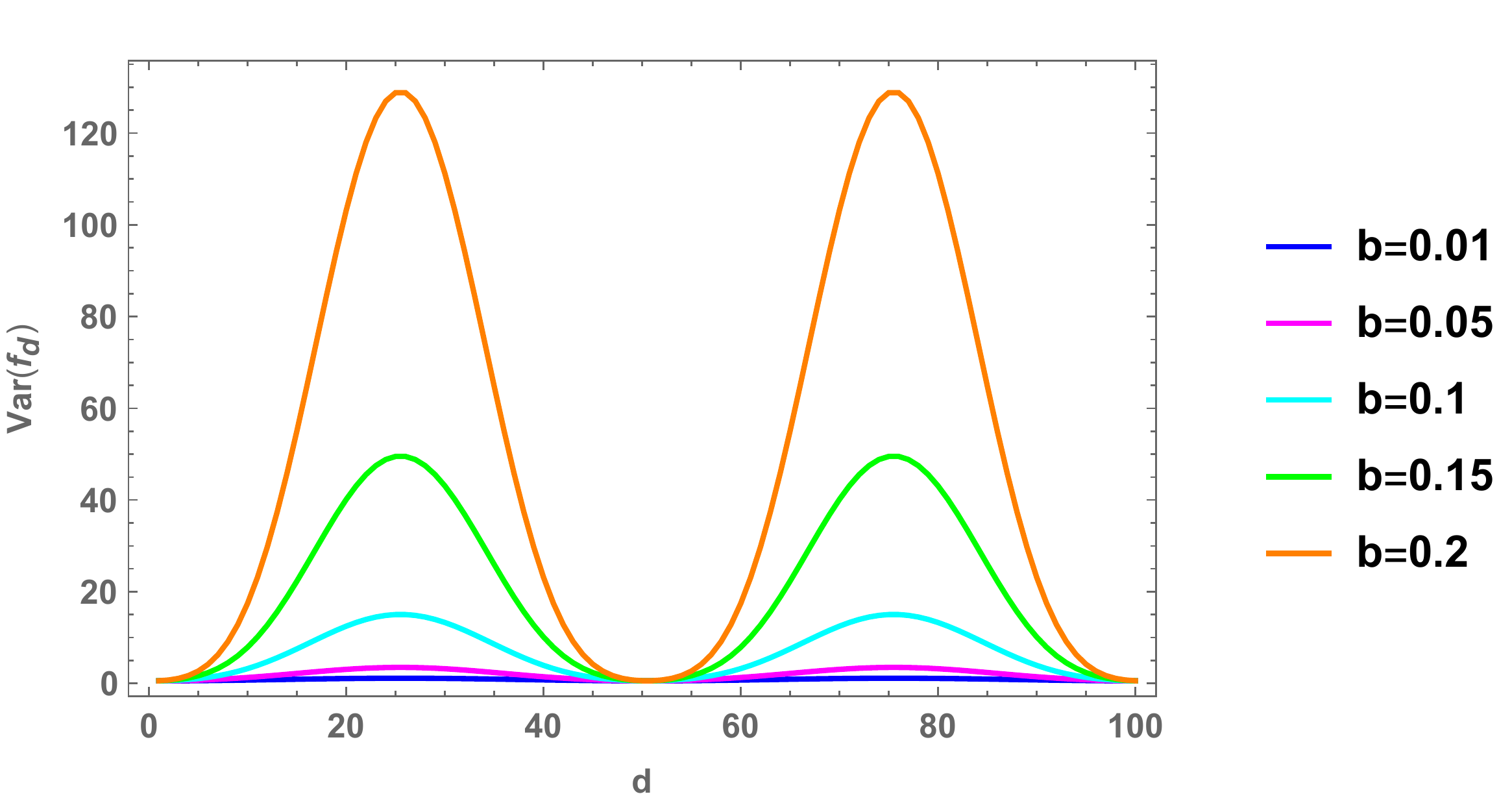}}
\subfloat[$(\amsbb{S}_d^a)$ for coordinate 1]{\includegraphics[width=3.5in]{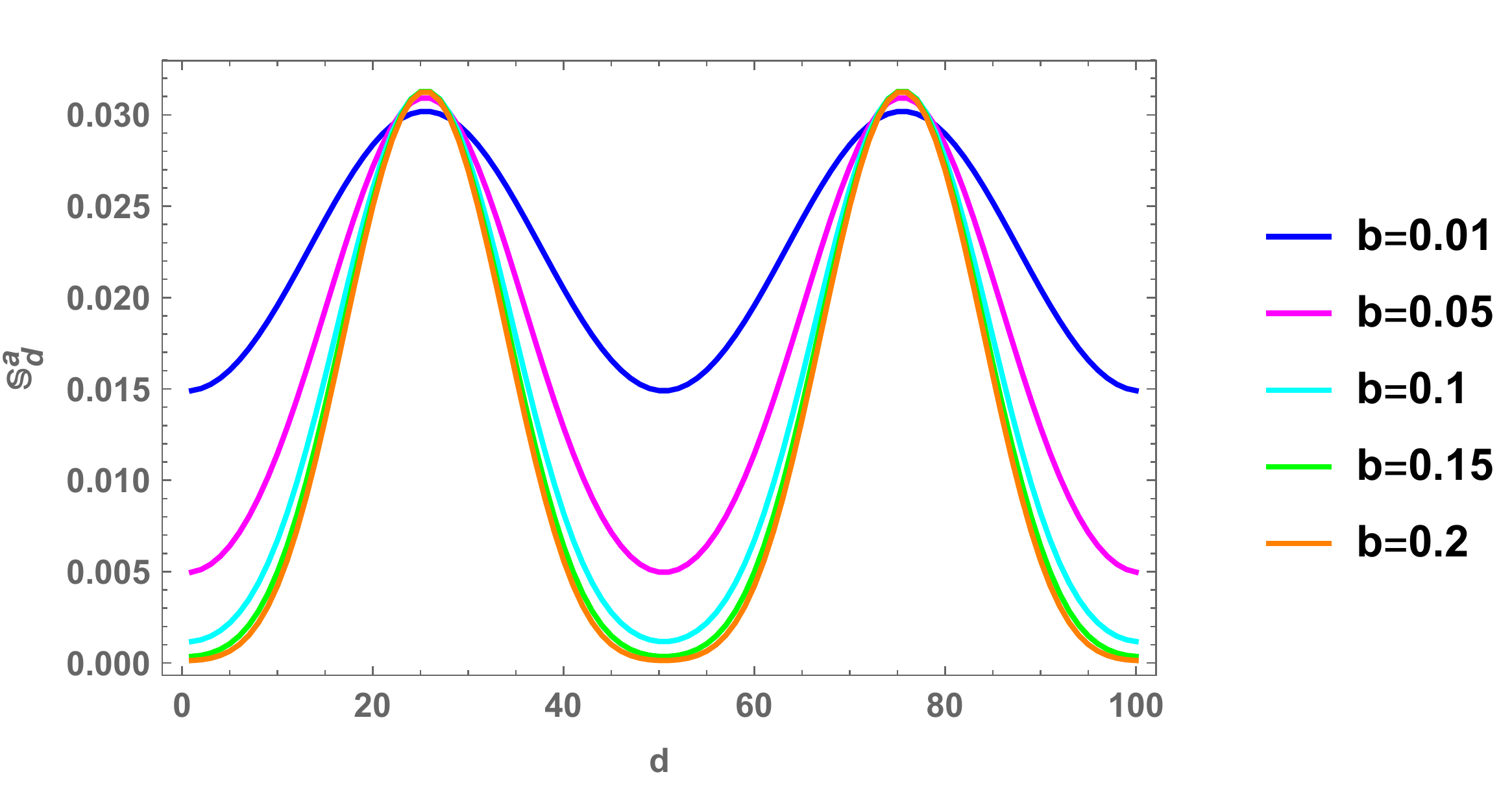}}\\
\subfloat[$(\Var f_d)$ for coordinate 2]{\includegraphics[width=3.5in]{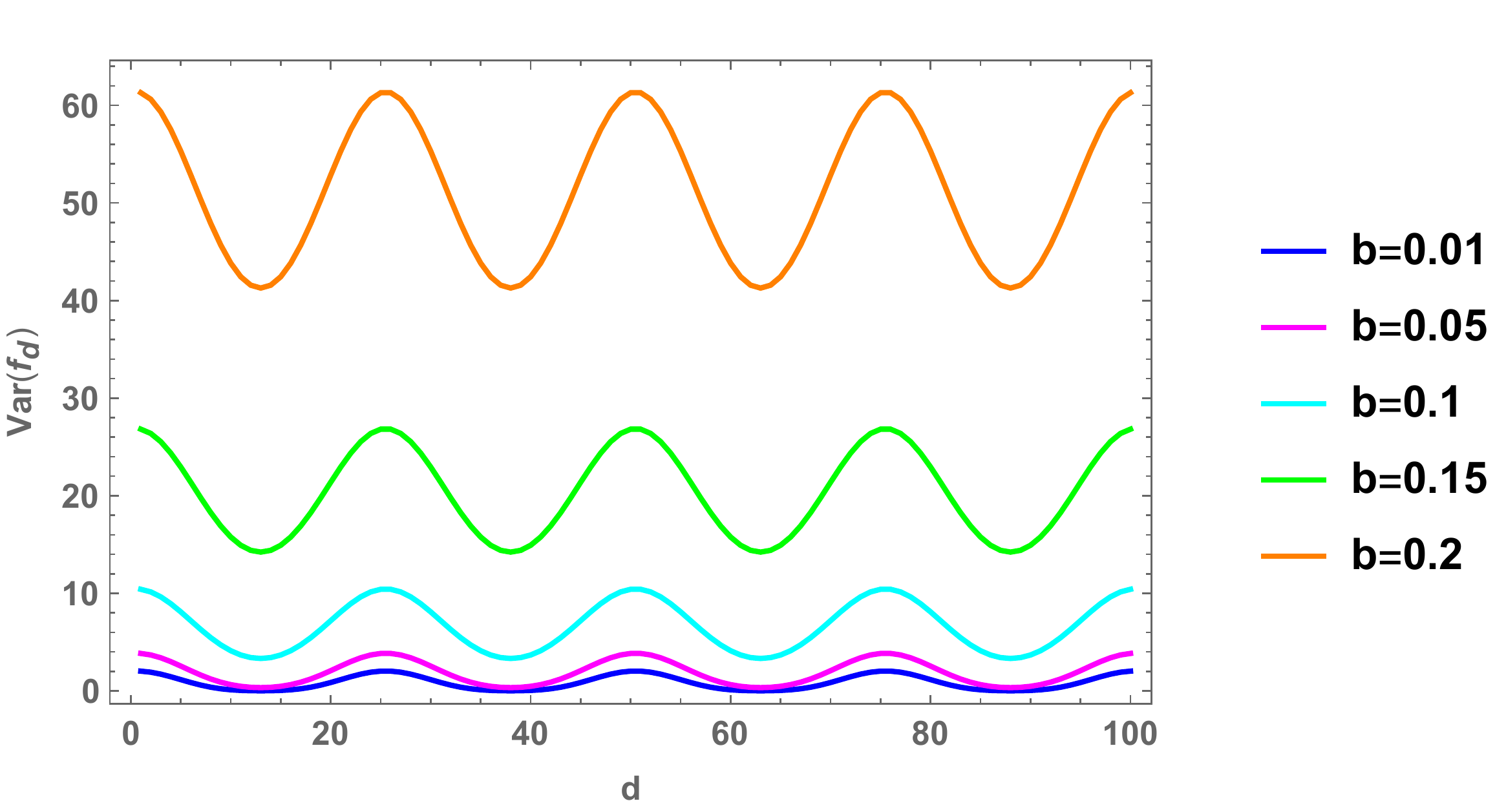}}
\subfloat[$(\amsbb{S}_d^a)$ for coordinate 2]{\includegraphics[width=3.5in]{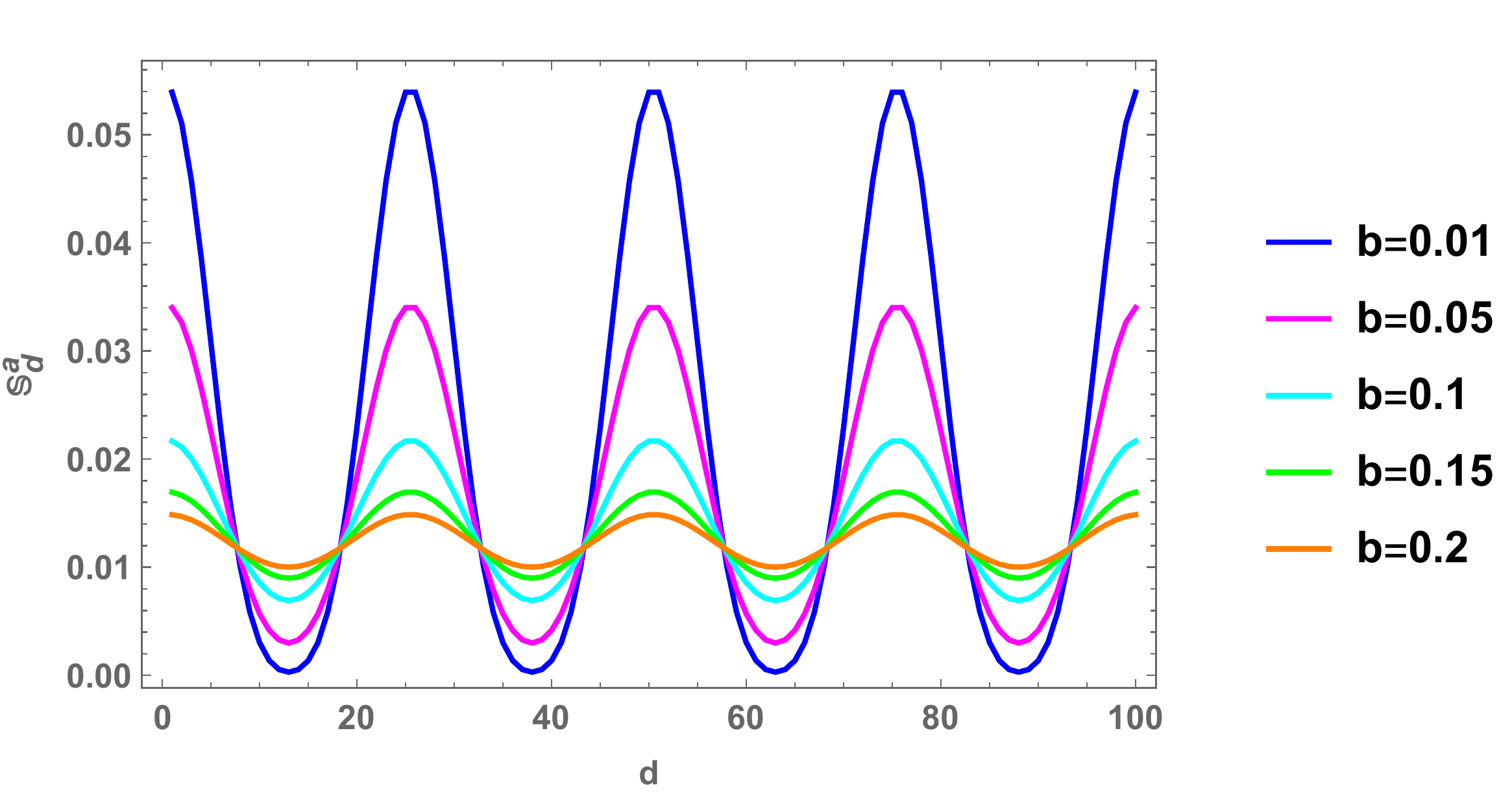}}\\
\subfloat[$(\Var f_d)$ for coordinate 3]{\includegraphics[width=3.5in]{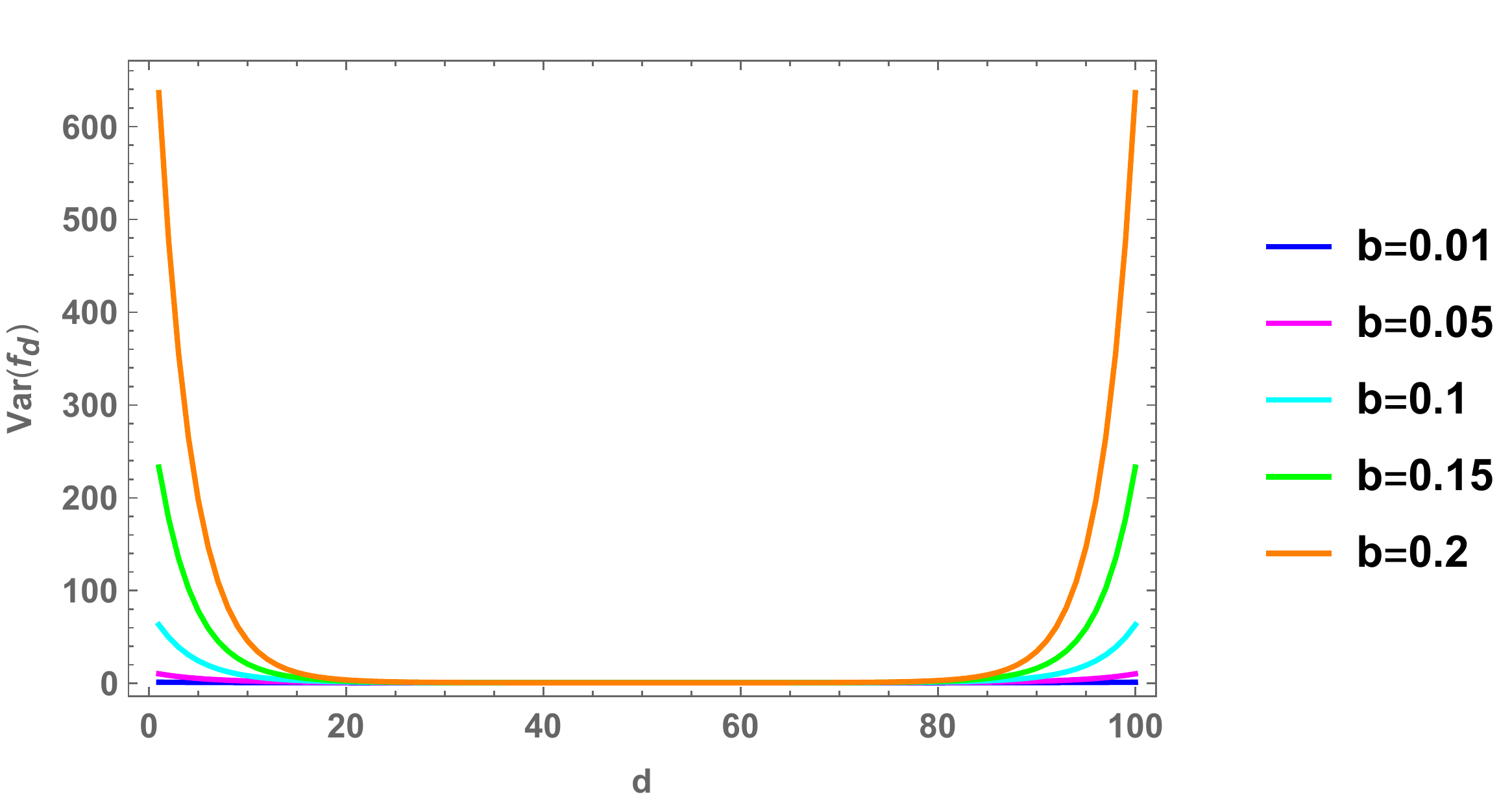}}
\subfloat[$(\amsbb{S}_d^a)$ for coordinate 3]{\includegraphics[width=3.5in]{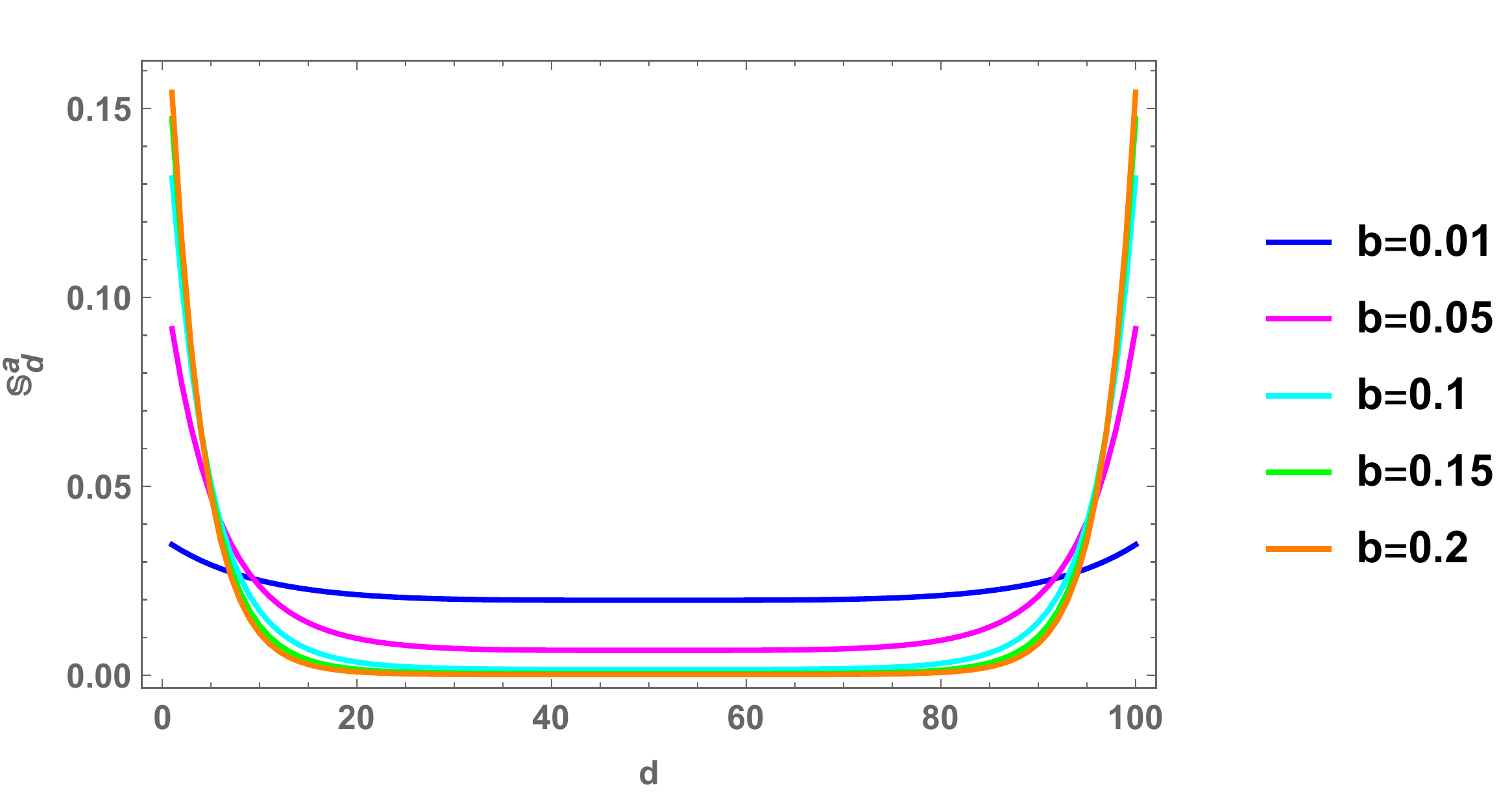}}\\
\endgroup
\caption{Ishigami function for Dirac $N$: variances $(\Var f_d)$ and structural sensitivity indices $(\amsbb{S}_d^a)$ for partitions by coordinates into $100$ intervals for $a=7$ and $b\in\{0.01,0.05,0.1,0.15,0.2\}$}\label{fig:4a}
\end{figure}

\FloatBarrier

\section{Fubini's Theorem} A foundational theorem of random measures is Fubini's Theorem. 
\begin{samepage}
\begin{thm}[Measure-kernel-function] Let $Q$ be a transition kernel from $(E,\mathcal{E})$ into $(F,\mathcal{F})$. Then \[Qf(x) = \int_F Q(x,\D y)f(y)\for x\in E\] defines a function $Qf$ that is in $\mathcal{E}_+$ for every function $f$ in $\mathcal{F}_+$; \[\nu Q(B)=\int_E\nu(\D x)Q(x,B)\for B\in\mathcal{F}\] defines a measure $\nu Q$ on $(F,\mathcal{F})$ for each measure $\nu$ on $(E,\mathcal{E})$; and \[(\nu Q)f = \nu(Qf)=\int_E\nu(\D x)\int_F Q(x,\D y)f(y)\] for every measure $\nu$ on $(E,\mathcal{E})$ and function $f$ in $\mathcal{F}_+$.  
\end{thm}
\end{samepage}

\begin{re}[Random measure] The random measure $N=(\kappa,\nu)$ on $(E,\mathcal{E})$ is a transition kernel from $(\Omega,\mathcal{H})$ into $(E,\mathcal{E})$. Then \[Nf(\omega) = \int_E N(\omega,\D x)f(x)\for \omega\in\Omega\] defines a positive random variable $Nf$ for every function $f$ in $\mathcal{E}_+$; \[c\nu(B)=\E N(B)=\int_\Omega\P(\D \omega)N(\omega,B)\for B\in\mathcal{E}\] defines a measure $c\nu=\E N$ on $(E,\mathcal{E})$ for each measure $\P$ on $(\Omega,\mathcal{H})$ called the mean or intensity measure; and \[\E Nf = c\nu f=\int_\Omega\P(\D\omega)\int_E N(\omega,\D x)f(x)\] for every measure $\P$ on $(\Omega,\mathcal{H})$ and function $f$ in $\mathcal{E}_+$.
\end{re}

\section{Additional examples}

\subsection{Symmetric polynomial with correlation}\label{sec:poly} As referenced in Table~\ref{tab:0}, here we consider correlation in a bivariate symmetric polynomial. Suppose we have $E=\R^2$ and $\nu=\text{Gaussian}([0,0],\begin{bmatrix}1&\rho\\\rho&1\end{bmatrix})$ and the polynomial $g(x_1,x_2)=1+x_1+x_2+x_1x_2$. This has mean and variance \begin{align*}\E g &= 1+ \rho\\\Var g &=3+2(\rho+1)\end{align*} This admits a HDMR analysis \begin{align*}g_1(x_1) &= x_1 + (\frac{\rho}{\rho^2+1})(x_1^2-1)\\g_2(x_2) &= x_2 + (\frac{\rho}{\rho^2+1})(x_2^2-1)\\g_{12}(x_1,x_2)&=\frac{\rho^2 x_1x_2-\rho \left(x_1^2+x_2^2-1\right)+x_1x_2-\rho^3}{\rho^2+1}\\\Var g_1 &=\Var g_2 = \frac{\rho^4+4 \rho^2+1}{\left(\rho^2+1\right)^2}\\\Cov(g_1,g_2)&=\frac{\rho\left(\rho^2(\rho (\rho+2)+2)+1\right)}{\left(\rho^2+1\right)^2}\\\Var g_{12}&=\frac{\left(1-\rho^2\right)^2}{\rho^2+1}\\\amsbb{S}_1^a &= \amsbb{S}_2^a =\frac{\rho^4+4 \rho^2+1}{\left(\rho^2+1\right)^2 (3+\rho (\rho+2))}\\\amsbb{S}_1^b &=\amsbb{S}_2^b = \frac{\rho \left(\rho^2(\rho (\rho+2)+2) +1\right)}{\left(\rho^2+1\right)^2 (3+\rho (\rho+2))}\\\amsbb{S}_{12}^a&=\frac{\left(1-\rho^2\right)^2}{\left(\rho^2+1\right) (3+\rho (\rho+2))}\\ ED&=\frac{2 \left(\rho^4+\rho^3+\rho^2+\rho+2\right)}{\left(\rho^2+1\right) (3+\rho (\rho+2))}\end{align*} with decomposition of variance \[\Var g = \Var g_1 + \Var g_2 + 2\Cov(g_1,g_2) + \Var g_{12}\] where we note that $\Cov(g_{12},g_1)=\Cov(g_{12},g_2)=0$ due to hierarchical orthogonality, so $\amsbb{S}_{12}^b=0$. 

\begin{figure}[h!]
\centering
\begingroup
\captionsetup[subfigure]{width=3in,font=normalsize}
\subfloat[$\amsbb{S}_1^a(\rho)+\amsbb{S}_2^a(\rho)+\amsbb{S}_1^b(\rho)+\amsbb{S}_2^b(\rho)$]{\includegraphics[width=3.5in]{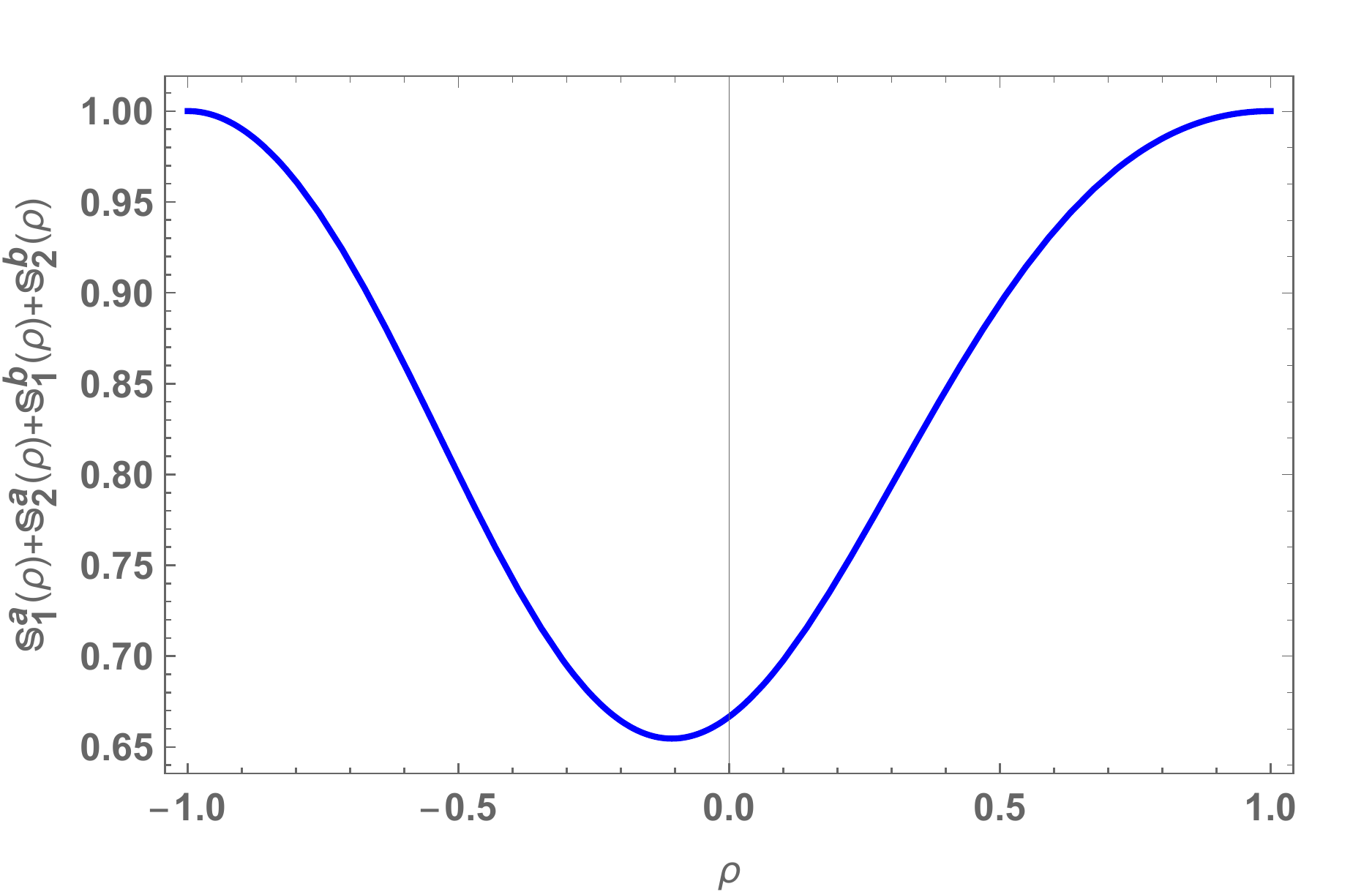}}
\subfloat[$\amsbb{S}_{12}^a(\rho)$]{\includegraphics[width=3.5in]{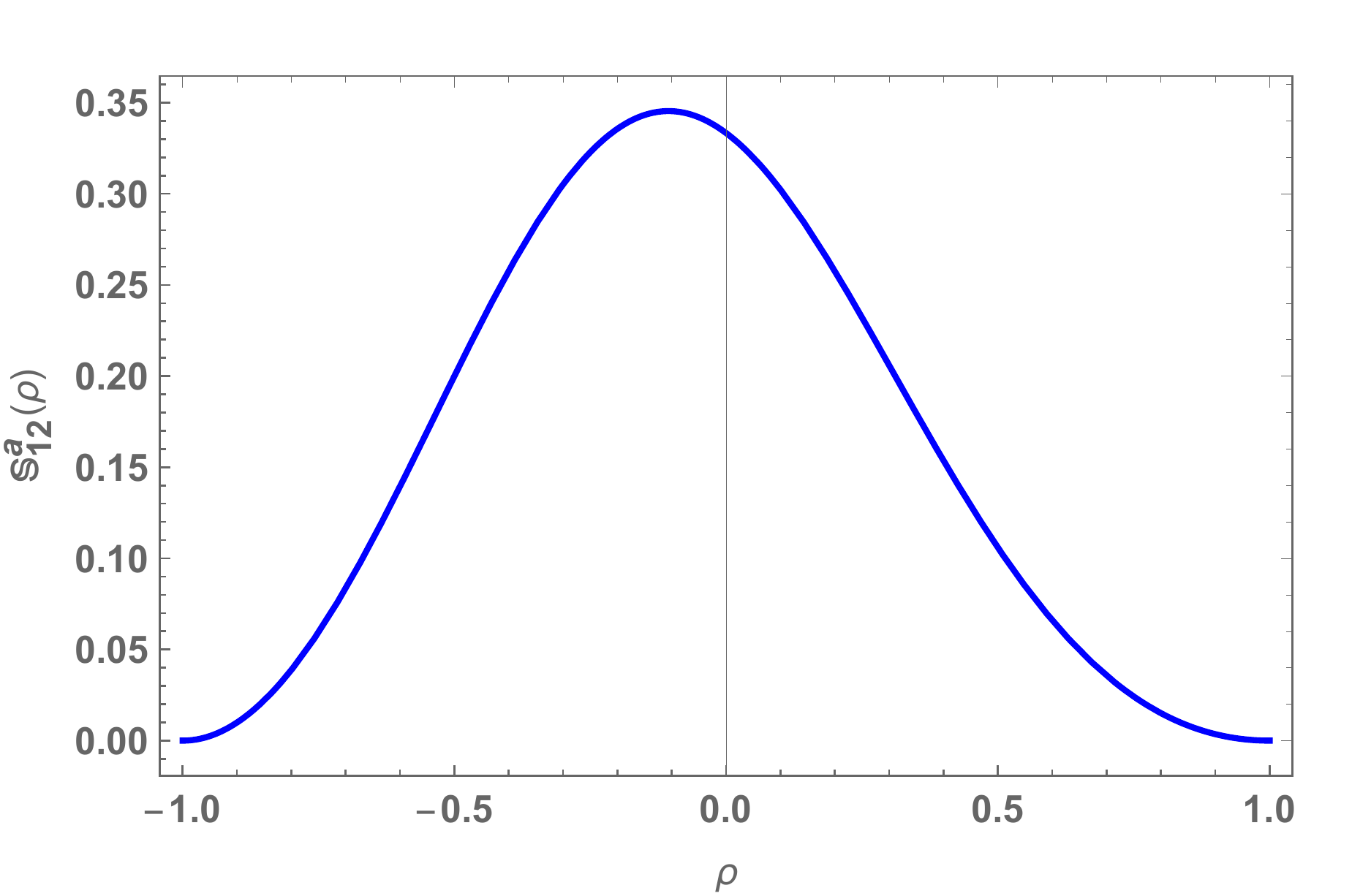}}\\
\subfloat[Effective dimension as a function of $\rho$]{\includegraphics[width=3.5in]{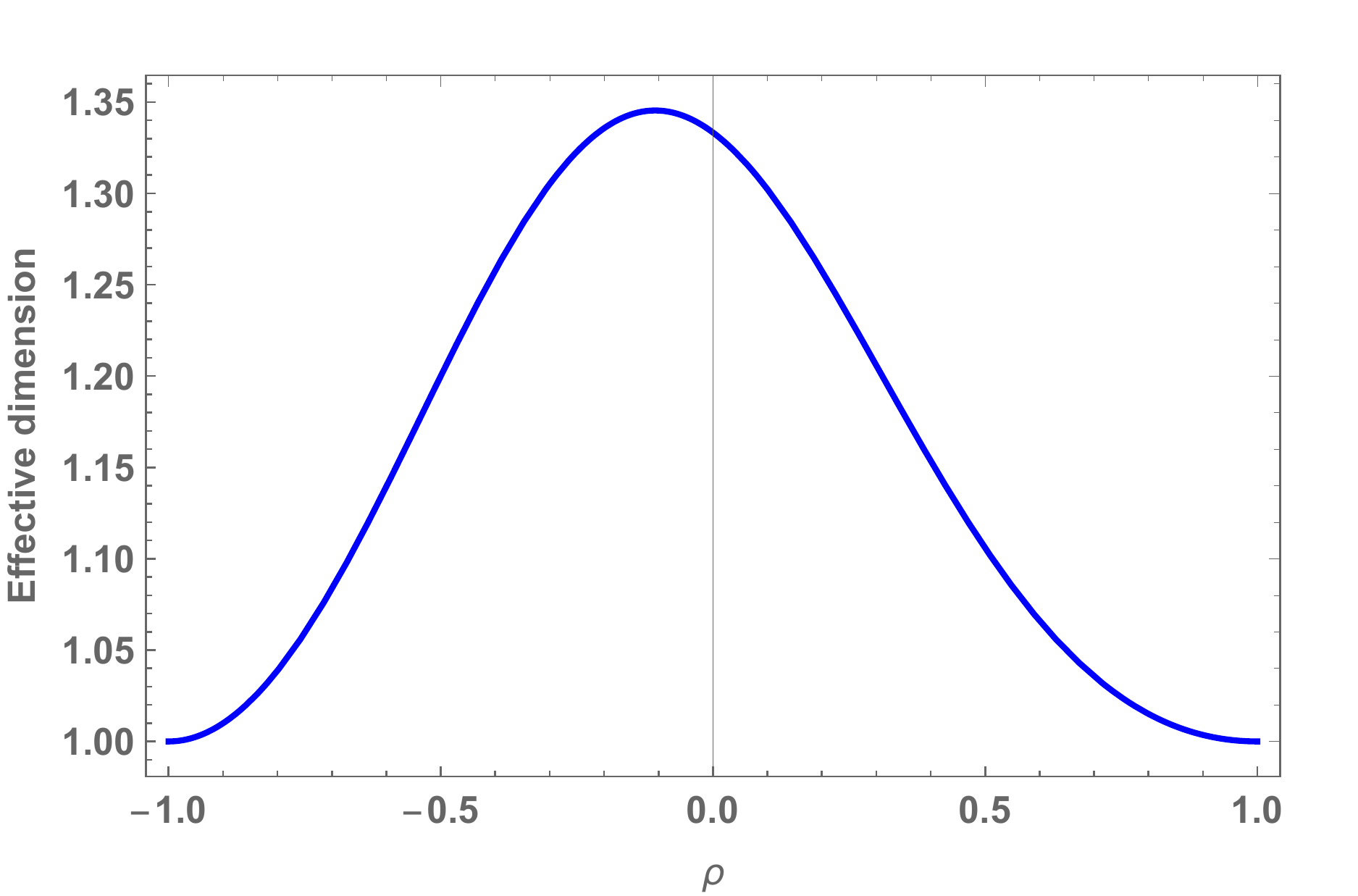}}
\subfloat[Entropy (on dimension) as a function of $\rho$]{\includegraphics[width=3.5in]{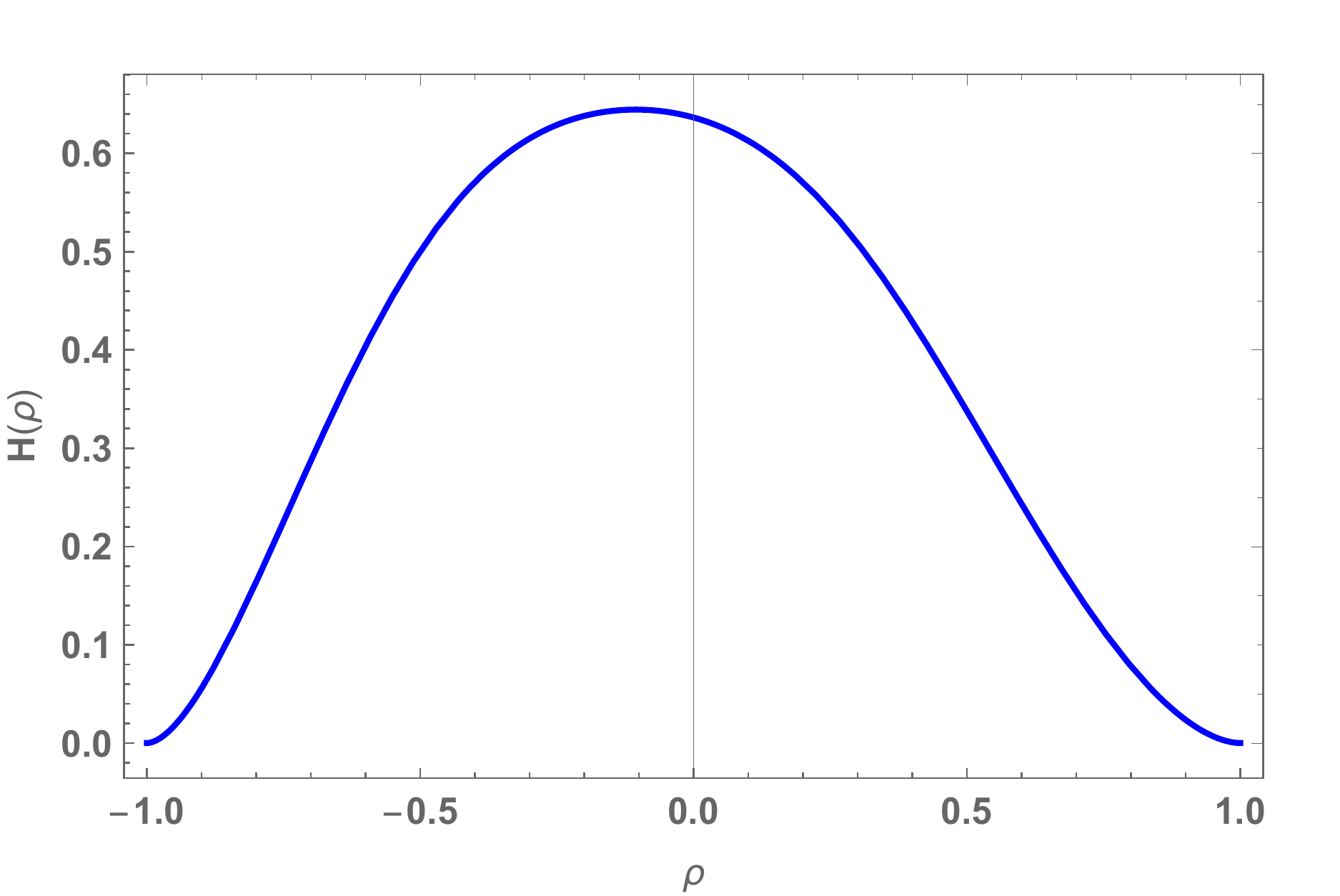}}\\ 
\endgroup
\caption{First and second order HDMR sensitivity indices for polynomial model with correlated inputs; effective dimension; entropy on subspace dimension; all as a function of correlation coefficient $\rho$}\label{fig:6}
\end{figure}

In Figure~\ref{fig:6} we show the first and second order HDMR sensitivity indices. The sensitivity indices show behavior whereby the second-order subspace vanishes and first-order achieve unity as $|\rho|\rightarrow1$. The effective dimension is unimodal, so too is the entropy. Interestingly these quantities all have local maxima and minima around $\rho\approx-0.106$ and none are symmetric, despite $g$ being symmetric.

\FloatBarrier

For the random measure $N=(\kappa,\nu)$ on $(E,\mathcal{E})$ and for the function $f=(g-\E g)^2\in\mathcal{E}_+$, we have the decomposition of $Nf$ of the mean \begin{equation}\E Nf = c\nu f = c \Var g = c\left(\Var g_1 + \Var g_2 + 2\Cov(g_1,g_2) + \Var g_{12}\right)\end{equation} and for the disjoint partition $\{A,\dotsb,B\}$ of $E$, putting $f_d=f\ind{D}$ for $D\in\{A,\dotsb,B\}$, the decomposition of the variance as \begin{align}\Var Nf &=c\nu f^2 + (\delta^2-c)(\nu f)^2\nonumber\\&= \sum_{D\in\{A,\dotsb,B\}} \Var Nf_d + \sum_{D_i\ne D_j\in\{A,\dotsb,B\}} \Cov(Nf_{d_i},Nf_{d_j})\nonumber\\&= \sum_{D\in\{A,\dotsb,B\}}(c\nu f_d^2 + (\delta^2-c)(\nu f_d)^2) + \sum_{D_i\ne D_j\in\{A,\dotsb,B\}}(\delta^2-c)\nu f_{d_i}\nu f_{d_j} \end{align} The second moment is \[\nu f^2 =3 \rho (\rho (3\rho (\rho+4)+46)+36)+57 \] The Laplace functional is defined through \[\nu e^{-f} = \int_{\R^2}\D\nu(x_1,x_2)e^{-f(x_1,x_2)}=\int_{\R^2}\D x_1\D x_2\frac{e^{-\left(\frac{1}{2(1-\rho^2)}(x_1^2+x_2^2-2 \rho x_1 x_2)+(x_1 x_2+x_1+x_2-\rho)^2\right)}}{2 \pi  \sqrt{1-p^2}}\] The sensitivity distribution of coordinate $i\in\{1,2\}$ is analytically computed.

\begin{figure}[h!]
\centering
\begingroup
\captionsetup[subfigure]{width=3in,font=normalsize}
\subfloat[$\amsbb{S}^a_i(x_i,\rho)$]{\includegraphics[width=3.5in]{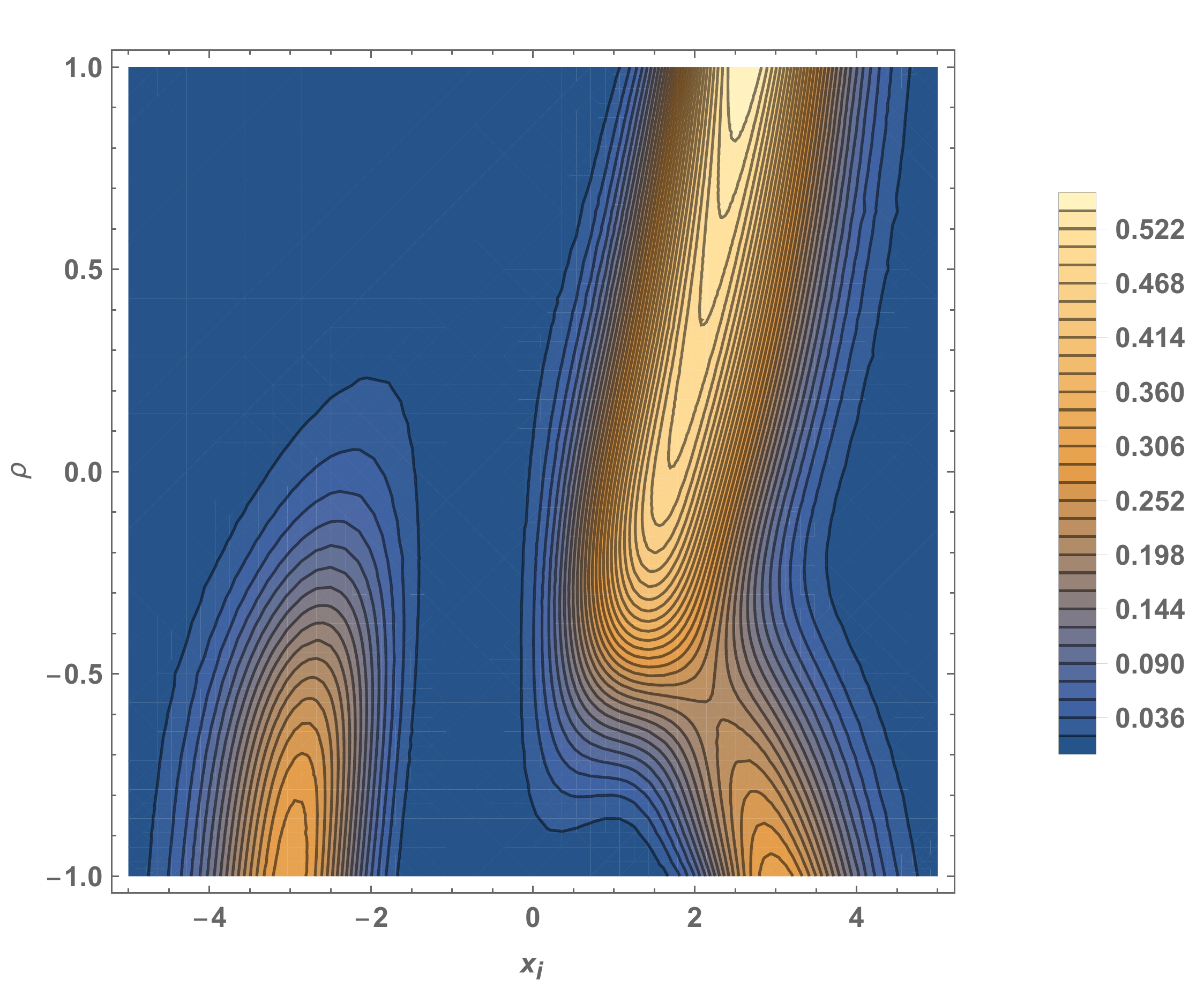}}
\subfloat[Entropy of $\amsbb{S}^a_i$ as a function of $\rho$]{\includegraphics[width=3.5in]{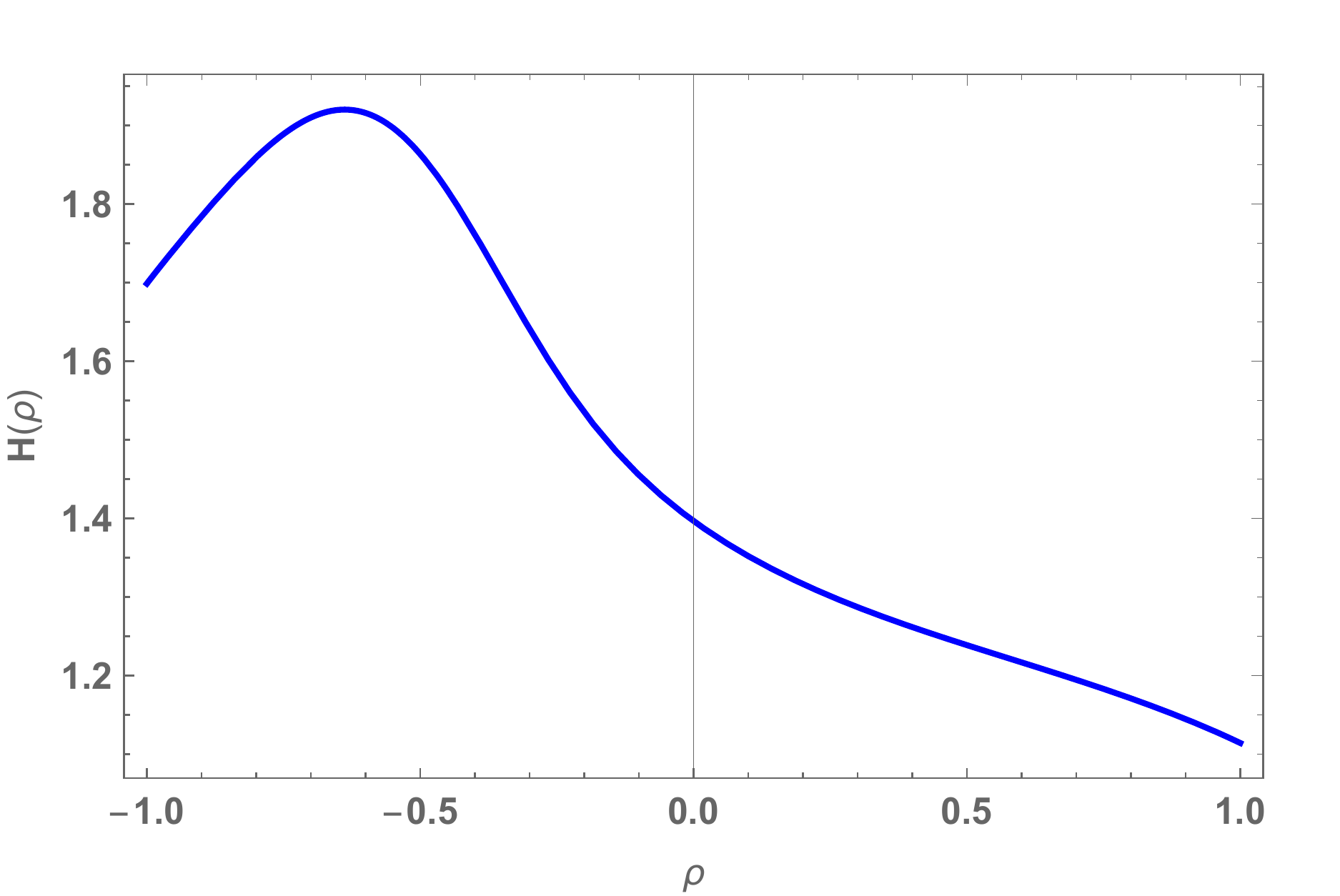}}\\
\subfloat[$\amsbb{S}^a_i$ and maximum entropy]{\includegraphics[width=4in]{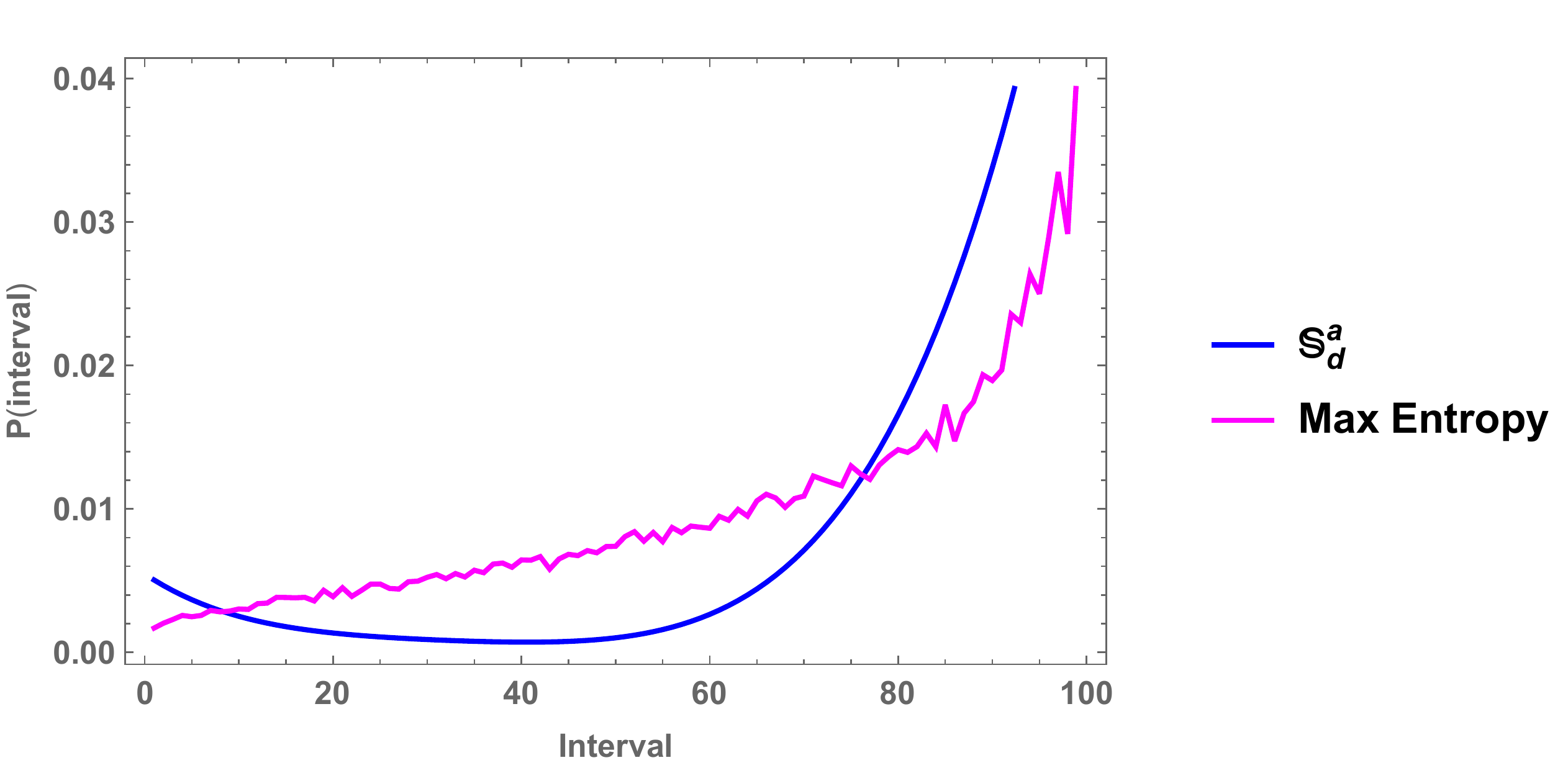}}\\
\endgroup
\caption{Bivariate polynomial for orthogonal $N$: structural sensitivity density $\amsbb{S}_i^a$ for coordinate $i$, entropy, and maximum entropy}\label{fig:7}
\end{figure}

In Figure~\ref{fig:7} we plot the density $\amsbb{S}^a_i$ for coordinate $i$ and its entropy. The density shows complex behavior, moving from bimodal to unimodal as $\rho$ ranges from $-1$ to $1$. Entropy has a maximum around $\rho\approx -0.64$. For the maximum entropy, we change the measure to a product of uniform distributions $\nu=\text{Uniform}[-1,1]\times\text{Uniform}[-1,1]$ so that $g$ takes values on $[0,4]$. The maximum entropy distribution is approximated using Metropolis-Hastings, where 125,000 samples were simulated, with the first 25,000 discarded as burn-in. We see that the coordinate uncertainty by sensitivity indices approximates the maximum entropy distribution in the coordinate.

\FloatBarrier 

\subsection{Graph property}\label{sec:graph} 

As referenced in Table~\ref{tab:0}, let $G(n,p)$ be the  Erd\H{o}s-Reny\'{i} (ER) graph space with $n$ vertices and edge probability $p$. We have $E=\{0,1\}^{\binom{n}{2}}$ with measure $\nu=\prod_i\nu_i$ where $\nu_i=\text{Bernoulli}(p)$ \[\nu\{x\} =\prod_i\nu_i\{x_i\}=\prod_i p^{x_i}(1-p)^{1-x_i}=p^{\sum_i x_i}(1-p)^{\binom{n}{2}-\sum_i x_i}\for x\in E\] Hence the law $\nu$ is defined in terms of the number of edges. Denote $E_i = \{0,1\}^i$ for $i\in\{1,\dotsb,\binom{n}{2}\}$. Suppose we partition $E$ into sets by the number of edges $\{A_0,\dotsb,A_{\binom{n}{2}}\}$ where $A_0=\{(0,\dotsb,0)\}$ with zero edges and $A_{\binom{n}{2}}=\{(1,\dotsb,1)\}$ with $\binom{n}{2}$ edges. Putting $y=\sum_ix_i$, we define \[\nu(A_y) =\binom{\binom{n}{2}}{y} p^y(1-p)^{\binom{n}{2}-y}=\text{Binomial}(\binom{n}{2},p)\for y\in\{0,1,\dotsb,\binom{n}{2}\}\] so the distribution of the number of edges is binomial. Suppose for each $x\in E$, we compute the spectral gap $g(x)\in\R_+$ of the graph Laplacian, where we define $g((0,\dotsb,0))=0$. We compute the mean and variance as \begin{align*}\E g &= \sum_{x\in E}\nu\{x\}g(x)\\\Var g&=\sum_{x\in E}\nu\{x\}g^2(x)-(\nu g)^2\end{align*}

\begin{figure}[h!]
\centering
\begingroup
\captionsetup[subfigure]{width=3in,font=normalsize}
\subfloat[$\E g$]{\includegraphics[width=3.5in]{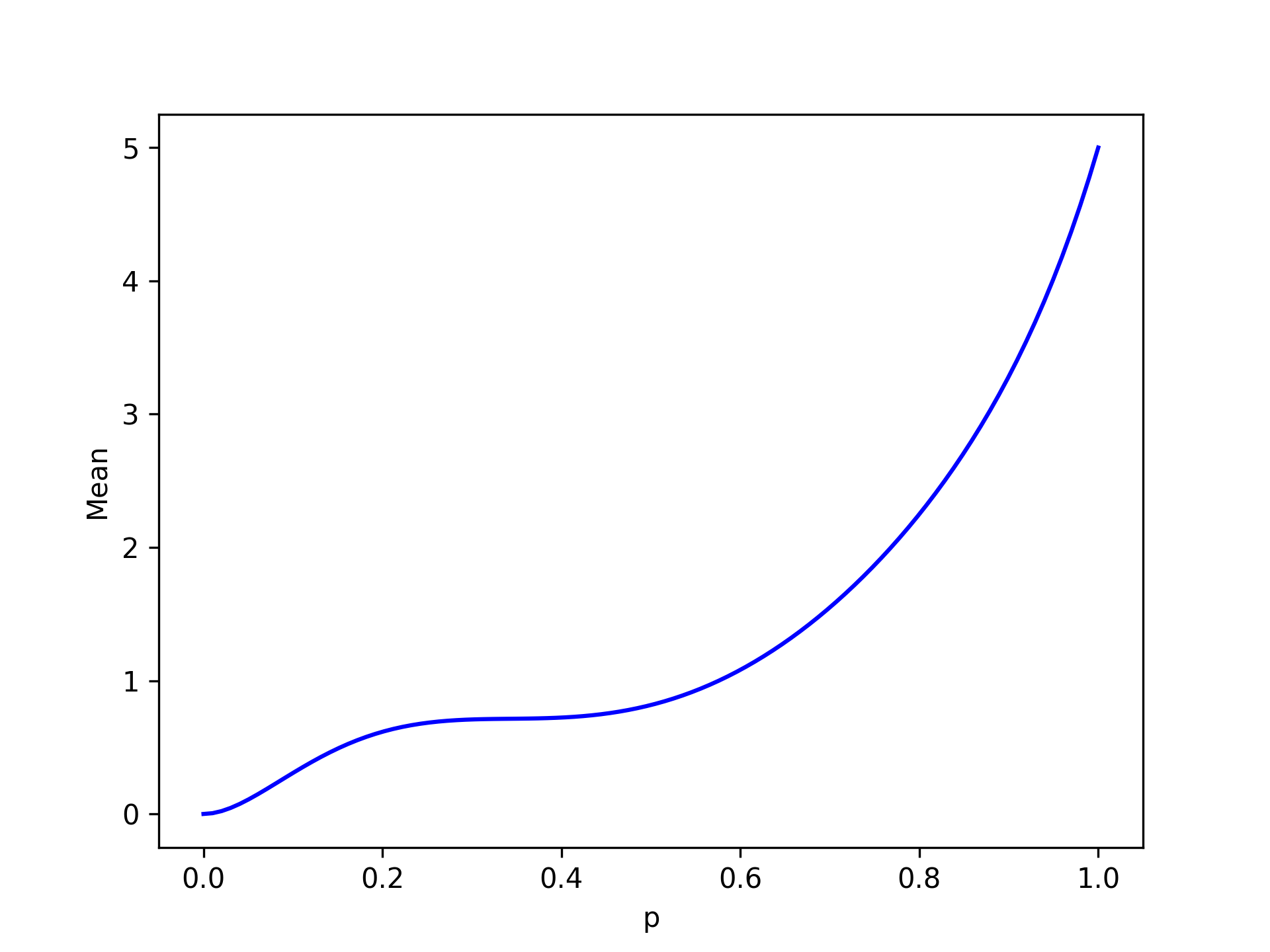}}
\subfloat[$\Var g$]{\includegraphics[width=3.5in]{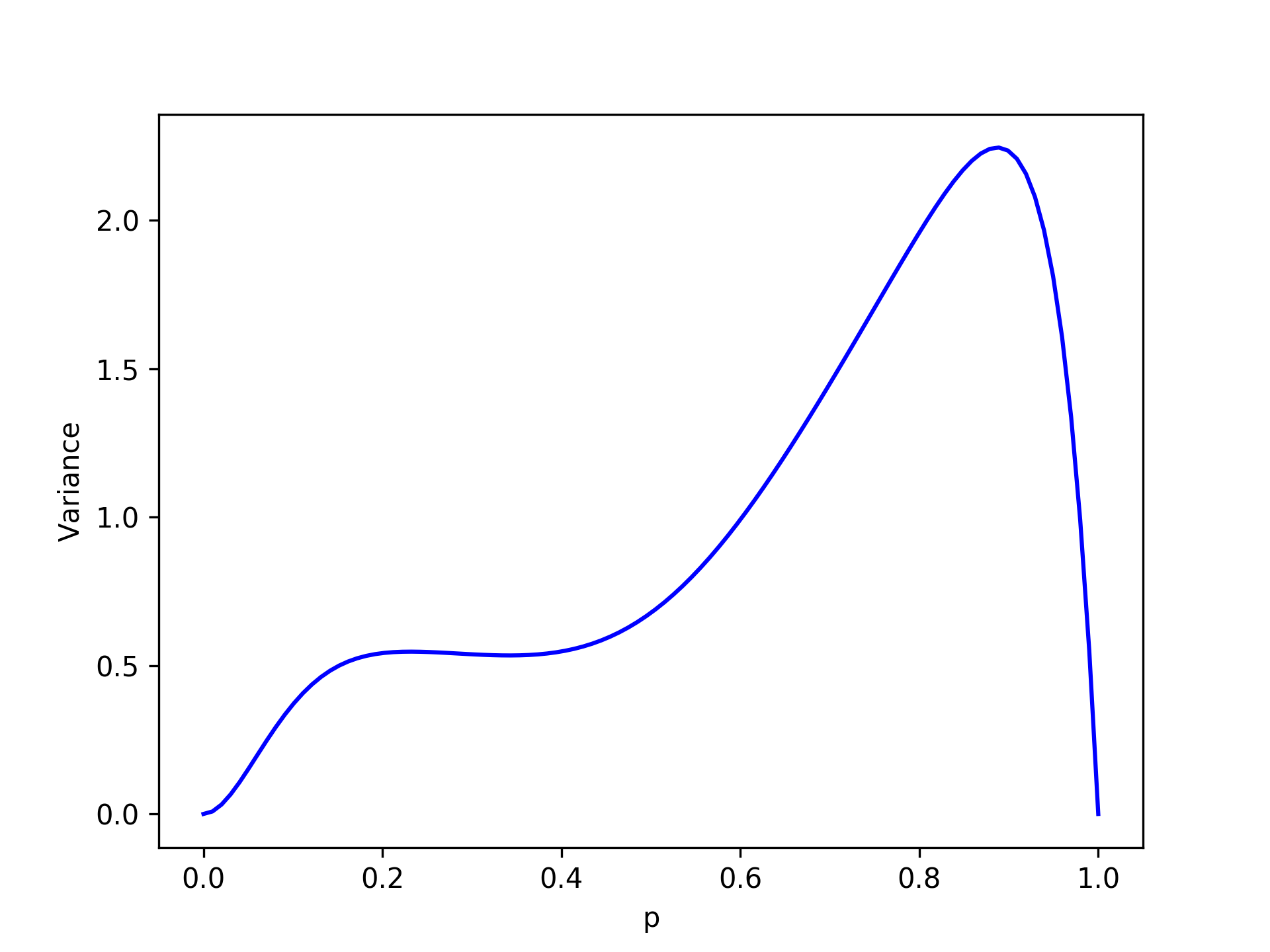}}\\
\subfloat[$(\amsbb{S}_u)$]{\includegraphics[width=3.5in]{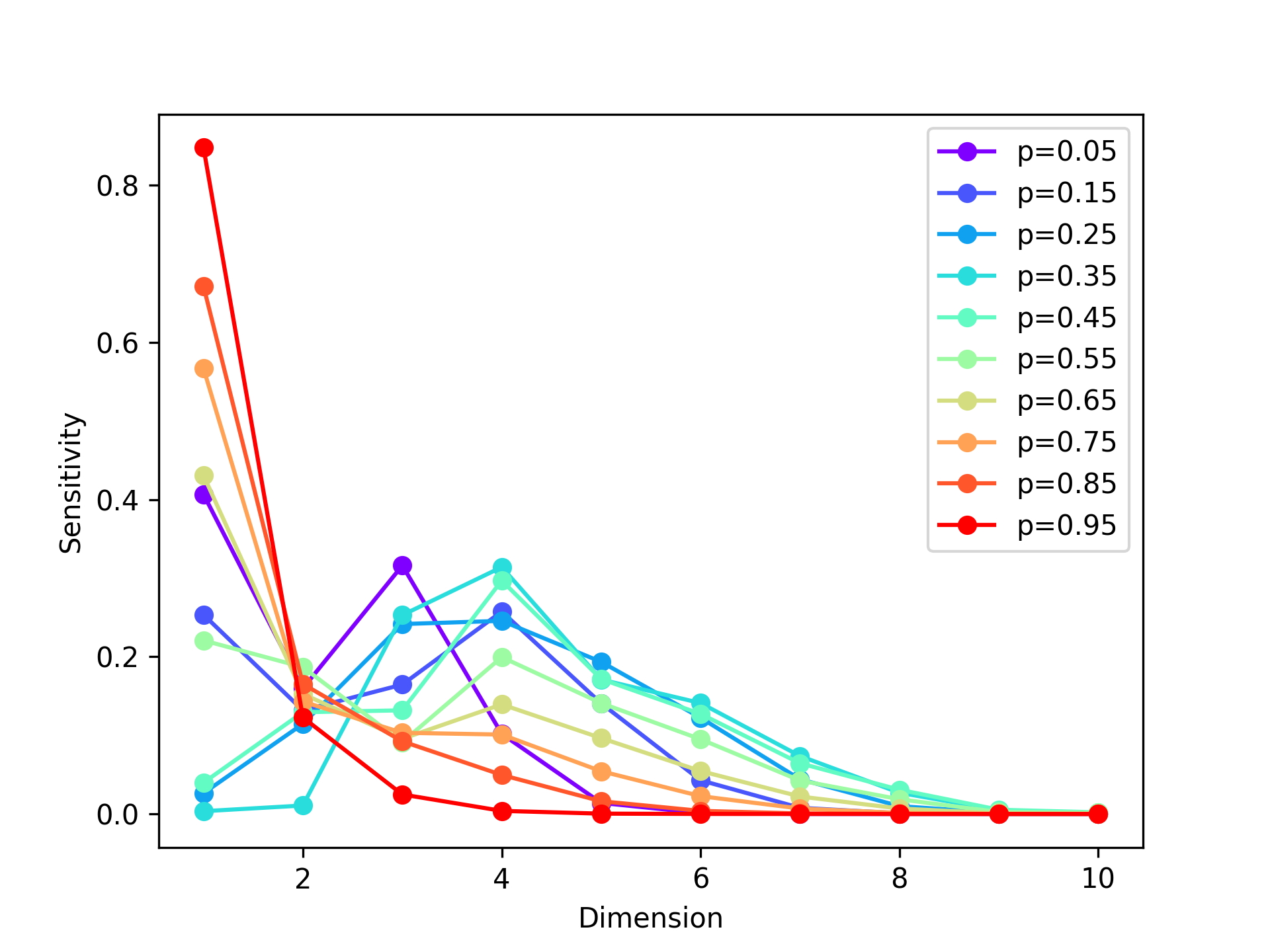}}
\subfloat[Entropy in $p$]{\includegraphics[width=3.5in]{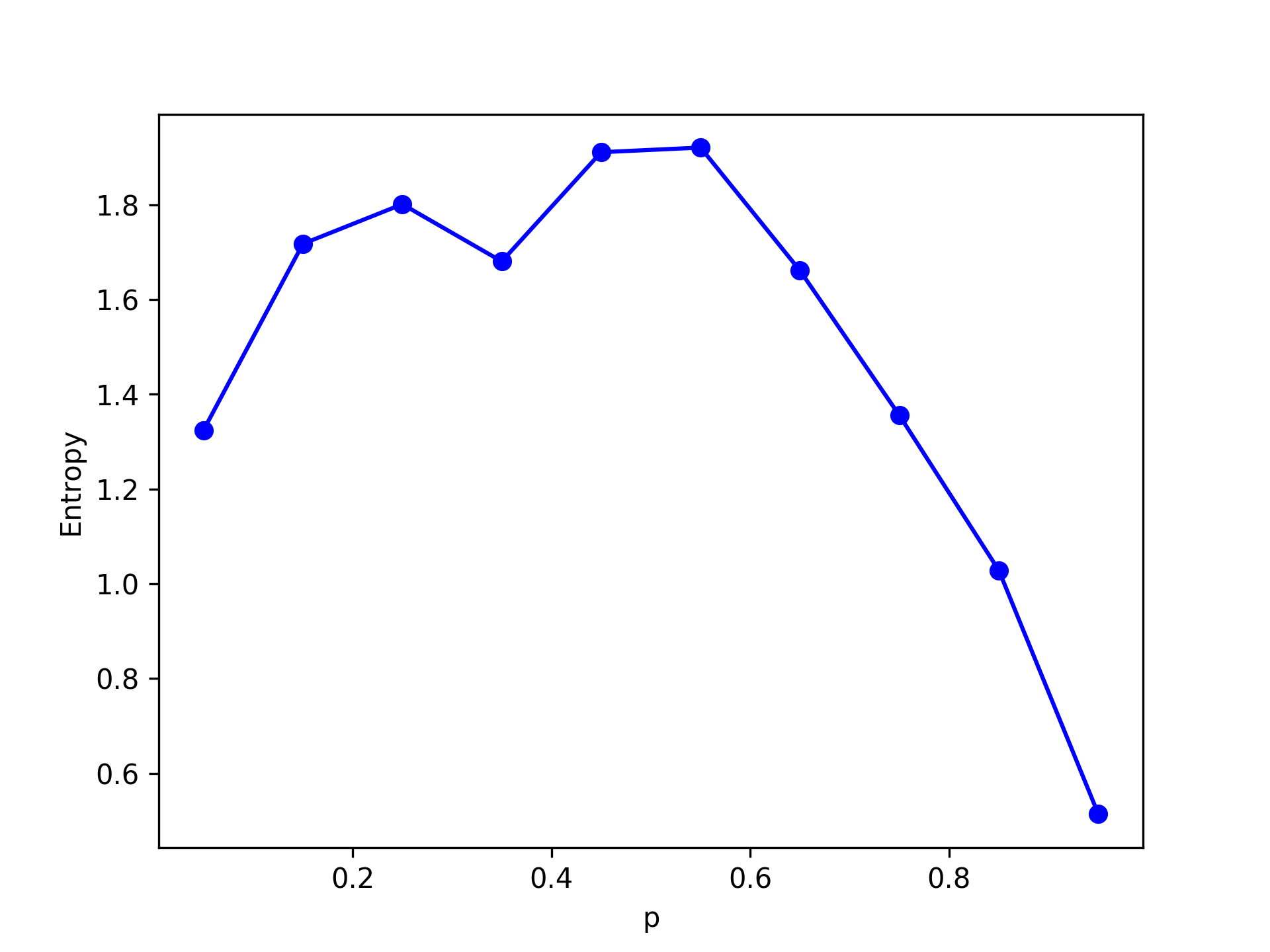}}\\
\endgroup
\caption{Spectral gap of graph Laplacian of Erd\H{o}s-Reny\'{i} $G(n=5,p)$: Mean, variance, sensitivity indices, and entropy as a function of $p$}\label{fig:8a}
\end{figure}

In possession of the component functions $\{g_u\}$, we compute the variances $(\Var g_u)$ and sensitivity indices $(\amsbb{S}_u^a)$. The sensitivity indices form a probability measure on the subspaces, and we compute its entropy. In Figure~\ref{fig:8a} for $n=5$ we plot the mean, variance, sensitivity indices, and entropy as a function of $p$. Mean and variance generally increase in $p$, whereas entropy increases then decreases, with a maximum near $p\approx 0.5$. 

\begin{figure}[h!]
\centering
\begingroup
\captionsetup[subfigure]{width=3in,font=normalsize}
\subfloat[Orthogonal $N$: $(\nu f_d^2)$ for coordinate 1]{\includegraphics[width=3.25in]{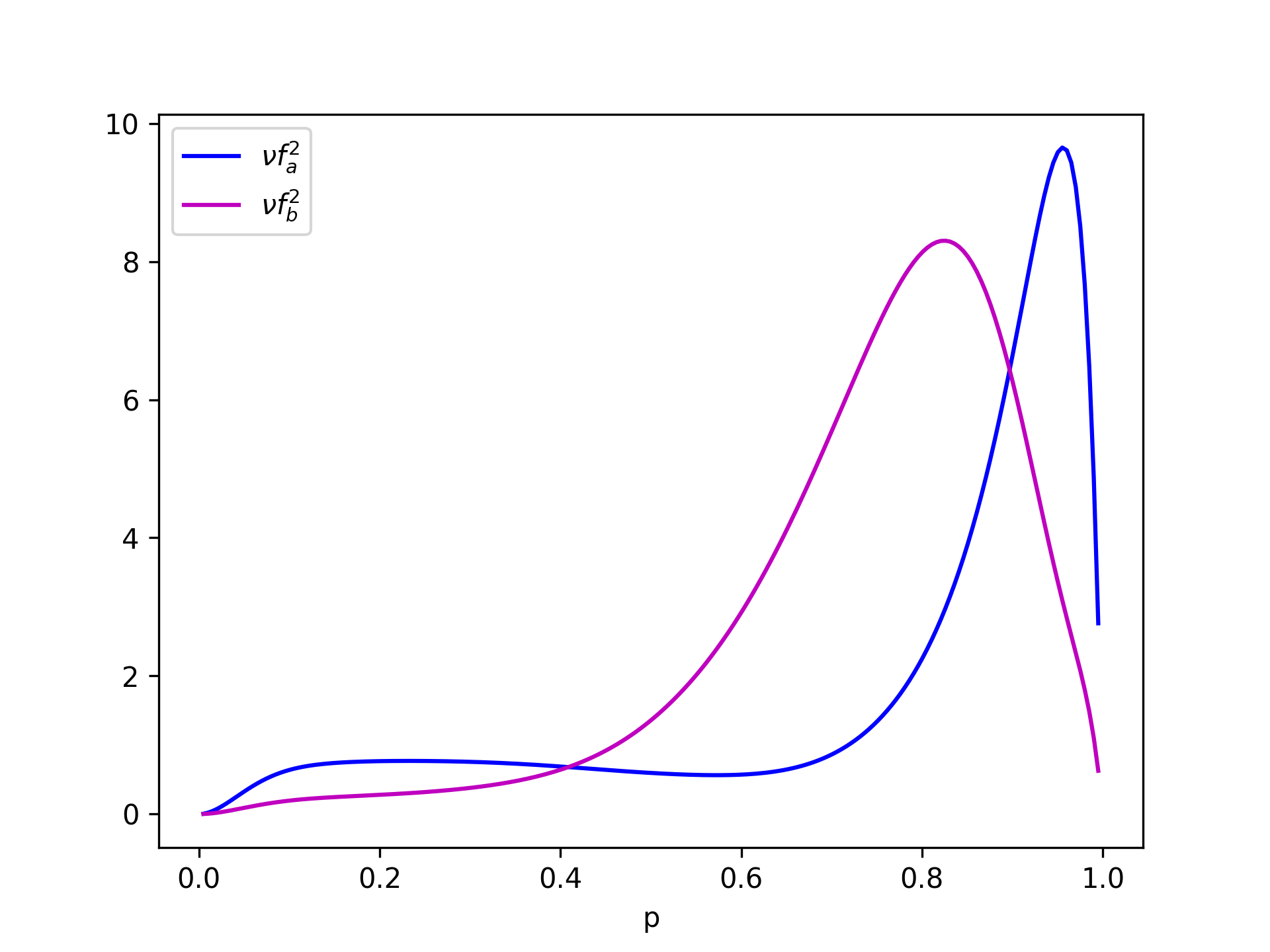}}
\subfloat[Orthogonal $N$: $(\amsbb{S}_d^a)$ for coordinate 1]{\includegraphics[width=3.25in]{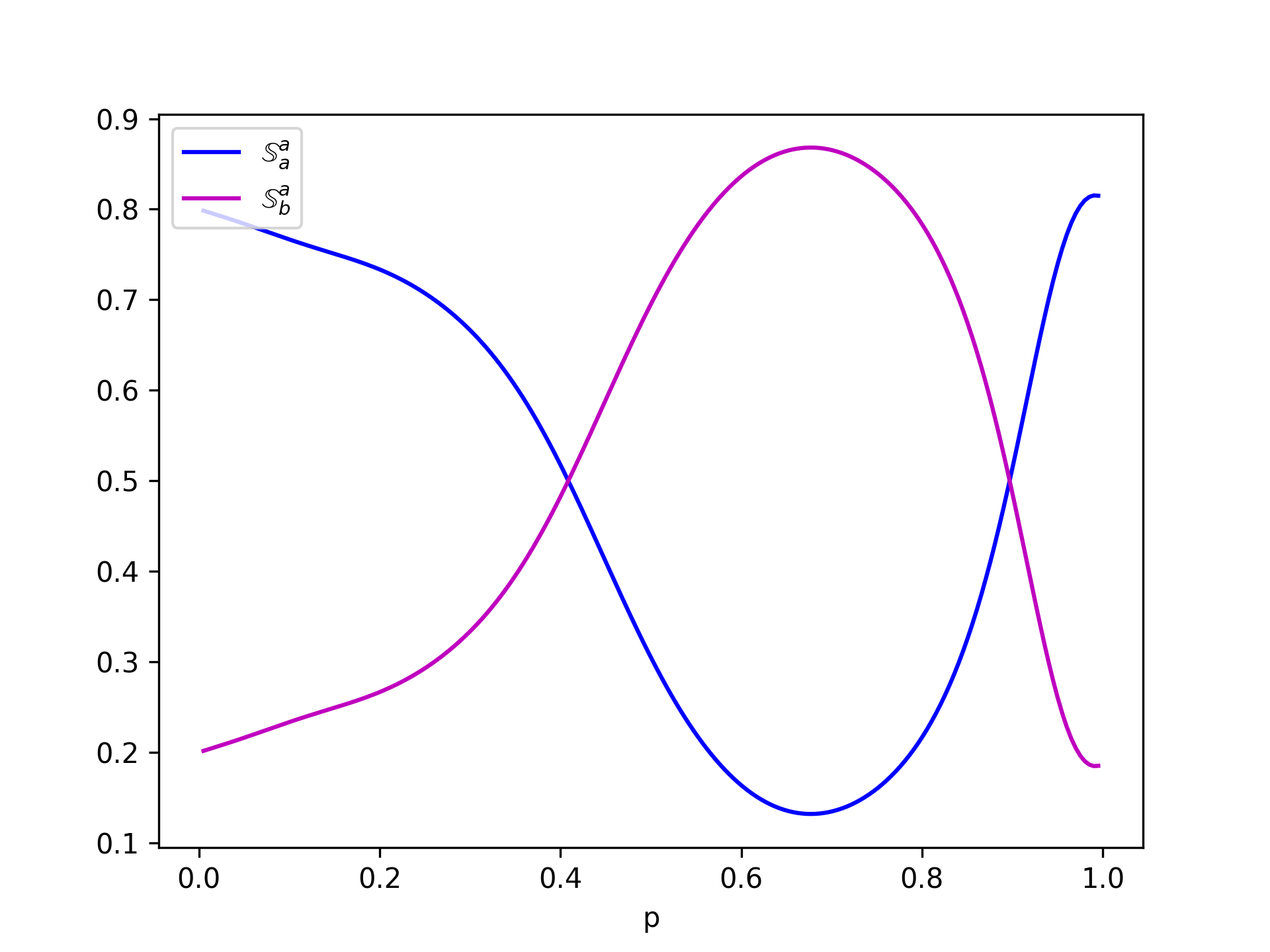}}\\
\subfloat[Dirac $N$: $(\Var f_d)$ for coordinate 1]{\includegraphics[width=3.25in]{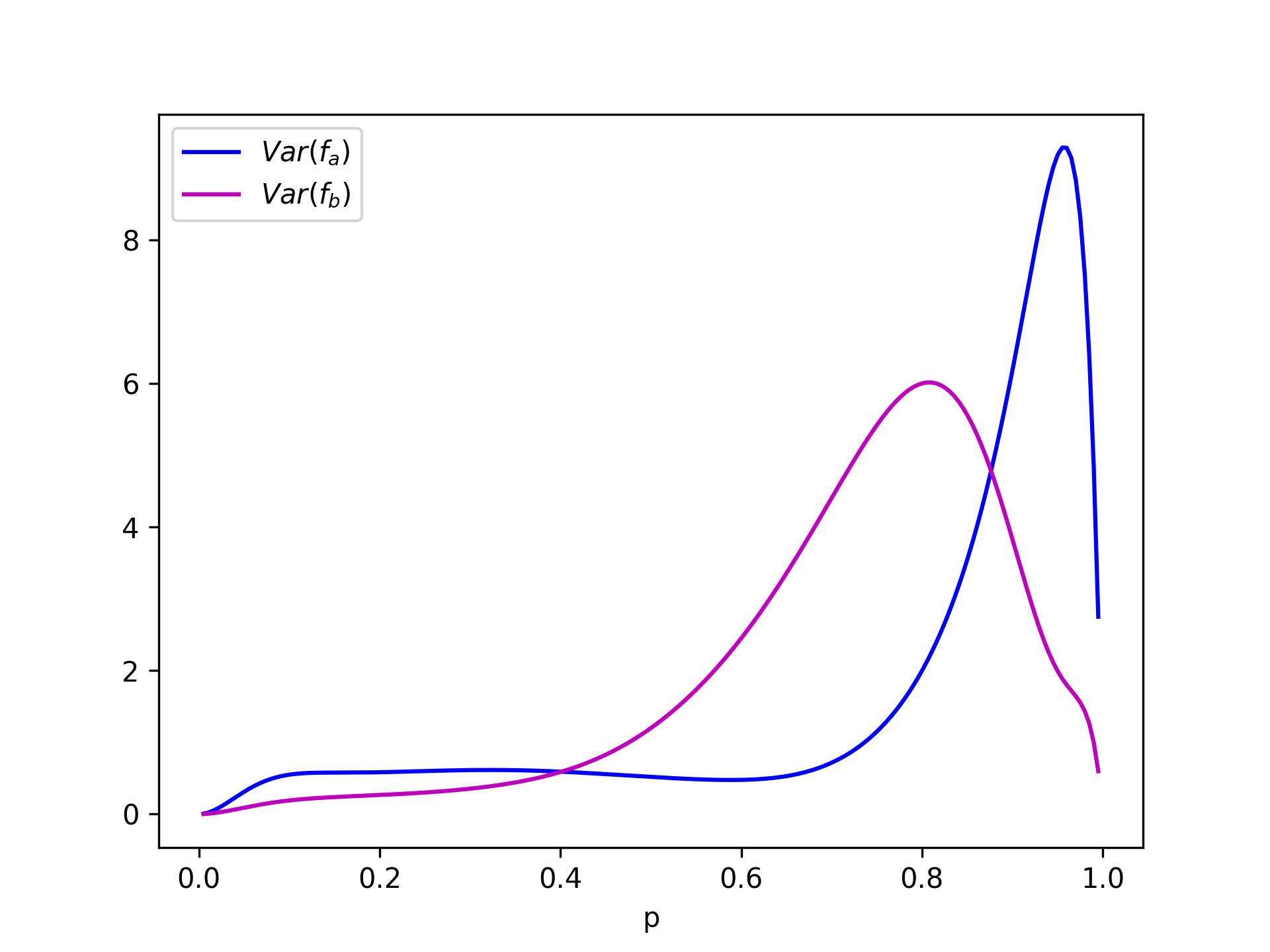}}
\subfloat[Dirac $N$: $(\amsbb{S}_d^a)$ for coordinate 1]{\includegraphics[width=3.25in]{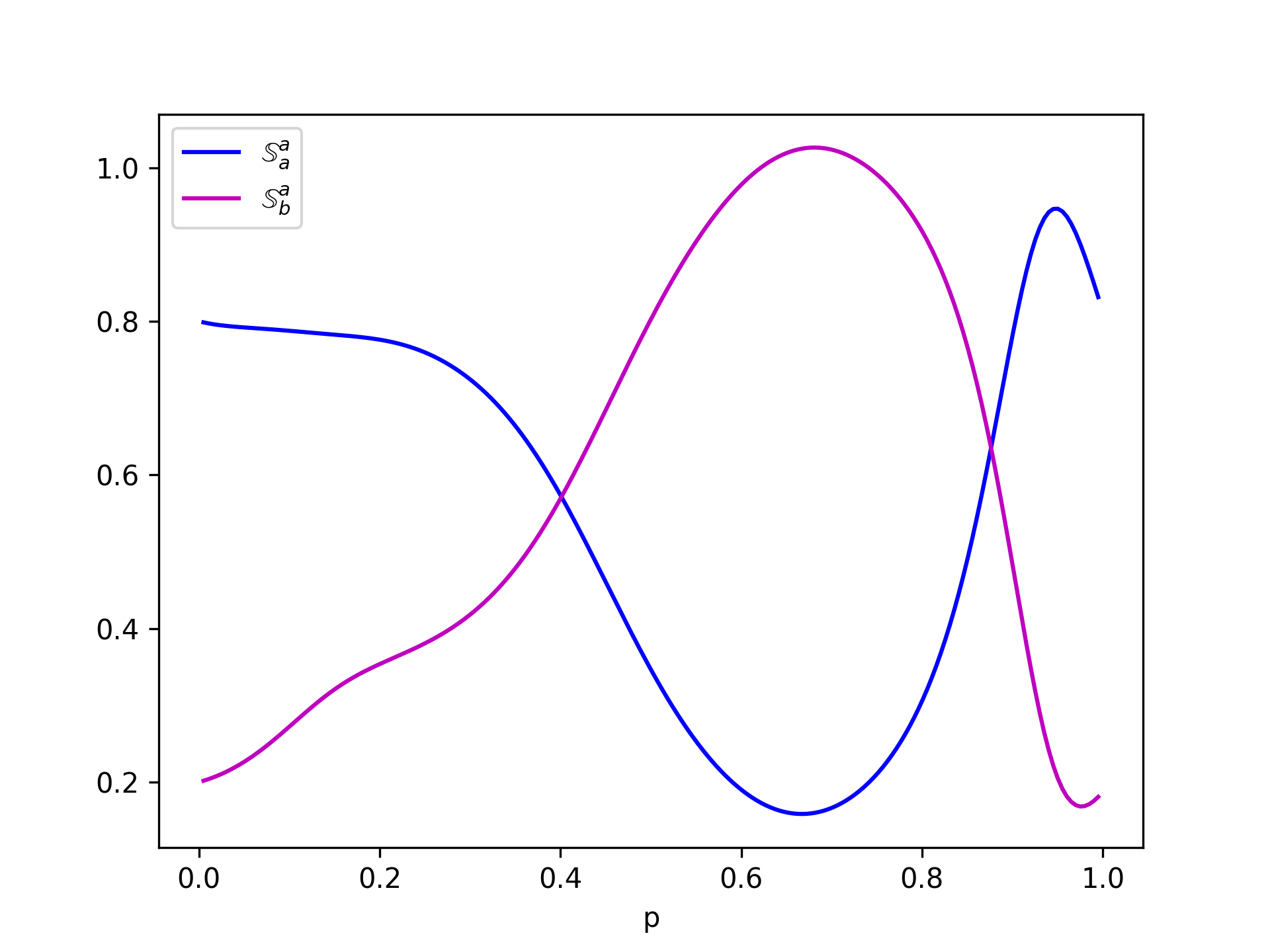}}\\
\subfloat[Entropy of $(\amsbb{S}_d^a)$ for coordinate 1 as a function of $p$]{\includegraphics[width=3.25in]{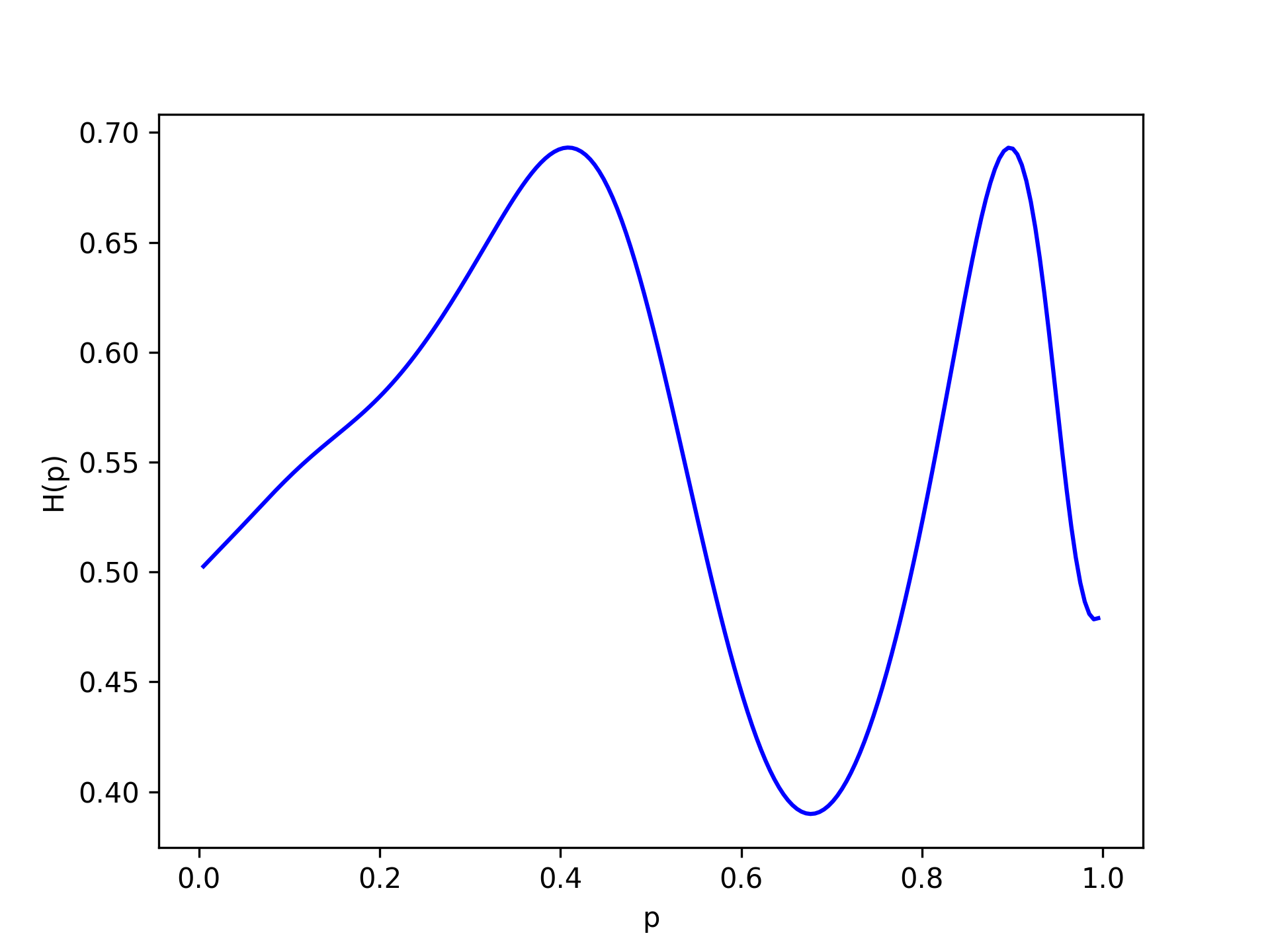}}
\subfloat[$\nu f^2$ and $\Var f$ as a function of $p$]{\includegraphics[width=3.25in]{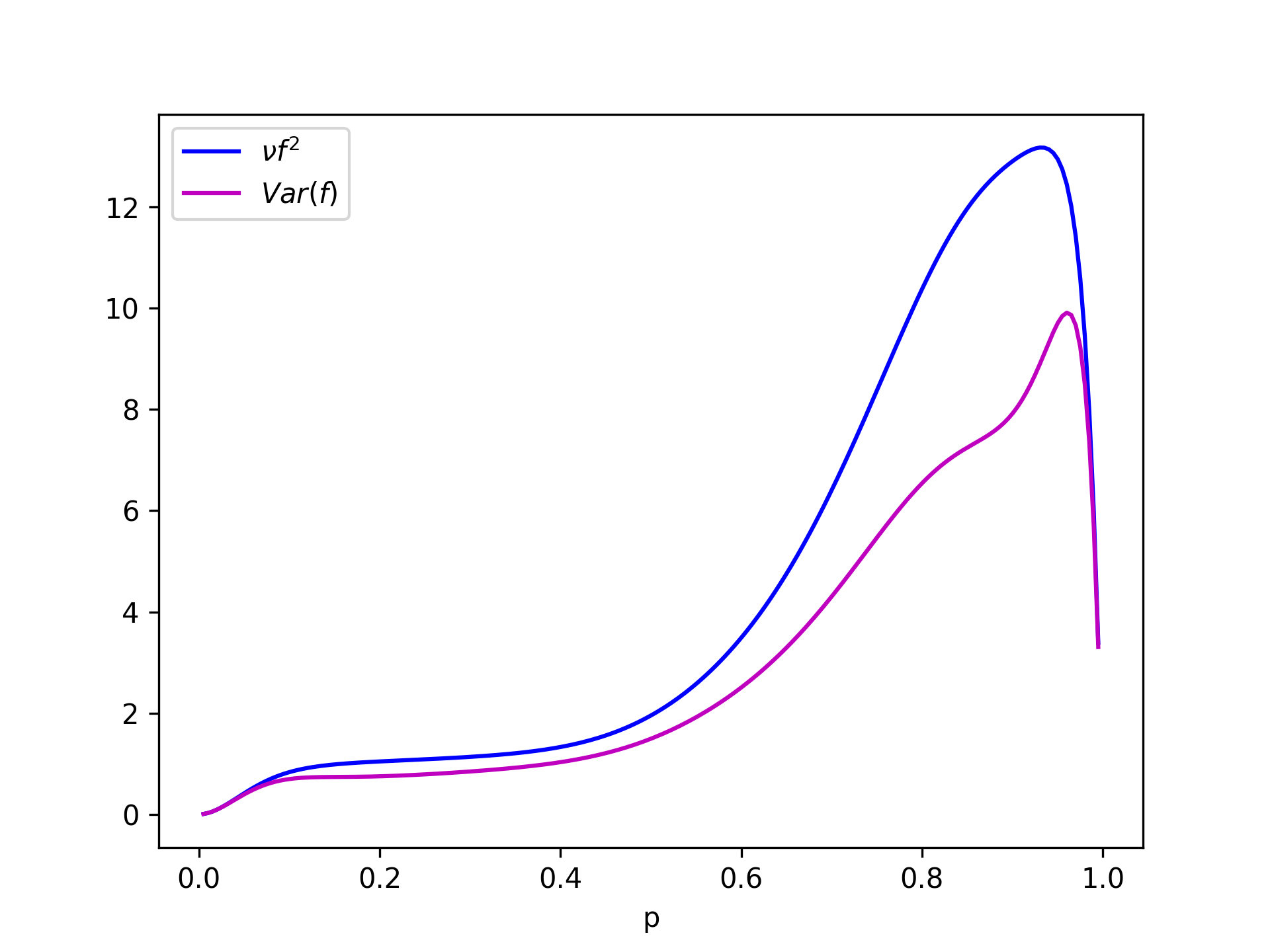}}\\
\endgroup
\caption{Spectral gap of graph Laplacian of Erd\H{o}s-Reny\'{i} $G(n=5,p)$ for orthogonal and Dirac $N$: second moments $(\nu f_d^2)$, variances $(\Var f_d)$, and structural sensitivity indices $(\amsbb{S}_d^a)$ for partitions by coordinates}\label{fig:8}
\end{figure}

Now consider the random measure $N=(\kappa,\nu)$ on $(E,\mathcal{E})$ and consider risk function $f=(g-\E g)^2\in\mathcal{E}_+$. The Laplace functional is defined through \[\nu e^{-f} = \int_E\nu(\D x)e^{-(g(x)-\E g)^2}\] Then the mean and variance are computed as \begin{align*}\E Nf &= c \nu f = c\Var g\\\Var Nf &= c\nu f^2 - (\delta^2-c)(\nu f)^2\end{align*}  As before we consider Dirac and orthogonal $N$. For each coordinate $i$, we take the partition $\{A_i,B_i\}$ of $E$, where $A_i=E_1\times\dotsb\times E_{i-1}\times\{0\}\times E_{i+1}\times\dotsb\times E_{\binom{n}{2}}$ and $B_i=E_1\times\dotsb\times E_{i-1}\times\{1\}\times E_{i+1}\times\dotsb\times E_{\binom{n}{2}}$ for the $i\in\{1,\dotsb,\binom{n}{2}\}$. We take $n=5$ and plot the spectral gap RM analysis in Figure~\ref{fig:8} of the variances, second moments, and sensitivity indices. The Dirac and orthogonal $N$ are similar, with non-linear surfaces. Entropy has two modes.


\FloatBarrier

\subsection{Ising random field}\label{sec:ising} As referenced in Table~\ref{tab:0}, we give an example for Remark~\ref{re:graph} for the nearest neighbor Ising model on a lattice. Let $N=(\kappa,\nu)$ be a random measure on $(E,\mathcal{E})$. Consider measurable space $(F,\mathcal{F})$, where $F=\{-1,1\}^{|V|}$ with dimension $2^{|V|}$ and discrete $\sigma$-algebra $\mathcal{F}=2^F$. 

Let $k:E\times F\mapsto\R_+$ be $\mathcal{E}\otimes\mathcal{F}$ measurable. Consider the random field $G$ on $(F,\mathcal{F})$ formed by $G(y)=Nf_y$ for $y\in F$ where $f_y(\cdot)=k(\cdot,y)\in\mathcal{E}_+$. Define the Hamiltonian function \[H(y) = -\sum_{i\sim j}y_iy_j\] and put \[Z_\beta = \sum_{y\in F}e^{-\beta H(y)}\] where $\beta\in\R_+$. The density is defined as \[\lambda_\beta\{y\} = \frac{e^{-\beta H(y)}}{Z_\beta}\for y\in F\] To define the Laplace transform of $G$, we define $f_\beta\in\mathcal{E}_+$ as \[f_\beta(x)=\sum_{y\in F}\lambda_\beta\{y\}k(x,y)\for x\in E\] and compute $\nu e^{-f_\beta}$. Let $|V|=n\in\N_{>0}$. 

\begin{figure}[h!]
\centering
\begingroup
\captionsetup[subfigure]{width=3in,font=normalsize}
\subfloat[Orthogonal $N$]{\includegraphics[width=3.5in]{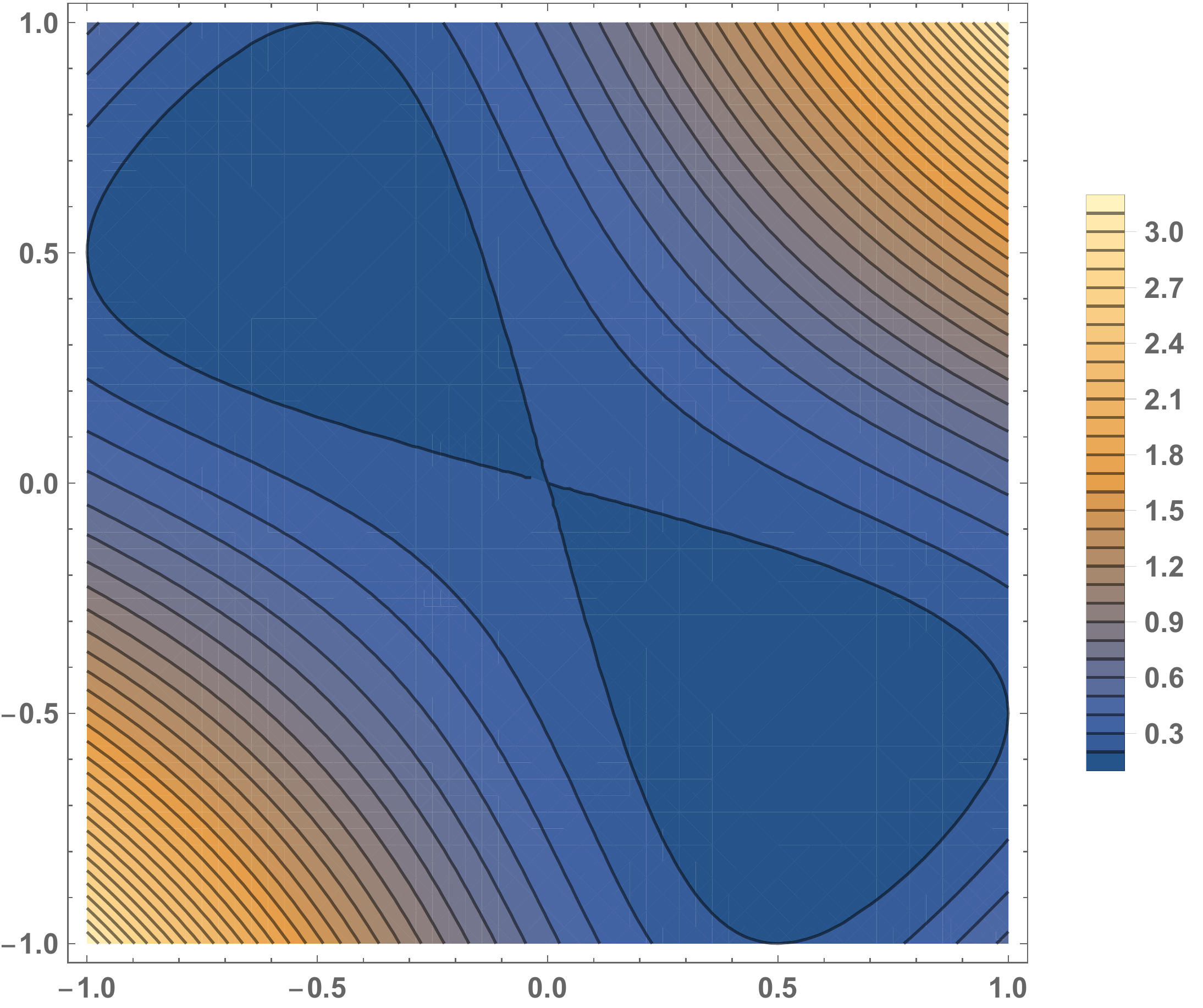}}
\subfloat[Dirac $N$]{\includegraphics[width=3.5in]{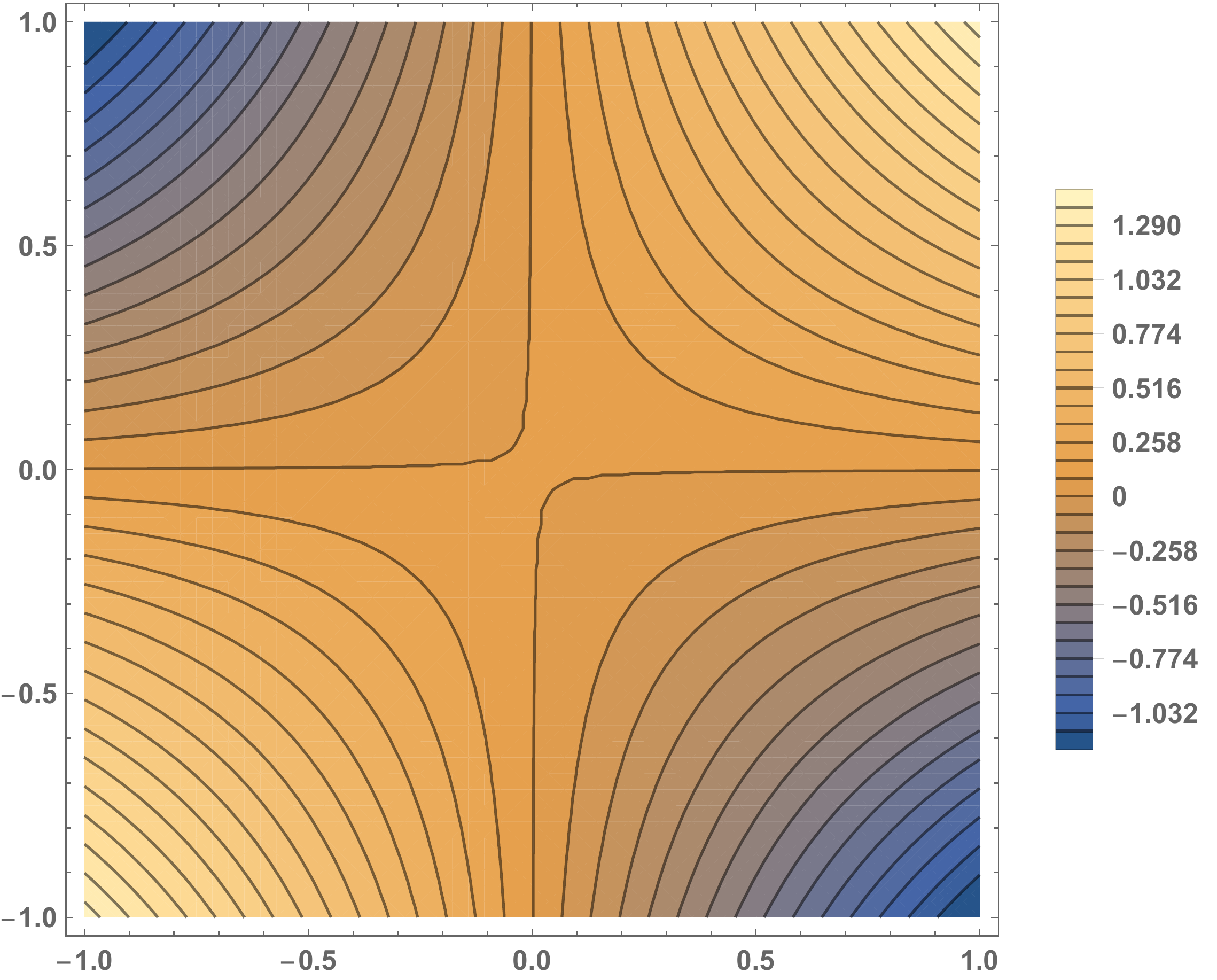}}\\
\endgroup
\caption{Random field kernel $K$ for orthogonal and Dirac $N$}\label{fig:kernel2}
\end{figure}

The average spin (magnetization) is \[g(y)=\frac{1}{n}\sum_{i}^n y_i\for y\in F\] Let $\nu=\text{Uniform}[-1,1]$ with $(E,\mathcal{E})=([-1,1],\mathcal{B}_{[-1,1]})$. Define $k$ as \[k(x,y) = (x-g(y))^2\for (x,y)\in E\times F\] The random field $G$ has the interpretation of magnetization distance on lattice points. We have \[\nu f_y = g^2(y)+\frac{1}{3}\for y\in F\] and \[\nu(f_yf_z) =\frac{1}{15} \left(5 g^2(y) \left(3 g^2(z)+1\right)+20 g(y) g(z)+5 g^2(z)+3\right)\for y,z\in F \] In Figure~\ref{fig:kernel2} we plot the kernel $K$ for orthogonal and Dirac $N$ and $c=1$.

\FloatBarrier

\end{document}